\documentclass[10pt, reqno]{amsart}
\pdfoutput=1
\usepackage{amssymb, amsmath, amsthm, graphics, latexsym, enumitem}

\usepackage{imakeidx}
\makeindex[title=Topics]
\makeindex[title=People,name=ppl]

\usepackage{lmodern}
\usepackage{newunicodechar}
\usepackage{ragged2e}
\usepackage{biblatex}
\usepackage{hyperref}
\usepackage{xurl}
\hypersetup{
    breaklinks=true,
    colorlinks=true,
    hyperindex=true,
    linkcolor=blue,
    filecolor=magenta,      
    urlcolor=cyan,
    pdftitle={Chronology},
    pdfpagemode=FullScreen}

\renewcommand{\descriptionlabel}[1]%
	{\hspace{\labelsep}\textsc{#1},}

\def\nats{\mathbb{N}}
\def\prls{\mathbb{R}_+}
\def\rls{\mathbb{R}}

\newunicodechar{ű}{\H{u}}
\newunicodechar{ő}{\H{o}}

\newcommand{\dK}{\ensuremath{\mathop\textrm{\LARGE{\lower.20ex\hbox{K}\normalsize}}}}
\newcommand{\tK}{\ensuremath{\mathop\textrm{\large{\lower.08ex\hbox{K}\normalsize}}}}

\def\surdex#1{\sqrt{\vphantom{#1}}#1}

\DeclareMathOperator{\arccot}{arccot}

\theoremstyle{plain}
\theoremstyle{definition}

\newtheorem*{theorem*}{Theorem}

\begin{document} 

\title[Continued square roots and continued compositions, through 2016]
{A chronology of continued square roots{\break}and other continued compositions,{\break}through the year 2016}
\author{Dixon J. Jones}
\address{Coralville, Iowa USA}
\email{d.j.jones.1798@gmail.com}
\thanks{The author is grateful to Z. Franco\index[ppl]{Franco, Z.}, S. G. Moreno\index[ppl]{Moreno, Samuel G\'{o}mez}, M. Somos\index[ppl]{Somos, Michael}, P. Trinh\index[ppl]{Trinh, P.}, and S. Wepster\index[ppl]{Wepster, S.} for structure and content suggestions; to M. P. Getz\index[ppl]{Getz, M. P.}, M. A. Jones\index[ppl]{Jones, M. A.}, J. P. Lambert\index[ppl]{Lambert, J. P.}, and W. Tape\index[ppl]{Tape, W.} for comments and advice; and to the Interlibrary Loan Department of the Elmer E. Rasmuson Library\index{Rasmuson Library, Elmer E.}, University of Alaska Fairbanks\index{University of Alaska Fairbanks}, for invaluable assistance. Any errors or omissions are entirely the author's.}

\date{\today}

\begin{abstract}
An infinite continued composition is an expression of the form
\begin{equation*}
\lim_{n\to\infty}t_0\circ t_1 \circ t_2 \circ \cdots \circ t_n(c)\;,
\end{equation*}
where the $t_i$ are maps from a set $D$ to itself, the initial value $c$ is a point in $D$, and the order of operations proceeds from right to left.\index{continued compositions!definition of}

This document is an annotated bibliography, in chronological order through the year 2016, of selected continued compositions whose primary sources have typically been obscure. In particular, we include continued square roots:
\begin{equation*}
a_0+\sqrt{a_1+\sqrt{a_2+\sqrt{\ldots}}}\;,
\end{equation*}
as well as continued powers, continued cotangents, and $f$-expansions. However, we do not include continued fractions, continued exponentials, or forms such as infinite sums and products in which the $t_i$ are linear functions, because the literature on these forms is extensive and well-summarized elsewhere.

\end{abstract}

\maketitle


\begin{RaggedRight}

\section{Introduction}\label{S:intro}

This is a bibliography, in chronological order through the year 2016, of selected \emph{continued function compositions} or, more briefly, \emph{continued compositions}\index{continued compositions!definition of}\footnote{In early drafts of this document, posted on a personal web site beginning in 2008, we used the term ``chain composition," which seemed truer to German \emph{Ketten-} constructions like \emph{Kettenbr\"{u}che}. However, ``continued" has long-entrenched connotations in English, and will probably be recognized more widely than ``chain."}, expressions of the form\index{continued compositions!definition of}
\begin{equation}\label{E:ccomp}
\lim_{n\to\infty}t_0\circ t_1 \circ t_2 \circ \cdots \circ t_n(c)\;,
\end{equation}
where the $t_i$ are maps from a set $D$ to itself, the \emph{initial value} $c$ is a point in $D$, and the order of operations proceeds from right to left. Under these conditions the \emph{$n$th approximant}\index{continued compositions!approximants of}
\begin{equation}\label{E:appx}
t_0\circ t_1 \circ t_2 \circ \cdots \circ t_n(c)
\end{equation}
is well-defined for each $n=0,1,2,\ldots$, and the continued composition converges if the sequence of approximants 
\begin{equation}\label{E:appx123}
t_0(c),\quad t_0\circ t_1(c), \quad t_0\circ t_1\circ t_2(c),
\end{equation}
and so on, has a unique, finite limit. Following \hyperref[Sch1992]{\textsc{Sch\"{o}nefuss 1992}}\index[ppl]{Schonefuss@Sch\"{o}nefuss, Lutz W.}, in emulation of a common notation for continued fractions we write expression \eqref{E:ccomp} more compactly as $\tK_{i=0}^\infty t_i(c)$.

Many mathematical objects can be thought of as continued compositions, but this chronology is limited to work in which one of the following appears. (We assume that $D$ is the complex plane here, but of course $D$ varies in the chronology from item to item.)

\begin{itemize}

\item\vspace{9pt} The \emph{continued, infinite, iterated,} or \emph{nested square root}\index{continued square roots} 
\begin{equation}\label{E:croot}
a_0+\sqrt{a_1+\sqrt{a_2+\sqrt{\ldots}}}\;,
\end{equation}
generated by $t_i(z)=a_i+\sqrt{z}$ with $c=0$. These generalize naturally to \emph{continued $r$th roots}\index{continued robb roots@continued $r$th roots} using $t_i(z)=a_i+\sqrt[r]{z}$ for $r>1$. 

\item\vspace{9pt} The \emph{continued $p$th power}\index{continued pobb powers@continued $p$th powers}
\begin{equation}\label{E:contpower}
a_0+(a_1+(a_2+(\ldots)^p)^p)^p\;,
\end{equation}
using $t_i(z)=a_i+z^p$, $p\in\mathbb{R}$, and $c=0$ for $p>0$; for $p<0$ we use $c=\infty$. In \hyperref[Dop1832b]{\textsc{Doppler 1832b}}\index[ppl]{Doppler, Christian} we find perhaps the earliest observation that this expression becomes a regular continued fraction\index{continued fractions} for $p=-1$; \hyperref[Her1935]{\textsc{Herschfeld 1935}}\index[ppl]{Herschfeld, Aaron} notes this as well, along with the fact that infinite sums\index{infinite sums} comprise the case $p=1$. An example of a continued square $(p=2)$\index{continued squares} is given without comment in \hyperref[Dix1878]{\textsc{Dixon 1878}}\index[ppl]{Dixon, T. S. E.}. The case
\[a_0+(a_1+(a_2+(\ldots)^{p_2})^{p_1})^{p_0}\;,\]
where the $p_i$ are in $(0,1)$ and may not all be equal (a form which could also be called a \emph{continued $r_i$th root}\index{continued rocc roots@continued $r_i$th roots}, where $r_i=\tfrac{1}{p_i}$), is mentioned in \hyperref[Her1935]{\textsc{Herschfeld 1935}}\index[ppl]{Herschfeld, Aaron}, \hyperref[And1985]{\textsc{Andrushkiw 1985}}\index[ppl]{Andrushkiw, R. L.}, and \hyperref[Muk2013]{\textsc{Mukherjee 2013}}\index[ppl]{Mukherjee, Soumendu Sundar}.

\item\vspace{9pt} Lehmer's\index[ppl]{Lehmer, D. H.} \emph{continued cotangent}\index{continued cotangents}, so-named because the continued composition using $t_i(x)=(a_ix+1)/(x-a_i)$ with $c=0$ simplifies to
\[\cot(\arccot a_0-\arccot a_1 +\arccot a_2 -\ldots)\;.\]
See \hyperref[Leh1938]{\textsc{Lehmer 1938}}\index[ppl]{Lehmer, D. H.}.

\item\vspace{9pt} The \emph{continued logarithm}\index{continued logarithms}, a name which has been applied to several disparate forms. Here we refer primarily to continued compositions in which $t_i(z)=a_i+\log_b(z)$. Other kinds of ``continued logarithms" use variants of $t_i(z)=t(z)=z/\ln(z)$, for instance in approximating the Lambert $W$ function\index{Lambert $W$ function}.\footnote{See e.g. Yanghua Wang\index[ppl]{Wang, Yanghua}, The Ricker wavelet\index{Ricker wavelet} and the Lambert $W$ function. \emph{Geophysical Journal International} (2015) \textbf{200}(1), 111--115, \href{https://doi.org/10.1093/gji/ggu384}{Oxford Academic}.} The ``binary (base 2) continued logarithm,"\index{continued logarithms!binary (base 2)} introduced by Gosper\index[ppl]{Gosper, Bill}\footnote{B. Gosper, Continued fraction arithmetic. Perl Paraphernalia, at \href{http://perl.plover.com/classes/cftalk/INFO/gosper.txt}{perl.plover.com}. See also J. M. Borwein\index[ppl]{Borwein, Jonathan M.}, K. G. Hare\index[ppl]{Hare, K. G.}, and J. G. Lynch\index[ppl]{Lynch, J. G.}, Generalized continued logarithms and related continued fractions, \href{https://arxiv.org/abs/1606.06984}, 22 June 2016}, is a continued fraction\index{continued fractions} of the form
\begin{equation*}
2^{a_0}+\cfrac{2^{a_0}}{2^{a_1} +\cfrac{2^{a_1}}{2^{a_2}+\cfrac{2^{a_2}}{\ddots}}}\;,
\end{equation*}
where one is interested in the integer sequence $\{a_n\}$.

\item\vspace{9pt} General continued compositions\index{continued compositions} of the form \eqref{E:ccomp}, including \emph{$f$-expansions}\index{f@$f$-expansions} (\hyperref[Bis1944]{\textsc{Bissinger 1944}})\index[ppl]{Bissinger, B. H.}, \emph{infinite processes}\index{infinite processes} (\hyperref[Thr1961]{\textsc{Thron 1961}})\index[ppl]{Thron, W. J.}, \emph{Kettenoperationen}\index{Kettenoperationen@\emph{Kettenoperationen}}\index{continued operations} or \emph{continued operations} (\hyperref[Lau1990]{\textsc{Laugwitz 1990}}\index[ppl]{Laugwitz, Detlef}), and others (see e.g. \hyperref[Kak1924]{\textsc{Kakeya 1924}}\index[ppl]{Kakeya, S\^{o}ichi}, \hyperref[Pau2013]{\textsc{Paulsen 2013}}\index[ppl]{Paulsen, William}). Because there is an extensive literature on $f$-expansions, and most of it dates from after 1950, only a few citations are included here. The book \hyperref[Sch2016]{\textsc{Schweiger 2016}}\index[ppl]{Schweiger, Fritz} gives a comprehensive overview of $f$-expansions and related topics.\index{f@$f$-expansions}
\end{itemize}


\bigskip Among the continued compositions not covered in this bibliography are:

\begin{itemize}

\item\vspace{9pt} The \emph{infinite sum}\index{infinite sums} $\sum a_i$, where $t_i(z)=a_i+z$ and $c=0$, and other examples for which the $t_i(z)$ are linear in $z$\index{continued compositions!of linear functions}, such as the \emph{Engel expansion}\index[ppl]{Engel, F.}\footnote{F. Engel, Entwicklung der Zahlen nach Stammbruechen, \emph{Verhandlungen der zweiundfünfzigsten Versammlung deutscher Philologen und Schulm\"{a}nner in Marburg}, 1913, 190--191, \href{https://babel.hathitrust.org/cgi/pt?id=uiug.30112109810074&seq=646}{HathiTrust}.}
\[\dfrac{1}{a_0}\Bigl(1+\dfrac{1}{a_1}\Bigl(1+\dfrac{1}{a_2}\Bigl(1+\ldots\Bigr)\Bigr)\Bigr)\;,\]
where $t_i(z)=(1+z)/a_i$ and $c=0$; and the \emph{Pierce expansion}\index[ppl]{Pierce, T. A.}\footnote{T. A. Pierce, On an algorithm and its use in approximating roots of algebraic equations, \emph{The American Mathematical Monthly} \textbf{36}(10), 1929, 523--525, \href{https://doi.org/10.1080/00029890.1929.11987017}{doi.org (Taylor \& Francis)}.}
\[a_0(1-a_1(1-a_2(1-\ldots)))\;,\]
using $t_i(z)=a_i(1-z)$ and $c=0$.

The \emph{infinite product}\index{infinite products} $\prod a_i$, generated by $t_i(z)=a_iz$ with $c=1$, is another continued composition of linear functions not covered here. However, we do include \emph{infinite products of finite continued compositions}, of which Vi\`{e}te's formula for $\tfrac{2}{\pi}$ is probably the most prominent example; see \hyperref[Vie1593]{\textsc{Vi\`{e}te 1593}}.\footnote{The exception being made for Vi\`{e}te's formula points up the fact that we are generally excluding sequences in linear terms of two or more indices, such as $\sum_i\prod_j a_{i,j}$.}

\item\vspace{9pt} Most types of \emph{continued fractions}\index{continued fractions}, such as the \emph{regular continued fraction}
\begin{equation*}
a_0+\cfrac{1}{a_1 +\cfrac{1}{a_2+\cfrac{1}{\ddots}}}\;,
\end{equation*}
generated by $t_i(z)=a_i+1/z$ with $c=\infty$. (It is more standard to use $t_i(z)=1/(a_i+z)$ with $c=0$, and tack on or ignore the term $a_0$.)

\item\vspace{9pt} The \emph{continued, infinite,} or \emph{iterated exponential}\index{continued exponentials}, or \emph{tower}
\[a_0e^{a_1e^{a_2e^{\cdot^{\cdot^{\cdot}}}}},\]
where $t_i(z)=a_ie^z$ and $c=1$. On the positive reals one can unambiguously define $t_i(x)=a_i^x$ to obtain
\[a_0^{a_1^{a_2^{\cdot^{\cdot^{\cdot}}}}}\;.\]

\end{itemize}

\noindent\bigskip The literature on these forms is too extensive to list here.\footnote{More than 6,000 unannotated references related to continued fractions\index{continued fractions} (including many of the sources cited in this chronology) are given in Claude Brezinski\index[ppl]{Brezinski, Claude}, \emph{A Bibliography on Continued Fractions, Pad\'{e} Approximation, Sequence Transformation and Related Subjects}, Ciencias [Sciences], 3. Prensas Universitarias de Zaragoza, Zaragoza, 1991. 348 pp. MR1148890, ISBN 847733238X.

A bibliography of 125 references on continued exponentials\index{continued exponentials} may be found in R. Arthur Knoebel\index[ppl]{Knoebel, R. Arthur}, Exponentials reiterated, \emph{The American Mathematical Monthly} \textbf{88}(4), 1981, 235--252, \href{https://mathscinet.ams.org/mathscinet/article?mr=610484}{MR0610484}, \href{https://doi.org/10.1080/00029890.1981.11995239}{doi.org (Taylor \& Francis)}.}


\section{The issue of associativity: continued versus iterated compositions}\label{S:assoc}\index{associativity of function compositions}

The superficially similar \emph{iterated composition}\index{iterated compositions}
\begin{equation}\label{E:iterated}
\lim_{n\to\infty} t_n\circ t_{n-1}\circ\cdots t_2\circ t_1\circ t_0(c)\;,
\end{equation}
where the order of operations proceeds from right to left, is sometimes confused with \eqref{E:ccomp}, particularly in the case where the $t_i$ are identical, that is, when $t_i=t$ for all $i=0,1,2\ldots$ In such a case, the approximants
\begin{equation}\label{E:fixpoint}
t(c)\;, \quad t\circ t(c)\;, \quad t\circ t\circ t(c)\;,
\end{equation}
and so on, can be viewed either as continued or as iterated compositions. 

Denoting an iterated composition's \emph{$n$th} iterate by
\[u_n(c)=t_n\circ t_{n-1}\circ\cdots t_2\circ t_1\circ t_0(c)\;,\]
one has the \emph{forward} recurrence relation\index{recurrence relation}
\[u_n(c)=t_n\circ u_{n-1}(c)\;,\]
but note that, in general, this does \emph{not} generate a continued composition. Indeed, using the abbreviation
\[\dK_j^n=t_{j}\circ t_{j+1}\circ\cdots\circ t_n(c)\;,\]
one has the \emph{backward} recurrence relation\index{recurrence relation}
\[\dK_j^n=t_j\circ \dK_{j+1}^{n}.\]

Up until 1911, most examples of continued nonlinear function compositions\footnote{Here we consider ``continued linear function compositions" to include continued fractions\index{continued fractions}.} were actually iterated compositions\index{iterated compositions} of the form \eqref{E:fixpoint}, constructed using the method of successive substitution.\index{successive substitution} In this procedure, one attempts to solve an equation of the form $x=t(x)$ by repeatedly substituting $t(x)$ for the argument of $t$:
\begin{align*}
x&=t(t(x))\;,\\
x&=t(t(t(x)))\;,\\
x&=t(t(t(t(x))))\;,
\end{align*}
and so on, which amounts to an algorithm for computing a fixed point of $t$. The conditions on $t$ under which this algorithm would converge began to be rigorously addressed in the second half of the 19th century; we cite only a few of the many papers on forward iteration from this period.

Ramanujan's continued square root identities from 1911 appear to be the first published examples of continued compositions \eqref{E:ccomp} in which the $t_i$ are nonlinear and not identical.\index[ppl]{Ramanujan, Srinivasa}\index{continued square roots!Ramanujan's}


\section{Terminology}\label{S:terminology}

To name the expression $\sqrt{a_0+\sqrt{a_1+\sqrt{a_2+\cdots}}}$, the works cited and quoted here have employed nearly every possible concatenation of one or more words from the set $\{$\emph{continued, infinite, nested, iterated}$\}$ with a word from the set $\{$\emph{root, radical, power}$\}$. However, the text not directly quoted from primary sources attempts to follow two guidelines in naming function compositions:\index{continued compositions!naming conventions for}

\begin{enumerate}

\item Adjectives should be specific about the associativity of the expression. Thus, as discussed in Section \ref{S:assoc}, \emph{continued} as a modifier of a named function composition indicates formation by a backward recurrence relation, while \emph{iterated} indicates a forward recurrence relation\index{recurrence relation}. On the other hand, the adjectives \emph{nested} and \emph{infinite} are ambiguous with regard to associativity.

\item Nouns should refer to functions rather than symbols. Thus we prefer \emph{continued square root} and \emph{continued $r$th root} to the less specific \emph{continued radical}.

\end{enumerate}


\section{About the format}\label{S:format} 

Citations are listed 1) in chronological order by year of publication, 2) alphabetically by author last name within a given year, and, if the author(s) published more than one paper in a year, 3a) chronologically when the order of publication is evident, or 3b) alphabetically by journal title. Undated references are listed last. Where possible, \emph{Mathematical Reviews} (MR) numbers, book ISBNs, and hyperlinks to online sources are provided. Hyperlinks and MR numbers between wedge symbols, e.g. $\langle${MR9999999}$\rangle$, indicate that the linked source post-dates the year 2016; for example, a preprint circulated before 31 Dec 2016 might have been published, and therefore reviewed, after that date. Most hyperlinks to external sources attempt to take the viewer directly to the article being annotated.


\section{About this version}\label{S:version} 

New or substantially modified listings are indicated with a dagger \dag. 

\begin{itemize}
\item \hyperref[Ber1692]{\textsc{Bernoulli 1692}} (page \pageref{Ber1692}),
\item \hyperref[Leg1816]{\textsc{Legendre 1816}} (page \pageref{Leg1816}),
\item \hyperref[Gal1830]{\textsc{Galois 1830}} (page \pageref{Gal1830}),
\item \hyperref[Fou1831]{\textsc{Fourier 1831}} (page \pageref{Fou1831}),
\item \hyperref[Pio1837]{\textsc{Pioch 1837}} (page \pageref{Pio1837}),
\item \hyperref[Hei1844]{\textsc{Heis 1844}} (page \pageref{Hei1844}),
\item \hyperref[Boc1895]{\textsc{Bochow 1895}} (page \pageref{Boc1895}),
\item \hyperref[Boc1896]{\textsc{Bochow 1896}} (page \pageref{Boc1896}),
\item \hyperref[Hey1898]{\textsc{Heymann 1898}} (page \pageref{Hey1898}),
\item \hyperref[Can1901]{\textsc{Cantor 1901}} (page \pageref{Can1901}),
\item \hyperref[Hey1904]{\textsc{Heymann 1904}} (page \pageref{Hey1904}),
\item \hyperref[Boc1905]{\textsc{Bochow 1905}} (page \pageref{Boc1905}),
\item \hyperref[Bal1920]{\textsc{Ball 1920}} (page \pageref{Bal1920}),
\item \hyperref[Gal2001]{\textsc{Galuzzi 2001}} (page \pageref{Gal2001}).
\end{itemize}

Minor changes have been made to 

\begin{itemize}
\item \hyperref[Arc-250]{\textsc{Archimedes $\sim$250 BCE}} (page \pageref{Arc-250}),
\item \hyperref[Ceu1596]{\textsc{Van Ceulen 1596}} (page \pageref{Ceu1596}),
\item \hyperref[MRS1830]{\textsc{``M. R. S." 1830}} (page \pageref{MRS1830}),
\item \hyperref[Dop1832b]{\textsc{Doppler 1832b}} (page \pageref{Dop1832b}),
\item \hyperref[Dix1878]{\textsc{Dixon 1878}} (page \pageref{Dix1878}),
\item \hyperref[Gun1880]{\textsc{G\"{u}nther 1880}} (page \pageref{Gun1880}).
\end{itemize}

Lastly, we have added an appendix of \hyperref[S:bios]{biographical sketches} for a few selected authors.
%


\section{Caveat lector}\label{S:caveat} 

This bibliography should be considered an inventory, not a comprehensive or interpretative history. Given the ever-increasing store of digitized source material, it should also be considered an incomplete work in progress. It was compiled and annotated by an interested layman who is not an authority on any of the mathematical topics it mentions, and who does not claim proficiency in languages other than English. Corrections, and suggestions for additional sources from 2016 and earlier, are welcome. 


\section{The Chronology}\label{S:chrono}

\begin{enumerate}[label={\arabic*}.]


\item\vspace{9pt} ($\sim$250 BCE) \label{Arc-250} Archimedes of Syracuse\index[ppl]{Archimedes|textbf}, \emph{Kuklou metr\={e}sis [Measurement of a circle]}. 

{\footnotesize For an English translation of the original Greek, see T. L. Heath\index[ppl]{Heath, T. L.}, \emph{The Works of Archimedes, ed. in Modern Notation, with Introductory Chapters}, C. J. Clay and Sons, Cambridge University Press Warehouse, London, 1897. Heath's translation was reprinted by Dover Publications, Inc., Mineola, New York, 2002, ISBN-13: 978-0486420844.}

Archimedes almost certainly did not manipulate continued square roots. Nonetheless, it is appropriate to cite this fragmentary manuscript as the inspiration for many of the continued square root expressions, listed below, which arise in the calculation of $\pi$.\index{pi@$\pi\;(3.14159\ldots)$!Archimedean algorithm for}\index{constants, named!pi@$\pi\;(3.14159\ldots)$} Archimedes' Proposition 2 is credited as the first known use of perimeters of regular polygons\index{polygons, regular} (specifically, 96-gons obtained iteratively by halving the sides of a hexagon), inscribed in and circumscribed around a circle, to compute lower and upper bounds, respectively, for $\pi$. Making rational approximations to irrational lengths at each step, Archimedes obtains the famous inequality
\[\dfrac{223}{71}<\pi<\dfrac{22}{7}\;.\]
Sometimes credited, sometimes not, the ``Archimedean algorithm"\index{Archimedean algorithm} of polygon construction will be employed by many authors, including \hyperref[Vie1593]{\textsc{Vi\`{e}te 1593}}\index[ppl]{Vi\`{e}te, Fran\c{c}ois}, \hyperref[Ceu1596]{\textsc{Van Ceulen 1596}}\index[ppl]{Ceulen, Ludolph van}, \hyperref[Eul1744]{\textsc{Euler 1744}}\index[ppl]{Euler, Leonhard}, \hyperref[Ens1799]{\textsc{Ensheim 1799}}\index[ppl]{Ensheim, Moses}, and \hyperref[Cat1842]{\textsc{Catalan 1842}}\index[ppl]{Catalan, E.}, to name just a few.


\item\vspace{9pt}(1593)\label{Vie1593} Fran\c{c}ois Vi\`{e}te\index[ppl]{Vi\`{e}te, Fran\c{c}ois|textbf}, \emph{Variorum de Rebus Mathematicis Responsorum, Liber VIII.}
\index{Vi\`{e}te's formula for $\tfrac{2}{\pi}$}

{\footnotesize Source: \href{https://archive.org/details/bub_gb_SM2Az-oKgoEC/page/n63/mode/2up}{Internet Archive}.}

{\footnotesize Reprinted, with an English translation of the Latin, in \emph{Pi: A Source Book, 3rd edition}, edited by Len Berggren\index[ppl]{Berggren, Len}, Jonathan Borwein\index[ppl]{Borwein, Jonathan M.}, and Peter Borwein\index[ppl]{Borwein, Peter}, Springer-Verlag, New York, 2004, 53--67 and 690--706. ISBN 978-0387205717.}

{\footnotesize Cited or referenced in
\hyperref[Rud1891]{\textsc{Rudio 1891}}\index[ppl]{Rudio, Ferdinand},
\hyperref[Bop1913]{\textsc{Bopp 1913}}\index[ppl]{Bopp, K.}, 
\hyperref[Leb1937]{\textsc{Lebesgue 1937}}\index[ppl]{Lebesgue, Henri},  
\hyperref[Osl1999]{\textsc{Osler 1999}}\index[ppl]{Osler, Thomas J.}, 
\hyperref[Ser2003]{\textsc{Servi 2003}}\index[ppl]{Servi, L. D.}, 
\hyperref[Lev2005]{\textsc{Levin 2005}} and \hyperref[Lev2006]{\textsc{2006}}\index[ppl]{Levin, Aaron},
\hyperref[Lim2007a]{\textsc{Lim 2007a}}\index[ppl]{Lim, Teik-Cheng}, 
\hyperref[Osl2007]{\textsc{Osler 2007}}\index[ppl]{Osler, Thomas J.}, 
\hyperref[Mor2013a]{\textsc{Moreno and Garc\'{i}a-Caballero 2013a}} and \hyperref[Mor2013b]{\textsc{2013b}}\index[ppl]{Moreno, Samuel G\'{o}mez}\index[ppl]{Garc\'{i}a-Caballero, Esther M.}, 
\hyperref[Gar2014b]{\textsc{Garc\'{i}a-Caballero, Moreno, and Prophet 2014b}}\index[ppl]{Garc\'{i}a-Caballero, Esther M.}\index[ppl]{Moreno, Samuel G\'{o}mez}\index[ppl]{Prophet, Michael P.}, 
\hyperref[Osl2016a]{\textsc{Osler 2016}}\index[ppl]{Osler, Thomas J.}, 
\hyperref[Osl2016b]{\textsc{Osler, Jacob, and Nishimura 2016}}\index[ppl]{Osler, Thomas J.}\index[ppl]{Jacob, Walter}\index[ppl]{Nishimura, Ryo}, and
\hyperref[Wei]{\textsc{Weisstein n.d.}}\index[ppl]{Weisstein, Eric W.}}

\vspace{4pt}

This work is generally considered the first to have expressed a number using nested square root expressions, namely
\begin{equation}\label{E:Vie1593-1}
\dfrac{2}{\pi}=
\sqrt{\tfrac{1}{2}}\cdot
\sqrt{\tfrac{1}{2}+\tfrac{1}{2}\sqrt{\tfrac{1}{2}}}\cdot
\sqrt{\tfrac{1}{2}+\tfrac{1}{2}\sqrt{\tfrac{1}{2}+\tfrac{1}{2}\sqrt{\tfrac{1}{2}}}}\cdot\cdots\;,
\end{equation}
a brute which has defeated a few subsequent authors and their typesetters. Indeed, H. W. Turnbull\index[ppl]{Turnbull, H. W.} gives a version\footnote{\emph{The World of Mathematics}, edited by James R. Newman\index[ppl]{Newman, James R.}, Simon and Schuster, New York, 1956, p. 121.} whose associativity\index{associativity of function compositions} is not easily parsed, which in turn led George B. Thomas, Jr.\index[ppl]{Thomas, George B., Jr.}, to present a misinterpreted and incorrect version\footnote{\emph{Calculus and Analytic Geometry, Alternate Edition}, Addison-Wesley Publishing Company, Inc, Reading, Massachusetts, 1972, p. 834.}. \hyperref[Her1935]{\textsc{Herschfeld 1935}}\index[ppl]{Herschfeld, Aaron} cites J. W. L. Glaisher\index[ppl]{Glaisher, J. W. L.} (\emph{Messenger of Mathematics} \textbf{2} (new series), 1873, p. 124) for the slightly cleaner
\begin{equation}\label{E:Vie1593-2}
\dfrac{2}{\pi}=\dfrac{\sqrt{2}}{2}\cdot
\dfrac{\sqrt{2+\sqrt{2}}}{2}\cdot
\dfrac{\sqrt{2+\sqrt{2+\sqrt{2}}}}{2}\cdot\cdots\;.
\end{equation}
but \hyperref[Cat1842]{\textsc{Catalan 1842}}\index[ppl]{Catalan, E.} derived the reciprocal of this form earlier (and apparently unaware of Vi\`{e}te)\index[ppl]{Vi\`{e}te, Fran\c{c}ois}. 

The formula is obtained by a variation on the Archimedean method\index{Archimedean algorithm}. Rather than approximating a circle's perimeter between inscribed and circumscribed regular polygons\index{polygons, regular}, Vi\`{e}te compares the areas of inscribed regular $2^n$-gons and $2^{n + 1}$-gons. A notably clear explanation of Vi\`{e}te's method is given in the appendix of \hyperref[Mor2013b]{\textsc{Moreno and Garc\'{i}a-Caballero 2013b}}\index[ppl]{Moreno, Samuel G\'{o}mez}\index[ppl]{Garc\'{i}a-Caballero, Esther M.}.

It is worth noting that the pages reprinted in \emph{Pi: A Source Book} are not from Vi\`{e}te's original 16th century work, but from \emph{Francisci Vietae, Opera Mathematica, in Unum Volumen Congesta}, a nicely typeset compilation of Vi\`{e}te's works published in 1646 by Frans van Schooten Junior\index[ppl]{Schooten, Frans van}, a professor at the Engineering School, Duytsche Mathematique, loosely connected to Leiden University. Vi\`{e}te did not use the ``modern" square root notation\index{notation!surd $\surd$} employed by Van Schooten.\footnote{\hyperref[Lef1897]{\textsc{Lefebvre 1897}}\index[ppl]{Lefebvre, B.}, \hyperref[Caj1928]{\textsc{Cajori 1928}}\index[ppl]{Cajori, Florian}, and others cite the 1634 \emph{Les Oeuvres Math\'{e}matiques} of Simon Stevin\index[ppl]{Stevin, Simon}, edited by Albert Girard\index[ppl]{Girard, Albert}, as one of the earliest sources of the  square root symbol as it is currently used. Cajori furthermore credits Descartes' 1637 \emph{G\'{e}om\'{e}trie} with introducing the \emph{vinculum} (the trailing ``over-bar") in conjunction with the surd symbol $\surd$ to create $\sqrt{\phantom{m}}$.} In the original \emph{Variorum}, equation \eqref{E:Vie1593-1} above is something like a diagram and description, almost a flow chart, in Latin, occupying about three quarters of the page that lies between ``pages" 30 and 31 (page numbers appear only on one of two facing pages). 


\item\vspace{9pt} \label{Ceu1596} (1596) Ludolph van Ceulen\index[ppl]{Ceulen, Ludolph van|textbf}, \emph{Vanden circkel}, Jan Andriesz, Delft.

{\footnotesize Source: \href{https://catalog.lindahall.org/discovery/delivery/01LINDAHALL_INST:LHL/1287123520005961}{Linda Hall Library}.}

{\footnotesize Cited in
\hyperref[Bos1910]{\textsc{Bosmans 1910.}}\index[ppl]{Bosmans, H., S.J.}}

\vspace{4pt}

In this book, his first, Van Ceulen computes $\pi$\index{pi@$\pi\;(3.14159\ldots)$!Archimedean algorithm for}\index{constants, named!pi@$\pi\;(3.14159\ldots)$} to 20 decimal places, using the Archimedean method\index{Archimedean algorithm} of regular polygons\index{polygons, regular} inscribed in and circumscribed around a circle. Van Ceulen's computations are based on his expansions of polygon side lengths as finite continued square roots of the constant term 2 (sometimes with other integers or fractions as the rightmost term). As noted in \hyperref[Caj1928]{\textsc{Cajori 1928}}, Van Ceulen\index[ppl]{Ceulen, Ludolph van} employed the ``cossic" square root notation\index{notation!cossic}, used as early as 1525 by Christoff Rudolff\index[ppl]{Rudolff, Christoff} (in his \emph{Behend vnnd Hubsch Rechnung durch die kunstreichen regeln Algebre so gemeincklich die Coss genent werden}): a dot is printed immediately to the right of the surd symbol to indicate that the radical continues to the right. Thus, in cossic notation\index{notation!cossic},
\[\surd.2+\surd.2+\surd.2+\surd.2\]
is equivalent to the modern
\[\sqrt{2+\sqrt{2+\sqrt{2+\sqrt{2}}}}\;.\]
One could view Van Ceulen's\index[ppl]{Ceulen, Ludolph van} finite expansions, more than Vi\`{e}te's\index[ppl]{Vi\`{e}te, Fran\c{c}ois} infinite one, as perhaps the earliest precursors of the continued square root form.

When Van Ceulen\index[ppl]{Ceulen, Ludolph van} died in 1610, his widow Adrienne Simons\index[ppl]{Simons, Adrienne} undertook to publish Latin versions of many of his mathematical works (originally written in Dutch) with the help of Willebrordus Snellius\index[ppl]{Snellius, Willebrordus (\emph{aka} Snell, Willebrord)} (whose un-Latinized name is variously given as Willebrord Snellius or Snell or Snel or Snel van Royen, and who is more famous for the property of wave refraction known as Snell's Law). The \emph{Vanden circkel} translation appears, along with three of Van Ceulen's pamphlets, in Snellius's \emph{Lvdolphi \`{a} Cevlen De Circvlo \& Adscriptis Liber}, Lvgd. Batav. Apud Iodocum \`{a} Colster, 1619. However, \hyperref[Bos1910]{\textsc{Bosmans 1910}}\index[ppl]{Bosmans, H., S.J.} comments that this translation has few pages in common with the original text.


\item\vspace{9pt} \label{Kep1621} (1621) Johannes Kepler\index[ppl]{Kepler, Johannes|textbf}, \emph{Epitome astronomiae Copernicanae. Libri V, VI, VII, Doctrina theorica.} Gottfried Tambach, Linz.
\index{root-finding approximations}
\index{iteration!17th century examples}

{\footnotesize Source: \href{https://www.digitale-sammlungen.de/en/view/bsb10626228?page=36,37}{M\"{u}nchener DigitalisierungsZentrum}.} 

{\footnotesize For an English translation of the original Latin, see \emph{Epitome of Copernican Astronomy; and, Harmonies of the World}, translated by Charles Glenn Wallis\index[ppl]{Wallis, Charles Glenn}, The St. John's Bookstore, Annapolis, 1939; reprinted by Prometheus Books, Amherst, NY, 1995, ISBN 978-1-57392-036-0.}

{\footnotesize Alluded to or mentioned in 
\hyperref[Pio1837]{\textsc{Pioch 1837}}\index[ppl]{Pioch, A.},
\hyperref[Hey1898]{\textsc{Heymann 1898}}\index[ppl]{Heymann, W.}, and 
\hyperref[Hey1904]{\textsc{1904}}. 
Cited in
\hyperref[Gau1809]{\textsc{Gauss 1809}}\index[ppl]{Gauss, Carl Friedrich} and
\hyperref[Swe2000]{\textsc{Swerdlow 2000}}\index[ppl]{Swerdlow, N. M.}. }

\vspace{4pt}

The historical development of continued compositions is intertwined with that of iterative methods, and an early instance is in section 5.2.4 of this book. Kepler attempts to get an approximate solution for $E$, given the quantities $M$ and $e$, in what has become known as ``Kepler's equation" or ``the Kepler problem"
\begin{equation}\label{E:Kep1621}
M=E+e\sin{E}\;.
\end{equation}
\hyperref[Swe2000]{\textsc{Swerdlow 2000}}\index[ppl]{Swerdlow, N. M.} gives a good explication of Kepler's geometric approach. The equation is also dealt with iteratively, but more efficiently, in \hyperref[Gau1809]{\textsc{Gauss 1809}}. Iterative solutions are also given in \hyperref[Pio1837]{\textsc{Pioch 1837}}\index[ppl]{Pioch, A.}, \hyperref[Hey1898]{\textsc{Heymann 1898}} and \hyperref[Hey1904]{\textsc{Heymann 1904}}\index[ppl]{Heymann, W.}.


\item\vspace{9pt} \label{Sne1621} (1621) Willebrordus Snellius\index[ppl]{Snellius, Willebrordus (\emph{aka} Snell, Willebrord)|textbf}, \emph{Cyclometricus: de circuli dimensione secundum logistarum abacos, \& ad mechanicem accuratissima; atque omnium parabilissima}. Elzevir, Lugdunum Batavorum [Leiden]. 

{\footnotesize Source: \href{https://catalog.lindahall.org/discovery/delivery/01LINDAHALL_INST:LHL/1287485760005961}{Linda Hall Library}.}

{\footnotesize Cited in
\hyperref[Bos1910]{\textsc{Bosmans 1910.}}\index[ppl]{Bosmans, H., S.J.}}

\vspace{4pt}

According to \hyperref[Bos1910]{\textsc{Bosmans 1910}}\index[ppl]{Bosmans, H., S.J.}, Snellius\index[ppl]{Snellius, Willebrordus (\emph{aka} Snell, Willebrord)} enthusiastically began translating \hyperref[Ceu1596]{\textsc{Van Ceulen 1596}} from Dutch into Latin, but became bored with Van Ceulen's enormous calculations. This drove him to find a more efficient version of the Archimedean algorithm\index{Archimedean algorithm}. The crux of Snellius' improvement is the observation that the perimeter of an inscribed $n$-gon converges to $\pi$ twice as fast as that of a circumscribed $n$-gon. Like Van Ceulen\index[ppl]{Ceulen, Ludolph van}, Snellius\index[ppl]{Snellius, Willebrordus (\emph{aka} Snell, Willebrord)} uses the cossic notation\index{notation!cossic} for his finite continued square roots. (See \hyperref[Bal1920]{\textsc{Ball 1920}}\index[ppl]{Ball, W. W. Rouse} for more details about the calculations of Van Ceulen and Snellius.)


\item\vspace{9pt} \label{New1669} (1669) Isaac Newton\index[ppl]{Newton, Isaac|textbf}, \emph{De analysi per aequationes numero terminorum infinitas}.
\index{iteration!of functions}
\index{root-finding approximations}
\index{Newton's or Newton-Raphson method}

{\footnotesize Source: \href{https://www.e-rara.ch/zut/content/zoom/2553722}{ETH-Bibliothek}.}

{\footnotesize For an English translation from the original Latin, see \emph{The Mathematical Papers of Isaac Newton, Vol. 2, 1667--1670,} edited by D. T. Whiteside\index[ppl]{Whiteside, D. T.}, Cambridge University Press, London, 1968, 221--227, ISBN-13 978-0521045964.}

{\footnotesize Alluded to or mentioned in
\hyperref[Fou1831]{\textsc{Fourier 1831}}\index[ppl]{Fourier, Joseph},
\hyperref[Dop1832b]{\textsc{Doppler 1832}}\index[ppl]{Doppler, Christian},
\hyperref[Pio1837]{\textsc{Pioch 1837}}\index[ppl]{Pioch, A.},
\hyperref[Bou1862]{\textsc{Bouch\'{e} 1862}}\index[ppl]{Bouch\'{e}, A.},
\hyperref[Ise1888]{\textsc{Isenkrahe 1888}}\index[ppl]{Isenkrahe, C.},
\hyperref[Hey1894a]{\textsc{Heymann 1894a}} and \hyperref[Hey1894b]{\textsc{1894b}}\index[ppl]{Heymann, W.}, 
\hyperref[Ise1897]{\textsc{Isenkrahe 1897}}\index[ppl]{Isenkrahe, C.}, and
\hyperref[Hey1898]{\textsc{Heymann 1898}}, \hyperref[Hey1901]{\textsc{1901}}, and \hyperref[Hey1904]{\textsc{1904}}\index[ppl]{Heymann, W.}. 
Cited in
\hyperref[Tou1996]{\textsc{Tourn\`{e}s 1996}}\index[ppl]{Tourn\`{e}s, Dominique}.}

\vspace{4pt}

Newton's iterative process for finding roots of equations (later known as Newton's method, or the Newton-Raphson method\index{Newton's or Newton-Raphson method}), was first revealed here. In this chronology, only \hyperref[Tou1996]{\textsc{Tourn\`{e}s 1996}}\index[ppl]{Tourn\`{e}s, Dominique} cites this work directly. However, several authors take Newton's method to be the inspiration for, or a special case of, the application of iterated and/or continued compositions in solving equations. \hyperref[Tou1996]{\textsc{Tourn\`{e}s 1996}} remarks, ``Without knowing whether there was transmission or rediscovery, the [method of successive substitution]\index{successive substitution!17th century examples} is found in the West in Vi\`{e}te\index[ppl]{Vi\`{e}te, Fran\c{c}ois} (1600), Kepler\index[ppl]{Kepler, Johannes} (1618), Harriot\index[ppl]{Harriot, Thomas} (1631), Oughtred\index[ppl]{Oughtred, William} (1652) and others, before Newton would reveal it systematically in \emph{De analysi per aequationes number terminorum infinitas} from 1669, then in the \emph{Methodus fluxionum serierum infinitarum} from 1671 and various subsequent writings$\,\ldots\,$ Newton, considerably further than his predecessors, gradually extended the technique of obtaining roots by successive extractions\,\ldots"


\item\vspace{9pt} \label{Ber1692} \dag(1692) Jacob Bernoulli\index[ppl]{Bernoulli, Jacob|textbf}, \emph{Positionum arithmeticarum de seriebus infinitis earumque summa finita pars altera}. Basel, Johann Konrad von Mechel. Reprinted in \emph{Jacobi Bernoulli, Basileensis, Opera. Tomus primus}, Sumptibus Hæredum Cramer \& Fratrum Philibert, Geneva (1744), 536--539.
\index{iteration!of functions}

{\footnotesize Source: \href{https://archive.org/details/jacobibernoulli00conggoog/page/n670/mode/2up}{Internet Archive}.}

{\footnotesize Mentioned in 
\hyperref[Gun1880]{\textsc{G\"{u}nther 1880}}\index[ppl]{Gunther@G\"{u}nther, Siegmund}. 
Cited in
\hyperref[Can1901]{\textsc{Cantor 1901}}\index[ppl]{Cantor, Moritz} and
\hyperref[Hey1904]{\textsc{Heymann 1904}}\index[ppl]{Heymann, W.}.}

\vspace{4pt}

According to \hyperref[Can1901]{\textsc{Cantor 1901}}\index[ppl]{Cantor, Moritz}, the five parts of Jacob Bernoulli's ``Arithmetical propositions on infinite series and their finite sum'' comprise his own work; over a period of about fifteen years, he assigned five of his students the task of defending his results as part of their dissertations. The last of the five was his nephew Nicolaus (I) Bernoulli\index[ppl]{Bernoulli, Nicolaus (I)}, who defended in 1704. The second, in 1692, was Hieronymus Beck\index[ppl]{Beck, Hieronymus}, who Cantor dismisses as one of ``[t]he other{\ldots}figures [who] have not acquired such merit that their names could claim special recognition.'' In the part defended by Beck, Bernoulli gives nine examples of identities derived by formal manipulation of infinite successive substitution. The following list reproduces the results of Propositions 27 through 35 as they are printed, with one notational exception: we write $x^2$ rather than $x\,x$ (Bernoulli uses exponent notation only for products of three or more factors).

\begin{equation*}
\begin{array}{c l} 
\text{Proposition}  & \text{Successive substitution identity}\\[2pt]
\hline
\phantom{\Big|}27  &  \surd(a\surd(a\surd(a\surd(a\surd(a\&c.)))))=a\;.\\
\phantom{\big|}28  &  \surd(a+\surd(a+\surd(a+\surd(a+\&c.))))=\tfrac{1}{2}+\surd(\tfrac{1}{4}+a)\;.\\[3pt]
\phantom{\big|}29  &  \surd(a\surd(b\surd(a\surd(b\surd(a\surd(b\&c.))))))=\;^3\!\!\!\surd a^2b\;.\\[3pt]
\phantom{\big|}30  &  \surd(a\;^3\!\!\!\surd(ab\surd(a\;^3\!\!\!\surd(ab,\&c.))))=\;^5\!\!\!\surd a^4b\;.\\[3pt]
		&  \;^3\!\!\!\surd(ab\surd(a\;^3\!\!\!\surd(ab\surd(a\&c.))))=\;^5\!\!\!\surd a^3b^2\;.\\[3pt]
\phantom{\big|}31  &  \surd(a\surd(a\surd(b \surd(a\surd(a\surd(b\&c.))))))=\;^7\!\!\!\surd a^6b\;.\\[3pt]
		&  \surd(a\surd(b\surd(a \surd(a\surd(b\surd(a\surd(a\&c.)))))))=\;^7\!\!\!\surd a^5b^2\;.\\[3pt]
		&  \surd(b\surd(a\surd(a \surd(b\surd(a\surd(a\&c.))))))=\;^7\!\!\!\surd a^3b^4\;.\\[3pt]
\phantom{\big|}32  &  \surd(-p+\surd(p^2+q\surd(-p+\surd(p^2+q\surd(-p+\&c.)))))\text{ is the}\\[1pt]
		&\text{root of the cubic equation }x^3=-2px+q\;.\\[3pt]
\phantom{\big|}33  &  \surd(+p+\surd(p^2+q\surd(+p+\surd(p^2+q\surd(\&c.)))))\text{ is the}\\[1pt]
		&\text{root of the cubic equation }x^3=+2px+q\;.\\[3pt]
\phantom{\big|}34  &  \surd(p\pm\surd(p^2-q\surd(p\pm\surd(p^2-q\surd(\&c.)))))\text{ is the}\\[1pt]
		&\text{root of the equation }x^3=+2px-q\;.\\[3pt]
\phantom{\big|}35  &  \surd(-p+\surd(p^2+r+q\surd(-p+\surd(p^2+r+q \&c.))))\text{ is the}\\[1pt]
		&\text{root of the biquadratic equation }x^4=-2px^2+qx+r\;.
\end{array}
\end{equation*}

Thus far, Bernoulli's\index[ppl]{Bernoulli, Jacob} Proposition 28 is the earliest to give 
\begin{equation}\label{E:Ber1692-1}
\dfrac{1+\sqrt{1+4a}}{2}
\end{equation}
as the limit of the continued square root\index{continued square roots!of terms $a_n=c$}
\begin{equation}\label{E:Ber1692-2}
\sqrt{a+\sqrt{a+\sqrt{a+\cdots}}}\;.
\end{equation}
predating the previous record-holder, \hyperref[Dop1832b]{\textsc{Doppler 1832b}}\index[ppl]{Doppler, Christian}, by $140$ years. Proofs are provided for the first three identities; the proof for Proposition 28 goes as follows: ``For if $x=\surd(a+\surd(a+\surd(a+\&c.)))$ then $x^2=a+\surd(a+\surd(a+\&c.))$ \& $x^2-a=\surd(a+\surd(a+\surd(a+\&c.)))=x$: hence $x^2=x+a$, \& $x=\tfrac{1}{2}+\surd(\tfrac{1}{4}+a)$. Q.E.D.'' After Proposition 35, Bernoulli says, ``All these Props. are demonstrated in the same way as Props. XXVII, XXVIII, \& XXIX.''


\item\vspace{9pt} \label{Eul1744} (1744) Leonhard Euler\index[ppl]{Euler, Leonhard|textbf}, De variis modis circuli quadraturam numeris proxime exprimendi. \emph{Commentarii Academiae Scientiarum Imperialis Petropolitanae} \textbf{9}, for the year 1737, 222--238. Enestr\"{o}m number 74.
\index{pi@$\pi\;(3.14159\ldots)$!continued square root expressions for}\index{constants, named!pi@$\pi\;(3.14159\ldots)$}
\index{infinite products!of secant functions}
\index{Archimedean algorithm} 

{\footnotesize Source: \href{https://scholarlycommons.pacific.edu/euler-works/74/}{The Euler Archive}.}
 
{\footnotesize Mentioned in 
\hyperref[Leb1937]{\textsc{Lebesgue 1937}}\index[ppl]{Lebesgue, Henri}.
Cited in 
\hyperref[Bop1913]{\textsc{Bopp 1913}}\index[ppl]{Bopp, K.}.}

\vspace{4pt}

In the first of nineteen sections, Euler uses inscribed and circumscribed 96-gons in the unit circle to establish the following bounds on $\pi$:
\begin{align}\label{E:Eul1744-1}
96&\sqrt{2-\sqrt{2+\sqrt{2+\sqrt{2+\sqrt{3}}}}}\\
&<\pi<\dfrac{192\sqrt{2-\sqrt{2+\sqrt{2+\sqrt{2+\sqrt{3}}}}}}{\sqrt{2+\sqrt{2+\sqrt{2+\sqrt{3}}}}}\;.\notag
\end{align}
In the last section he represents an arc length $A$ on the unit circle by the expansion
\begin{equation}\label{E:Eul1744-2}
A=\dfrac{\sin A}{\cos\tfrac{A}{2}\cdot\cos\tfrac{A}{4}\cdot\cos\tfrac{A}{8}\cdot\cos\tfrac{A}{16}\cdots}\;.
\end{equation}
The special case $A=\tfrac{\pi}{2}$ will appear in \hyperref[Eul1763]{\textsc{Euler 1763}} (cited in \hyperref[Cat1842]{\textsc{Catalan 1842}} and \hyperref[Cat1872]{\textsc{1872}}, and in \hyperref[Can1908]{\textsc{Candido 1908}}), while a reprise of \eqref{E:Eul1744-2} in \hyperref[Eul1783]{\textsc{Euler 1783}} will be cited by Rudio in 1890 (see \hyperref[Rud1891]{\textsc{Rudio 1891}}). However, none of these later writers seem to have been aware of this early paper. A note downloaded from an old (and now defunct) version of the Euler Archive stated that this paper ``was presented to the St. Petersburg Academy on February 20, 1738."


\item\vspace{9pt} \label{Eul1763} (1763) Leonhard Euler\index[ppl]{Euler, Leonhard|textbf}, Annotationes in locum quendam Cartesii ad circuli quadraturam spectantem. \emph{Novi commentarii academiae scientiarum imperialis Petropolitanae, pro Annis MDCCLX. et MDCCLXI.} \textbf{8}, 24--27 (summary) and 157--168. Enestr\"{o}m number 275.
\index{Vi\`{e}te's formula for $\tfrac{2}{\pi}$}
\index{infinite products!of secant functions}

{\footnotesize Source: \href{https://scholarlycommons.pacific.edu/euler-works/275/}{The Euler Archive}.}

{\footnotesize For an English translation of the original Latin, see ``Annotations to a certain passage of Descartes for finding the quadrature of the circle", translated by Jordan Bell\index[ppl]{Bell, Jordan}, 23 May 2007 (\href{https://arxiv.org/abs/0705.3423v1}{arXiv.org}).}

{\footnotesize Cited in
\hyperref[Cat1842]{\textsc{Catalan 1842}} and \hyperref[Cat1872]{\textsc{1872}}\index[ppl]{Catalan, E.}, and
\hyperref[Can1908]{\textsc{Candido 1908}}\index[ppl]{Candido, G.}.}

\vspace{4pt}

This article contains the expansion
\begin{equation}\label{E:Eul1763-1}
\dfrac{\pi}{2}=\lim_{n\to\infty}\sec\dfrac{\pi}{2^2}\sec\dfrac{\pi}{2^4}\cdots\sec\dfrac{\pi}{2^n}\;,
\end{equation}
a special case of identity \eqref{E:Eul1744-2} above in \hyperref[Eul1744]{\textsc{Euler 1744}}. In 1842, Catalan\index[ppl]{Catalan, E.} will independently rediscover Vi\`{e}te's\index[ppl]{Vi\`{e}te, Fran\c{c}ois} continued square root product formula for $\tfrac{2}{\pi}$, and will cite this paper and equation \eqref{E:Eul1763-1} above in his proof. 

\hyperref[Cat1872]{\textsc{Catalan 1872}}\index[ppl]{Catalan, E.} and \hyperref[Can1908]{\textsc{Candido 1908}}\index[ppl]{Candido, G.} claim that equation \eqref{E:Eul1763-1} dates from 1760 and/or 1761. These citations allude to the years covered by the journal volume, not its publication date. A note at \href{https://scholarlycommons.pacific.edu/euler-works/275/}{The Euler Archive} states: ``According to Jacobi\index[ppl]{Jacobi, C. G. J}, a treatise with this title was read to the Berlin Academy on July 20, 1758. According to St. Petersburg Academy records, it was presented to that Academy on October 15, 1759."


\item\vspace{9pt} \label{Eul1783} (1783) Leonhard Euler\index[ppl]{Euler, Leonhard|textbf}, Variae observationes circa angulos in progressione geometrica progredientes. \emph{Opuscula analytica, Vol. I}, Typis academiae imperialis scientiarum, Petropoli [St.~Petersburg], 345--352. Enestr\"{o}m number 561.
\index{infinite products!of secant functions}

{\footnotesize Source: \href{https://scholarlycommons.pacific.edu/euler-works/561/}{The Euler Archive}.}  

{\footnotesize For an English translation of the original Latin, see ``Various observations on angles proceeding in geometric progression," translated by Jordan Bell\index[ppl]{Bell, Jordan}, 8 September 2010 (\href{https://arxiv.org/abs/1009.1439}{arXiv.org}).}

{\footnotesize Cited in
\hyperref[Rud1891]{\textsc{Rudio 1891}}\index[ppl]{Rudio, Ferdinand}.}

\vspace{4pt}

Alluding indirectly to his paper from 1744, Euler here uses the identity $\sin s=2\sin\tfrac{s}{2}\cos\tfrac{s}{2}$ to derive the expansion \eqref{E:Eul1744-2} above for the length of an arc on the unit circle (using $s$ rather than $A$ to represent arc length). \hyperref[Rud1891]{\textsc{Rudio 1891}} will take this expansion as a starting point for developing Vi\`{e}te's formula for $\tfrac{2}{\pi}$\index{Vi\`{e}te's formula for $\tfrac{2}{\pi}$}\index[ppl]{Vi\`{e}te, Fran\c{c}ois}. A note downloaded from an old (and now defunct) version of the Euler Archive stated that this paper ``was presented to the St. Petersburg Academy on November 15, 1773."


\item\vspace{9pt} \label{Ens1799} (1799) ``Citoyen Ensheim" [Moses Ensheim]\index[ppl]{Ensheim, Moses|textbf}, \emph{Recherches sur les calculs diff\'{e}rentiel et int\'{e}gral}. De l'imprimerie d'Agasse, rue de Poitevins, no. 13, an VII, Paris. 28 pages.
\index{successive substitution!18th century examples}

{\footnotesize Source: \href{https://gallica.bnf.fr/ark:/12148/bpt6k3181761/f19.item}{Gallica}.}

{\footnotesize Cited in
\hyperref[Bop1913]{\textsc{Bopp 1913}}\index[ppl]{Bopp, K.}.}

\vspace{4pt}

On pages 17 and 18 of this rare pamphlet, Ensheim approximates $\pi$ by the Archimedean method\index{Archimedean algorithm} of inscribed and circumscribed regular polygons\index{polygons, regular} in the unit circle\index{pi@$\pi\;(3.14159\ldots)$!continued square root expressions for}\index{constants, named!pi@$\pi\;(3.14159\ldots)$}. He uses successive substitution\index{successive substitution!18th century examples} to develop a nested expression for $\cos\tfrac{m}{2^k}$, and independently derives the identity \eqref{E:Eul1744-2} above from \hyperref[Eul1744]{\textsc{Euler 1744}}\index{infinite products!of secant functions}. Ensheim then applies his general formulas to the sequence of $3\cdot 2^k$-gons, $k=1,2,3,\ldots$, to generate the following bounds on $\pi$ (of which Euler's inequality \eqref{E:Eul1744-1} above is a special case):
\begin{align*}
3\cdot 2^q&(2-(2+(2+(2\cdots+(2+\sqrt{3})^{\frac{1}{2}}\cdots)^{\frac{1}{2}})^{\frac{1}{2}})^{\frac{1}{2}})^{\frac{1}{2}}\\
&<\pi<3\cdot 2^{q+1}\dfrac{(2-(2+(2+(2\cdots+(2+\sqrt{3})^{\frac{1}{2}}\cdots)^{\frac{1}{2}})^{\frac{1}{2}})^{\frac{1}{2}})^{\frac{1}{2}}}{(2+(2+(2+(2\cdots+(2+\sqrt{3})^{\frac{1}{2}}\cdots)^{\frac{1}{2}})^{\frac{1}{2}})^{\frac{1}{2}})^{\frac{1}{2}}}\;,
\end{align*}
where the left-hand expression and the numerator and denominator on the right each contain $q$ square roots (including $\sqrt{3}$). In a footnote to his discussion of the sequence 
\[(2+\sqrt{3})^{\frac{1}{2}},\; (2+(2+\sqrt{3})^{\frac{1}{2}})^{\frac{1}{2}},\ldots\;,\]
Ensheim expresses amazement that his fellow mathematicians ``hardly treat this species of infinite radical."


\item\vspace{9pt} \label{Gau1809} (1809) Carl Friedrich Gauss\index[ppl]{Gauss, Carl Friedrich|textbf}, \emph{Theoria motus corporum coelestium in sectionibus conicis solem ambientium}, Frid. Perthes \& I. H. Besser, Hamburg.
\index{Kepler problem}
\index{successive substitution!19th century examples}
\index{iteration!of functions}

{\footnotesize Source: \href{https://archive.org/details/bub_gb_ORUOAAAAQAAJ/page/n25/mode/2up}{Internet Archive}.}

{\footnotesize For an English translation of the original Latin, see \emph{Theory of the motion of the heavenly bodies moving about the sun in conic sections: a translation of Gauss's ``Theoria motus." With an appendix}, translated by Charles Henry Davis\index[ppl]{Davis, Charles Henry}, Little, Brown, and Co., Boston, 1857 (\href{https://archive.org/details/theoryofmotionof00gausrich/theoryofmotionof00gausrich/page/n33/mode/2up}{Internet Archive}).}

{\footnotesize Cited in 
\hyperref[Ise1888]{\textsc{Isenkrahe 1888}} and \hyperref[Ise1897]{\textsc{1897}}\index[ppl]{Isenkrahe, C.}.}

\vspace{4pt}

The ``Kepler problem" (see \hyperref[Kep1621]{\textsc{Kepler 1621}} and \hyperref[Swe2000]{\textsc{Swerdlow 2000}}) is to solve for $E$ in equation \eqref{E:Kep1621}, given the quantities $M$ and $e$. In article 11 of Section 1, Gauss proposes an iterative solution to the modified equation\footnote{Curtis Wilson (in ``Carl Friedrich Gauss, book on celestial mechanics (1809)", chapter 23 of \emph{Landmark Writings in Western Mathematics 1640--1940}, edited by I. Grattan-Guinness, Elsevier (2005), ISBN: 0-444-50871-6, \href{https://books.google.com/books?id=UdGBy8iLpocC\&lpg=PP1\&pg=PA324\#v=onepage\&q\&f=false}{Google Books}, p. 324) explains why Gauss used \eqref{E:Gau1809-0} instead of \eqref{E:Kep1621}: ``Following a policy adopted about 1800 by the French \emph{Bureau des Longitudes}, Gauss{\ldots}measures anomalies in elliptical orbits from perihelion, as is necessarily done in parabolic and hyperbolic orbits. `Kepler's equation' takes the form [\eqref{E:Gau1809-0}] with a minus sign on the right rather than the plus sign used by Kepler."}
\begin{equation}\label{E:Gau1809-0}
M=E-e\sin E\;;
\end{equation}
but not, however, the obvious regimen of successive substitution\index{successive substitution!19th century examples}, which would produce the expansion
\[E=M+e\sin(M+e\sin(M+\cdots+e\sin M)\cdots)\;.\]
Instead, he makes some simplifying assumptions which yield
\begin{equation}\label{E:Gau1809-1}
E=M+e\sin\epsilon\pm\dfrac{\lambda}{\mu\mp\lambda}(M+e\sin\epsilon-\epsilon)\;,
\end{equation}
where $\epsilon$ is an approximate value or initial guess for $E$. Gauss then posits (in Davis's translation): ``If the assumed value $\epsilon$ differs too much from the truth$\,\ldots\,$at least a much more suitable value will be found by this method, with which the same operation can be repeated, once, or several times if it should appear necessary. It is very apparent, that if the difference of the first value $\epsilon$ from the truth is regarded as a quantity of the first order, the error of the new value would be referred to the second order, and if the operation were further repeated, it would be reduced to the fourth order, the eighth order, etc. Moreover, the less the eccentricity [$e$], the more rapidly will the successive corrections converge."




\item\vspace{9pt} \label{Leg1816} \dag(1816) Adrien-Marie Legendre\index[ppl]{Legendre, Adrien-Marie|textbf}, \emph{Suppl\'{e}ment \`{a} l'essai sur la th\'{e}orie des nombres, 2. \'{E}d. [No. 1]-2.}

{\footnotesize Source: \href{https://gallica.bnf.fr/ark:/12148/bpt6k62826k/f574.item}{Gallica}.}

{\footnotesize Cited in 
\hyperref[Pio1837]{\textsc{Pioch 1837}}\index[ppl]{Pioch, A.} and \hyperref[Hey1904]{\textsc{Heymann 1904}}\index[ppl]{Heymann, W.}.}
\index{iteration!of functions}
\index{iteration!graphical methods}
\index{root-finding approximations}

The \emph{Supplement} (1816) to the 1808 edition of Legendre's \emph{Essay on the Theory of Numbers} is divided into three chapters, the last of which concerns two methods for approximating the roots of an equation of the form $F(x)=0$, where $F$ is a degree-$n$ polynomial. The first method is an iterative process involving an auxiliary function 
\[\phi(x)=\frac{ax^{n-k}+bx^{n-k-1}+cx^{n-k-2}+\text{etc.}}{1+\frac{f}{x}+\frac{g}{x^2}+\frac{h}{x^3}++\text{etc.}}\;,\]
which is demonstrated to be monotonic and greater than $F$. The iteration generates approximations 
\[\alpha\;,\quad\alpha'=\sqrt[n]{\phi(\alpha)}\;,\quad\alpha''=\sqrt[n]{\phi(\alpha')}\;,\ldots\]
\hyperref[Gal2001]{\textsc{Galuzzi 2001}}\index[ppl]{Galuzzi, Massimo} gives a concise explication of this essentially geometrical algorithm and of Galois' subsequent improvements to it (\hyperref[Gal1830]{\textsc{Galois 1830}}).\index[ppl]{Galois, \'{E}variste}


\item\vspace{9pt} \label{Sch1821} (1821) Henri Gerner Schmidten\index[ppl]{Schmidten, Henri Gerner|textbf}, M\'{e}moire sur l'int\'{e}gration des \'{e}quations lin\'{e}aires. \emph{Annales de Math\'{e}matiques Pures et Appliqu\'{e}es} \textbf{11}, 269--316.  

{\footnotesize Source: \href{http://www.numdam.org/item/AMPA_1820-1821__11__269_0.pdf}{NUMDAM}.}

{\footnotesize Cited in 
\hyperref[MRS1830]{\textsc{``M.R.S." 1830}}\index[ppl]{M. R. S.},
\hyperref[Pio1837]{\textsc{Pioch 1837}}\index[ppl]{Pioch, A.}, and 
\hyperref[Tou1996]{\textsc{Tourn\`{e}s 1996}}\index[ppl]{Tourn\`{e}s, Dominique}.}
\index{successive substitution!19th century examples}

\vspace{4pt}

The author's first name is given in other sources as Henrik, and his surname as von Schmidten and de Schmidten. He begins his paper by assuming that a solution $y$ exists for the implicit differential equation\index{differential equations!implicit}
\begin{equation}\label{E:Sch1821-1}
\varphi. y=f. y\;,
\end{equation}
where $\varphi. y$ is ``a function which contains the differential coefficients or the highest order differences in the proposed equation," and $f.y$ is ``any other function of the independent variables of the differential or differential coefficients." If $X$ is a solution of the homogenous equation $\varphi.X=0$, equation \eqref{E:Sch1821-1} can be written as $\varphi.y=\varphi.X+f.y$; from this, Schmidten isolates $y$ on the left by applying the inverse $\tfrac{1}{\varphi}$; the result is 
\[y=X+\dfrac{1}{\varphi}f.y\;.\]
He continues, ``By means of this implicit relation we shall easily find the explicit value of $y$ by successive substitutions; this will be\index{successive substitution!19th century examples}
\[y=X+\tfrac{1}{\varphi}f\bigl(X+\tfrac{1}{\varphi}f\bigl(X+\tfrac{1}{\varphi}f\bigl(X+\ldots\;.\]

Some difficulties are immediately apprehended: ``It will be seen, however, that the value of $y$ will in general remain very complicated, unless $\varphi.y$ and $f.y$ are linear with respect to $y$, which embraces a very extended and very important class of equations: that of \emph{linear equations}$\ldots$ In this case we have
\[
y = X + \tfrac{1}{\varphi}fX+\tfrac{1}{\varphi}f\bigl(\tfrac{1}{\varphi}fX\bigr)+\tfrac{1}{\varphi}f\bigl(\tfrac{1}{\varphi}f\bigl(\tfrac{1}{\varphi}fX\bigr)\bigr)+\ldots
\]
and I propose to state its principal consequences$\ldots$"

\hyperref[Tou1996]{\textsc{Tourn\`{e}s 1996}}\index[ppl]{Tourn\`{e}s, Dominique} observes, ``[T]he beginning of the text contains an abstract exposition of the method of successive substitutions, of great formal beauty and surprisingly modern."


\item\vspace{9pt} \label{Sch1826} (1826) ``Herrn Prof. Schmidten" [Henri Gerner Schmidten],\index[ppl]{Schmidten, Henri Gerner|textbf} Versuch \"{u}ber die Integration der Differential-Gleichungen. \emph{Journal f\"{u}r die reine und angewandte Mathematik} \textbf{1}, 137--151.
\index{successive substitution!19th century examples} 

{\footnotesize Source: \href{https://gdz.sub.uni-goettingen.de/id/PPN243919689_0001?tify=\%7B\%22pages\%22\%3A\%5B143\%5D\%2C\%22pan\%22\%3A\%7B\%22x\%22\%3A0.527\%2C\%22y\%22\%3A0.667\%7D\%2C\%22view\%22\%3A\%22info\%22\%2C\%22zoom\%22\%3A0.484\%7D}{G\"{o}ttinger Digitalisierungszentrum}.}

\vspace{4pt}

In this ``[e]xcerpt from an essay read by the author, in the Royal Danish Academy of Sciences\index{Royal Danish Academy of Sciences}, Copenhagen, in the Danish language," Schmidten furthers his applications of successive substitution. Beginning with a differential equation\index{differential equations} of the form
\[F\left(x,y,y',y'',\cdots,y^{(n)}\right) = 0\;,\]
he shows how, under helpful conditions, one may iteratively rearrange and integrate to produce a reduced equation
\[f_n(x,y) = c_n + \psi y\;,\]
where $c_n$ is a constant of integration. Schmidten then assumes that $f_n$ can be inverted (with inverse $P$), by which one obtains
\[y = P(x,\; c_n + \psi y)\;.\]
This, in turn, is the setup he needs to generate a solution by successive substitution\index{successive substitution!19th century examples}:
\[y = P(x, c_n + \psi(P(x,c_n+\psi(\cdots \psi y\cdots))))\;.\]
Some formidable expressions arise from his algorithm. One of his solutions is
\[y=\dfrac{\alpha}{\int dx\left(\gamma+\dfrac{\delta\alpha}{\int dx\left(\gamma+\dfrac{\delta\alpha}{\int dx\left(\cdots\right)}\right)}\right)}\;.\]


\item\vspace{9pt} \label{MRS1830} (1830) {``M. R. S."},\index[ppl]{M. R. S.|textbf} Note sur quelques expressions alg\'{e}briques peu connues. \emph{Annales de Math\'{e}matiques Pures et Appliqu\'{e}es} \textbf{20}, 352--366.  
\index{successive substitution!19th century examples}

{\footnotesize Source: \href{http://www.numdam.org/item/AMPA_1829-1830__20__352_1.pdf}{NUMDAM}.}

{\footnotesize Cited in
\hyperref[Pio1837]{\textsc{Pioch 1837}}\index[ppl]{Pioch, A.}.}

\vspace{4pt}

The author of this paper, hiding behind an acronym, further develops the key idea in \hyperref[Sch1821]{\textsc{Schmidten 1821}} about successive substitution. After several specific examples, he proposes that the general equation $F(x)=0$ can be made to take the form $x=A+f(x)$, and this in turn can be used with the method of successive substitution to obtain
\begin{equation*}\label{E:MRS1830-1}
x=A+f(A+f(A+f(A+\ldots)))\;.
\end{equation*}
He goes on to use $x=Af(x)$ for
\[x=Af(Af(Af(A\ldots)))\;,\]
as well as $x=A^{f(x)}$ to generate
\[x=A^{f(A^{f(A^{f(A^{\ldots})})})}\;,\]
and $x=A+f(B+\varphi(x))$ to get
\[x=A+f(B+\varphi(A+f(B+\varphi(A+\ldots))))\;.\]
However, Monsieur R. S. then attempts to reverse the process, using the divergent sum\index{infinite sums!divergent}
\[x=1-2+4-8+16-32+64-\ldots\]
to deduce the ``equation" $x=1-2x$, from which he concludes that $x=\tfrac{1}{3}$. His argument that 
\[a-a+a-a+a-a+\ldots\]
sums to $\tfrac{1}{2} a$ was one of many in a debate, reaching back to the early 18th century, about how to deal with divergent infinite sums\index{infinite sums!divergent}; see for instance the \href{https://en.wikipedia.org/wiki/History_of_Grandi's_series}{Wikipedia} entry on Grandi's series\index{Grandi's series}\index{series!Grandi's}.\footnote{It seems possible that this paper is the work of Henri Gerner Schmidten\index[ppl]{Schmidten, Henri Gerner}. As shown by two of the three papers credited to him here, Schmidten appears diffident about revealing his full name. He also seems to be a prominent proponent of successive substitution\index{successive substitution!19th century examples} during this time period. The work of {``M.~R.~S."} fits nicely into the sequence of papers by Schmidten in the decade preceding and including 1830.} 


\item\vspace{9pt} \label{Sch1830} (1830) ``Mr. de Schmidten, prof. des math\'{e}m. \`{a} Copenhague" [Henri Gerner Schmidten]\index[ppl]{Schmidten, Henri Gerner|textbf}, Sur un principe g\'{e}n\'{e}ral dans la th\'{e}orie des s\'{e}ries. \emph{Journal f\"{u}r die reine und angewandte Mathematik} \textbf{5}, 388--396.  
\index{successive substitution!19th century examples}

{\footnotesize Source: \href{https://doi.org/10.1515/crll.1830.5.388}{doi.org (De Gruyter)}.}

\vspace{4pt}

Schmidten broadens his conception of infinite nested function compositions by defining
\[F(x)=\varphi[A_1+f_1(x)\varphi(A_2+f_2(x)\varphi(A_3+\cdots))]\;,\]
where $\varphi$ is an invertible function with inverse $\psi(x)$, the $f_i$ are functions of $x$, and $A_1, A_2, A_3\ldots$ are quantities which depend on $F(x)$. He furthermore assumes that $x_1, x_2, x_3,\ldots$ are values of $x$ for which $f_i(x_i)=0$ for $i=1,2,3,\ldots$, and that these values exist, although he admits that $x_i$ may not \emph{uniquely} solve $f_i(x)=0$. Under these assumptions, he is able to unwrap the finite approximants\index{approximants!of continued compositions} of $F$ and solve for the $A_i$ in sequential order:
\begin{align*}
A_1&=\psi(F(x_1))\;,\quad
A_2=\psi\!\left(\dfrac{\psi(F(x_2))-A_1}{f_1(x_2)}\right),\\
A_3&=\psi\!\left(\raisebox{-6pt}{\textrm{$\dfrac{\psi\!\left(\dfrac{\psi(F(x_3))-A_1}{f_1(x_3)}\right)}{f_2(x_3)}$}}\right),
\end{align*}
and so on. He gives two examples in which these calculations are simplified: first, setting $\varphi(x)=\tfrac{1}{x}$, and second, using $\varphi(x)=x$, which turns $F$ into an infinite sum:
\[F(x) = A_1 + A_2f_1(x) + A_3f_1(x)f_2(x)+\cdots\;.\]
Furthermore, if for each $i=1,2,3\ldots$ one has $f_i(x)=x$, then the Taylor series\index{Taylor series}\index{series!Taylor} for $F$ is generated. Schmidten also considers functions defined implicitly, for instance
\[y=\varphi(x, y)\;,\]
which expands using successive substitution\index{successive substitution!19th century examples} into
\[y = \varphi [x, \varphi (x, \varphi (x,\ldots))]\;.\]
(Such an expansion underlies the continued cotangent\index{continued cotangents} developed in \hyperref[Leh1938]{\textsc{Lehmer 1938}}\index[ppl]{Lehmer, D. H.}.) On page 391, there is a tantalizing allusion to the ``sum" $\tfrac{1}{2}=1-1+1-1+\cdots$\index{series!Grandi's}\index{Grandi's series}\index{infinite sums!divergent} mentioned in \hyperref[MRS1830]{\textsc{``M.R.S" 1830}}\index[ppl]{M. R. S.} --- perhaps another clue in support of the idea that Schmidten was ``M.R.S." 


\item\vspace{9pt} \label{Gal1830} \dag(1830) \'{E}variste Galois\index[ppl]{Galois, \'{E}variste|textbf}, Note sur la r\'{e}solution des \'{e}quations numeriques. \emph{Bulletin des sciences math\'{e}matiques (Bulletin de F\'{e}russac)}, \textbf{13}, 413--414. In \emph{{\OE}uvres math\'{e}matiques d'\'{E}variste Galois}, Gauthier-Villars et Fils, Paris, 1897, 13--14.
\index{iteration!of functions}
\index{root-finding approximations}

{\footnotesize Source: \href{https://books.google.com/books?id=0jGzsXc2TbgC&newbks=1&newbks_redir=0&dq=Note\%20sur\%20la\%20r\%C3\%A9solution\%20des\%20\%C3\%A9quation\%20num\%C3\%A9riques\%20Galois&pg=PA13#v=onepage&q=Note\%20sur\%20la\%20r\%C3\%A9solution\%20des\%20\%C3\%A9quation\%20num\%C3\%A9riques\%20Galois&f=false}{Google Books}.}

\vspace{4pt}

\hyperref[Gal2001]{\textsc{Galuzzi 2001}}\index[ppl]{Galuzzi, Massimo} says, ``Galois rewords the geometrical argument of [\hyperref[Leg1816]{\textsc{Legendre 1816}}\index[ppl]{Legendre, Adrien-Marie}]{\ldots}it is clear that he thinks in terms of the sequence
\[a\;,\;\phi(a)\;,\;\phi(\phi(a))\;,\;\phi(\phi(\phi(a)))\;,\;\ldots\]
whose convergence is granted by the conditions posed.'' Galois' motivation is that “[Legendre's] method has the drawback of requiring the calculation of the $n$th root at every step.''


\item\vspace{9pt} \label{Fou1831} \dag(1831) Joseph Fourier\index[ppl]{Fourier, Joseph|textbf}, \emph{Analyse des \'{e}quations d\'{e}termin\'{e}es}, Vol. 1. Edited by Claude-Louis Navier\index[ppl]{Navier, Claude-Louis}. Chex Firmin Didot Fr\`{e}res, Paris.
\index{successive substitution!19th century examples}
\index{root-finding approximations}

{\footnotesize Source: \href{https://gallica.bnf.fr/ark:/12148/bpt6k1057816b/f221.item}{Gallica}.}

{\footnotesize Cited in 
\hyperref[Pio1837]{\textsc{Pioch 1837}}\index[ppl]{Pioch, A.}.}

In the ``publisher's note'', editor Navier describes this posthumously-printed book's subjects as ``the resolution of symbolic equations, the separation of the roots of numerical equations, and the distinction between the two cases where there exists between two very close limits either a pair of imaginary roots, or two almost equal real roots.'' In Book 2,\\``Method for calculating the values of roots whose limits are known, and various remarks on the convergence of approximations and on the distinction of roots'', Sections 20 and 21, Fourier introduces and demonstrates an algorithm for performing division of two (presumably large) integers, which he calls ``ordered division''\index{ordered division}  (\emph{division ordonn\'{e}e}), and which is now known as Fourier division or cross division.\footnote{Ronald W. Doerfler\index[ppl]{Doerfler, Ronald W.}, \emph{Dead Reckoning: Calculating Without Instruments}, Houston, Gulf Publishing Company, 1993, pp. 28--33.} Uspensky\footnote{J. V. Uspensky\index[ppl]{Upensky, J. V.}, \emph{Theory of Equations}, New York, McGraw-Hill Book Co., 1948, pp. 159--164.} describes it as a method in which ``a group or small number of first digits of the divisor consisting, say, of two or three digits is selected as an `abridged' divisor, and all the divisions are made with this abridged divisor, the remaining digits of the divisor being taken into consideration gradually to make the so-called `corrections.'" In Section 22, Fourier shows how ordered division can be adapted to solve a quadratic equation, essentially by turning the quadratic into a successive substitution scheme disguised as a division problem. He rearranges the example $x^2+765432x=123456$ into the form
\[x=\frac{123456}{765432+x}\;,\]
and shows how a modified ordered division process produces increasingly accurate approximations. In Section 23 he details the steps for finding a root of the cubic $x^3+345x=12$ using ordered division on
\[x=\frac{12}{345+x^2}\;.\]
Underlying these processes is Fourier's variation on Newton's method\index{Newton's or Newton-Raphson method}, which he develops in Book 1 and the first sections of Book 2.


\item\vspace{9pt} \label{Bol1832} (1832) Farkas Bolyai\index[ppl]{Bolyai, Farkas (\emph{aka} Wolfgang)|textbf}, \emph{Tentamen juventutem studiosam in elementa matheseos purae, elementaris ac sublimioris, methodo intuitiva, evidentiaque huic propria, introducendi: cum appendice triplici}, Vol. 1. Maros Vasarhelyini: Typis Collegii Reformatorum per Josephum et Simeonem Kali de fels\H{o} Vist.  
\index{successive substitution!19th century examples}
\index{continued robb roots@continued $r$th roots!of constant nonnegative real terms}
\index{trinomial equations}
\index{iterated square roots!of constant nonnegative real terms}

{\footnotesize Source: \href{https://archive.org/details/tentamenjuventu1boly/page/412/mode/2up}{Internet Archive}.}

{\footnotesize Cited in 
\hyperref[Far1881]{\textsc{Farkas 1881}}\index[ppl]{Farkas, Gyul\'{a}t\'{o}l (\emph{aka} Gyula or Jules)} and 
\hyperref[Sza2010]{\textsc{Szab\'{o} 2010}}\index[ppl]{Szab\'{o}, P\'{e}ter G\'{a}bor}.}

\vspace{4pt}

Known as ``Wolfgang" in Germany, Bolyai was the father of J\'{a}nos\index[ppl]{Bolyai, J\'{a}nos}, whose groundbreaking work on non-euclidean geometry appears as an appendix in this volume (which is written in Latin). On page 412, following a discussion of continued fraction solutions to some simple trinomial equations, the senior Bolyai proposes to solve $x^2-x=a$ by successive substitution, using $x=\sqrt{a+x}$ to obtain
\[x=\sqrt{a+\sqrt{a+\sqrt{a+\cdots}}}\;.\]
He lets $a=1$ to get the golden ratio\index{golden ratio ($\tfrac{1+\sqrt{5}}{2}=1.61803\ldots$)}\index{constants, named!golden ratio ($\tfrac{1+\sqrt{5}}{2}=1.61803\ldots$)} $x=\tfrac{1+\sqrt{5}}{2}$; he exhibits $x=2$ as the limit when $a=2$, and also as the limit of a pair of periodic continued fractions. He then proves that the iterated $m$th root
\[x=\sqrt[m]{a+\sqrt[m]{a+\sqrt[m]{a+\cdots}}}\;.\]
converges (assuming $a>0$) by showing that the finite expansions are increasing and bounded above.


\item\vspace{9pt} \label{Dop1832a}  (1832a) Christian Doppler\index[ppl]{Doppler, Christian|textbf}, \"{U}ber die Konvergenz einer unendlichen Logarithmenfolge. \emph{Jahrb\"{u}cher des kaiserlichen k\"{o}niglichen polytechnischen Institutes in Wien} \textbf{17}, 172--174.  
\index{continued logarithms}
\index{successive substitution!19th century examples}

{\footnotesize Source: \href{https://books.google.com/books?id=lC6iVzC52kEC&newbks=1&newbks_redir=0&dq=\%22\%C3\%9Cber\%20die\%20Konvergenz\%20einer\%20unendlichen\%20Logarithmenfolge\%22\&pg=PA172#v=onepage\&q\&f=false}{Google Books}.}

\vspace{4pt}

The author is well-known for the Doppler effect in physics; this paper and the one immediately below predate his wave propagation work by about ten years.

In this note, comprising two pages of text and one (printed sideways) of calculations, Doppler asserts that the \emph{logarithmische Gr\"{o}{\ss}enfolge}
\[s=\log(a+\log(a+\log(a+\log(a+\cdots))))\]
converges for $a>1$, by an increasing-and-bounded-sequence argument. (His logarithms are to base 10.) He also claims that ``if $a<1$, then [the expression $s$] also converges as long as it does not become imaginary." For $a=5$ he computes $s=5.7604576$ (versus the machine-generated $5.7604567835$).


\item\vspace{9pt} \label{Dop1832b} (1832b) Christian Doppler\index[ppl]{Doppler, Christian|textbf}, \"{U}ber Kettenwurzeln und deren Konvergenz. \emph{Jahrb\"{u}cher des kaiserlichen k\"{o}niglichen polytechnischen Institutes in Wien} \textbf{17}, 175--200.
\index{successive substitution!19th century examples}
\index{continued square roots!of periodic real terms}
\index{continued pocc powers@continued $p_i$th powers}
\index{Kettenwurzeln@\emph{Kettenwurzeln}} 

{\footnotesize Source: \href{https://books.google.com/books?id=lC6iVzC52kEC&newbks=1&newbks_redir=0&pg=PA175#v=onepage&q&f=false}{Google Books}.}

{\footnotesize Cited in 
\hyperref[Jon2015]{\textsc{Jones 2015}}\index[ppl]{Jones, Dixon J.}.}

\vspace{4pt}

The second of Doppler's\index[ppl]{Doppler, Christian} back-to-back papers in this journal issue contains perhaps the earliest use of the German \emph{Kettenwurzel}\index{Kettenwurzeln@\emph{Kettenwurzeln}}, which he defines more generally than do most subsequent authors as 

\begin{equation}\label{E:Dop1832-1}
A\sqrt[a]{\phantom{|}}[\alpha+B\sqrt[b]{\phantom{|}}(\beta+C\sqrt[c]{\phantom{|}}(\gamma+D\sqrt[d]{\phantom{|}}(\delta+E\sqrt[e]{\phantom{|}}(\epsilon+\ldots))))]\;.
\end{equation}
However, he goes on to say that ``the exponent values $a, b, c, d, e,$ etc., as well as the values and signs of the multipliers $A, B, C, D, E,$ etc. are assumed to be periodic in their infinite progression." This assumption is crucial for applying the method of successive substitution\index{successive substitution!19th century examples} to obtain convergence results. Doppler\index[ppl]{Doppler, Christian} recognizes that for special values of the exponents $a, b, c,\ldots$ (which are not assumed to be positive integers) one obtains the continued fraction\index{continued fractions}

\begin{equation}\label{E:Dop1832-3}
\cfrac[l]{A}{\alpha+\cfrac[l]{B}{\beta+\cfrac[l]{C}{\gamma+\cfrac[l]{D}{\delta+\cfrac[l]{E}{\epsilon+\ldots}}}}}
\end{equation}
\vspace{6pt}

\noindent as well as what are called continued reciprocal powers:\index{continued reciprocal powers}

\begin{equation}\label{E:Dop1832-4}
\cfrac[l]{A}{{\phantom{\biggl |}}^a\biggl(\alpha+\cfrac[l]{B}{{\phantom{\biggl |}}^b\biggl(\beta+\cfrac[l]{C}{{\phantom{\biggl |}}^c\biggl(\gamma+\cfrac{D}{{\phantom{\biggl |}}^d\biggl(\delta+\ldots}}}}
\end{equation}
\medskip and continued reciprocal roots\index{continued reciprocal roots}:

\begin{equation}\label{E:Dop1832-5}
\cfrac[l]{A}{\sqrt[a]{\alpha+\cfrac[l]{B}{\sqrt[b]{\beta+\cfrac[l]{C}{\sqrt[c]{\gamma+\cfrac[l]{D}{\sqrt[d]{\delta+\ldots}}}}}}}}\;.
\end{equation}

\bigskip\noindent Before 1900 we find the form \eqref{E:Dop1832-5} in \hyperref[Pio1837]{\textsc{Pioch 1837}}\index[ppl]{Pioch, A.}, \hyperref[Dix1878]{\textsc{Dixon 1878}}\index[ppl]{Dixon, T. S. E.} and  \hyperref[Gun1880]{\textsc{G\"{u}nther 1880}}\index[ppl]{Gunther@G\"{u}nther, Siegmund}. A continued reciprocal square\index{continued reciprocal squares} appears in \hyperref[Pio1837]{\textsc{Pioch 1837}}\index[ppl]{Pioch, A.}, and equation \eqref{E:Dix1878-1} (see \hyperref[Dix1878]{\textsc{Dixon 1878}}) can be interpreted as a continued reciprocal $\tfrac{m}{n}$ power if $n<m$, but the general continued reciprocal power \eqref{E:Dop1832-4} does not appear again in the papers listed here until the 1990s (\hyperref[Lau1990]{\textsc{Laugwitz 1990}}, \hyperref[Sch1992]{\textsc{Sch\"{o}nefuss 1992}}).\index[ppl]{Laugwitz, Detlef}\index[ppl]{Schonefuss@Sch\"{o}nefuss, Lutz W.}


\item\vspace{9pt} \label{Eic1834} (1834) C. F. Eichhorn\index[ppl]{Eichhorn, Christian Friedrich|textbf}, \emph{Principien einer allgemeinen Functionenrechnung}, Helwingschen Hof-Buchhandlung, Hannover.
\index{successive substitution!19th century examples}

{\footnotesize Source: \href{https://gdz.sub.uni-goettingen.de/id/PPN587693088}{G\"{o}ttinger Digitalisierungszentrum}.}

\vspace{4pt}

This unusual book begins by restricting its ``general function calculations'' to linear functions $f$, for which $f(ca+cb)=cf(a)+cf(b)$. The author employs modern notation like $f^n$ for the composition of $f$ with itself $n$ times, including $f^{-1}$ for the inverse of $f$; on the other hand, he also turns $y+a\theta(y)=\psi(x)$ into the dubious
\[y(1+a\theta(\;\;))=\psi(x)\;,\quad\textrm{and consequently}\quad y=\dfrac{\psi(x)}{(1+a\theta(\;\;))}\]
(p. 81). Pages 81--98 and 138--139 develop some infinitely nested expressions based on successive substitution. A review of this book in \emph{Repertorium der gesammten deutschen Literatur} \textbf{2}, 1834, 536--537, complains, ``[T]his work will be intelligible and useful only to those analysts who have studied Euler\index[ppl]{Euler, Leonhard}, Legendre\index[ppl]{Legendre, Adrien-Marie}, Lacroix\index[ppl]{Lacroix, S. F.}, Burg, and others, written with greater clarity and elegance$\ldots$the errors of printing and calculation take the reader's attention much more than is already the case, due to [the general difficulty of the subject]." 


\item\vspace{9pt} \label{Pio1837} \dag(1837) A. Pioch\index[ppl]{Pioch, A.|textbf}, \emph{M\'{e}moire sur la r\'{e}solution des \'{e}quations, suivi de notes sur l'\'{e}valuation des fonctions sym\'{e}triques et sur la d\'{e}termination des tangentes et des plans tangents}, Leroux, Brussels.
\index{successive substitution!19th century examples}
\index{continued reciprocal roots}
\index{continued reciprocal squares}

{\footnotesize Cited in 
\hyperref[Rea1877]{\textsc{Realis 1877}}\index[ppl]{Realis, S.}.}

\vspace{4pt}

A footnote in \hyperref[Rea1877]{\textsc{Realis 1877}}\index[ppl]{Realis, S.}, credited to ``E. C.'' (likely E. Catalan\index[ppl]{Catalan, E.}, the journal editor), states, ``An instructor in the military academy of Brussels, named Pioch, who died very young, had considered equations of the form $x = f(x)$. His 1837 \emph{Dissertation on the resolution of equations} contains some remarkable ideas." Pioch's principle goal in this book is to prove that ``[t]he root of an equation of any degree can always be expressed by a finite number of non-symmetric rational functions of all the roots.'' He goes on to say, ``Abel\index[ppl]{Abel, Niels Henrik} stated it differently{\ldots}he says it is only true when the root of an equation can be expressed algebraically. I will show that the reasoning which led him to this theorem is inaccurate, and that the theorem itself, as he states it, implies a contradiction, and that the author [Abel] has consequently turned it into a vicious circle.'' Such statements have deterred us from devoting much time to the first three of this work's four chapters; furthermore, Pioch is dismissive of Galois\index[ppl]{Galois, \'{E}variste}, whose work was poorly understood at the time. 

Pioch's Chapter 4 is a fairly detailed, rigorous, and general discussion of successive substitution. He credits \hyperref[Sch1821]{\textsc{Schmidten 1821}} and \hyperref[MRS1830]{\textsc{``M. R. S.'' 1830}}\index[ppl]{M. R. S.} for inspiring his interest in the topic, and points out that \hyperref[Fou1831]{\textsc{Fourier 1831}} makes use of this method; he also claims that material in the appendix to Legendre's \emph{Theory of Numbers} comprises a special case of his (Pioch's) theory. He is critical of ``M. R. S.'', who ``has not shown, in the particular cases he considered, whether the successive values deduced from [successive substitution] formulas converge towards the effective value of the root sought.'' Pioch admits, however, that \hyperref[MRS1830]{\textsc{``M. R. S.'' 1830}} ``contains fruitful ideas, [suggesting] research that is perhaps far from being exhausted, and which deserves the attention of analysts.''

Pioch considers iteration of the equation $x_n=f(x_{n-1})$ using the three principle examples in \hyperref[MRS1830]{\textsc{``M. R. S.'' 1830}}: $f(x)=A+\phi(x)$, $f(x)=A\cdot\phi(x)$, and $f(x)=A^{\phi(x)}$. He is careful to distinguish between monotonic and alternating behavior of the iterates in convergence and divergence. He considers various cases in which $\phi(x)$ is a rational function, including $\phi(x)=P_m(x)/Q_n(x)$, where $P$ and $Q$ are polynomials of degrees $m$ and $n$ respectively, and $\phi(x)=A_0+A_1/x$, which produces a continued fraction. At the end of the first section in Chapter 4, he concludes that the last two of his three versions of $f(x)$ can be considered special cases of the first (transforming the last case into the first using logarithms). 

In the chapter's second section, Pioch calculates approximate roots of various sample equations. Here it is notable that he recognizes that a given polynomial equation may be rearranged in more than one way to create an iteration scheme. For instance, from $x^2-ax-b=0$ he derives $x = a+\frac{b}{x}$ and $x=-\frac{b}{a-x}$ (compare \hyperref[Hof1881]{\textsc{Hoffmann 1881}}\index[ppl]{Hoffmann, K. E.}, who resolved this kind of issue arising in \hyperref[Gun1880]{\textsc{G\"{u}nther 1880}}\index[ppl]{Gunther@G\"{u}nther, Siegmund}). Rearranging the cubic equations
\[x^3-a x^2-c=0\quad\text{and}\quad x^3-bx-c=0\]
in the respective forms
\[x= a+\frac{c}{x^2}\quad\text{and}\quad x=\sqrt{b+\frac{c}{x}}\;,\]
Pioch derives the continued reciprocal square\index{continued reciprocal squares}
\begin{align*}
x &= a+
\cfrac{c}{\left(\raisebox{3.6ex}{$a+
\cfrac{c}{\left(\raisebox{1.6ex}{$a+
\cfrac{c}{\left(a+\ldots\right)\!\raisebox{1.6ex}{$^2$}}$}
  \right)^2}$}
  \right)^2}
\intertext{and the continued reciprocal square root\index{continued reciprocal square root}}
x&=\sqrt{b+\cfrac{c}{\sqrt{b+\cfrac{c}{\sqrt{b+\cfrac{c}{\sqrt{b+\cdots}}}}}}}
\end{align*}
(but does not assign names to these forms).

In the third and last of the chapter's sections, titled ``Application of the preceding theories to the resolution of transcendental equations,'' he asserts that the method of successive substitution applies when $\phi(x)$ is a trigonometric or other transcendental function. He specifically mentions ``Kepler's equation'' $x = p + q \sin x$ \index{Kepler problem}\index[ppl]{Kepler, Johannes}, and offers a numerical example of an iterated solution when $p=4$ and $q=\tfrac{5}{4}$. The levels of detail and rigor are lower in this section than in the previous two.

A curious tic of the author's is his consistent spelling of Fourier\index[ppl]{Fourier, Joseph} as ``Fourrier."


\item\vspace{9pt} \label{Cat1842} (1842) E. Catalan\index[ppl]{Catalan, E.|textbf}, Note sur le rapport de la circonf\'{e}rence au diam\`{e}tre. \emph{Nouvelles Annales de Math\'{e}matiques} 1re s\'{e}rie. \textbf{1}, 190--196.
\index{Vi\`{e}te's formula for $\tfrac{2}{\pi}$}
\index[ppl]{Vi\`{e}te, Fran\c{c}ois}

{\footnotesize Source: \href{http://www.numdam.org/item/NAM_1842_1_1__190_1.pdf}{NUMDAM}.} 

{\footnotesize Cited in
\hyperref[Lef1897]{\textsc{Lefebvre 1897}}\index[ppl]{Lefebvre, B.},  
\hyperref[Wie1904b]{\textsc{Wiernsberger 1904b}} and \hyperref[Wie1905]{\textsc{1905}}\index[ppl]{Wiernsberger, Paul}, and
\hyperref[Bop1913]{\textsc{Bopp 1913}}\index[ppl]{Bopp, K.}.}

\vspace{4pt}

Catalan\index[ppl]{Catalan, E.} rediscovers Vi\`{e}te's\index[ppl]{Vi\`{e}te, Fran\c{c}ois} formula (although Vi\`{e}te's name is not mentioned), expressed as the reciprocal of equation \eqref{E:Vie1593-2} above. He also derives the formula\index{pi@$\pi\;(3.14159\ldots)$!continued square root expressions for}\index{constants, named!pi@$\pi\;(3.14159\ldots)$}
\begin{equation}\label{E:Cat1842-1}
\pi=\lim_{n\to\infty}2^n\sqrt{2-\sqrt{2+\sqrt{2+\sqrt{2+\cdots+\sqrt{2}}}}}\;,
\end{equation}
where there are $n$ 2s under the radical. This expansion will repeatedly be rediscovered, for instance in \hyperref[Fan1850]{\textsc{Fanien 1850}}\index[ppl]{Fanien, A.}, \hyperref[Did1872]{\textsc{Didion 1872}}\index[ppl]{Didion, Isidore}, \hyperref[Pie1891]{\textsc{Pierce 1891}}\index[ppl]{Pierce, George Winslow}, \hyperref[Can1908]{\textsc{Candido 1908}}\index[ppl]{Candido, G.}, \hyperref[Hau2003]{\textsc{Hauser 2003}}\index[ppl]{Hauser, Clemens}, and \hyperref[Cha2016]{\textsc{Chang and Chang 2016}}\index[ppl]{Chang, Mu-Ling}\index[ppl]{Chang, Chia-Chin (Cristi)}. 


\item\vspace{9pt} \label{Hei1844} \dag(1844) Eduard Heis\index[ppl]{Heis, Eduard|textbf}, \emph{Sammlung von Beispielen und Aufgaben aus der allgemeinen Arithmetik und Algebra. Für Gymnasien, höhere Bürgerschulen und Gewerbschulen in systematischer Folge bearbeitet}, first edition. M. Du Mont-Schauberg, K\"{o}ln.  

{\footnotesize Source: \href{https://books.google.com/books/download/Sammlung_von_Beispielen_und_Aufgaben_aus.pdf?id=Bxd1pAulrAIC&hl=en&capid=AFLRE73TSxoEi0wRS4yFEjQ_578TPkpbTjATP0w0hfgpFjftQwWUchvN_9XZjSbPMcCQljwddcIM7ly3jA7BnhgDJpF2mO_GJQ&continue=https://books.google.com/books/download/Sammlung_von_Beispielen_und_Aufgaben_aus.pdf\%3Fid\%3DBxd1pAulrAIC\%26output\%3Dpdf\%26hl\%3Den}{Google Books}.}
\index{continued cube roots!as solutions to cubic equations}
\index{continued reciprocal roots}



{\footnotesize Cited in
\hyperref[Hey1898]{\textsc{Heymann 1898}.}\index[ppl]{Heymann, W.}}

\vspace{4pt}

The first edition of Heis's \emph{Collection of examples and exercises from general arithmetic and algebra. Presented in a systematic sequence for grammar schools, higher secondary schools, and vocational schools}, published in 1837, comprised six chapters ranging from elementary arithmetic to geometric and arithmetic progressions, combinatorics, and probability. In the third edition, Chapters 7 and 8 were added; in Chapter 7, ``Equations of Higher Degrees and Transcendental Equations'', Section 105, ``Applications of higher-degree equations'', we find the following two problems and their answers:

\vspace{4pt}

\begin{quote}
\textbf{24}) Which equation needs to be solved to determine the value of the infinite continued fraction [\emph{unendlichen Kettenbruches}] $a+\surdex{\dfrac{b}{a}}\raisebox{-7pt}{$\;+$}\raisebox{-14pt}{$\surdex{\dfrac{b}{a}}$}\,\raisebox{-21pt}{$\ldots$}$? 

\vspace{4pt}

Answer: $x^3 - 2ax^2 + a^2x - b = 0$.

\vspace{4pt}

\textbf{25}) What equation gives the value of the infinite series [\emph{unendlichen Reihe}] $^3\!\!\!\!\surd(a + \;^3\!\!\!\!\surd(a + \;^3\!\!\!\!\surd(a + \ldots)$? 

\vspace{4pt}

Answer: $x^3 - x = a$.

\end{quote}

\vspace{4pt}

The notation and terminology leave much to be desired. Heis's book saw 107 editions in various languages, and was in print until the early 1900s, about 30 years after his death.


\item\vspace{9pt} \label{Fan1850} (1850) A. Fanien\index[ppl]{Fanien, A.|textbf}, Sur le calcul de $pi$. \emph{Nouvelles Annales de Math\'{e}matiques} 1re s\'{e}rie. \textbf{9}, 190--192.
\index{Vi\`{e}te's formula for $\tfrac{2}{\pi}$}
\index[ppl]{Vi\`{e}te, Fran\c{c}ois}

{\footnotesize Source: \href{http://www.numdam.org/item/NAM_1850_1_9__190_0.pdf}{NUMDAM}.}

{\footnotesize Cited in
\hyperref[Fam1854]{\textsc{Famin 1854}.}\index[ppl]{Famin, A.}}

\vspace{4pt}

This is another derivation, via regular polygons\index{polygons, regular} inscribed in the unit circle, of the formula for $\pi$ given in \hyperref[Cat1842]{\textsc{Catalan 1842}}\index[ppl]{Catalan, E.}, equation \eqref{E:Cat1842-1} above\index{pi@$\pi\;(3.14159\ldots)$!continued square root expressions for}\index{constants, named!pi@$\pi\;(3.14159\ldots)$}. In a footnote, Fanien claims that \eqref{E:Cat1842-1} was known to Vi\`{e}te\index[ppl]{Vi\`{e}te, Fran\c{c}ois}. 


\item\vspace{9pt} \label{Fam1854} (1854) A. Famin, Programme du Baccalaur\'{e}at \`{e}s-sciences---G\'{e}om\'{e}trie, no. 31, Math\'{e}matiques appliqu\'{e}es, no. 55. \emph{Revue de L'enseignement Chr\'{e}tien; Recueil P\'{e}riodique Publi\'{e} par les Professeurs de L'Assomption}, \textbf{3}, (January), L. Giraud, Libraire, Nimes.
\index{pi@$\pi\;(3.14159\ldots)$!continued square root expressions for}\index{constants, named!pi@$\pi\;(3.14159\ldots)$}

{\footnotesize Source: \href{https://books.google.com/books?id=RI3-rfDt5-4C&newbks=1&newbks_redir=0&dq=Programme\%20du\%20Baccalaur\%C3\%A9at\%C3\%A8s-sciences---G\%C3\%A9om\%C3\%A9trie\%2C\%20no.\%2031\%2C\%20Math\%C3\%A9matiques\%20appliqu\%C3\%A9es\%2C\%20no.\%2055\&pg=PA379\#v=onepage\&q\&f=false}{Google Books}.}

\vspace{4pt}

Yet another derivation of identity \eqref{E:Cat1842-1}.


\item\vspace{9pt} \label{Ami1855} (1855) A. Amiot\index[ppl]{Amiot, A.|textbf}, \emph{\'{E}l\'{e}ments de G\'{e}om\'{e}trie: r\'{e}dig\'{e}s d'apr\`{e}s le nouveau programme de l'enseignement scientifique des Lyc\'{e}es}, Paris, Dezobry, E. Magdeleine et Cie, Lib.-\'{E}diteurs.
\index{Vi\`{e}te's formula for $\tfrac{2}{\pi}$}

{\footnotesize Source: \href{https://gallica.bnf.fr/ark:/12148/bpt6k6471737v/f129.item.texteImage}{Gallica}.}

\vspace{4pt}

This geometry textbook was reprinted many times between 1855 and 1875; sometime after 1860 the firm of Ch. Delagrave et Cie, Libraires-\'{E}diteurs became the publisher. Pages 119 and 120 of the first edition develop the formula \eqref{E:Cat1842-1} above from \hyperref[Cat1842]{\textsc{Catalan 1842}}.\index{pi@$\pi\;(3.14159\ldots)$!continued square root expressions for}\index{constants, named!pi@$\pi\;(3.14159\ldots)$} No references are cited.


\item\vspace{9pt} \label{Bou1862} (1862) A. Bouch\'{e}\index[ppl]{Bouch\'{e}, A.|textbf}, Premier essai sur la th\'{e}orie des radicaux continus, et sur ses applications \`{a} l'alg\'{e}bre et au calcul infinit\'{e}simal. \emph{M\'{e}moires de la Soci\'{e}t\'{e} Acad\'{e}mique de Maine et Loire} \textbf{12}, Mallet-Bachelier, Paris, 81--151.
\index{successive substitution!19th century examples}
\index{continued radicals}
\index{radicaux continus@\emph{radicaux continus}}
\index{continued robb roots@continued $r$th roots!as solutions to polynomial equations}

{\footnotesize Source: \href{https://books.google.com/books?id=dbk-AQAAMAAJ&lpg=RA1-PA81&ots=oggtxwhLqN&dq=\%22Premier\%20essai\%20sur\%20la\%20th\%C3\%A9orie\%20des\%20radicaux\%20continus\%2C\%20et\%20sur\%20ses\%20applications\%20\%C3\%A0\%20l'alg\%C3\%A9bre\%20et\%20au\%20calcul\%20infinit\%C3\%A9simal\%22\&pg=RA1-PA81\#v=onepage\&q\&f=false}{Google Books}.}

{\footnotesize Cited in 
\hyperref[Rea1877]{\textsc{Realis 1877}}\index[ppl]{Realis, S.}.}

\vspace{4pt}

In this journal's Vol. 9 (1861), the minutes of a meeting held on 9 January 1861 state: ``M. Bouch\'{e}\index[ppl]{Bouch\'{e}, A.} summarizes his theory of continued radicals; this work is referred to the editorial board." This long paper is apparently the result. Using successive substitution, Bouch\'{e} develops continued $r$th root solutions to certain polynomial and other functional equations. The first chapter introduces the article's primary target, the trinomial equation\index{trinomial equations} $y^m-y-P=0$, where $m>1$ and $P>0$; successive substitution is accomplished with $y=\sqrt[m]{P+y}$. The example $y^3-y-7=0$ is computed in detail in Chapter 2. In Chapter 3 the author considers the consequences of substituting $y=z/k+t$ in his trinomial $m$th degree equation, although focusing primarily on the case $m=2$. The case $m=3$ consumes Chapter 4; the polynomial $y^m-Ay^n-By^p-Cyq-D=0$ is the subject of Chapter 5, along with some formal observations on properties of a continued radical solution. The author extends his methods to $x^m-Ax^{m-1}-Bx^{m-2}-\ldots-Gx-H=0$ in Chapter 6, to $F(x)=x^m+\sum_{i=1}^m P_ix^{m-i}=0$ in Chapter 7, and to other functional equations like $y^y=e$ in Chapter 8.

\hyperref[Rea1877]{\textsc{Realis 1877}} remarks that ``the theory of continued radicals is far from new, and Mr. A. Bouch\'{e} is not the first to have studied it.''


\item\vspace{9pt} \label{San1862} (1862) L\'{e}on Sancery\index[ppl]{Sancery, L\'{e}on|textbf}, De la m\'{e}thode des substitutions successives pour le calcul des racines des \'{e}quations. \emph{Nouvelles annales de math\'{e}matiques, 2nd series} \textbf{1}, 305--312 and 384--400.
\index{successive substitution!19th century examples}
\index{iteration!of functions}
\index{iteration!convergence conditions for}

{\footnotesize Source: \href{http://www.numdam.org/item/NAM_1862_2_1__384_1.pdf}{NUMDAM}.}

\vspace{4pt}

Sancery supposes that $F(x)=0$ is an algebraic or transcendental equation which may be expressed in the form $x=\varphi(x)$, and for which $x=\alpha$ is a solution. By successive substitution from an initial value $x_1$, he forms the sequence $x_n=\varphi(x_{n-1})$, and asks what conditions must be met by $\varphi(x)$ so that the $x_n$ converge to $\alpha$. In short order he deduces the condition $0<\varphi'(x)<1$ if the sequence $\{x_n\}$ is to converge monotonically, and $-1<\varphi'(x)<0$ if the sequence is alternating. He applies his tests to the quadratic equation $x^2+px+q=0$, recast as $x=-\tfrac{q}{p}-\tfrac{x^2}{p}$, and to the cubic trinomial\index{trinomial equations}\index{continued square roots!as solutions to quadratic equations} $x^3-px+q=0$\index{continued cube roots!as solutions to cubic equations}, where again the linear term is isolated on the left. The limitations of these examples lead him to consider the case where $\varphi(x)$ is invertible, and ultimately to the case in which $F(x)=0$ can only be expressed as $\psi(x)=\varphi(x)$, where $\psi(x)$ is not a linear function.




\item\vspace{9pt} \label{Did1872} (1872) ``M. le g\'{e}n\'{e}ral Didion" [Isidore Didion]\index[ppl]{Didion, Isidore|textbf}, Expression du rapport de la circonf\'{e}rence au diam\`{e}tre et nouvelle fonction. \emph{Comptes rendus hebdomadaires des s\'{e}ances de l'Acad\'{e}mie des Sciences} \textbf{74}, 36--39.
\index{Vi\`{e}te's formula for $\tfrac{2}{\pi}$}
\index{pi@$\pi\;(3.14159\ldots)$!continued square root expressions for}\index{constants, named!pi@$\pi\;(3.14159\ldots)$}

{\footnotesize Source: \href{https://gallica.bnf.fr/ark:/12148/bpt6k3031q/f36.item}{Gallica}.}

{\footnotesize Cited in 
\hyperref[Cat1872]{\textsc{Catalan 1872}}\index[ppl]{Catalan, E.} and
\hyperref[Wie1904b]{\textsc{Wiernsberger 1904b}}\index[ppl]{Wiernsberger, Paul}.}

\vspace{4pt}

General Didion independently rediscovers formula \eqref{E:Cat1842-1} above,\index{pi@$\pi\;(3.14159\ldots)$!continued square root expressions for}\index{constants, named!pi@$\pi\;(3.14159\ldots)$} along with related expressions arising from regular polygons\index{polygons, regular} inscribed in the unit circle. He seems to sense that the associativity of continued square roots\index{continued square roots}\index{associativity of function compositions}\index{function compositions!associativity of} of constant terms is ambiguous, and instead of
\begin{equation}\label{E:Did1872-1}
\sqrt{2-\sqrt{2+\sqrt{2\ldots+\sqrt{2+\sqrt{2+\sqrt{4-C^2}}}}}}
\end{equation}
he prefers to write
\begin{equation}\label{E:Did1872-2}
\sqrt{2-\left[\left\{\ldots\left[\left(\sqrt{4-C^2}+2\right)^{\frac{1}{2}}+2\right]^{\frac{1}{2}}\ldots+2\right\}^{\frac{1}{2}}+2\right]^{\frac{1}{2}}}\;,
\end{equation}
explaining that ``replacing the $\sqrt{\phantom{m}}$ sign by elevation to the power $\tfrac{1}{2}$, the expression \eqref{E:Did1872-1} for the [polygon] side will become \eqref{E:Did1872-2}, in which the indication of successive operations follows the natural order of writing, from left to right, and approaches that of series and continued fractions."


\item\vspace{9pt} \label{Cat1872} (1872) E. Catalan\index[ppl]{Catalan, E.|textbf}, Sur une communication r\'{e}cente de M. le g\'{e}n\'{e}ral Didion, concernant une expression du rapport de la circonf\'{e}rence au diam\`{e}tre. \emph{Comptes rendus hebdomadaires des s\'{e}ances de l'Acad\'{e}mie des Sciences} \textbf{74}, 177.
\index{Vi\`{e}te's formula for $\tfrac{2}{\pi}$}
\index{pi@$\pi\;(3.14159\ldots)$!continued square root expressions for}\index{constants, named!pi@$\pi\;(3.14159\ldots)$}
\index{letters to editors}

{\footnotesize Source: \href{https://gallica.bnf.fr/ark:/12148/bpt6k3031q/f177.item}{Gallica}.}

{\footnotesize Cited in
\hyperref[Wie1904b]{\textsc{Wiernsberger 1904b}}\index[ppl]{Wiernsberger, Paul}.}

\vspace{4pt}

Catalan immediately responds to \hyperref[Did1872]{\textsc{Didion 1872}}\index[ppl]{Didion, Isidore} with a curt letter stating, first, that Catalan himself had established the general's continued square root\index{continued square roots} formulas thirty years earlier (\hyperref[Cat1842]{\textsc{Catalan 1842}}), and, second, that ``[t]he true author of these various formulas is, if I am not mistaken, Euler\index[ppl]{Euler, Leonhard}. As early as 1760, this great Geometer gave this curious relation in the \emph{Nouveaux commentaires de P\'{e}tersbourg}:\index{infinite products!of secant functions}
\[\dfrac{\pi}{2} = \sec\dfrac{\pi}{4}\sec\dfrac{\pi}{8}\sec\dfrac{\pi}{16}\ldots\,.\]
It is easy to see that this does not differ, in substance, from the main formulas in question." Neither Catalan nor Didion seems to have been aware of Vi\`{e}te's work.\index[ppl]{Vi\`{e}te, Fran\c{c}ois}


\item\vspace{9pt} \label{Ast1877} (1877) J. J. \r{A}strand\index[ppl]{Astrand@\r{A}strand, J. J.|textbf}, Neue, einfache Transformation und Aufl\"{o}sung der Gleichungen von der Form $x^n-ax\pm b = 0$. \emph{Astronomische Nachrichten} \textbf{89}, 347--350.
\index{continued robb roots@continued $r$th roots!of terms $a_n=c$}
\index{continued robb roots@continued $r$th roots!as solutions to polynomial equations}
\index{trinomial equations}

{\footnotesize Source: \href{https://doi.org/10.1002/asna.18770892203}{doi.org (Wiley Online)}.}

{\footnotesize Cited in 
\hyperref[Gun1880]{\textsc{G\"{u}nther 1880}} and \hyperref[Gun1881]{\textsc{1881}}\index[ppl]{Gunther@G\"{u}nther, Siegmund}.}

\vspace{4pt}

In this note, the author proposes to solve the equation in the title by substituting $y\sqrt[n-1]{a}$ for $x$, which yields the transformed equation $y^n-y\pm c=0$, where $c=b/(a\sqrt[n-1]{a})$. Successive substitution\index{successive substitution!19th century examples} produces the solution\index{continued robb roots@continued $r$th roots!of terms $a_n=\pm c$}
\[\sqrt[n-1]{a}\sqrt[n]{\mp c+\sqrt[n]{\mp c +\sqrt[n]{\mp c+\cdots}}}\;.\]
\r{A}strand then performs the arithmetic to solve $x^3-7x+7=0$, which he cites as an equation from Lagrange's\index[ppl]{Lagrange, Joseph-Louis} \emph{Trait\'{e} de r\'{e}solution des equations num\'{e}riques} (1808). 


\item\vspace{9pt} \label{Rea1877} (1877) S. Realis\index[ppl]{Realis, S.|textbf}, Sur quelques questions propos\'{e}es dans la nouvelle correspondance, Question 142. \emph{Nouvelle correspondance math\'{e}matique}, \textbf{3}, 193--194.
\index{continued square roots!of terms $a_n=\pm2$}

{\footnotesize Source: \href{https://books.google.com/books?id=llhOAAAAYAAJ\&pg=PA193\#v=onepage\&q\&f=false}{Google Books}.}

\vspace{4pt}

The question, attributed to ``\'{E}. Lucas," asks for an equation whose roots are
\begin{equation}\label{E:Catalan}
x=\pm\sqrt{2\pm\sqrt{2\pm\sqrt{2\ldots}}}\;.
\end{equation}
In reply, Realis offers $x^2-x-2=0$, and mentions its relation to the more general\index{continued robb roots@continued $r$th roots!of terms $a_n=c$}
\[x=\sqrt[m]{p+\sqrt[m]{p+\sqrt[m]{p+\ldots}}};.\]
However, Realis seems not to have appreciated the subtlety of M. Lucas's $\pm$ signs, which, if assumed to have some periodic pattern, would produce various trigonometric identities\index{continued square roots!and trigonometric functions}. (If ``\'{E}. Lucas" is in fact Edouard Lucas\index[ppl]{Lucas, Edouard}, then this question may foreshadow some of the results in \hyperref[Luc1878]{\textsc{Lucas 1878}} concerning the connection between expressions of the form \eqref{E:Catalan} and trigonometric half-angle formulas.) 


\item\vspace{9pt} \label{Dix1878} (1878) T. S. E. Dixon\index[ppl]{Dixon, T. S. E.|textbf}, Continued roots. \emph{The Analyst} \textbf{5}(1), 20--21. 

{\footnotesize Source: \href{https://doi.org/10.2307/2635558}{doi.org (JSTOR)}.}

{\footnotesize Cited in
\hyperref[Jon2008]{\textsc{Jones 2008}} and \hyperref[Jon2015]{\textsc{2015}}\index[ppl]{Jones, Dixon J.}.}
\index{continued powers}
\index{continued robb roots@continued $r$th roots!of constant nonnegative real terms}

\vspace{4pt}

This note, barely more than a page long, contains the earliest appearance yet of the English-language terms ``continued root" and ``continued power". The author observes that the limit of \index{continued robb roots@continued $r$th roots!as solutions to polynomial equations}
\[\sqrt[n]{q+p\sqrt[n]{q+p\sqrt[n]{q+p\sqrt[n]{\cdots}}}}\]
is a root of the equation $x^n-px=q$\index{trinomial equations}; that \index{continued robb roots@continued $r$th roots!as solutions to polynomial equations}
\[\sqrt[m]{q-p\sqrt[m/n]{q-p\sqrt[m/n]{q-\ldots}}}\;,\]
``solves" the equation $x^m+px^n=q$; and that 
\[\sqrt{q-p\sqrt{q-p\sqrt{q-p\sqrt{\cdots}}}}\quad\textrm{and}\quad\dfrac{q}{p}-{\dfrac{1}{p}\;}^2\negmedspace\Biggl(\dfrac{q}{p}-{\dfrac{1}{p}\;}^2\negmedspace\Biggl(\dfrac{q}{p}-\cdots\Biggr)\Biggr)\]
are both roots of $x^2+px=q$.\index{continued squares} (The latter notation, placing the exponent to the left of the left parenthesis,\index{notation!left-justified exponent $^p(\cdot)$} was independently proposed more than a century later in \hyperref[Jon1991]{\textsc{Jones 1991}}\index[ppl]{Jones, Dixon J.}.) A continued reciprocal root solution\index{continued reciprocal roots} to the equation $x^{m+n}+px^n=q$ is exhibited:
\begin{equation}\label{E:Dix1878-1}
\begin{array}{llll}
x=\sqrt[n]{q + \quad\quad}\\
&\hspace*{-42pt}\overline{p+\sqrt[n/m]{q+\quad\quad}}\\
&&\hspace*{-42pt}\overline{p+\sqrt[n/m]{q+\quad\quad}}\\
&&&\hspace*{-42pt}\overline{p+\sqrt[n/m]{\cdots}}\;.
\end{array}
\end{equation}
Nothing is specified about $m$ or $n$, but if $n<m$ is assumed, then \eqref{E:Dix1878-1} could be construed as a continued reciprocal $\tfrac{m}{n}$ power\index{continued reciprocal powers}. The results are formal, with no consideration given to convergence, nor are any limitations on the values of $p$ or $q$ mentioned. Most of these observations are covered or anticipated in \hyperref[Dop1832b]{\textsc{Doppler 1832b}}\index[ppl]{Doppler, Christian}.


\item\vspace{9pt} \label{Luc1878} (1878) Edouard Lucas\index[ppl]{Lucas, Edouard|textbf}, Th\'{e}orie des fonctions num\'{e}riques simplement p\'{e}riodiques. \emph{American Journal of Mathematics} \textbf{1}, 184--240.
\index{continued square roots!of terms $a_n=c$}
\index{continued square roots!and trigonometric functions}

{\footnotesize Source: \href{https://doi.org/10.2307/2369308}{doi.org (JSTOR)}.}

{\footnotesize For an English translation, see \emph{The Theory of Simply Periodic Numerical Functions}, translated by Sidney Kravitz\index[ppl]{Kravitz, Sidney}, edited by Douglas Lind\index[ppl]{Lind, Douglas}, Fibonacci Association (1969).}

{\footnotesize Cited in 
\hyperref[Gun1881]{\textsc{G\"{u}nther 1881}}\index[ppl]{Gunther@G\"{u}nther, Siegmund}, and
\hyperref[Vel2016c]{\textsc{Vellucci and Bersani 2016c}}\index[ppl]{Vellucci, Pierluigi}\index[ppl]{Bersani, Alberto Maria}.}

\vspace{4pt}

Section XV of this paper deals with the ``relation of the functions $U_n$ and $V_n$ with continued radicals," where\index{continued square roots!of periodic real terms} 
\[U_n=\dfrac{a^n-b^n}{a-b}\;,\quad V_n=a^n+b^n\;,\] 
and $a$ and $b$ are the two roots of the equation $x^2=Px-Q$ for $P$ and $Q$ relatively prime integers. By successive substitution\index{successive substitution!19th century examples} one may write
\[a=\sqrt{-Q+P\sqrt{-Q+P\sqrt{-Q+\cdots}}}\]
(Lucas's Equation 88). Lucas also builds iterated formulas for $V_n$, and notes their resemblance to expressions for $\cos\tfrac{\pi}{2^n}$, $\cos\tfrac{\pi}{3\cdot2^n}$, and $\cos\tfrac{\pi}{5\cdot2^n}$.


\item\vspace{9pt} \label{Gun1880} (1880) S. G\"{u}nther\index[ppl]{Gunther@G\"{u}nther, Siegmund|textbf}, Eine didaktisch wichtige Aufl\"{o}sung der trinomischer Gleichungen, \emph{Zeitschrift f\"{u}r mathematischen und naturwissenschaftlichen Unterricht [Hoffmann Z.]} \textbf{11} 68--72.
\index{trinomial equations}
\index{continued reciprocal roots!as solutions to polynomial equations}
\index{pension problem}

{\footnotesize Source: \href{https://books.google.com/books?id=ii0nAYhDKREC&newbks=1&newbks_redir=0&vq=G\%C3\%BCnther\&pg=PA68\#v=onepage\&q\&f=false}{Google Books}.}

{\footnotesize Originally published with the same title in \emph{Verhandlungen der vierunddreissigsten Versammlung deutscher Philologen und Schulm\"{a}nner in Trier, vom 24. bis 27. September 1879}. B.~G. Teubner, Leipzig, 187--190.}

{\footnotesize Cited in 
\hyperref[Sch1880]{\textsc{Schaewen 1880}}\index[ppl]{Schaewen, Paul von},
\hyperref[Gun1881]{\textsc{G\"{u}nther 1881}}\index[ppl]{Gunther@G\"{u}nther, Siegmund},
\hyperref[Hof1881]{\textsc{Hoffmann 1881}}\index[ppl]{Hoffmann, K. E.}, 
\hyperref[Ise1888]{\textsc{Isenkrahe 1888}}\index[ppl]{Isenkrahe, C.}, 
\hyperref[Hey1894a]{\textsc{Heymann 1894a}}\index[ppl]{Heymann, W.}, 
\hyperref[Ise1897]{\textsc{Isenkrahe 1897}}\index[ppl]{Isenkrahe, C.}, 
\hyperref[Hey1898]{\textsc{Heymann 1898}}\index[ppl]{Heymann, W.}, 
\hyperref[Hey1904]{\textsc{Heymann 1904}}\index[ppl]{Heymann, W.},
\hyperref[Gol1911]{\textsc{Goldziher 1911}}\index[ppl]{Goldziher, Karl}, and 
\hyperref[Jon2015]{\textsc{Jones 2015}}\index[ppl]{Jones, Dixon J.}.}
\index{continued reciprocal roots} 

\vspace{4pt}

In his introduction, G\"{u}nther recounts the efforts of Lambert\index[ppl]{Lambert, Johann Heinrich}, Malfatti\index[ppl]{Malfatti, Gian Francesco}, Gauss\index[ppl]{Gauss, Carl Friedrich}, and others to solve the trinomial equation $x^p\pm\alpha x=\beta$; mentions \hyperref[Ast1877]{\textsc{\r{A}strand 1877}}\index[ppl]{Astrand@\r{A}strand, J. J.} in this regard; and credits Jacob Bernoulli\index[ppl]{Bernoulli, Jacob} (meaning, one supposes, \hyperref[Ber1692]{\textsc{Bernoulli 1692}}\index[ppl]{Bernoulli, Jacob}) with the ideas concerning successive substitution that underlie \r{A}strand's method. G\"{u}nther himself proposes a continued reciprocal root solution, similar to expression \eqref{E:Dop1832-3} above in \hyperref[Dop1832b]{\textsc{Doppler 1832b}}\index[ppl]{Doppler, Christian}. One of his goals is to solve the ``pension problem"\index{pension problem}: finding the rate of interest $q$ in the equation
\[r\dfrac{q^n-1}{q-1}=aq^n\;,\]
where $a$ is the single premium, $r$ the pension rate, and $n$ the time. G\"{u}nther rewrites this as $q^{n+1}-\tfrac{a+r}{a}q^n=-r$ and invokes successive substitution\index{successive substitution!19th century examples} to get
\begin{equation*}\label{E:Gun1880-1}
\begin{array}{llll}
q=\sqrt[n]{\dfrac{-r}{\phantom{xxxxxxxxxxxxxxxxxx}}}\\
&\hspace*{-104pt}{-\tfrac{a+r}{a}+\sqrt[n]{\dfrac{-r}{-\tfrac{a+r}{a}+\ldots}}}
\end{array}
\end{equation*}

A review in \emph{Jahrbuch \"{u}ber die Fortschritte der Mathematik} \textbf{12} (1882), p. 76, states: ``After historical notes about the solution of the trinomial equation, the author gives the solution of the equation 
\[x^{m+n}+ax^m=b\]
in the following form
\begin{equation}\label{E:Gun1880-2}
\begin{array}{llll}
x=\sqrt[m]{b + \quad\quad}\\
&\hspace*{-42pt}\overline{a+\sqrt[m/n]{b+\quad\quad}}\\
&&\hspace*{-42pt}\overline{a+\sqrt[m/n]{b+\quad\quad}}\\
&&&\hspace*{-42pt}\overline{a+\sqrt[m/n]{\cdots}}\;.
\end{array}
\end{equation}
He then discusses the forms this structure takes in the special cases where $m/n = 1$, $m=1$, $n=1$, and carefully puts an end to defects of this solution, namely the lack of a readily detectable convergence, among others." However, it emerges immediately with \hyperref[Sch1880]{\textsc{Schaewen 1880}}, and subsequently in \hyperref[Hof1881]{\textsc{Hoffmann 1881}}\index[ppl]{Hoffmann, K. E.}, \hyperref[Net1887]{\textsc{Netto 1887}}, and \hyperref[Ise1888]{\textsc{Isenkrahe 1888}}\index[ppl]{Isenkrahe, C.}, that G\"{u}nther's\index[ppl]{Gunther@G\"{u}nther, Siegmund} solution has defects which are not resolved. 


\item\vspace{9pt} \label{Sch1880} (1880) ``Gymnasialleher v. Schaewen" [Paul von Schaewen]\index[ppl]{Schaewen, Paul von|textbf}, Zur L\"{o}sung trinomischer Gleichungen. \emph{Zeitschrift f\"{u}r mathematischen und naturwissenschaftlichen Unterricht [Hoffmann Z.]} \textbf{11}, 264--267.
\index{continued reciprocal roots!as solutions to polynomial equations}
\index{trinomial equations}

{\footnotesize Source: \href{https://books.google.com/books?id=ii0nAYhDKREC&newbks=1&newbks_redir=0&pg=PA264#v=onepage&q&f=false}{Google Books}.}

{\footnotesize Cited in 
\hyperref[Hey1894a]{\textsc{Heymann 1894a}} and \hyperref[Hey1894b]{\textsc{1894b}}\index[ppl]{Heymann, W.}, 
\hyperref[Hey1898]{\textsc{Heymann 1898}}\index[ppl]{Heymann, W.}, and in
\hyperref[Gol1911]{\textsc{Goldziher 1911}}\index[ppl]{Goldziher, Karl}.}

\vspace{4pt}

In this note, a high school teacher in Saarbr\"{u}cken takes issue with several claims made in \hyperref[Gun1880]{\textsc{G\"{u}nther 1880}}\index[ppl]{Gunther@G\"{u}nther, Siegmund}. After offering corrections to a few typographical errors in G\"{u}nther's paper, Schaewen works several examples showing that G\"{u}nther's continued reciprocal root\index{continued reciprocal roots} algorithm converges to the intended root of the initial trinomial equation\index{trinomial equations} only in special cases, and gives some empirical deductions based on these calculations. In a gracious postscript, Dr. G\"{u}nther admits that ``[von Schaewen's] remark that the [continued reciprocal roots] are not suitable for numerical computation appears to be generally true. For the time being, as long as the nature of these new forms is still little studied, we can only conclude that in some instances we have found [his observations] confirmed, and that in others there was a fairly rapid convergence." The paper is (thus far) unique in using the term \emph{Kettenradical}\index{Kettenradical@\emph{Kettenradical}}, in this case in reference to \eqref{E:Gun1880-2}.


\item\vspace{9pt} \label{Far1881} (1881) Gyul\'{a}t\'{o}l Farkas\index[ppl]{Farkas, Gyul\'{a}t\'{o}l (\emph{aka} Gyula or Jules)|textbf}, A Bolyai-f\'{e}le Algorithmus, \emph{\'{E}rtekez\'{e}sek a Mathematikai Tudom\'{a}nyok K\"{o}r\'{e}b\H{o}l}, \textbf{8}, 1--8.

{\footnotesize Source: \href{https://real-eod.mtak.hu/2740/1/Matematikai_ertekezesek_8_kotet_3_szam.pdf}{REAL-EOD}.}

{\footnotesize Cited in 
\hyperref[Sza2010]{\textsc{Szab\'{o} 2010}}\index[ppl]{Szab\'{o}, P\'{e}ter G\'{a}bor}.}
\index{successive substitution!19th century examples}
\index{trinomial equations}
\index{Bolyai's algorithm}

\vspace{4pt}

The author (whose anglicized name is Gyula Farkas) was a Hungarian mathematician and physicist. In this article he calls the successive substitution technique in \hyperref[Bol1832]{\textsc{Bolyai 1832}} ``Bolyai's algorithm", and extends it to solve the trinomial equation $x^m=ax+b$. \hyperref[Sza2010]{\textsc{Szab\'{o} 2010}} and other sources suggest that, because of this article, successive substitution using $m$th roots  became known in some circles as ``Bolyai's algorithm".


\item\vspace{9pt} \label{Gun1881} (1881) Siegmund G\"{u}nther\index[ppl]{Gunther@G\"{u}nther, Siegmund|textbf}, \emph{Die Lehre von den gew\"{o}hnlichen und verallgeminerten Hyperbelfunktionen}, Verlag von Louis Nebert, Halle a/S.
\index{root-finding approximations}
\index{continued square roots!and trigonometric functions}

{\footnotesize Source: \href{https://gdz.sub.uni-goettingen.de/id/PPN578420767}{G\"{o}ttinger Digitalisierungszentrum}.}

{\footnotesize Cited in 
\hyperref[Hey1894b]{\textsc{Heymann 1894b}}\index[ppl]{Heymann, W.}.}

\vspace{4pt}

Pages 75--82 give an overview of, and some additional detail about, the functions $U_n$ and $V_n$ defined in \hyperref[Luc1878]{\textsc{Lucas 1878}}\index[ppl]{Lucas, Edouard}, which can be used to generate continued square root representations of trigonometric functions of special angles. 


\item\vspace{9pt} \label{Hof1881} (1881) K. E. Hoffmann\index[ppl]{Hoffmann, K. E.|textbf}, Ueber die Aufl\"{o}sung der trinomischen Gleichungen durch kettenbruch\"{a}hnliche Algorithmen, \emph{Archiv der Mathematik und Physik}, \textbf{66}, 33--45.
\index{trinomial equations}
\index{successive substitution!19th century examples}
\index{iteration!convergence conditions for}

{\footnotesize Source: \href{https://books.google.com/books?id=YjyK4OnRfygC\&pg=PA33\#v=onepage\&q\&f=false}{Google Books}.}

{\footnotesize Cited in
\hyperref[Ise1888]{\textsc{Isenkrahe 1888}} and \hyperref[Ise1897]{\textsc{1897}}\index[ppl]{Isenkrahe, C.}, 
\hyperref[Gol1911]{\textsc{Goldziher 1911}}\index[ppl]{Goldziher, Karl}, 
\hyperref[Thr1961]{\textsc{Thron 1961}}\index[ppl]{Thron, W. J.}, and 
\hyperref[Jon2015]{\textsc{Jones 2015}}\index[ppl]{Jones, Dixon J.}.}

\vspace{4pt}

Hoffmann addresses some of the difficulties arising in \hyperref[Gun1880]{\textsc{G\"{u}nther 1880}}\index[ppl]{Gunther@G\"{u}nther, Siegmund}. He observes that the trinomial equation
\[f(x) = x ^ m + px ^ n + q = 0\]
can be rearranged, on the one hand, in the form 
\[x^n(x^{m-n}+p)=-q\] 
to obtain ``Algorithm A": 
\[x_{i+1}=\sqrt[n]{\dfrac{-q}{p+x_{i}^{m-n}}}\]
as the basis for successive substitution\index{successive substitution!19th century examples}, starting at $x_0=0$. On the other hand, the trinomial\index{trinomial equations} may also take the form 
\[x^{m-n}=-p-\dfrac{q}{x^n}\]
which yields ``Algorithm B": 
\[x_{i+1}=\sqrt[m-n]{-p-\dfrac{q}{x_{i}^n}}\;,\]
with initial value $x_0 = \infty$. He shows that some of G\"{u}nther's calculations mistakenly assume that iteration of each form yields the same limiting value. One of Hoffmann's principle results is the following: ``All real roots of trinomial equations\index{trinomial equations} can be found by [G\"{u}nther's] continued fraction-like algorithms, as follows: Algorithm A supplies roots lying in the interval
\[-\left(\!\sqrt[m-n]{\tfrac{np}{m}}\,\right)\quad\textrm{ to }+\left(\!\sqrt[m-n]{\tfrac{np}{m}}\,\right)\;,\]
while one finds the roots lying outside this interval using Algorithm B." 


\item\vspace{9pt} \label{Far1884} (1884) Jules (Gyul\'{a}t\'{o}l) Farkas\index[ppl]{Farkas, Gyul\'{a}t\'{o}l (\emph{aka} Gyula or Jules)|textbf}, Sur les fonctions it\'{e}ratives, \emph{Journal de math\'{e}matiques pures et appliqu\'{e}es} 3rd series, \textbf{10}, p. 101--108, \textbf{8}, 101--108.
\index{successive substitution!19th century examples}
\index{iteration!of functions}
\index{iteration!convergence conditions for}

{\footnotesize Source: \href{http://www.numdam.org/item/JMPA_1884_3_10__101_0.pdf}{NUMDAM}.}

{\footnotesize Cited in \hyperref[Sza2010]{\textsc{Szab\'{o} 2010}}}\index[ppl]{Szab\'{o}, P\'{e}ter G\'{a}bor}.

\vspace{4pt}

Following his previous work on iterated functions (\hyperref[Far1881]{\textsc{Farkas 1881}}), Farkas writes, ``If the function $f(z)$ satisfies certain conditions, its iterates are convergent, and, in my first Note, I describe such conditions. The problem of analytical iteration aims to express the iterates of an analytical function by an analytical function of the index $k$ considered as an independent variable\ldots. [I]n my second Note, I establish and treat the general form of the analytical function of two variables $F(k,z)$ defined by the expression $F(k,z)=f^k(z)$." He gives convergence conditions for $f(z)$ defined by $a+\sqrt[n]{z}$, $a+b(\sin z)^n$, $a+b(\log z)^n$, $\cos(a-bz)$, and $\log(a+bz)$.


\item\vspace{9pt} \label{Net1887} (1887) Eugen Netto\index[ppl]{Netto, Eugen|textbf}, Ueber einen Algorithmus zur Aufl\"{o}sung numerischer algebraischer Gleichungen, \emph{Mathematische Annalen} \textbf{29} 141--147.
\index{successive substitution!19th century examples}
\index{iteration!of functions}
\index{periodic points}

{\footnotesize Source: \href{https://gdz.sub.uni-goettingen.de/id/PPN235181684_0029?tify=\%7B\%22pages\%22\%3A\%5B147\%5D\%2C\%22pan\%22\%3A\%7B\%22x\%22\%3A0.51\%2C\%22y\%22\%3A0.805\%7D\%2C\%22view\%22\%3A\%22info\%22\%2C\%22zoom\%22\%3A0.375\%7D}{G\"{o}ttinger Digitalisierungszentrum}.}

{\footnotesize Cited in 
\hyperref[Ise1888]{\textsc{Isenkrahe 1888}} and \hyperref[Ise1897]{\textsc{1897}}\index[ppl]{Isenkrahe, C.}, and 
\hyperref[Gol1911]{\textsc{Goldziher 1911}}\index[ppl]{Goldziher, Karl}.}

\vspace{4pt}

The paper's notation has been somewhat simplified in what follows. The author notes that iteration of
\begin{equation}\label{E:Net1897-1}
x_{k+1}=\sqrt[n]{x_k+a}
\end{equation}
has been used ``many times" to solve the trinomial equation\index{trinomial equations}\index{continued robb roots@continued $r$th roots!as solutions to polynomial equations}\index{continued robb roots@continued $r$th roots!as solutions to polynomial equations} $x^n-x-a=0$, but no references are cited. He sets $f(x)=ax^{n-1}+bx^{n-2}+\cdots+c$, and lists $f$'s real roots as $r_1, r_2, \ldots, r_i,\ldots$ in descending order (assuming initially that none are duplicated). From an arbitrary initial value $x_0$, he inductively defines $x_{k+1}=\sqrt[n]{f(x_k)}$. For an arbitrary root $r_i$, he establishes a sequence of distances $\delta_k=x_k-r_i$, and deduces that for $\delta_k$ sufficiently small, 
\[\delta_{k+1}=\dfrac{f'(r_i)}{nr_{i}^{n-1}}\delta_k=q(r_i)\delta_k\]
``with arbitrary precision." He identifies six cases, accordingly as $q(r_i)$ is positive or negative and of magnitude greater than, equal to, or less than $1$; he meticulously analyzes the sequence of iterates $x_i$ in each case; and closes with a detailed example using the trinomial $x^n-x-a=0$.\index{trinomial equations}

Perhaps the most interesting aspect of this paper is its glimpse into the intricate structure of periodic points\index{periodic points} arising out of algorithm \eqref{E:Net1897-1} under certain conditions. In the long paragraph beginning at the top of page 144, Netto discovers the ``remarkable circumstance" that ``\emph{if $r_{2i+1}$ is negative and $f'(r_{2i+1})$ is greater than $nr_{2i+1}^n$, there is between $r_{2i}$ and $r_{2i+1}$, and also between $r_{2i+2}$ and $r_{2i+1}$, an odd number of values $x_0$ such that two values $x_0, x_1$ periodically repeat so that the algorithm does not converge.} In the same way, it is also possible that $x_\lambda = x_0$ without the equality $x_\mu = x_0$ even for $\mu<\lambda$. Indeed, there may be an infinite number of values $x_0$ with this property$\ldots$ \emph{If we let $\lambda$ pass through the sequence of primes, the existence of infinitely many such $x_0$ follows.}" (Italics are in the original.) The resemblance to Sharkovskii's\index[ppl]{Sharkovskii, O. M.} results from 1964 is a bit uncanny. 




\item\vspace{9pt} \label{Ise1888} (1888) C. Isenkrahe\index[ppl]{Isenkrahe, C.|textbf}, Ueber die Anwendung iterirter Functionen zur Darstellung der Wurzeln algebraischer und transcendenter Gleichungen, \emph{Mathematische Annalen} \textbf{31} 309--317.
\index{successive substitution!19th century examples}
\index{iteration!convergence conditions for}

{\footnotesize Source: \href{https://doi.org/10.1007/BF01206217}{doi.org (Springer)}.} 

{\footnotesize Cited in
\hyperref[Ise1897]{\textsc{Isenkrahe 1897}}\index[ppl]{Isenkrahe, C.}, 
\hyperref[Gol1911]{\textsc{Goldziher 1911}}\index[ppl]{Goldziher, Karl}, 
\hyperref[Thr1961]{\textsc{Thron 1961}}\index[ppl]{Thron, W. J.}, 
and \hyperref[Jon2015]{\textsc{Jones 2015}}\index[ppl]{Jones, Dixon J.}.}

\vspace{4pt}

After reviewing the results of \hyperref[Hof1881]{\textsc{Hoffmann 1881}}\index[ppl]{Hoffmann, K. E.} and \hyperref[Net1887]{\textsc{Netto 1887}} concerning the solution of polynomial equations of the form $x=f(x)$ by iterated function algorithms\index{iterated compositions!as solutions to polynomial equations}\index{continued robb roots@continued $r$th roots!as solutions to polynomial equations}, Isenkrahe establishes general principles governing the convergence of such algorithms. Some of his main results are:

\begin{enumerate}

\item\vspace{9pt} Iteration of $f(x)$ can only converge to a fixed value $\xi_a$ when $|f'(\xi_a)|<1$.

\item\vspace{9pt} If the iteration of $f(x)$ converges to $\xi_a$, this convergence is always more rapid the smaller the modulus of $f'(x)$, and the slower, the nearer this modulus is to 1.

\item\vspace{9pt} For real values, the iteration approaches the limit $\xi_a$ only from one side if $f'(\xi_a)$ is positive, while the value oscillates around $\xi_a$ if $f'(\xi_a)$ is negative.

\item\vspace{9pt} Iteration of the function $\tfrac{f(x)-xf'(x)}{1-f'(x)}$ yields the complete real and complex roots of the equation $x=f(x)$, depending on the starting point.

\end{enumerate}

\vspace{9pt}

\noindent(Though elegantly stated, observations (a)--(c) were known at least as far back as \hyperref[San1862]{\textsc{Sancery 1862}}.) Isenkrahe points out that Newton's method\index{Newton's or Newton-Raphson method} for approximating roots of equations is an application of (d); he shows that Gauss's\index[ppl]{Gauss, Carl Friedrich} iterative solution to the Kepler problem (see \hyperref[Gau1809]{\textsc{Gauss 1809}})\index{Kepler problem} is an example of Newton's method\index{Newton's or Newton-Raphson method}\index[ppl]{Newton, Isaac}, justifies its convergence, and offers an elegant simplification of formula \eqref{E:Gau1809-1} above. Isenkrahe closes by considering the consequences of his analysis for complex roots. Near the beginning of the paper, and also at the end, he intimates that a more extensive work on this subject is forthcoming; see \hyperref[Ise1897]{\textsc{Isenkrahe 1897}}.


\item\vspace{9pt} \label{Fou1890} (1890) Maurice Fouch\'{e}\index[ppl]{Fouch\'{e}, Maurice|textbf}, Remarques sur la m\'{e}thode des p\'{e}rim\`{e}tres pour calculer le nombre $\pi$. \emph{Bulletin de la Soci\`{e}t\`{e} Math\'{e}matique de France} \textbf{18}, 135--138.
\index{Vi\`{e}te's formula for $\tfrac{2}{\pi}$}
\index{pi@$\pi\;(3.14159\ldots)$!continued square root expressions for}\index{constants, named!pi@$\pi\;(3.14159\ldots)$}

{\footnotesize Source: \href{http://www.numdam.org/item/10.24033/bsmf.411.pdf}{NUMDAM}.}

\vspace{4pt}

Presents a relatively short derivation of the reciprocal of Vi\`{e}te's\index[ppl]{Vi\`{e}te, Fran\c{c}ois} formula \eqref{E:Vie1593-2} above. Curiously, the author calls this only ``the known formula" for $\pi/2$; Vi\`{e}te's name is not mentioned. Compare \hyperref[Cat1842]{\textsc{Catalan 1842}}.


\item\vspace{9pt} \label{Pie1891} (1891) George Winslow Pierce\index[ppl]{Pierce, George Winslow|textbf}, \emph{The Life-Romance of an Algebraist}, J. G. Cupples, Boston, 1891, 18.
\index{pi@$\pi\;(3.14159\ldots)$!continued square root expressions for}\index{constants, named!pi@$\pi\;(3.14159\ldots)$}

{\footnotesize Source: \href{https://babel.hathitrust.org/cgi/pt?id=uc2.ark:/13960/t6542wp0v&seq=48}{HathiTrust}.}

{\footnotesize Cited in
\hyperref[Her1935]{\textsc{Herschfeld 1935}}\index[ppl]{Herschfeld, Aaron}.}

\vspace{4pt}

Mentions in passing the formula\index{pi@$\pi\;(3.14159\ldots)$!continued square root expressions for}\index{constants, named!pi@$\pi\;(3.14159\ldots)$}
\begin{equation}\label{E:Pie1891-1}
\pi=\lim_{n\to\infty}2^n\underbrace{\sqrt{2-\sqrt{2+\sqrt{2+\sqrt{2+\cdots}}}}}_{n\textrm{ square roots}}\;,
\end{equation}
given in \hyperref[Cat1842]{\textsc{Catalan 1842}} and many other sources, and rigorously interpreted in \hyperref[Her1935]{\textsc{Herschfeld 1935}}\index[ppl]{Herschfeld, Aaron}.


\item\vspace{9pt} \label{Rud1891} (1891) Ferdinand Rudio\index[ppl]{Rudio, Ferdinand|textbf}, \"{U}ber die Konvergenz einer von Vieta herr\"{u}hrenden eigent\"{u}mlichen Producktentwicklung. \emph{Zeitschrift f\"{u}r Mathematik und Physik, historisch-literarische Abteilung} \textbf{36}, 139--140.
\index{Vi\`{e}te's formula for $\tfrac{2}{\pi}$}
\index[ppl]{Vi\`{e}te, Fran\c{c}ois}

{\footnotesize Source: \href{https://gdz.sub.uni-goettingen.de/id/PPN599415665_0036?tify=\%7B\%22pages\%22\%3A\%5B557\%5D\%2C\%22view\%22\%3A\%22info\%22\%7D}{G\"{o}ttinger Digitalisierungszentrum}.}

{\footnotesize Cited in 
\hyperref[Bop1913]{\textsc{Bopp 1913}}\index[ppl]{Bopp, K.}, 
\hyperref[Sch1961]{\textsc{Schuske and Thron 1961}}\index[ppl]{Schuske, Georgellen}\index[ppl]{Thron, W. J.}, 
\hyperref[Siz1986]{\textsc{Sizer 1986}}\index[ppl]{Sizer, Walter S.}, and
\hyperref[Cas1988]{\textsc{Castellanos 1988}}\index[ppl]{Castellanos, Dario}.}

\vspace{4pt}

\hyperref[Cas1988]{\textsc{Castellanos 1988}}\index[ppl]{Castellanos, Dario} states that the convergence of Vi\`{e}te's\index[ppl]{Vi\`{e}te, Fran\c{c}ois} expression \eqref{E:Vie1593-1} above was proved for the first time in this short note. In fact, the note is a letter to an editor who more or less goaded Rudio to provide such a proof. In 1890, in the journal \emph{Vierteljahrsschrift der Naturforschenden Gesellschaft in Z\"{u}rich} \textbf{35}, Rudio published a long article on the history of the circle quadrature problem. Mentioning (the reciprocal of) Vi\`{e}te's formula for $\tfrac{2}{\pi}$, he added a footnote quoting the infinite product of secants\index{infinite products!of secant functions} from \hyperref[Eul1783]{\textsc{Euler 1783}}\index[ppl]{Euler, Leonhard}.\footnote{In the footnote on page 17 of his 1890 quadrature paper, Rudio expresses doubt that Euler was aware of Vi\`{e}te's formula.} The following year, in \emph{Zeitschrift f\"{u}r Mathematik und Physik, historisch-literarische Abteilung} \textbf{36}, editor Moritz Cantor\index[ppl]{Cantor, Moritz} reviewed Rudio's article. Cantor pointed out that a rigorous investigation of convergence had not yet been done for Vi\`{e}te's formula, and suggested that Rudio should have done so. Rudio's page-and-a-half response to Cantor uses the cosine half-angle formula to replace the secants in Euler's infinite product with finite continued square roots, and carefully establishes the convergence of the resulting expansion.


\item\vspace{9pt} \label{Hey1894a} (1894a) W. Heymann\index[ppl]{Heymann, W.|textbf}, Ueber die Aufl\"{o}sung der Gleichungen vom f\"{u}nften Grade. \emph{Zeitschrift f\"{u}r Mathematik und Physik} \textbf{39}, 162--182, 193--202, 257--272, 321--354.
\index{iteration!of functions}
\index{iteration!graphical methods}
\index{Kettenfunctionen@\emph{Kettenfunctionen}}

{\footnotesize Source: \href{https://gdz.sub.uni-goettingen.de/id/PPN599415665_0039?tify=\%7B\%22pages\%22\%3A\%5B325\%5D\%2C\%22view\%22\%3A\%22info\%22\%7D}{G\"{o}ttinger Digitalisierungszentrum}.}

{\footnotesize Cited in 
\hyperref[Hey1894b]{\textsc{Heymann 1894b}}\index[ppl]{Heymann, W.} and 
\hyperref[Gol1911]{\textsc{Goldziher 1911}}\index[ppl]{Goldziher, Karl}; plagiarized in \hyperref[Rab1911]{\textsc{Rabinowitsch 1911}}\index[ppl]{Rabinowitsch, Izko-Ewna}.}

\vspace{4pt}

This paper spans more than 90 journal pages in four installments; it possesses two main ``parts" (which are not quite synchronized with the installments), and offers a potpourri of section styles. Part I gives an overview of approaches to the quintic equation, citing the work of Klein\index[ppl]{Klein, Felix}, Abel\index[ppl]{Abel, Niels Henrik}, Paul Gordan\index[ppl]{Gordan, Paul}, and others; many special cases are examined in detail. Our interest here is the fourth installment: Section C (titled ``Ueber Kettenfunctionen")\index{Kettenfunctionen@\emph{Kettenfunctionen}} of Part II.\footnote{The sprawling organization of this paper has led to erroneous citations. Some later authors have cited Section C as a stand-alone paper called ``Ueber Kettenfunctionen."} The author defines a \emph{Kettenfunction} by ``a system of equations
\[x = \varphi_1x_1,\quad x_1=\varphi_2x_2,\quad x_2=\varphi_3x_3,\,\ldots,\,x_{n-1}=\varphi_n x_n\]
or by an expression
\[x = \varphi_1\varphi_2\varphi_3\cdots\varphi_n x_n\]
where the $\varphi$ are given functions. For
\[\varphi_k = a_k + \dfrac{1}{x_k}\quad(k=1,2,\ldots,n)\]
an ordinary continued fraction\index{continued fractions} arises; for
\[\varphi_k=\sqrt[p]{a_k + x_k}\]
one obtains a continued root\index{continued robb roots@continued $r$th roots!of arbitrary nonnegative real terms}\index{continued robb roots@continued $r$th roots!as solutions to polynomial equations}
\[x = \sqrt[p]{a_1+\sqrt[p]{a_2+\sqrt[p]{a_3+\cdots\sqrt[p]{a_n+x_n}}}}\]
and in the same way one can speak of a continued power\index{continued powers} [\emph{Kettenpotenz}\index{Kettenpotenz@\emph{Kettenpotenz}}], a continued logarithm\index{continued logarithms} [\emph{Kettenlogarithmus}\index{Kettenlogarithmus@\emph{Kettenlogarithmus}}], and so on." Heymann immediately follows this bold definition by limiting the discussion to periodic continued compositions, which, as we have seen in \hyperref[Dop1832b]{\textsc{Doppler 1832b}}\index[ppl]{Doppler, Christian} (not cited in this paper), are handled using successive substitution\index{successive substitution!19th century examples} and fixed-point iteration\index{fixed-point iteration}: ``The special subject of our paper does not allow us to enter into general questions about continued compositions. On the contrary, we restrict ourselves to$\,\ldots\,$continued compositions of period 2$\ldots$ Let
\begin{equation}\label{E:1}
y=\varphi(x)\quad\textrm{and}\quad y=\psi(x)
\end{equation}
be any two functions, but assume that they can easily be inverted." The author's interest in \emph{graphical} iteration is exquisitely demonstrated in 25 figures, all on a separate oversized sheet. 

In a postscript, Heymann reports having recently found \hyperref[Gun1880]{\textsc{G\"{u}nther 1880}}\index[ppl]{Gunther@G\"{u}nther, Siegmund} and \hyperref[Sch1880]{\textsc{Schaewen 1880}}\index[ppl]{Schaewen, Paul von}, and comments that his approach to continued compositions resolves G\"{u}nther's difficulties. The date ``June 1893" at the end indicates that this paper precedes \hyperref[Hey1894b]{\textsc{Heymann 1894b}}, which is dated ``December 1893." 


\item\vspace{9pt} \label{Hey1894b} (1894b) W. Heymann\index[ppl]{Heymann, W.|textbf}, Theorie der An- und Uml\"{a}ufe und Aufl\"{o}sung der Gleichungen vom vierten, f\"{u}nften und sechsten Grade mittelst goniometrischer und hyperbolischer Functionen. \emph{Journal f\"{u}r die reine und angewandte Mathematik} \textbf{113}, 267--302.
\index{successive substitution!19th century examples}
\index{iteration!of functions}
\index{iteration!graphical methods}
\index{iteration!cyclic orbits}
\index{Kettenfunctionen@\emph{Kettenfunctionen}}

{\footnotesize Source: \href{https://gdz.sub.uni-goettingen.de/id/PPN243919689_0113?tify=\%7B\%22view\%22:\%22info\%22,\%22pages\%22:\%5B271\%5D\%7D}{G\"{o}ttinger Digitalisierungszentrum}.}

{\footnotesize Cited in
\hyperref[Hey1901]{\textsc{Heymann 1901}}\index[ppl]{Heymann, W.} and \hyperref[Hey1904]{\textsc{1904}}.}

\vspace{4pt}

The prolific Heymann continues his research, initiated in \hyperref[Hey1894a]{\textsc{Heymann 1894a}}, into the solution of equations by iterated functions. In Section I, he begins to document the existence of cyclic orbits (\emph{indifferente Uml\"{a}ufe}), lays the groundwork for function iteration in polar coordinates, and briefly considers how iteration would produce imaginary solutions. Section II is devoted to a detailed investigation of $n$th-degree trinomial equations\index{trinomial equations}, $n\ge 1$, which are reduced to four cases corresponding to the possible arrangements of $+$ and $-$ signs between the three terms. His iterative solutions thus involve the functions $\varphi(x)=x^r$ and $\psi(x)=\pm x \pm c$. He makes the calculations (independently duplicated a century later in \hyperref[Jon1991]{\textsc{Jones 1991}}\index[ppl]{Jones, Dixon J.} and \hyperref[Sch1992]{\textsc{Sch\"{o}nefuss 1992}}) for the ranges of $c$ yielding an attracting fixed point. In Section III, Heymann tackles polynomials of the fourth, fifth, and sixth degree, simplified using circular and hyperbolic trig substitutions; he works out many numerical examples\index{continued robb roots@continued $r$th roots!as solutions to polynomial equations}. Finally, in Section IV he considers the question of reformulating the pair of functions $\varphi(x)$ and $\psi(x)$ so that their intersection is unchanged, while the rate of convergence of the iterative algorithm is increased. 


\item\vspace{9pt} \label{Boc1895} \dag(1895) Karl Bochow\index[ppl]{Bochow, Karl|textbf}, Eine einheitliche Theorie der regelm\"{a}ssigen Vielecke: Teil 1. In \emph{F\"{u}nfter Jahresbericht \"{u}ber die St\"{a}dtische Realschule zu Magdeburg. Ostern 1894 bis Ostern 1895}. Progr. No. 278. E. Baensch jun., Magdeburg, 1--14.
\index{continued square roots!and trigonometric functions}
\index{continued square roots!geometric interpretations of}
\index{continued square roots!of terms $a_n=\pm2$}

{\footnotesize Source: \href{https://www.google.com/books/edition/Eine_einheitliche_Theorie_der_regelm\%C3\%A4ss/G4AWFF3GunMC?hl=en\&gbpv=1\&dq=\%22Eine+einheitliche+theorie+der+regelm\%C3\%A4ssigen+vielecke\%22&pg=PA1&printsec=frontcover}{Google Books}.}

{\footnotesize Cited in 
\hyperref[Boc1896]{\textsc{Bochow 1896}}, 
\hyperref[Wie1904b]{\textsc{Wiernsberger 1904b}}\index[ppl]{Wiernsberger, Paul},
\hyperref[Boc1905]{\textsc{Bochow 1905}}\index[ppl]{Bochow, Karl}, and
\hyperref[Boc1910]{\textsc{Bochow 1910}}\index[ppl]{Bochow, Karl}.}

\vspace{4pt}

Published in the annual report of the high school in which the author was a teacher, this is the first part of a two-part work on the continued square roots of terms of magnitude $2$. He says, ``In what follows, I develop a general principle that yields a series of equations for the side and every diagonal of any regular polygon\index{polygons, regular}, always according to the same rule: it is the derivation of circular chords whose central angles are commensurate with each other and with the extended polygon[.] The theorem of angles and exterior angles of an isosceles triangle and the similarity of triangles are needed, so this theory is completely elementary. For any polygon, it gives a series of second-degree equations in which the sides and diagonals are of the same unknown; and this system of equations yields a common expression for \emph{all chords whose central angles are commensurable to the stretched one}: to my knowledge, such an expression \emph{is put forward here for the first time}: it is the ``periodic root of two'' (emphasis in the original). Bochow continues, ``[T]he method developed here is also useful for teaching{\ldots}for me there is no doubt that the theory of regular polygons could also be presented in this way in school; what is simpler scientifically is also simpler for teaching.'' The publication has a three-panel oversized sheet of figures inserted.


\item\vspace{9pt} \label{Boc1896} \dag(1896) Karl Bochow\index[ppl]{Bochow, Karl|textbf}, Eine einheitliche Theorie der regelm\"{a}ssigen Vielecke: Teil 2. In \emph{Sechster Jahresbericht \"{u}ber die St\"{a}dtische Realschule zu Magdeburg. Ostern 1895 bis Ostern 1896}. Progr. No. 279. E. Baensch jun., Magdeburg, 1--20.
\index{continued square roots!and trigonometric functions}
\index{continued square roots!of terms $a_n=\pm2$}
\index{continued square roots!geometric interpretations of}

{\footnotesize Source: \href{https://play.google.com/store/books/details/Karl_Bochow_Eine_einheitliche_Theorie_der_regelm\%C3\%A4s?id=2ACtBZ1ytjcC\&pli=1}{Google Play Store}.}

{\footnotesize Cited in 
\hyperref[Wie1904b]{\textsc{Wiernsberger 1904b}}\index[ppl]{Wiernsberger, Paul},
\hyperref[Boc1905]{\textsc{Bochow 1905}}\index[ppl]{Bochow, Karl}, and
\hyperref[Boc1910]{\textsc{Bochow 1910}}\index[ppl]{Bochow, Karl}.}

\vspace{4pt}

Part two of the work started in \hyperref[Boc1896]{\textsc{Bochow 1896}}. The continued square root expressions have up to four terms.


\item\vspace{9pt} \label{Ise1897} (1897) C. Isenkrahe\index[ppl]{Isenkrahe, C.|textbf}, \emph{Das Verfahren der Funktionswiederholung, seine Geometrische Veranschaulichung und Algebraische Anwendung}, B.~G. Teubner, Leipzig.
\index{iteration!of functions}
\index{iteration!graphical methods}

{\footnotesize Source: \href{https://books.google.com/books?id=PQRVAAAAYAAJ\&printsec=frontcover\&source=gbs_ge_summary_r\&cad=0\#v=onepage\&q\&f=false}{Google Books}.}

{\footnotesize Cited in
\hyperref[Hey1904]{\textsc{Heymann 1904}}\index[ppl]{Heymann, W.} and 
\hyperref[Gol1911]{\textsc{Goldziher 1911}}\index[ppl]{Goldziher, Karl}.}

\vspace{4pt}

At the beginning and close of \hyperref[Ise1888]{\textsc{Isenkrahe 1888}}, the author vowed to say more on the topic of functional iteration, and in this booklet of 114 pages he apparently fulfills that promise. It was published as a ``scientific supplement" to the annual report of the Kaiser-Wilhelm-Gymnasium, the high school in Trier, Germany, where Isenkrahe was a teacher and, eventually, director. Drawing upon his earlier paper, and citing a textbook by Eugen Netto\index[ppl]{Netto, Eugen} (\emph{Vorlesungen \"{u}ber Algebra, Vol. 1}, B.~G. Teubner, Leipzig, 1896) several times, the book gives a concise development of iterated functions as solutions to equations of the form $F(x)=0$. The geometrical interpretation of fixed-point iteration\index{fixed-point iteration} is emphasized in nearly 80 graphs.\footnote{The parallels in approach and subject matter between this book and \hyperref[Hey1894a]{\textsc{Heymann 1894a}}\index[ppl]{Heymann, W.} and \hyperref[Hey1894a]{\textsc{Heymann 1894b}} are intriguing. Indeed, Heymann in \hyperref[Hey1904]{\textsc{Heymann 1904}} makes a point of claiming priority for the aspects of graphical analysis found in Isenkrahe's book.}

An interesting aspect of Isenkrahe's approach is his handling of the associativity of infinite function compositions\index{associativity of function compositions}. Perhaps recognizing the difference between iterated\index{iterated compositions} and continued compositions\index{continued compositions}, he invents the terms \emph{Bruchkette}\index{Bruchkette@\emph{Bruchkette}} and \emph{Wurzelkette}\index{Wurzelkette@\emph{Wurzelkette}} as the iterated equivalents of \emph{Kettenbr\"{u}chen}\index{Kettenbr\"{u}chen@\emph{Kettenbr\"{u}chen}} (continued fractions) and \emph{Kettenwurzeln}\index{Kettenwurzeln@\emph{Kettenwurzeln}} (continued roots). However, the expression of the iterated/continued dichotomy in mathematical notation is not so clear. For certain values of $a_0$ and $a_1$, he develops what we might call an \emph{iterated square} solution to the quadratic equation $a_0+a_1x+x^2=0$:
\begin{equation}\label{E:Ise1897-1}
x=-\dfrac{a_0}{a_1}-\dfrac{1}{a_1}\left(-\dfrac{a_0}{a_1}\cdots-\dfrac{1}{a_1}\left(-\dfrac{a_0}{a_1}-\dfrac{1}{a_1}\left(-\dfrac{a_0}{a_1}\right)^2\right)^2\cdots\right)^2\;.
\end{equation}
Perhaps concerned that this looks too unconventional, he assures us that \eqref{E:Ise1897-1} ``becomes somewhat clearer when the powers are recorded as radicals with a fractional index," that is,
\begin{equation}\label{E:Ise1897-2}
x=\sqrt[1/2]{-\dfrac{a_0}{a_1}-\dfrac{1}{a_1}\sqrt[1/2]{-\dfrac{a_0}{a_1}\cdots-\dfrac{1}{a_1}\sqrt[1/2]{-\dfrac{a_0}{a_1}-\cdots}}}\;,
\end{equation}
where $\sqrt[1/2]{x}=x^2$. The trick of writing $\sqrt[m/n]{x}$ for $x^{n/m}$ is repeated later, when Isenkrahe reprises his 1888 analysis of \hyperref[Hof1881]{\textsc{Hoffmann 1881}} and the trinomial $x^m - px^n+q=0$.\index{trinomial equations}


\item\vspace{9pt} \label{Lef1897} (1897) B. Lefebvre, S. J.\index[ppl]{Lefebvre, B.|textbf}, \emph{Cours D\'{e}velopp\'{e} d'Alg\`{e}bre \'{E}l\'{e}mentaire pr\'{e}c\'{e}d\'{e} d'un Aper\c{c}u Historique sur les Origines des Math\'{e}matiques \'{E}l\'{e}mentaires et suivi d'un Recueil d'Exercises et de Probl\`{e}mes, Tome I, Calcul Alg\'{e}brique}, Namur Librairie Classique de Ad. Wesmael-Charlier, \'{E}diteur, Belgium.
\index{notation!surd $\surd$}
\index{continued square roots}
\index{history!of square root notation}

{\footnotesize Source: \href{https://gallica.bnf.fr/ark:/12148/bpt6k207267n/f355.item.r=B}{Gallica}.}

\vspace{4pt}

This school textbook includes some problems (numbers 13, 23, and 30 on pages 306, 308, and 309, respectively) derived from \hyperref[Cat1842]{\textsc{Catalan 1842}}, and an interesting example of a telescoping continued square root of 6 terms. Incidentally, a long footnote spanning pages 105 to 107 contains a surprisingly detailed history of square root notation. 


\item\vspace{9pt} \label{Hey1898} \dag(1898) W. Heymann, \"{U}ber die elementare Aufl\"{o}sung transcendenter Gleichungen. Mit Beitr\"{a}gen zur Ingenieur-Mathematik. \emph{Zeitschrift f\"{u}r mathematischen und naturwissenschaftlichen Unterricht} \textbf{29}, 1--15.
\index{iteration!of functions}
\index{iteration!graphical methods}
\index{trinomial equations}
\index{continued robb roots@continued $r$th roots!as solutions to polynomial equations}

{\footnotesize Source: \href{https://digital.slub-dresden.de/werkansicht/dlf/102439/29}{Dresden State and University Library (SLUB)}.}

{\footnotesize Cited in
\hyperref[Gol1911]{\textsc{Goldziher 1911}}\index[ppl]{Goldziher, Karl}.}

Heymann offers an overview of his graphical iteration methods (\hyperref[Hey1894a]{\textsc{Heymann 1894a}} and \hyperref[Hey1894b]{\textsc{1894a}}), focusing on solutions to practical problems. The second section, titled ``Infinite radicals'', alludes to \hyperref[Gun1880]{\textsc{G\"{u}nther 1880}}\index[ppl]{Gunther@G\"{u}nther, Siegmund} and \hyperref[Sch1880]{\textsc{Schaewen 1880}}\index[ppl]{Schaewen, Paul von} by way of introducing two ways to create an iteration scheme from the trinomial equation $x^n-x-c=0$:
\[x=\sqrt[n]{y}\;,\;\;y=c+x\;,\quad\text{and}\quad x=-c+y\;,\;\;y=x^n=\sqrt[1/n]{x}\;.\]
 Initialized at $x_0$, the first pair yields
\[x=\sqrt[n]{c+\sqrt[n]{c+\sqrt[n]{c+\cdots\sqrt[n]{c+x_0}}}}\]
and the second
\[x = -c +\sqrt[1/n]{-c+\sqrt[1/n]{-c+\sqrt[1/n]{-c+\cdots\sqrt[1/n]{x_0}}}}\;,\]
(compare \hyperref[Ise1897]{\textsc{Isenkrahe 1897}}\index[ppl]{Isenkrahe, C.} for a similar notation). The third and final section, ``Exercises'', offers six applications of Heymann's iterative solutions, including the Kepler problem (\hyperref[Kep1621]{\textsc{Kepler 1621}}) and a problem involving the catenary curve. The fifth exercise, to solve for the pressure value in a hydrodynamics problem, leads to the trinomial equation $x^5=ax+b$ and the expansions
\[x=\sqrt[5]{b+a\sqrt[5]{b}}\;,\quad x=\sqrt[5]{b+a\sqrt[5]{b+a\sqrt[5]{b}}}\;,\quad\text{etc.}\;.\]
Heymann closes with a reprise on the interest rate problem posed in\\ \hyperref[Gun1880]{\textsc{G\"{u}nther 1880}}\index[ppl]{Gunther@G\"{u}nther, Siegmund}.


\item\vspace{9pt} \label{Boc1899} (1899) Karl Bochow\index[ppl]{Bochow, Karl|textbf}, Problem 1740. \emph{Zeitschrift f\"{u}r mathematischen und naturwissenschaftlichen Unterricht} \textbf{30}, 24--28. Solution by the proposer, \emph{ibid.}, \textbf{31}, 191--192.
\index{continued square roots!and trigonometric functions}
\index{continued square roots!of terms $a_n=\pm2$}

{\footnotesize Problem source: \href{https://books.google.com/books?id=ut5LAAAAMAAJ&pg=PA27#v=onepage&q&f=false}{Google Books}.}

{\footnotesize Solution source: \href{https://books.google.com/books?id=fd9LAAAAMAAJ&newbks=1&newbks_redir=0&dq=\%22Zeitschrift\%20f\%C3\%BCr\%20mathematischen\%20und\%20naturwissenschaftlichen\%20Unterricht\%22\&pg=PA191\#v=onepage\&q\&f=false}{Google Books}.}

{\footnotesize Cited in
\hyperref[Boc1905]{\textsc{Bochow 1905}} and
\hyperref[Boc1910]{\textsc{Bochow 1910}}\index[ppl]{Bochow, Karl}.}

\vspace{4pt}

The problem reads roughly as follows: If $a$ is a number between 0 and $\tfrac{1}{2}$, one can write
\begin{equation}\label{E:Boc1899-1}
2\sin a\pi =\sqrt{2+\lambda_1\sqrt{2+\lambda_2\sqrt{2+\lambda_3\sqrt{2+\lambda_4\sqrt{2+\ldots\textrm{in inf.}}}}}}
\end{equation}
and at the same time $a$ can be written as a series\index{real numbers!binary representations of}
\begin{equation}\label{E:Boc1899-2}
a = \dfrac{1}{2^2}+\dfrac{\lambda_1}{2^3}+\dfrac{\lambda_1\lambda_2}{2^4}+\dfrac{\lambda_1\lambda_2\lambda_3}{2^5}+\dfrac{\lambda_1\lambda_2\lambda_3\lambda_4}{2^6}+\cdots
\end{equation}
in which each $\lambda_i$ is either $+1$ or $-1$.

Similarly, if every $\lambda_i$ is equal to $\pm1$, 
\[a = \dfrac{1}{2^2}+\dfrac{\lambda_1}{2^3}+\dfrac{\lambda_2}{2^4}+\dfrac{\lambda_3}{2^5}+\cdots\]
and if, simultaneously,
\[2\sin a\pi = \sqrt{2+\lambda_1\sqrt{2+\lambda_1\lambda_2\sqrt{2+\lambda_2\lambda_3\sqrt{2+\cdots}}}}\;,\]
then the double sine, the continued square root\index{continued square roots}, the value of $a$, and its series expansion ``determine each other unambiguously so that this interrelation contains the general solution of the problem: for a given angle to find its [sine] function, and for a given function value its corresponding angle." 


\item\vspace{9pt} \label{Cwo1899} (1899) Kasimir Cwojdzinski\index[ppl]{Cwojdzinski, Kasimir|textbf}, Kettenwurzeln. \emph{Archiv der Mathematik und Physik} \textbf{17}, 29--35.
\index{Kettenwurzeln@\emph{Kettenwurzeln}}
\index{trinomial equations}

{\footnotesize Source: \href{https://archive.org/details/archivdermathem82unkngoog/page/n42/mode/1up}{Internet Archive}.}

\vspace{4pt}

The paper begins with a discussion remarkably similar to \hyperref[Dix1878]{\textsc{Dixon 1878}}\index[ppl]{Dixon, T. S. E.}, but less compressed, giving many of the same formal continued $r$th root constructions of solutions to trinomial equations\index{trinomial equations}\index{continued robb roots@continued $r$th roots!as solutions to polynomial equations} in one variable. The author defines the continued $n$th root
\[\pm a\sqrt[n]{\pm b\pm a\sqrt[n]{\pm b\pm a\sqrt[n]{\pm b \pm\ldots}}}\;,\]
as the limit of the sequence of approximants using the forward recurrence relation\index{recurrence relation}
\[w_{n+1}=\pm a\sqrt[n]{\pm b\pm w_n}\;.\]
He then discusses some particular cases. When $w_{n+1}=a\sqrt[n]{b-w_n}$, he argues that the even and odd approximants respectively increase and decrease toward a limit, provided that $b>a\sqrt[n]{b}$. The cases $w_{n+1}=a\sqrt[n]{b+w_n}$ and $w_{n+1}=-a\sqrt[n]{-b-w_n}$ are claimed respectively to increase and decrease towards their limits. However, although more is offered than in \hyperref[Dix1878]{\textsc{Dixon 1878}}\index[ppl]{Dixon, T. S. E.} to justify the formal expansions, the arguments for convergence are not rigorous. The author gives a numerical evaluation of 
\[\sqrt[\pi/\psi]{e-\sqrt[\pi/\psi]{e-\cdots}}\]
where $\psi=\tfrac{1}{2}(\sqrt{5}-1)$, but his list of calculated approximants is incorrect. The paper is plagued with nearly two dozen typographical errors in its seven pages.


\item\vspace{9pt} \label{Can1901} \dag(1901) Moritz Cantor\index[ppl]{Cantor, Moritz|textbf}, \emph{Vorlesungen \"{U}ber Geschichte der Mathematik}, Vol. 3., from 1668--1758. Second Edition. B. G. Teubner, Leipzig.
\index{successive substitution!17th century examples}
\index{history!of successive substitution}

{\footnotesize Source: \href{https://books.google.com/books/about/Vorlesungen_\%C3\%BCber_Geschichte_der_Mathema.html?id=9SA3ZSOFMWMC}{Google Books}.}

{\footnotesize Cited in 
\hyperref[Hey1904]{\textsc{Heymann 1904}}\index[ppl]{Heymann, W.}.}

\vspace{4pt}

In discussing Bernoulli's work on infinite series, Cantor mentions \hyperref[Ber1692]{\textsc{Bernoulli 1692}}\index[ppl]{Bernoulli, Jacob} on page 96, and reprints Bernoulli's expressions
\[\sqrt{a\sqrt{a\sqrt{a\sqrt{\ldots}}}}\quad\text{and}\quad\sqrt{a+\sqrt{a+\sqrt{a+\ldots}}}\;.\]


\item\vspace{9pt} \label{Hey1901} (1901) W. Heymann\index[ppl]{Heymann, W.|textbf}, \"{U}ber Wurzelgruppen, welche durch Uml\"{a}ufe ausgeschnitten werden. \emph{Zeitschrift f\"{u}r Mathematik und Physik} \textbf{46}, 265--299.
\index{iteration!of functions}
\index{iteration!graphical methods}

{\footnotesize Source: \href{https://archive.org/details/zeitschriftfrma12runggoog/page/265/mode/1up}{Internet Archive}.}

\vspace{4pt}

Continuing his investigations of fixed-point iteration from the previous decade, Heymann considers the question of periodic points\index{periodic points}, especially points of periods 2 and 3. Many examples are exhibited, including the continued square root\index{continued square roots!of periodic real terms}
\[\sqrt{a+\sqrt{b+\sqrt{a+\sqrt{b+\cdots}}}}\;.\]
He also applies his graphical theory to continued fractions and geometric series.


\item\vspace{9pt} \label{Boc1902} (1902) Karl Bochow\index[ppl]{Bochow, Karl|textbf}, Zur Behandlung der regelm\"{a}ssigen Vielecke. \emph{Unterrichtsbl\"{a}tter f\"{u}r Mathematik und Naturwissenschaften} \textbf{8}(5), 109--113.
\index{continued square roots!and trigonometric functions}
\index{continued square roots!of terms $a_n=\pm2$}
\index{continued square roots!geometric interpretations of}

{\footnotesize Source: \href{https://babel.hathitrust.org/cgi/pt?id=njp.32101079819924&seq=389}{HathiTrust}.}

{\footnotesize Cited in
\hyperref[Wie1904b]{\textsc{Wiernsberger 1904b}}\index[ppl]{Wiernsberger, Paul}, 
\hyperref[Boc1904]{\textsc{Bochow 1904}}\index[ppl]{Bochow, Karl}, and
\hyperref[Boc1905]{\textsc{Bochow 1905}}\index[ppl]{Bochow, Karl}.}

\vspace{4pt}

This paper and \hyperref[Boc1904]{\textsc{Bochow 1904}} comprise Bochow's lecture at the Eleventh Annual General Meeting of the Association for the Advancement of Teaching in Mathematics and Natural Sciences at Dusseldorf during Whitsun Week 1902. The topic is the calculation of side lengths of regular polygons\index{polygons, regular}, resulting in the continued square root expansions\index{continued square roots!and trigonometric functions}\index{continued square roots!of terms $a_n=\pm2$}
\begin{align}\label{E:Boc1902-1}
2\cos\left(\dfrac{p}{2^q}\pi\right)&=\sqrt{2\pm\sqrt{2\pm\sqrt{2\pm\cdots\pm\sqrt{2\pm{\sqrt{2}}}}}}\\
2\cos\left(\dfrac{p}{3\cdot2^q}\pi\right)&=\sqrt{2\pm\sqrt{2\pm\sqrt{2\pm\cdots\pm\sqrt{2\pm{\sqrt{3}}}}}}\;,
\end{align}
where $p$, $q$, and the periodic sequence of the signs are given in a table. (Compare \hyperref[Luc1878]{\textsc{Lucas 1878}} and \hyperref[Boc1899]{\textsc{Bochow 1899}}).

Bochow's lecture was followed by a discussion involving three other teachers. Notes on this discussion, written up by J. E. B\"{o}ttcher\index[ppl]{Bottcher@B\"{o}ttcher, J. E.} of Leipzig, offer some suggestions on how Bochow's methods could be implemented in the classroom.


\item\vspace{9pt} \label{Wie1903} (1903) Paul Wiernsberger, \index[ppl]{Wiernsberger, Paul|textbf} Convergence des radicaux superpos\'{e}s p\'{e}riodiques. \emph{Comptes rendus hebdomadaires des s\'{e}ances de l'Acad\'{e}mie des Sciences} \textbf{137}, 1233--1234.
\index{continued square roots!of terms $a_n=\pm2$}
\index{continued square roots!and trigonometric functions}
\index{radicaux superpos\'{e}s@\emph{radicaux superpos\'{e}s}}
\index{continued square roots!geometric interpretations of}

{\footnotesize Source: \href{https://gallica.bnf.fr/ark:/12148/bpt6k62622724/f37.item}{Gallica}.} 

{\footnotesize Cited in 
\hyperref[Wie1904a]{\textsc{Wiernsberger 1904a}} and \hyperref[Wie1904b]{\textsc{1904b}}.}

\vspace{4pt}

This page-and-a-half research report exhibits, with little motivation or explanation, the continued radical $\sqrt{2\pm\sqrt{2\pm\ldots\pm\sqrt{2}}}$ (having signs in periodic sequence) as it pertains to the lengths of sides of regular polygons\index{polygons, regular} in the unit circle. The material is expanded upon in \hyperref[Wie1904b]{\textsc{Wiernsberger 1904b}}. 


\item\vspace{9pt} \label{Boc1904} (1904) Karl Bochow\index[ppl]{Bochow, Karl|textbf}, Zur Behandlung der regelm\"{a}ssigen Vielecke. (Fortsetzung). Das regul\"{a}re Achtzehneck als Beispiel f\"{u}r ein N\"{a}herungsverfahren zur Konstruktion und zur Berechnung von regelm\"{a}ssigen Vielecken und von Winkelfunktionen. \emph{Unterrichtsbl\"{a}tter f\"{u}r Mathematik und Naturwissenschaften} \textbf{10}(1), 12--16.
\index{continued square roots!of terms $a_n=\pm2$}
\index{continued square roots!and trigonometric functions}

{\footnotesize Source: \href{https://books.google.com/books?id=_bs7AQAAMAAJ&pg=PA12}{Google Books}.}

{\footnotesize Cited in
\hyperref[Wie1904b]{\textsc{Wiernsberger 1904b}}\index[ppl]{Wiernsberger, Paul},
\hyperref[Boc1905]{\textsc{Bochow 1905}}\index[ppl]{Bochow, Karl}, and 
\hyperref[Boc1910]{\textsc{Bochow 1910}}\index[ppl]{Bochow, Karl}.}

\vspace{4pt}

In this continuation of the lecture \hyperref[Boc1902]{\textsc{Bochow 1902}}, the author works many examples demonstrating continued square root expansions of trigonometric functions in a regular 18-gon.


\item\vspace{9pt} \label{Hey1904} \dag(1904) W. Heymann\index[ppl]{Heymann, W.|textbf}, \"{U}ber die Aufl\"{o}sung der Gleichungen durch Iteration auf geometrischer Grundlage. \emph{Jahresbericht der technischen Staatslehranstalten in Chemnitz f\"{u}r die Zeit von Ostern 1903 bis Ostern 1904}, J. C. F. Pickenhahn \& S., Chemnitz.
\index{iteration!of functions}
\index{iteration!graphical methods}

{\footnotesize Cited in 
\hyperref[Gol1911]{\textsc{Goldziher 1911}}\index[ppl]{Goldziher, Karl}.}

\vspace{4pt}

In 45 pages, Heymann gives a further recounting of his graphics-oriented work on iterated functions, perhaps for the benefit of the audiences of the State Technical School in Chemnitz\footnote{Now the Technische Universit\"{a}t Chemnitz\index{Technische Universit\"{a}t Chemnitz}.} (where he was a professor), rather than for mathematicians per se. About half of the paper is given to the working of examples of root-finding by iteration, most of them seen in his previous work.

A valuable inclusion in this paper is Section 6 of Part II, in which Heymann provides authors' names and some citations for previous work involving iteration and successive substitution. Regarding the ``geometric interpretation of the iteration process,'' Heymann says, ``[h]ere I could only mention Legendre [\hyperref[Leg1816]{\textsc{Legendre 1816}}\index[ppl]{Legendre, Adrien-Marie}], in whose work on number theory (appendix to the second volume) the iteration method is applied to numerical equations and the analytical process is traced using intersecting curves. However, one will find that the geometric aspect recedes considerably in Legendre's work; finer distinctions are entirely absent and probably not intended.'' Heymann then continues, ``A book by Isenkrahe [\hyperref[Ise1897]{\textsc{Isenkrahe 1897}}\index[ppl]{Isenkrahe, C.}] also includes geometric considerations\ldots. To maintain my priority, insofar as it extends to geometric theory, I would therefore like to expressly emphasize that my relevant works were already published in 1894.''

\hyperref[Gol1911]{\textsc{Goldziher 1911}}\index[ppl]{Goldziher, Karl} states that the method of iterated functions was ``applied to the pension problem\index{pension problem} in the 1880s [e.g. in \hyperref[Gun1880]{\textsc{G\"{u}nther 1880}}\index[ppl]{Gunther@G\"{u}nther, Siegmund}], and after investigations of the process by Netto\index[ppl]{Netto, Eugen}, Isenkrahe\index[ppl]{Isenkrahe, C.}, and Heymann, Heymann\index[ppl]{Heymann, W.} [in the work cited here and in \hyperref[Hey1894a]{\textsc{Heymann 1894a}}\index[ppl]{Heymann, W.}] succeeded in giving the final formula for the practical case.''


\item\vspace{9pt} \label{Wie1904a} (1904a) Paul Wiernsberger,\index[ppl]{Wiernsberger, Paul|textbf} Sur les expressions form\'{e}es de radicaux superpos\'{e}s. \emph{Comptes rendus hebdomadaires des s\'{e}ances de l'Acad\'{e}mie des sciences.} \textbf{138}, 1401--1403.
\index{continued square roots!geometric interpretations of}
\index{continued square roots!of terms $a_n=\pm2$}
\index{continued square roots!and trigonometric functions}
\index{radicaux superpos\'{e}s@\emph{radicaux superpos\'{e}s}}

{\footnotesize Source: \href{https://gallica.bnf.fr/ark:/12148/bpt6k3092p/f1491.item}{Gallica}.} 

{\footnotesize Cited in 
\hyperref[Wie1904b]{\textsc{Wiernsberger 1904b}}\index[ppl]{Wiernsberger, Paul}.}

\vspace{4pt}

Revealing more of the work that will eventually comprise \hyperref[Wie1904b]{\textsc{Wiernsberger 1904b}}, the author looks not only at 
\begin{equation}\label{E:Wie1904a-1}
\sqrt{2+\epsilon_1\sqrt{2+\epsilon_3\sqrt{2+\ldots}}}\;,\quad |\epsilon_h|=1\;,
\end{equation}
but also discusses convergence criteria for the more general\index{continued robb roots@continued $r$th roots!convergence conditions for}\index{continued robb roots@continued $r$th roots!of arbitrary nonnegative real terms}
\begin{equation}\label{E:Wie1904a-2}
\sqrt[n]{a_1+\epsilon_1\sqrt[n]{a_2+\cdots+\epsilon_{h-1}\sqrt[n]{a_h+\cdots}}}\;,\quad |\epsilon_h|=1\;.
\end{equation}

Of particular note is this result: ``If all $\epsilon_h$ are positive, the necessary and sufficient condition for the expression \eqref{E:Wie1904a-2} to be convergent is that the numbers $a_h$ all satisfy the inequality
\begin{equation}\label{E:Wie1904a-3}
a_h^{n^h}<A\quad(h=1,2,3,\ldots)\;,
\end{equation}
where $A$ is a finite number." (On page 1402, equation \eqref{E:Wie1904a-3} above is misprinted as $a^{n^h}<A$, with no subscript on the $a$. Other misprints include a missing $n$ in a nest of $n$th roots on page 1403.) As of this writing, Wiernsberger's is the earliest statement of this important convergence condition, predating \hyperref[Pol1916]{\textsc{P\'{o}lya 1916}}, Vijayaraghavan\index[ppl]{Vijayaraghavan, T.} in \hyperref[Ram1927]{\textsc{Ramanujan 1927}}\index[ppl]{Ramanujan, Srinivasa}, and \hyperref[Her1935]{\textsc{Herschfeld 1935}}\index[ppl]{Herschfeld, Aaron} (all of which more rigorously invoke the $\limsup$). Wiernsberger also attempts convergence criteria for the case when the $\epsilon_h$ are not all positive, but with less definitive results. All the proofs are sketched.


\item\vspace{9pt} \label{Wie1904b} (1904b) Paul Wiernsberger,\index[ppl]{Wiernsberger, Paul|textbf} \emph{Recherches diverses sur des polygones r\'{e}guliers et les radicaux superpos\'{e}s}, A. Rey, Lyons.

{\footnotesize Cited in 
\hyperref[Boc1905]{\textsc{Bochow 1905}}\index[ppl]{Bochow, Karl} and
\hyperref[Boc1910]{\textsc{Bochow 1910}}\index[ppl]{Bochow, Karl}.}
\index{continued square roots!of terms $a_n=\pm2$}
\index{continued square roots!and trigonometric functions}
\index{continued square roots!of arbitrary nonnegative real terms}
\index{radicaux superpos\'{e}s@\emph{radicaux superpos\'{e}s}}
\index{continued square roots!geometric interpretations of}

\vspace{4pt}

Wiernsberger's\index[ppl]{Wiernsberger, Paul} doctoral thesis gives an overview of 19th-century geometry of regular polygons\index{polygons, regular}, leading to detailed calculations involving the continued square root
\[\sqrt{2\pm\sqrt{2\pm\ldots\pm\sqrt{2}}}\;,\]
where the signs occur in periodic sequence. Among the trigonometric formulas eventually derived are
\begin{align}
\sin\tfrac{\pi}{2}x=&\tfrac{1}{2}\sqrt{2+\epsilon_1\sqrt{2+\cdots+\epsilon_h\sqrt{2+\cdots}}}\label{E:Wie1904b-1}\\
\intertext{and}
\cos\tfrac{\pi}{2}x=&\tfrac{1}{2}\sqrt{2-\epsilon_1\sqrt{2+\cdots+\epsilon_h\sqrt{2+\cdots}}}\;,\label{E:Wie1904b-2}
\end{align}
where the $\epsilon_i=\pm 1$ form a sequence unique to $x$ (compare \hyperref[Boc1899]{\textsc{Bochow 1899}}\index[ppl]{Bochow, Karl}, which Wiernsberger cites).

In the book's last chapter, Wiernsberger undertakes what is, as of this writing, the first rigorous discussion of convergence for the general continued $n$th root $\sqrt[n]{a_1+\sqrt[n]{a_2+\cdots+\sqrt[n]{a_h+\cdots}}}$, where the $a_h$ are positive, $n$ is greater than or equal to 2, and the positive root is always taken.\index{continued robb roots@continued $r$th roots!convergence conditions for}\index{continued robb roots@continued $r$th roots!of arbitrary nonnegative real terms} Here the typo in the convergence condition \eqref{E:Wie1904a-3} from \hyperref[Wie1904a]{\textsc{Wiernsberger 1904a}} is corrected; the proof is essentially the same as that given over thirty years later in \hyperref[Her1935]{\textsc{Herschfeld 1935}}\index[ppl]{Herschfeld, Aaron}. Wiernsberger independently rediscovers Bernoulli's\index[ppl]{Bernoulli, Jacob} limit \eqref{E:Ber1692-1} for a continued square root of constant terms $a>0$\index{continued square roots!of terms $a_n=c$}, and investigates special cases in which one negative square root appears in the infinite nest; as an example, Wiernsberger\index[ppl]{Wiernsberger, Paul} uses equation \eqref{E:Ber1692-1} above to show that $\sqrt{a-\sqrt{a+\sqrt{a+\cdots}}}$ converges for $a>(5+2\sqrt{5})/4$.

The bibliography contains 105 items in chronological order, most of which concern topics related to the properties of regular polygons\index{polygons, regular} inscribed in and/or circumscribed around a circle. The earliest source is Section VII of Gauss's\index[ppl]{Gauss, Carl Friedrich} \emph{Disquisitiones arithmetic\ae} (1801), while the latest is \hyperref[Wie1904a]{\textsc{Wiernsberger 1904a}}. In connection with continued square roots\index{continued square roots}, the author cites \hyperref[Cat1842]{\textsc{Catalan 1842}}, \hyperref[Did1872]{\textsc{Didion 1872}}, and several articles by Karl Bochow\index[ppl]{Bochow, Karl}; however, \hyperref[Vie1593]{\textsc{Vi\`{e}te 1593}}\index[ppl]{Vi\`{e}te, Fran\c{c}ois} is absent.


\item\vspace{9pt} \label{Boc1905} \dag(1905) Karl Bochow\index[ppl]{Bochow, Karl|textbf}, Die Funktionen rationaler Winkel, besonders die numerische Berechnung der Winkelfunktionen ohne Benutzung der\\trigonometrischen Reihen und ohne Kenntnis der Zahl $\pi$. In \emph{Funfzehnter Jahresbericht
\"{u}ber die St\"{a}dtische Realschule zu Magdeburg. Ostern 1904 bis Ostern 1905.} Progr. No. 317. E. Baensch, Magdeburg. 
\index{continued square roots!and trigonometric functions}
\index{continued square roots!of terms $a_n=\pm2$}
\index{continued square roots!geometric interpretations of}

{\footnotesize Source: \href{https://digital.ulb.hhu.de/ulbdsp/periodical/pageview/5782226}{University and State Library of D\"{u}sseldorf}.}
 
{\footnotesize Cited in 
\hyperref[Boc1910]{\textsc{Bochow 1910}}\index[ppl]{Bochow, Karl} and
\hyperref[Sze1930]{\textsc{Szeg\"{o} 1930}}\index[ppl]{Szeg\"{o}, G.}.}

\vspace{4pt}

Like \hyperref[Boc1905]{\textsc{Bochow 1895}} and \hyperref[Boc1896]{\textsc{Bochow 1896}}, this paper was published in the annual report of the Municipal Secondary School in Magdeburg, Germany, where the author was a senior teacher. Bochow begins by listing his published work on trigonometrically-derived continued square roots going back a dozen years; but this leads to something of a lament. ``I attempted to ground geometry in the concept of direction. I succeeded in defining the trigonometric functions before the parallels, before the triangle; I had recognized the dyadic continued root as a sufficient expression of them. {\ldots}Various circumstances prevented me from publishing these results, most of which I have possessed since the early 1890s. One journal found the material too elementary, another kept the manuscripts for years after promising prompt publication, and I lacked the time due to extensive other work. In 1904, I intended to use my school's curriculum to present the most important theorems on continued roots in a coherent manner, but I had to postpone this. Then, in February 1904, I received a letter from a mathematician living in Lyon, Mr. Paul Wiernsberger, informing me that he had read my `Unified Theory' [\hyperref[Boc1905]{\textsc{Bochow 1895}} and \hyperref[Boc1896]{\textsc{Bochow 1896}}] and continued working on the subject. I informed the gentleman of my works that had since been published and now resolved to contribute further material to this program. Mr. Wiernsberger has beaten me to it. His work [\hyperref[Wie1904b]{\textsc{Wiernsberger 1904b}}\index[ppl]{Wiernsberger, Paul}], sent to me in January 1905, contains much that I have long possessed but not published, but it also surpasses me in many important points. However, I think that the following elementary and goniometrically oriented presentation will still be able to stand alongside Wiernsberger's work, since the starting point and goal are different; my geometric derivations are also useful in teaching.''\footnote{Bochow goes on: ``To my regret, I did not inform Mr. Wiernsberger at the time that I had already published the general expression for the relationship between the continued root and the trigonometric function in 1899 [\hyperref[Boc1899]{\textsc{Bochow 1899}}\index[ppl]{Bochow, Karl}]. I had completely forgotten that it was printed in the problem repertory of Hoffmann's journal. Thus, this information is also missing from the extraordinarily detailed bibliography that Mr. Wiernsberger includes with his work, through my fault.''}

In the article's 40 pages, Bochow lays out his system, with numerous examples and tables.


\item\vspace{9pt} \label{Wie1905} (1905) Paul Wiernsberger,\index[ppl]{Wiernsberger, Paul|textbf} Sur les polygones r\'{e}guliers et les radicaux carr\'{e}s superpos\'{e}s. \emph{Journal f\"{u}r die reine und ungewandte Mathematik} \textbf{130}, 144--152.  MR1580680.
\index{continued square roots!and trigonometric functions}
\index{continued square roots!of terms $a_n=\pm2$}
\index{radicaux superpos\'{e}s@\emph{radicaux superpos\'{e}s}}
\index{continued square roots!geometric interpretations of}

{\footnotesize Source: \href{https://doi.org/10.1515/crll.1905.130.144}{doi.org (De Gruyter)}.} 

\vspace{4pt}

This is an overview and extension of \hyperref[Wie1904b]{\textsc{Wiernsberger 1904b}}. In addition to developing the formulas \eqref{E:Wie1904b-1}, \eqref{E:Wie1904b-2}, and similar expressions for $\sin\pi x$ and $\cos\pi x$,\index{continued square roots!and trigonometric functions} Wiernsberger\index[ppl]{Wiernsberger, Paul} derives an expression for $\pi$ as the infinite product of continued square roots\index{infinite products!of continued square roots}\index{pi@$\pi\;(3.14159\ldots)$!infinite products for}\index{constants, named!pi@$\pi\;(3.14159\ldots)$}
\begin{align}\label{E:Wie1905-1}
\pi=&\dfrac{2^{k-1}}{a}\times C_{2^k,a}\times\\
&\lim_{h=\infty}\left\{\dfrac{2}{C_{2^{k+1},2^k-a}}\cdot\dfrac{2}{C_{2^{k+2},2^{k+1}-a}}\cdots\dfrac{2}{C_{2^h,2^{h-1}-a}}\right\}\;,
\end{align}
where $C_{m,n}$ is a continued square root\index{continued square roots!of terms $a_n=\pm2$} representation of $2\sin\tfrac{n\pi}{m}$ and $4a\le2^{k+1}<2^h$; for $a=1$ this yields the \hyperref[Cat1842]{\textsc{Catalan 1842}} reciprocal version of Vi\`{e}te's\index[ppl]{Vi\`{e}te, Fran\c{c}ois} formula \eqref{E:Vie1593-1} above.\index{Vi\`{e}te's formula for $\tfrac{2}{\pi}$}


\item\vspace{9pt} \label{Ala1907} (1907) C. Alasia, Question 878. \emph{Supplemento al Periodico di Matematica} \textbf{11}(1), 16. Solution by the proposer, \emph{ibid.} \textbf{11}(3), 46--47.

{\footnotesize Cited in
\hyperref[Cip1908]{\textsc{Cipolla 1908}}\index[ppl]{Cipolla, Michele}.}
\index{continued square roots!and trigonometric functions}
\index{continued square roots!of terms $a_n=\pm2$}

\vspace{4pt}

The problem reads roughly as follows: Prove that 
\[\sin(45^{\circ}\pm x)=\sqrt{\dfrac{1\pm\sin2x}{2}}\] 
and deduce the relation
\[2\cos\left(60^{\circ}\pm\dfrac{30}{2^n}\right)=\sqrt{2-\sqrt{2-\sqrt{2-\cdots-\sqrt{2}}}}\;,\]
where $2$ is repeated $n$ times on the right, and where the $+$ sign or $-$ sign on the left must be taken accordingly as $n$ is odd or even. Of the 14 solutions received, 3 were incomplete and the rest were deemed ``too cumbersome."


\item\vspace{9pt} \label{Cip1908} (1908) Michele Cipolla\index[ppl]{Cipolla, Michele|textbf}, Intorno ad un radicale continuo. \emph{Periodico di Matematica per l'insegnamento Secondario,} Series 3, \textbf{5}, 179--185.  
\index{continued square roots!and trigonometric functions}
\index{continued square roots!of terms $a_n=\pm2$}

{\footnotesize Source: \href{https://books.google.com/books?id=nFakZLlfAU8C&newbks=1&newbks_redir=0&dq=\%22Intorno\%20ad\%20un\%20radicale\%20continuo\%22\&pg=PA179\#v=onepage\&q\&f=false}{Google Books}.}

{\footnotesize Cited in 
\hyperref[Her1935]{\textsc{Herschfeld 1935}}\index[ppl]{Herschfeld, Aaron}, 
\hyperref[Jon2008]{\textsc{Jones 2008}}\index[ppl]{Jones, Dixon J.}, 
\hyperref[Mor2012]{\textsc{Moreno and Garc\'{i}a-Caballero 2012}}\index[ppl]{Moreno, Samuel G\'{o}mez}\index[ppl]{Garc\'{i}a-Caballero, Esther M.}, 
\hyperref[Gar2014b]{\textsc{Garc\'{i}a-Caballero, Moreno, and Prophet 2014b}}\index[ppl]{Garc\'{i}a-Caballero, Esther M.}\index[ppl]{Moreno, Samuel G\'{o}mez}\index[ppl]{Prophet, Michael P.}, and
\hyperref[Vel2016c]{\textsc{Vellucci and Bersani 2016c}}\index[ppl]{Vellucci, Pierluigi}\index[ppl]{Bersani, Alberto Maria}.}

\vspace{4pt}

Uses the cosine half-angle formula and induction to prove: If $i_j$ is $+1$ or $-1$, and 
\[\epsilon_{j,n-1}=\dfrac{1-i_ji_{j+1}\cdots i_{n-1}}{2}\]
for $j=1,\ldots,n-1$, then 
\begin{multline*}
\qquad\sqrt{2+i_{n-1}\sqrt{2+i_{n-2}\sqrt{2+\cdots+i_1\sqrt{2}}}}\\
=2\cos (1+2\epsilon_{1,n-1}+2^2\epsilon_{2,n-1}+\cdots+2^{n-1}\epsilon_{n-1,n-1})\dfrac{\pi}{2^n+1}\;;
\end{multline*}
the left side is a ``left continued radical" (\emph{radicale continuo a sinistra}, or what we are calling an iterated square root)\index{continued radicals!left}\index{radicale continuo a sinistra@\emph{radicale continuo a sinistra}} in the limit as $n\to\infty$. A similar formulation using $i_{n-j}$ instead of $i_j$ gives a ``right continued radical" (\emph{radicale continuo a destra}, or a continued square root)\index{infinite radicals!right}\index{radicale continuo a destra@\emph{radicale continuo a destra}} in the limit: Setting 
\[\lambda_j = \dfrac{1-i_1i_2\ldots i_j}{2}\]
for $j=1,\ldots,n-1$, one has
\begin{multline}\label{E:Cip1908-1}
\qquad\sqrt{2+i_{1}\sqrt{2+i_{2}\sqrt{2+\cdots+i_{n-1}\sqrt{2}}}}\\
=2\cos (2^{n-1}\lambda_1+2^{n-2}\lambda_2+\cdots+2\lambda_{n-1}+1)\dfrac{\pi}{2^n+1}\;.
\end{multline}
The paper gives examples using various periodic sequences $\{i_j\}$ in the finite and infinite cases, including\index{continued square roots!of periodic real terms}\index{continued square roots!of terms $a_n=\pm2$}
\[\sqrt{2-\sqrt{2-\sqrt{2-\cdots}}}=1 =\cdots-\sqrt{2-\sqrt{2-\sqrt{2}}}\]
and
\[\sqrt{2-\sqrt{2+\sqrt{2-\sqrt{2+\cdots}}}}=\dfrac{\sqrt{5}-1}{2}\;.\]
Compare the work of Bochow\index[ppl]{Bochow, Karl} and Wiernsberger\index[ppl]{Wiernsberger, Paul}. Cipolla makes a clearer distinction between iterated and continued function compositions (\S\ref{S:assoc}) than \hyperref[Ise1897]{\textsc{Isenkrahe 1897}}.


\item\vspace{9pt} \label{Can1908} (1908) G. Candido,\index[ppl]{Candido, G.|textbf} Sul numero $\pi$. \emph{Supplemento al Periodico di Matematica} \textbf{11}, 113--115.
\index{pi@$\pi\;(3.14159\ldots)$!infinite products for}\index{constants, named!pi@$\pi\;(3.14159\ldots)$}
\index{continued square roots!of terms $a_n=\pm2$}

{\footnotesize Cited in 
\hyperref[Her1935]{\textsc{Herschfeld 1935}}\index[ppl]{Herschfeld, Aaron}.}

\vspace{4pt}

Uses trig identities\index{continued square roots!and trigonometric functions} to independently develop several formulas already in existence at the time, including the reciprocal of equation \eqref{E:Vie1593-2} above (Vi\`{e}te\index[ppl]{Vi\`{e}te, Fran\c{c}ois} is not mentioned; \hyperref[Cat1842]{\textsc{Catalan 1842}} has the same result\index{Vi\`{e}te's formula for $\tfrac{2}{\pi}$}), and equation \eqref{E:Pie1891-1}. The paper concludes with the formula
\[\dfrac{\pi}{2}=\lim_{n\to\infty} 2^{n-1}\dfrac{\sqrt{2-\sqrt{2+\cdots+\sqrt{2}}}}{\sqrt{2+\sqrt{2+\cdots+\sqrt{2}}}}\;,\]
where there are $n-1$ terms each in the numerator and denominator; this was essentially derived in \hyperref[Did1872]{\textsc{Didion 1872}}. 


\item\vspace{9pt} \label{Boc1910} (1910) Karl Bochow\index[ppl]{Bochow, Karl|textbf}, Kettenwurzeln und Winkelfunktionen. \emph{Zeitschrift f\"{u}r mathematischen und naturwissenschaftlichen Unterricht} \textbf{41}, 161--186.
\index{Kettenwurzeln@\emph{Kettenwurzeln}}
\index{continued square roots!of terms $a_n=\pm2$}
\index{continued square roots!and trigonometric functions}
\index{continued square roots!geometric interpretations of}

{\footnotesize Source: \href{https://books.google.com/books?id=AuPuAAAAMAAJ&lpg=PA161&dq=\%22Kettenwurzeln\%20und\%20Winkelfunktionen\%22&pg=PA161#v=onepage&q&f=false}{Google Books}.}  

{\footnotesize Cited in
\hyperref[Sze1930]{\textsc{Szeg\"{o} 1930}}\index[ppl]{Szeg\"{o}, G.}, 
\hyperref[Her1935]{\textsc{Herschfeld 1935}}\index[ppl]{Herschfeld, Aaron}, and 
\hyperref[Kom1951]{\textsc{Kommerell 1951}}\index[ppl]{Kommerell, Karl}.}

\vspace{4pt}

The author's intent is to ``present to a wider readership a lecture that I gave some time ago at the Mathematical Society in Halle on this subject, in the hope of arousing interest in the matter.'' Bochow gives a lengthy overview of his work over two decades on continued square roots with terms $\pm 2$ arising from trigonometric identities, interpreted geometrically as lengths of chords in a circle; see \hyperref[Boc1895]{\textsc{Bochow 1895}}, \hyperref[Boc1896]{\textsc{Bochow 1896}}, \hyperref[Boc1899]{\textsc{1899}}, \hyperref[Boc1902]{\textsc{1902}}, \hyperref[Boc1904]{\textsc{1904,}} and \hyperref[Boc1905]{\textsc{1905}}.\index{continued square roots!of terms $a_n=\pm2$} \hyperref[Wie1904b]{\textsc{Wiernsberger 1904b}} is also cited as a principal reference. Much of the paper is devoted to computing examples of finite continued square roots having up to six terms. 


\item\vspace{9pt} \label{Bos1910} (1910) H. Bosmans, S. J.\index[ppl]{Bosmans, H., S.J.|textbf}, Un \'{e}mule de Vi\`{e}te: Ludolphe van Ceulen. Analyse de son ``Trait\'{e} du cercle." \emph{Annales de la soci\'{e}t\'{e} scientifique de Bruxelles}, \textbf{34}, second partie, 88--139.
\index{history!of continued square roots}  

{\footnotesize Source: \href{https://www.biodiversitylibrary.org/item/158901#page/254/mode/2up}{Biodiversity Heritage Library}.}

\vspace{4pt}

In this valuable overview and explication of \hyperref[Ceu1596]{\textsc{Van Ceulen 1596}}\index[ppl]{Ceulen, Ludolph van}, Bosmans shows in modern language and notation the geometric and trigonometric tools Van Ceulen employed to obtain his continued square root approximations to $\pi$\index{pi@$\pi\;(3.14159\ldots)$!continued square root expressions for}\index{constants, named!pi@$\pi\;(3.14159\ldots)$}.


\item\vspace{9pt} \label{Gol1911} (1911) Karl Goldziher\index[ppl]{Goldziher, Karl|textbf}, Beitr\"{a}ge zur Praxis der f\"{u}r die Berechnung des Rentenzinsfusses verwendbaren speziellen trinomischen Gleichung. \emph{Zeitschrift f\"{u}r Mathematik und Physik} \textbf{59}, 410--431.  

{\footnotesize Source: \href{https://books.google.com/books?id=rN1oWLc63_4C&pg=PA410#v=onepage&q&f=false}{Google Books}.}

\vspace{4pt}

This is a survey of practical methods for solving trinomial equations\index{trinomial equations} arising in the computation of pension interest rates\index{pension problem}. The paper's second section draws on the work of G\"{u}nther\index[ppl]{Gunther@G\"{u}nther, Siegmund}, Hoffmann\index[ppl]{Hoffmann, K. E.}, Netto, and Isenkrahe\index[ppl]{Isenkrahe, C.} from the 1880s concerning fixed point iteration, by which one generates pension equation solutions expressible as continued $r$th roots\index{continued robb roots@continued $r$th roots}. 


\item\vspace{9pt} \label{Rab1911} (1911) Izko-Ewna Rabinowitsch\index[ppl]{Rabinowitsch, Izko-Ewna|textbf}, Beitr\"{a}ge zur aufl\"{o}sung der algebraischen gleichungen 5. grades. Buchdruckerei K. J. Wyss, Bern. 32 pp.
\index{plagiarism}

{\footnotesize Source: \href{https://gdz.sub.uni-goettingen.de/id/PPN31492888X?tify=\%7B\%22view\%22\%3A\%22info\%22\%7D}{G\"{o}ttinger Digitalisierungszentrum}.}

\vspace{4pt}

This doctoral dissertation is noted here primarily for its blatant plagiarism, on page 27, of the discussion of ``Kettenfunktion,"\index{Kettenfunctionen@\emph{Kettenfunctionen}} ``Kettenwurzel,"\index{Kettenwurzeln@\emph{Kettenwurzeln}} ``Kettenbruch,"\index{Kettenbr\"{u}chen@\emph{Kettenbr\"{u}chen}} and ``Kettenlogarithmus"\index{Kettenlogarithmus@\emph{Kettenlogarithmus}} given on page 322 of \hyperref[Hey1894a]{\textsc{Heymann 1894a}}. The nearly identical copy is marred only by its failure to record the index $p$ in defining a $p$th root function. Nowhere does Rabinowitsch cite or credit Heymann.


\item\vspace{9pt} \label{Ram1911} (1911) Srinivasa Ramanujan\index[ppl]{Ramanujan, Srinivasa|textbf}, Question 289. \emph{Journal of the Indian Mathematical Society} \textbf{3}(2), 90. Solution by the proposer, \emph{ibid.} \textbf{4}(6), 1912, 226. The problem and its solution are reprinted in \hyperref[Ram1927]{\textsc{Ramanujan 1927}}\index[ppl]{Ramanujan, Srinivasa}, p. 323.
\index{continued square roots!Ramanujan's}


{\footnotesize Problem source: \href{https://babel.hathitrust.org/cgi/pt?id=umn.319510002394343&seq=642}{HathiTrust}.}

{\footnotesize Solution source: \href{https://babel.hathitrust.org/cgi/pt?id=chi.63388476&seq=252}{HathiTrust}.}

{\footnotesize Cited in 
\hyperref[Her1935]{\textsc{Herschfeld 1935}}\index[ppl]{Herschfeld, Aaron}, 
the solution to \hyperref[Han1955]{\textsc{Hanisch 1955}}\index[ppl]{Hanisch, Herman}, 
\hyperref[Sch1961]{\textsc{Schuske and Thron 1961}}\index[ppl]{Schuske, Georgellen}\index[ppl]{Thron, W. J.}, 
\hyperref[Bor1991]{\textsc{Borwein and de Barra 1991}}\index[ppl]{Borwein, Jonathan M.}\index[ppl]{de Barra, G.}, 
\hyperref[Ber1999]{\textsc{Berndt, Choi, and Kang 1999}}\index[ppl]{Berndt, Bruce C.}\index[ppl]{Choi, Youn-Seo}\index[ppl]{Kang, Soon-Yi}, 
\hyperref[Rao2005]{\textsc{Rao and Vanden Berghe 2005}}\index[ppl]{Rao, K. Srinivasa}\index[ppl]{Vanden Berghe, G.}, 
\hyperref[Muk2013]{\textsc{Mukherjee 2013}}\index[ppl]{Mukherjee, Soumendu Sundar}, 
\hyperref[Lyn2014]{\textsc{Lynd 2014}}\index[ppl]{Lynd, Chris D.}, and
\hyperref[Wei]{\textsc{Weisstein n.d.}}\index[ppl]{Weisstein, Eric W.}}

\vspace{4pt}

The problem is to prove the formulas
\begin{equation}\label{E:Ram1911-1}
3=\sqrt{1+2\sqrt{1+3\sqrt{1+\cdots}}}
\end{equation}
and
\begin{equation}\label{E:Ram1911-2}
4=\sqrt{6+2\sqrt{7+3\sqrt{8+\cdots}}}\;.
\end{equation}
These identities are notable because they do not represent the iteration of a forward recurrence relation\index{recurrence relation} to find a fixed point (although, in his notebooks, Ramanujan derived the more general precursors to these formulas by successive substitution\index{successive substitution!20th century examples}; see \hyperref[Ber1989]{\textsc{Berndt 1989}}\index[ppl]{Berndt, Bruce C.}). On p. 348 of \hyperref[Ram1927]{\textsc{Ramanujan 1927}}\index[ppl]{Ramanujan, Srinivasa}, and in \hyperref[Her1935]{\textsc{Herschfeld 1935}}\index[ppl]{Herschfeld, Aaron}, it is noted that Ramanujan's solution is incomplete, and rigorous convergence arguments are supplied.


\item\vspace{9pt} \label{Bop1913} (1913) K. Bopp\index[ppl]{Bopp, K.|textbf}, Eine Schrift von Ensheim ``Recherches sur les calculs diff\'{e}rentiel et int\'{e}gral" mit einem sich darauf beziehenden, nicht in die ``Oeuvres" \"{u}bergegangenen Brief. \emph{Sitzungsberichte der Heidelberger Akademie der Wissenschaften, Mathematisch-naturwissenschaftliche Klasse, Abteilung A. Mathematisch-physikalische Wissenschaften}, essay 7, 49 pages. 

{\footnotesize Source: \href{https://archiv.ub.uni-heidelberg.de/volltextserver/12423/}{heiDOK}.}
\index[ppl]{Ensheim, Moses}

April 10, 1913 was the one hundredth anniversary of the death of Joseph-Louis Lagrange\index[ppl]{Lagrange, Joseph-Louis}. In connection with this centenary, Karl Bopp's essay reveals the existence of the pamphlet \hyperref[Ens1799]{\textsc{Ensheim 1799}} (considered rare even in 1913), in the back of which was printed a previously unknown letter by Lagrange, written to Ensheim, dated 23 Nov. 1799.\footnote{Apparently Ensheim had sent some of his work to Lagrange, and presumed to print Lagrange's reply on the inside back cover of his pamphlet. Following perfunctory excuses about its tardiness, Lagrange's short letter notes the similarity of Ensheim's calculus methods to John Landen's\index[ppl]{Landen, John} some thirty years earlier, and closes with a generally complimentary flourish.} Bopp explicates the content of \hyperref[Ens1799]{\textsc{Ensheim 1799}}, including Ensheim's independent rediscovery of the method of inscribed and circumscribed regular polygons\index{polygons, regular} in a circle to generate continued square root expressions as upper and lower bounds for $\pi$.\index{pi@$\pi\;(3.14159\ldots)$!Archimedean algorithm for}\index{constants, named!pi@$\pi\;(3.14159\ldots)$}\footnote{Notable in this connection is that Bopp expands Ensheim's inequalities to an excessive degree, creating some of the most elaborate continued square root expressions to be found in any of the sources reviewed here.} Bopp cites \hyperref[Eul1744]{\textsc{Euler 1744}} (presumably unknown to Ensheim) for an earlier example of this polygonal approximation, and also notes the exchange between M. Cantor and F. Rudio which resulted in \hyperref[Rud1891]{\textsc{Rudio 1891}}.


\item\vspace{9pt} \label{Ram1914} (1914) Srinivasa Ramanujan\index[ppl]{Ramanujan, Srinivasa|textbf}, Question 507. \emph{Journal of the Indian Mathematical Society} \textbf{5}, 240. Solution by the proposer, \emph{ibid.} \textbf{6}, 1914, 74--77. The problem and its solution are reprinted in \hyperref[Ram1927]{\textsc{Ramanujan 1927}}\index[ppl]{Ramanujan, Srinivasa}, 327--329.
\index{continued square roots!Ramanujan's}

{\footnotesize Problem source: \href{https://babel.hathitrust.org/cgi/pt?id=chi.63388476&seq=526}{HathiTrust}.}

{\footnotesize Solution source: \href{https://babel.hathitrust.org/cgi/pt?id=chi.63388476&seq=614}{HathiTrust}.}

{\footnotesize Cited in 
\hyperref[Ber1993]{\textsc{Berndt and Bhargava 1993}}\index[ppl]{Berndt, Bruce C.}\index[ppl]{Bhargava, S.}, 
\hyperref[Ber1994]{\textsc{Berndt 1994}}\index[ppl]{Berndt, Bruce C.}, and 
and \hyperref[Ber1999]{\textsc{Berndt, Choi, and Kang 1999}}\index[ppl]{Berndt, Bruce C.}\index[ppl]{Choi, Youn-Seo}\index[ppl]{Kang, Soon-Yi}.}

\vspace{4pt}

The problem states: ``Solve completely
\begin{equation*}
x^2=y+a\;,\quad y^2=z+a\;,\quad z^2=x+a\;;
\end{equation*}
and hence show that 
\begin{table}[ht]
\begin{center}
\begin{tabular}{crcl}
($a$)&$\surd[8-\surd\{8+\surd(8-\ldots)\}]\negthickspace\negthickspace\negthickspace$&=&$\negthickspace\negthickspace\negthickspace1+2\surd3\sin 20^{\circ}$\\
($b$)&$\surd[11-2\surd\{11+2\surd(-\ldots)\}]\negthickspace\negthickspace\negthickspace$&=&$\negthickspace\negthickspace\negthickspace1+4\sin 10^{\circ}$\\
($c$)&$\surd[23-2\surd\{23+2\surd(23+2\surd 23-\ldots)\}]\negthickspace\negthickspace\negthickspace$&=&$\negthickspace\negthickspace\negthickspace1+4\surd 3 \sin 20^{\circ}\;."$
\end{tabular}
\end{center}
\end{table}


\item\vspace{9pt} \label{Ram1915} (1915) Srinivasa Ramanujan\index[ppl]{Ramanujan, Srinivasa|textbf}, Question 722. \emph{Journal of the Indian Mathematical Society} \textbf{8}, 240. Solution in M. B. Rao, Cyclic equations, \emph{ibid.} \textbf{16}, 1925, 139--154. Solution by G. N. Watson, \emph{ibid.} \textbf{18}, 1929, 113--117. The problem is reprinted in \hyperref[Ram1927]{\textsc{Ramanujan 1927}}\index[ppl]{Ramanujan, Srinivasa}, 332.

{\footnotesize Problem source: \href{https://babel.hathitrust.org/cgi/pt?id=chi.63388483&seq=252}{HathiTrust}.}

{\footnotesize Cited in
\hyperref[Ber1993]{\textsc{Berndt and Bhargava 1993}}\index[ppl]{Berndt, Bruce C.}\index[ppl]{Bhargava, S.},
\hyperref[Ber1994]{\textsc{Berndt 1994}}\index[ppl]{Berndt, Bruce C.}, and 
\hyperref[Ber1999]{\textsc{Berndt, Choi, and Kang 1999}}\index[ppl]{Berndt, Bruce C.}\index[ppl]{Choi, Youn-Seo}\index[ppl]{Kang, Soon-Yi}.}
\index{continued square roots!Ramanujan's}

\vspace{4pt}

The problem states: ``Solve completely
\begin{equation*}
x^2=a+y\;,\quad y^2=a+z\;,\quad z^2=a+u\;, \quad u^2=a+x\;;
\end{equation*}
and deduce that, if
\[x=\surd[5+\surd\{5+\surd(5-\surd(5+x))\}]\;,\]
then
\[x=\tfrac{1}{2}\{2+\surd5+\surd(15-6\surd5)\}\;;\]
and that, if
\[x=\surd[5+\surd\{5-\surd(5-\surd(5+x))\}]\;,\]
then
\[x=\tfrac{1}{4}[\surd5-2+\surd(13-4\surd5)+\surd\{50+12\surd5-2\surd(65-20\surd5)\}]."\]
\hyperref[Ber1993]{\textsc{Berndt and Bhargava 1993}}\index[ppl]{Berndt, Bruce C.}\index[ppl]{Bhargava, S.} note that the solution by G. N. Watson\index[ppl]{Watson, G. N.} appeared in the journal 14 years after the problem.


\item\vspace{9pt} \label{Gin1916} (1916) J. J. Ginsburg\index[ppl]{Ginsburg, J. J.|textbf}, Problem 460. \emph{The American Mathematical Monthly} \textbf{23}(6), 209. Solution by Nathan Altshiller\index[ppl]{Altshiller, Nathan}, \emph{ibid.} \textbf{24}(1), 1917, 32--33
\index{continued square roots!of terms $a_n=1$}
\index{continued square roots!of terms $a_n=c$}

{\footnotesize Solution source: \href{https://doi.org/10.2307/2972669}{doi.org (JSTOR)}.} 

{\footnotesize Cited in 
\hyperref[Her1935]{\textsc{Herschfeld 1935}}\index[ppl]{Herschfeld, Aaron}, 
\hyperref[Her1935]{\textsc{Sizer 1986}}\index[ppl]{Sizer, Walter S.}, 
\hyperref[Bor1991]{\textsc{Borwein and de Barra 1991}}\index[ppl]{Borwein, Jonathan M.}\index[ppl]{de Barra, G.}, and
\hyperref[Jon2008]{\textsc{Jones 2008}}\index[ppl]{Jones, Dixon J.}.}

\vspace{4pt}

The problem asks for the limit of the continued square root\index{continued square roots}
\[\sqrt{1+\sqrt{1+\sqrt{1+\cdots}}}\;.\]
The solution considers the more general form
\[\sqrt{a+\sqrt{a+\sqrt{a+\cdots}}}\]
and gives the limit as
\[\dfrac{1+\sqrt{1+4a}}{2}\;,\]
a result shown previously in \hyperref[Ber1692]{\textsc{Bernoulli 1692}}\index[ppl]{Bernoulli, Jacob}.


\item\vspace{9pt} \label{Pol1916} (1916) G. P\'{o}lya\index[ppl]{Polya@P\'{o}lya, G.|textbf}, Problem proposal 501. \emph{Archiv der Mathematik und Physik, Series 3} \textbf{24}, 84. Solution by G. Szeg\"{o}\index[ppl]{Szeg\"{o}, G.|textbf}, \emph{ibid.} \textbf{25}, 1917, 88--89.
\index{continued square roots!of arbitrary nonnegative real terms}
\index{continued square roots!convergence conditions for}

{\footnotesize Solution source: \href{https://books.google.com/books?id=-4QbAQAAMAAJ\&newbks=1\&newbks_redir=0\&pg=PA88\#v=onepage\&q\&f=false}{Google Books}.}

{\footnotesize Cited in
\hyperref[Her1935]{\textsc{Herschfeld 1935}}\index[ppl]{Herschfeld, Aaron}, 
\hyperref[Pol1925]{\textsc{P\'{o}lya and Szeg\"{o} 1925}}\index[ppl]{Polya@P\'{o}lya, G.}\index[ppl]{Szeg\"{o}, G.}, and
\hyperref[Sch1961]{\textsc{Schuske and Thron 1961}}\index[ppl]{Schuske, Georgellen}\index[ppl]{Thron, W. J.}.}

\vspace{4pt}

The problem is to prove that for positive real $a_i$,
\[\lim_{n\to\infty} \sqrt{a_1+\sqrt{a_2+\sqrt{\cdots+\sqrt{a_n}}}}\]
converges or diverges accordingly as 
\[\limsup_{n\to\infty} \dfrac{\log\log a_n}{n}\]
is less than or greater than $\log 2$. The problem reappears in \hyperref[Pol1925]{\textsc{P\'{o}lya and Szeg\"{o} 1925}}.


\item\vspace{9pt} \label{Pin1918a} (1918a) Salvatore Pincherle\index[ppl]{Pincherle, Salvatore|textbf}, Sulle catene di radicali quadratici. \emph{Rendiconto delle sessioni della Reale Accademia delle Scienze dell'instituto di Bologna. Classe di scienze fisiche. Nuova Serie} \textbf{22} (1917--18), 35--55.
\index{continued square roots!of periodic real terms}

{\footnotesize Source: \href{https://babel.hathitrust.org/cgi/pt?id=osu.32435061153979&seq=39}{HathiTrust.}}

{\footnotesize Cited in 
\hyperref[Pin1918b]{\textsc{Pincherle 1918b}} and \hyperref[Pin1918c]{\textsc{1918c}}\index[ppl]{Pincherle, Salvatore}, 
\hyperref[Sca1920]{\textsc{Scarpis 1920}}\index[ppl]{Scarpis, Umberto}, and 
\hyperref[Pol1925]{\textsc{P\'{o}lya and Szeg\"{o} 1925}}\index[ppl]{Polya@P\'{o}lya, G.}\index[ppl]{Szeg\"{o}, G.}.}

\vspace{4pt}

Initiating a somewhat confusing record of continued square root\index{continued square roots} publications in one year, this paper is the first of two parts, both parts having the same title, but the second (\hyperref[Pin1918b]{\textsc{Pincherle 1918b}}) appearing in a different journal. The author considers the continued square root\index{continued square roots!of periodic real terms}
\begin{equation}\label{E:Pin1918a-1}
\sqrt{a\pm\sqrt{a\pm\sqrt{a\pm\ldots\pm\sqrt{a}}}}\;,
\end{equation}
where the signs occur in a periodic order. He introduces definitions and notation, then investigates convergence in the cases $a>2$ and $a=2$. 


\item\vspace{9pt} \label{Pin1918b} (1918b) Salvatore Pincherle\index[ppl]{Pincherle, Salvatore|textbf}, Sulle catene di radicali quadratici. \emph{Atti della Rendiconti Accademia della Scienze di Torino. Classe di scienze fisiche, matematiche e naturale.} \textbf{53} (1917--1918), 437/745--455/763.
\index{continued square roots!of periodic real terms}

{\footnotesize Source: \href{https://babel.hathitrust.org/cgi/pt?id=osu.32435061422614&seq=503}{HathiTrust.}}

{\footnotesize Cited in 
\hyperref[Pin1918c]{\textsc{Pincherle 1918c}}\index[ppl]{Pincherle, Salvatore}, 
\hyperref[Sca1920]{\textsc{Scarpis 1920}}\index[ppl]{Scarpis, Umberto},
\hyperref[Pol1925]{\textsc{P\'{o}lya and Szeg\"{o} 1925}}\index[ppl]{Polya@P\'{o}lya, G.}\index[ppl]{Szeg\"{o}, G.}, and 
\hyperref[Kub2002]{\textsc{Kuba and Schoissengeier 2002}}\index[ppl]{Kuba, G.}\index[ppl]{Schoissengeier, J.}.} 

\vspace{4pt}

Presented in Turin, this paper is the second part of \hyperref[Pin1918a]{\textsc{Pincherle 1918a}}, which was read in Bologna. Here, Pincherle repeats most of his notation and definitions for a new audience, then wraps up his initial investigation into the continued square root\index{continued square roots} \eqref{E:Pin1918a-1} above by considering the cases in which $0\le a<2$. 


\item\vspace{9pt} \label{Pin1918c} (1918c) Salvatore Pincherle\index[ppl]{Pincherle, Salvatore|textbf}, Sulle radici reali delle equazioni iterate di una equazione quadratica. \emph{Atti dell' Accademia Nazionale dei Lincei Rendiconti. Classe di scienze fisiche, matematiche e naturali, Roma}. Series 5, \textbf{27}, 2nd sem., 177--183. 
\index{continued square roots!of periodic real terms} 

{\footnotesize Source: \href{http://operedigitali.lincei.it/rendicontiFMN/rol/pdf/S5V27T2A1918P177_183.pdf}{Accademia Nazionale dei Lincei}.}

{\footnotesize Cited in
\hyperref[Sca1920]{\textsc{Scarpis 1920}}\index[ppl]{Scarpis, Umberto},
\hyperref[Pol1925]{\textsc{P\'{o}lya and Szeg\"{o} 1925}}\index[ppl]{Polya@P\'{o}lya, G.}\index[ppl]{Szeg\"{o}, G.}, and 
\hyperref[Kub2002]{\textsc{Kuba and Schoissengeier 2002}}\index[ppl]{Kuba, G.}\index[ppl]{Schoissengeier, J.}.} 

\vspace{4pt}

The author inductively defines $\alpha_1(x)=x^2-a$ and $\alpha_n(x)=\alpha_{n-1}^2(x)-a$. He observes that, to each arrangement of the signs $+$ or $-$ in the terms of the continued square root \eqref{E:Pin1918a-1} above, having $n$ radicals, there corresponds a different solution of $\alpha_n(x)=0$. He then gives a criterion for these solutions to be real. See also \hyperref[Vel2016a]{\textsc{Vellucci and Bersani 2016a}}\index[ppl]{Vellucci, Pierluigi}\index[ppl]{Bersani, Alberto Maria}.


\item\vspace{9pt} \label{Bal1920} \dag(1920) W. W. Rouse Ball\index[ppl]{Ball, W. W. Rouse|textbf}, \emph{Mathematical Recreations and Essays}, 9th edition, Macmillan and Co., Limited, London.
\index{pi@$\pi\;(3.14159\ldots)$!continued square root expressions for}\index{constants, named!pi@$\pi\;(3.14159\ldots)$}

{\footnotesize Source: \href{https://archive.org/details/mathematicalrecr1920ball/page/301/mode/1up}{Internet Archive.}}

On pp. 301--302, Ball offers further details about the work of Van Ceulen\index[ppl]{Ceulen, Ludolph van} and Snellius\index[ppl]{Snellius, Willebrordus (\emph{aka} Snell, Willebrord)} in calculating the digits of $\pi$. He comments that van Ceulen's approximation ``was calculated by finding the perimeters of the inscribed and circumscribed regular polygons\index{polygons, regular} of $60\times 2^{33}$ sides, obtained by the repeated use of a theorem of his discovery equivalent to the formula $1-\cos A=2\sin^2\tfrac{1}{2}A.$" Ball also explains that Snell's improvements on Van Ceulen's methods allowed him to obtain ``from a hexagon an approximation as correct as that for which Archimedes\index[ppl]{Archimedes} had required a polygon of $96$ sides, while from a polygon of $96$ sides he determined the value of $\pi$ correct to seven decimal places instead of the two places obtained by Archimedes.''


\item\vspace{9pt} \label{Sca1920} (1920) Umberto Scarpis\index[ppl]{Scarpis, Umberto|textbf}, Catene periodiche di radicali. \emph{Giornale di matematiche di Battaglini} \textbf{58}, 1--13.  
\index{continued square roots!of periodic real terms}

{\footnotesize Source: \href{https://books.google.com/books?id=stgSAQAAMAAJ&newbks=1&newbks_redir=0&pg=RA1-PA1#v=onepage&q&f=false}{Google Books}.}

\vspace{4pt}

Writing in a journal intended for university students, the author presents a relatively accessible account of some results in \hyperref[Pin1918a]{\textsc{Pincherle 1918a}}, \hyperref[Pin1918b]{\textsc{1918b}}, and \hyperref[Pin1918c]{\textsc{1918c}} concerning
\[+\sqrt{a\pm \sqrt{a\pm \sqrt{a\ldots}}}\;,\]
where the signs are taken to be in a periodic sequence. 


\item\vspace{9pt} \label{Kak1924} (1924) S. Kakeya\index[ppl]{Kakeya, S\^{o}ichi|textbf}, On a generalized scale of notations, \emph{Japanese Journal of Mathematics} \textbf{1} 95--108. 
\index{f@$f$-expansions}
\index{real numbers!g@$g$-adic expansions of}

{\footnotesize Source: \href{https://doi.org/10.4099/jjm1924.1.0_95}{doi.org (J-STAGE)}.}

{\footnotesize Cited in 
\hyperref[Thr1961]{\textsc{Thron 1961}}\index[ppl]{Thron, W. J.} and 
\hyperref[Sch2016]{\textsc{Schweiger 2016}}\index[ppl]{Schweiger, Fritz}.}

\vspace{4pt}

This paper anticipates the $f$-expansion of \hyperref[Bis1944]{\textsc{Bissinger 1944}}\index[ppl]{Bissinger, B. H.}. The author writes, ``We have usually two methods of denoting the [real] numbers, namely the method of decimal fractions\index{real numbers!decimal expansions of} and of continued fractions\index{continued fractions}\index{real numbers!continued fraction expansions of}. These two methods can be included in the same consideration, and a generalized method of notation can be obtained, to which the usual methods belong as special cases.

``We consider a function $f(x)$ of the following properties:

``1. \emph{$f(x)$ is continuous and monotonic (havin[g] no constant part) in the interval $0 \le x \le 1$.}

``We distinguish the two cases when $f(x)$ is increasing and decreasing by the names \emph{case} A and \emph{case} B.

``2. \emph{$f(0)=\alpha, f(1)=\beta$, where $\alpha$ and $\beta$ are integers.}

``$\alpha$ is less or greater than $\beta$ according as the case is A or B.

``$f(x)$ is differentiable almost everywhere. If $Df(x)$ denotes one of the derivatives of $f(x)$, then $Df(x)$, if it is not zero, is positive or negative according as the case is A or B.

``3. \emph{The measure of the set of $x$ in the interval $0 \le x \le 1$, for which $|Df(x)|\le 1$ is zero.} Namely $|Df(x)|>1$ almost everywhere.

``Then evidently we have the \emph{inverse function} $g(x)$ of $f(x)\,\ldots$"

Ultimately, expansions of the form $a_1+g(a_2+g(a_3 +\ldots+g(a_n)\ldots))$ are developed for arbitrary real numbers. Kakeya gives examples of functions $g$ which produce decimal\index{real numbers!decimal expansions of}, continued fraction\index{continued fractions}, continued square\index{continued squares}, and continued logarithm\index{continued logarithms} expansions. No references are cited.


\item\vspace{9pt} \label{Poi1925} (1925) Ossian Poirier\index[ppl]{Poirier, Ossian|textbf}, \emph{Angles et sinus: th\'{e}orie des radicaux carr\'{e}s sur 2 superpos\'{e}s.} J. Hermann, Paris.
\index{continued square roots!and trigonometric functions}
\index{continued square roots!of terms $a_n=\pm2$}

\vspace{4pt}

In the \emph{Avertissement} at the front of this eccentric book, the author admits that he is \emph{\'{e}tranger aux math\'{e}matiques.} Over 146 pages, he slowly and elaborately develops the trigonometric theory of continued square roots of $\pm2$, and invents yet another notation system intended to simplify the typography. Poirier's idea is to write
\[\surd a\pm \surd b \pm \surd c \pm \surd d | | | |\]
to represent
\[\sqrt{a\pm\sqrt{b\pm\sqrt{c\pm\sqrt{d}}}}\;.\]
He seems to have worked in isolation; no sources are cited. The book is mysterious. Almost nothing can be discovered about its author. In the dozens of databases consulted for this bibliography, the name ``Ossian Poirier" appears independently from this book only a few times, in connection with the military, law enforcement, or society page notices; it cannot even be determined whether these few disparate references are about the same person, much less whether the writer was an amateur mathematician or was using a pseudonym.


\item\vspace{9pt} \label{Sze1925} (1925)\footnote{Wherever this problem is cited, the date 1924 is given to Vol. 33. The problem itself was received on April 30, 1924; however, the date on the journal's title page for Vol. 33 is 1925.} G. Szeg\"{o}\index[ppl]{Szeg\"{o}, G.|textbf}, Aufgabe 18. \emph{Jahresbericht der Deutschen Mathematiker-Vereinigung} \textbf{33}, 69. Solution (received October 6, 1924) by N. Obreschkoff\index[ppl]{Obreschkoff, N.}, \emph{ibid.}, 117--118.
\index{continued square roots!of terms $a_n=\pm2$}

{\footnotesize Solution source: \href{https://gdz.sub.uni-goettingen.de/id/PPN37721857X_0033?tify=\%7B\%22pages\%22\%3A\%5B382\%5D\%2C\%22view\%22\%3A\%22info\%22\%7D}{G\"{o}ttinger Digitalisierungszentrum}.}

{\footnotesize Cited in
\hyperref[Sze1930]{\textsc{Szeg\"{o} 1930}}\index[ppl]{Szeg\"{o}, G.} and
\hyperref[Her1935]{\textsc{Herschfeld 1935}}\index[ppl]{Herschfeld, Aaron}.}

\vspace{4pt}

The problem is to prove that
\begin{equation}\label{E:Sze1925-1}
\epsilon_0\sqrt{2+\epsilon_1\sqrt{2+\epsilon_2\sqrt{2+\cdots}}}=2\sin\left(\dfrac{\pi}{4}\sum_{n=0}^\infty\dfrac{\epsilon_0\epsilon_1\epsilon_2\cdots\epsilon_n}{2^n}\right)\;,
\end{equation}
where the $\epsilon_i$ assume one of the values $-1,0,1$. (The formula was discovered at least 25 years earlier; see \hyperref[Boc1899]{\textsc{Bochow 1899}}, and also \hyperref[Sze1930]{\textsc{Szeg\"{o} 1930}}.)


\item\vspace{9pt} \label{Pol1925} (1925) G. P\'{o}lya\index[ppl]{Polya@P\'{o}lya, G.|textbf} and G. Szeg\"{o}\index[ppl]{Szeg\"{o}, G.|textbf}, \emph{Aufgaben und Lehrs\"{a}tze aus der Analysis}, Vol. 1, Julius Springer, Berlin. Reprinted as \emph{Problems and Theorems in Analysis}, Vol. I, translated by D. Aeppli\index[ppl]{Aeppli, D.}, Springer-Verlag New York, 1972.
\index{continued square roots!of terms $a_n=\pm2$}
\index{continued fractions!simple}
\index{continued square roots!of arbitrary nonnegative real terms}
\index{continued square roots!convergence conditions for}

{\footnotesize Cited in 
\hyperref[Her1935]{\textsc{Herschfeld 1935}}\index[ppl]{Herschfeld, Aaron}, 
\hyperref[And1985]{\textsc{Andrushkiw 1985}}\index[ppl]{Andrushkiw, R. L.}, 
\hyperref[Jon1991]{\textsc{Jones 1991}} and \hyperref[Jon1995]{\textsc{1995}}\index[ppl]{Jones, Dixon J.}, 
\hyperref[Ben1999]{\textsc{Bencze 1999}}\index[ppl]{Bencze, Mih\'{a}ly}, 
\hyperref[Ber1999]{\textsc{Berndt, Choi, and Kang 1999}}\index[ppl]{Berndt, Bruce C.}\index[ppl]{Choi, Youn-Seo}\index[ppl]{Kang, Soon-Yi}, 
\hyperref[Kub2002]{\textsc{Kuba and Schoissengeier 2002}}\index[ppl]{Kuba, G.}\index[ppl]{Schoissengeier, J.}, 
\hyperref[Jon2008]{\textsc{Jones 2008}}\index[ppl]{Jones, Dixon J.}, 
\hyperref[Mor2012]{\textsc{Moreno and Garc\'{i}a-Caballero 2012}} and \hyperref[Mor2013b]{\textsc{2013b}}\index[ppl]{Moreno, Samuel G\'{o}mez}\index[ppl]{Garc\'{i}a-Caballero, Esther M.}, 
\hyperref[Sen2013]{\textsc{Senadheera 2013}}\index[ppl]{Senadheera, Jayantha}, and
\hyperref[Wei]{\textsc{Weisstein n.d.}}\index[ppl]{Weisstein, Eric W.}}

\vspace{4pt}

Problem 161 asks for justification for the equation\index{continued square roots!of terms $a_n=1$}\index{continued fractions}
\[\sqrt{1+\sqrt{1+\sqrt{1+\cdots}}}=1+\cfrac{1}{1+\cfrac{1}{1+\ddots}}\;.\]

Problem 162 restates \hyperref[Pol1916]{\textsc{P\'{o}lya 1916}}: for positive real $a_i$,
\[\lim_{n\to\infty} \sqrt{a_1+\sqrt{a_2+\sqrt{\cdots+\sqrt{a_n}}}}\]
converges or diverges accordingly as 
\[\limsup_{n\to\infty} \dfrac{\log\log a_n}{n}\]
is less than or greater than $\log 2$. (The 1972 version of this book in English mistakenly gives the value 2 instead of $\log 2$.)\index{continued square roots!of arbitrary nonnegative real terms}\index{continued square roots!convergence conditions for} 

Problem 163 asks for proof that 
\[\lim_{n\to\infty} \sqrt{a_1+\sqrt{a_2+\sqrt{\cdots+\sqrt{a_n}}}}\]
converges if the series $\sum_{n=1}^\infty 2^{-n} a_n(a_1a_2\cdots a_n)^{-1/2}$ converges.

Problem 183 restates \hyperref[Sze1925]{\textsc{Szeg\"{o} 1925}}: prove that 
\begin{equation}\label{E:Pol1925-1}
\epsilon_0\sqrt{2+\epsilon_1\sqrt{2+\epsilon_2\sqrt{2+\cdots}}}=2\sin\left(\dfrac{\pi}{4}\sum_{n=0}^\infty\dfrac{\epsilon_0\epsilon_1\epsilon_2\cdots\epsilon_n}{2^n}\right)\;,
\end{equation}
where the $\epsilon_i$ assume one of the values $-1,0,1$ (see also \hyperref[Boc1899]{\textsc{Bochow 1899}}). Problem 184 asks for proof that every $x\in[-2,2]$ can be written in the form
\[x=\epsilon_0\sqrt{2+\epsilon_1\sqrt{2+\epsilon_2\sqrt{2+\cdots}}}\;,\]
where the $\epsilon_i$ assume either of the values $-1$ or $1$. Problem 185 asks to show that a real number $x$ is of the form $x=2\cos k\pi$, $k$ rational, if and only if the sequence $\epsilon_i$ is periodic after a certain term.\index{continued square roots!and trigonometric functions}\index{continued square roots!of periodic real terms}\index{continued square roots!of terms $a_n=\pm2$} 


\item\vspace{9pt} \label{Ram1927} (1927) Srinivasa Ramanujan, \emph{Collected Papers of Srinivasa Ramanujan}\index[ppl]{Ramanujan, Srinivasa}, G. H. Hardy\index[ppl]{Hardy, G. H.|textbf}, P. V. Seshu Aiyar\index[ppl]{Seshu Aiyar, P. V.|textbf}, and B. M. Wilson\index[ppl]{Wilson, B. M.|textbf}, eds. The University Press, Cambridge, England. 
\index{continued square roots!Ramanujan's}
\index{continued square roots!convergence conditions for}
\index{continued square roots!of arbitrary nonnegative real terms}

{\footnotesize Cited in
\hyperref[Her1935]{\textsc{Herschfeld 1935}}\index[ppl]{Herschfeld, Aaron}, 
the solution to \hyperref[Han1955]{\textsc{Hanisch 1955}}\index[ppl]{Hanisch, Herman}, 
\hyperref[Siz1986]{\textsc{Sizer 1986}}\index[ppl]{Sizer, Walter S.}, 
\hyperref[Bor1991]{\textsc{Borwein and de Barra 1991}}\index[ppl]{Borwein, Jonathan M.}\index[ppl]{de Barra, G.}, 
\hyperref[Pic1991]{\textsc{Pickover and Lakhtakia 1991}}\index[ppl]{Lakhtakia, A.}, 
\hyperref[Ber1993]{\textsc{Berndt and Bhargava 1993}}\index[ppl]{Berndt, Bruce C.}\index[ppl]{Bhargava, S.}, 
\hyperref[Ber1994]{\textsc{Berndt 1994}}\index[ppl]{Berndt, Bruce C.}, 
\hyperref[Ber1999]{\textsc{Berndt, Choi, and Kang 1999}}\index[ppl]{Berndt, Bruce C.}\index[ppl]{Choi, Youn-Seo}\index[ppl]{Kang, Soon-Yi}, 
\hyperref[Rao2005]{\textsc{Rao and Vanden Berghe 2005}}\index[ppl]{Rao, K. Srinivasa}\index[ppl]{Vanden Berghe, G.}, 
\hyperref[Son2007]{\textsc{Sondow and Hadjicostas 2007}}\index[ppl]{Sondow, Jonathan}\index[ppl]{Hadjicostas, Petros}, 
\hyperref[Vel2016c]{\textsc{Vellucci and Bersani 2016c}}\index[ppl]{Vellucci, Pierluigi}\index[ppl]{Bersani, Alberto Maria}, and 
\hyperref[Wei]{\textsc{Weisstein n.d.}}\index[ppl]{Weisstein, Eric W.}}

\vspace{4pt}

In a note on page 348 discussing \hyperref[Ram1911]{\textsc{Ramanujan 1911}}, T. Vijayaraghavan\index[ppl]{Vijayaraghavan, T.|textbf} is credited with the following result: If
\begin{equation}
a_n\ge 0,\quad T_n=\sqrt{a_1+\sqrt{a_2+\sqrt{a_3+\cdots+\sqrt{a_n}}}}\;,
\end{equation}
then a necessary and sufficient condition for the existence of $\lim_{n\to\infty}T_n$ is that
\begin{equation}\label{E:Ram1927-1}
\limsup_{n\to\infty}\dfrac{\log a_n}{2^n}<\infty\;.
\end{equation}
As pointed out in \hyperref[Son2007]{\textsc{Sondow and Hadjicostas 2007}}\index[ppl]{Hadjicostas, Petros}, this is essentially equivalent to \eqref{E:Her1935-1} in \hyperref[Her1935]{\textsc{Herschfeld 1935}}\index[ppl]{Herschfeld, Aaron}, although it is rather more reminiscent of \hyperref[Pol1916]{\textsc{P\'{o}lya 1916}}\index[ppl]{Polya@P\'{o}lya, G.} in its formulation; indeed, Problem 162 in \hyperref[Pol1925]{\textsc{P\'{o}lya and Szeg\"{o} 1925}} is cited by the \emph{Collected Papers} editors and/or T. Vijayaraghavan\index[ppl]{Vijayaraghavan, T.} as ``a less precise form of the convergence criterion."

The earliest statement of this convergence criterion found thus far is in \hyperref[Wie1904a]{\textsc{Wiernsberger 1904a}}\index[ppl]{Wiernsberger, Paul}.

\hyperref[Ber1999]{\textsc{Berndt, Choi, and Kang 1999}}\index[ppl]{Berndt, Bruce C.}\index[ppl]{Choi, Youn-Seo}\index[ppl]{Kang, Soon-Yi} further comment that the note cited here ``was considerably amplified in a letter from Vijayaraghavan\index[ppl]{Vijayaraghavan, T.} to B. M. Wilson\index[ppl]{Wilson, B. M.} on January 4, 1928." That letter is reproduced in \hyperref[Ber1995]{\textsc{Berndt and Rankin 1995}}. 


\item\vspace{9pt} \label{Caj1928} (1928-1929) Florian Cajori\index[ppl]{Cajori, Florian|textbf}, \emph{A History of Mathematical Notations, Vols. I and II}, Open Court Publishing Co., Chicago. Reprinted in one volume by Dover Publications, Inc., Mineola, New York, 1993, ISBN-13: 978-0486677668.
\index{notation!for continued square roots}
\index{continued square roots!notation for}
\index{history!of square root notation}

\vspace{4pt}

Paragraphs 332 and 342--350 contain interesting remarks about the history of continued square root notation, including an expression by Christophorus Dibuadius\index[ppl]{Dibuadius, Christophorus} from 1605 for the length of a side of a regular 128-gon.


\item\vspace{9pt} \label{Sze1930} (1930) G. Szeg\"{o}\index[ppl]{Szeg\"{o}, G.|textbf}, Bemerkung zur Aufgabe 18. \emph{Jahresbericht der Deutschen Mathematiker-Vereinigung} \textbf{39}, 6.
\index{mea culpas}
\index{letters to editors}

{\footnotesize Source: \href{https://gdz.sub.uni-goettingen.de/id/PPN37721857X_0039?tify=\%7B\%22pages\%22\%3A\%5B301\%5D\%2C\%22view\%22\%3A\%22info\%22\%7D}{G\"{o}ttinger Digitalisierungszentrum}.}

{\footnotesize Cited in 
\hyperref[Her1935]{\textsc{Herschfeld 1935}}\index[ppl]{Herschfeld, Aaron}.}

\vspace{4pt}

In this one-paragraph letter, Szeg\"{o} acknowledges having been kindly informed by Dr. K. Bochow that the problem posed in \hyperref[Sze1925]{\textsc{Szeg\"{o} 1925}} (and reprinted as Problem 183 in \hyperref[Pol1925]{\textsc{P\'{o}lya and Szeg\"{o} 1925}}\index[ppl]{Polya@P\'{o}lya, G.}\index[ppl]{Szeg\"{o}, G.}) had been posed and solved in \hyperref[Boc1905]{\textsc{Bochow 1905}} and \hyperref[Boc1910]{\textsc{Bochow 1910}}. (See also \hyperref[Boc1899]{\textsc{Bochow 1899}}.)


\item\vspace{9pt} \label{Bol1935a} (1935a) Alexander Y. Boldyreff\index[ppl]{Boldyreff, Alexander Y.|textbf}, Problem 73. \emph{National Mathematics Magazine} \textbf{9}(4), 117. Solution by Dewey C. Duncan\index[ppl]{Duncan, Dewey C.|textbf}, \emph{ibid.} \textbf{9}(6), 1935, 177--178.

{\footnotesize Solution source: \href{https://doi.org/10.2307/3028823}{doi.org (JSTOR)}.}

\vspace{4pt}

It is proposed that\index{continued square roots!of terms $a_n=c$}
\[\dfrac{\sqrt{20+\sqrt{20+\sqrt{\cdots}}}+\sqrt{12+\sqrt{12+\sqrt{\cdots}}}}{\sqrt{6+\sqrt{6+\sqrt{\cdots}}}}=3\;.\]
The ``proof" by formal manipulation uses the formula
\[n+1=\sqrt{n(n+1)+\sqrt{n(n+1)+\sqrt{\cdots}}}\;.\]


\item\vspace{9pt} \label{Bol1935b} (1935b) Alexander Y. Boldyreff\index[ppl]{Boldyreff, Alexander Y.|textbf}, Problem 75. \emph{National Mathematics Magazine} \textbf{9}(4), 118. Solution by Dewey C. Duncan\index[ppl]{Duncan, Dewey C.|textbf}, \emph{ibid.} \textbf{9}(7), 1935, 208--209. Solution by Theodore Bennett\index[ppl]{Bennett, Theodore}, \emph{ibid.} \textbf{9}(8), 247--248.

{\footnotesize Solution source: \href{https://doi.org/10.2307/3028022}{doi.org (JSTOR)}.}

{\footnotesize Cited in 
\hyperref[Her1935]{\textsc{Herschfeld 1935}}\index[ppl]{Herschfeld, Aaron},
the solution to \hyperref[Den1983]{\textsc{Dence 1983}}\index[ppl]{Dence, Thomas P.}, and
\hyperref[Siz1986]{\textsc{Sizer 1986}}\index[ppl]{Sizer, Walter S.}.}

\vspace{4pt}

It is proposed that\index{continued robb roots@continued $r$th roots!of constant nonnegative real terms}
\[a=\sqrt[n]{ab+(a^{n-1}-b)\sqrt[n]{ab+(a^{n-1}-b)\sqrt[n]{\cdots}}}\]
where $a,b$ are positive real numbers and $n$ is a positive integer. The first solution involves only formal manipulations; the second justifies the existence of the limit as well.







\item\vspace{9pt} \label{Har1935} (1935) Vincent C. Harris\index[ppl]{Harris, Vincent C.|textbf}, Problem 78. \emph{National Mathematics Magazine} \textbf{9}(6), 180. Solution by Theodore Bennett\index[ppl]{Bennett, Theodore}, \emph{ibid.} \textbf{9}(8), 1935, 251--252.

{\footnotesize Solution source: \href{https://doi.org/10.2307/3028022}{doi.org (JSTOR)}.} 

{\footnotesize Cited in 
\hyperref[Her1935]{\textsc{Herschfeld 1935}}\index[ppl]{Herschfeld, Aaron}, 
\hyperref[Den1983]{\textsc{Dence 1983}}\index[ppl]{Dence, Thomas P.} and
\hyperref[Siz1986]{\textsc{Sizer 1986}}\index[ppl]{Sizer, Walter S.}.}

\vspace{4pt}

The value of
\[\sqrt{12+\sqrt{48+\sqrt{768+\sqrt{196608+\sqrt{\cdots}}}}}\]
is determined to be $1+\sqrt{13}$, using a method similar to Altshiller's\index[ppl]{Altshiller, Nathan} general solution of \hyperref[Gin1916]{\textsc{Ginsburg 1916}}\index[ppl]{Ginsburg, J. J.}. Aaron Herschfeld's critique of the formal manipulation of the infinite expressions in the above Problems 75 and 78 is printed on page 252.


\item\vspace{9pt} \label{Her1935} (1935) Aaron Herschfeld\index[ppl]{Herschfeld, Aaron|textbf}, On infinite radicals. \emph{The American Mathematical Monthly} \textbf{42}(7), 419--429.
\index{infinite radicals}

{\footnotesize \href{https://mathscinet.ams.org/mathscinet/article?mr=1523428}{MR1523428 (citation only)}}

{\footnotesize Source: \href{https://doi.org/10.1080/00029890.1935.11987745}{doi.org (Taylor \& Francis)}.}

{\footnotesize Cited in 
the solution to \hyperref[Ogi1949]{\textsc{Ogilvy 1949}}\index[ppl]{Ogilvy, C. S.}, 
the solution to \hyperref[Han1955]{\textsc{Hanisch 1955}}\index[ppl]{Hanisch, Herman}, 
\hyperref[Sch1961]{\textsc{Schuske and Thron 1961}}\index[ppl]{Schuske, Georgellen}\index[ppl]{Thron, W. J.}, 
\hyperref[Thr1961]{\textsc{Thron 1961}}\index[ppl]{Thron, W. J.}, 
\hyperref[Sch1962]{\textsc{Schuske and Thron 1962}}\index[ppl]{Schuske, Georgellen}\index[ppl]{Thron, W. J.}, 
\hyperref[Ogi1970]{\textsc{Ogilvy 1970}}\index[ppl]{Ogilvy, C. S.}, 
\hyperref[Roh1974]{\textsc{Rohde 1974}}\index[ppl]{Rohde, Hanns-Walter}, 
\hyperref[Jon1991]{\textsc{Jones 1991}} and \hyperref[Jon1995]{\textsc{1995}}\index[ppl]{Jones, Dixon J.}, 
\hyperref[Ber1999]{\textsc{Berndt, Choi, and Kang 1999}}\index[ppl]{Berndt, Bruce C.}\index[ppl]{Choi, Youn-Seo}\index[ppl]{Kang, Soon-Yi}, 
\hyperref[Lau1999]{\textsc{Laugwitz and Sch\"{o}nefuss 1999}}\index[ppl]{Laugwitz, Detlef}\index[ppl]{Schonefuss@Sch\"{o}nefuss, Lutz W.}, 
\hyperref[Lev2005]{\textsc{Levin 2005}}\index[ppl]{Levin, Aaron},
\hyperref[Rao2005]{\textsc{Rao and Vanden Berghe 2005}}\index[ppl]{Rao, K. Srinivasa}\index[ppl]{Vanden Berghe, G.}, 
\hyperref[Hum2007]{\textsc{Humphries 2007}}\index[ppl]{Humphries, Peter J.}, 
\hyperref[Son2007]{\textsc{Sondow and Hadjicostas 2007}}\index[ppl]{Sondow, Jonathan}\index[ppl]{Hadjicostas, Petros},
\hyperref[Zim2008a]{\textsc{Zimmerman and Ho 2008a}}\index[ppl]{Zimmerman, Seth}\index[ppl]{Ho, Chungwu}, 
\hyperref[Eft2012]{\textsc{Efthimiou 2012}}\index[ppl]{Efthimiou, Costas J.}, 
\hyperref[Glu2012]{\textsc{Gluzman and Yukalov 2012}}\index[ppl]{Gluzman, S.}\index[ppl]{Yukalov, V. I.},
\hyperref[Muk2013]{\textsc{Mukherjee 2013}}\index[ppl]{Mukherjee, Soumendu Sundar}, 
\hyperref[Sen2013]{\textsc{Senadheera 2013}}\index[ppl]{Senadheera, Jayantha}, 
\hyperref[Xi2013]{\textsc{Xi and Qi 2013}}\index[ppl]{Xi, Bo-Yan}\index[ppl]{Qi, Feng},
\hyperref[Cla2014]{\textsc{Clark and Richmond 2014}}\index[ppl]{Clark, Tyler}\index[ppl]{Richmond, Tom},
\hyperref[Lyn2014]{\textsc{Lynd 2014}}\index[ppl]{Lynd, Chris D.}, 
\hyperref[Jon2015]{\textsc{Jones 2015}}\index[ppl]{Jones, Dixon J.},
\hyperref[Les2016]{\textsc{Lesher and Lynd 2016}}\index[ppl]{Lesher, Devyn A.}\index[ppl]{Lynd, Chris D.}, 
\hyperref[Vel2016c]{\textsc{Vellucci and Bersani 2016c}}\index[ppl]{Vellucci, Pierluigi}\index[ppl]{Bersani, Alberto Maria}, and
\hyperref[Wei]{\textsc{Weisstein n.d.}}\index[ppl]{Weisstein, Eric W.}}

\vspace{4pt}

Herschfeld proves that the ``right infinite radical"\index{continued square roots!of arbitrary nonnegative real terms}\index{infinite radicals!right}
\[\sqrt{a_1+\sqrt{a_2+\sqrt{a_3+\cdots}}}\]
converges if
\begin{equation}\label{E:Her1935-1}
\limsup_{n\to\infty} a_n^{2^{-n}}<+\infty\;,
\end{equation}
where the $a_n$ are nonnegative reals\index{continued square roots!convergence conditions for}; compare \hyperref[Wie1904b]{\textsc{Wiernsberger 1904b}}, \hyperref[Pol1916]{\textsc{P\'{o}lya 1916}}\index[ppl]{Polya@P\'{o}lya, G.}, and Vijayaraghavan\index[ppl]{Vijayaraghavan, T.} in \hyperref[Ram1927]{\textsc{Ramanujan 1927}}\index[ppl]{Ramanujan, Srinivasa} for essentially the same result. The author calculates the approximate value\index{Kasner number $(1.75793\ldots)$}\index{constants, named!Kasner number $(1.75793\ldots)$}\footnote{Herschfeld calls this the Kasner number, in honor of Edward Kasner\index[ppl]{Kasner, Edward}, who ``for approximately twenty-five years$\,\ldots\,$has periodically suggested to his classes at Columbia University\index{Columbia University} the investigation of the problem of `infinite radicals.'"\index{infinite radicals} The web site \href{http://mathworld.wolfram.com/NestedRadicalConstant.html}{Wolfram MathWorld}\index{Wolfram MathWorld} calls this the ``nested radical constant."\index{nested radical constant|seeonly{Kasner number}}}
\[\sqrt{1+\sqrt{2+\sqrt{3+\cdots}}}\approx 1.757933\;;\]
and gives the formula
\[x(2+x)=x\sqrt{2^2+x\sqrt{2^4+x\sqrt{2^6+\cdots}}}\;.\]
Herschfeld follows \hyperref[Cip1908]{\textsc{Cipolla 1908}}\index[ppl]{Cipolla, Michele} in distinguishing between these ``right" infinite radicals (which are continued compositions)\index{infinite radicals!right} and ``left" infinite radicals\index{infinite radicals!left}\index{iterated square roots}, which are iterated compositions of the form
\[\cdots\sqrt{a_3+\sqrt{a_2+\sqrt{a_1}}}\;.\]
Herschfeld shows that a left infinite radical converges if and only if the sequence $a_n$ has a limit $a$. The limit of the left infinite radical is then\index{iterated square roots!convergence conditions for}
\[\dfrac{1+\sqrt{1+4a}}{2}\]
when $a>0$ (compare \hyperref[Ber1692]{\textsc{Bernoulli 1692}}\index[ppl]{Bernoulli, Jacob}; and \hyperref[Gin1916]{\textsc{Ginsburg 1916}}\index[ppl]{Ginsburg, J. J.}, which Herschfeld cites). Lastly, his Theorem III states without proof that for $a_n\ge 0$ and $r_n\in(0,1]$, suppose the series $\sum_{n=1}^\infty r_1r_2\cdots r_n$ converges; then the right infinite radical\index{continued rocc roots@continued $r_i$th roots!convergence conditions for}
\[(a_1+(a_2+(a_3+\cdots)^{r_3})^{r_2})^{r_1}\]
converges if and only if 
\[\limsup_{n\to\infty} a_n^{r_1r_2\cdots r_n}<+\infty\;.\]
%


\item\vspace{9pt} \label{Leb1937} (1937) Henri Lebesgue\index[ppl]{Lebesgue, Henri|textbf}, Sur certaines expressions irrationnelles illimit\'{e}s. \emph{Bulletin of the Calcutta Mathematical Society} \textbf{29}, 17--28.
\index{continued square roots!and trigonometric functions}
\index{continued square roots!of terms $a_n=\pm2$}

{\footnotesize Cited in 
\hyperref[Aok2016]{\textsc{Aoki and Kojima 2016}}\index[ppl]{Aoki, Noboru}\index[ppl]{Kojima, Shota}.}

\vspace{4pt}

Lebesgue claims the ideas for this paper occurred to him while lecturing at the Coll\`{e}ge de France in 1930; however, his insights unwittingly duplicate those of Rudio\index[ppl]{Rudio, Ferdinand}, Bochow\index[ppl]{Bochow, Karl}, Wiernsberger\index[ppl]{Wiernsberger, Paul}, Cipolla\index[ppl]{Cipolla, Michele}, P\'{o}lya\index[ppl]{Polya@P\'{o}lya, G.} and Szeg\"{o}\index[ppl]{Szeg\"{o}, G.}, Poirier\index[ppl]{Poirier, Ossian}, and others in explicating the continued square roots of $\pm2$\index{continued square roots!and trigonometric functions}\index{continued square roots!of terms $a_n=\pm2$} arising from trigonometric half-angle formulas. He attributes formula \eqref{E:Eul1744-2} to Euler in this regard, and realizes that Vi\`{e}te had obtained essentially the same result\index{Vi\`{e}te's formula for $\tfrac{2}{\pi}$}. Fortunately, he finds out and acknowledges that his work comprises independent rediscoveries of at least one author's work (see the next item below). 


\item\vspace{9pt} \label{Leb1938} (1938) Henri Lebesgue\index[ppl]{Lebesgue, Henri|textbf}, Sur certaines expressions irrationnelles illimit\'{e}s. \emph{Bulletin of the Calcutta Mathematical Society} \textbf{30}, 9--10.
\index{mea culpas}
\index{letters to editors}

{\footnotesize Cited in 
\hyperref[Aok2016]{\textsc{Aoki and Kojima 2016}}\index[ppl]{Aoki, Noboru}\index[ppl]{Kojima, Shota}.}

\vspace{4pt}

In this courteous note, about a page long, Lebesgue reports having been informed by his colleague Hadamard that the principal theorem proved in \hyperref[Leb1937]{\textsc{Lebesgue 1937}} was previously proved, using the same methods, in \hyperref[Wie1904b]{\textsc{Wiernsberger 1904b}}\index[ppl]{Wiernsberger, Paul}. Lebesgue makes a complimentary remark about the extent of Wiernsberger's thesis bibliography, and directs readers to Wiernsberger's papers from 1903, 1904, and 1905. 


\item\vspace{9pt} \label{Leh1938} (1938) D. H. Lehmer\index[ppl]{Lehmer, D. H.|textbf}, A cotangent analogue of continued fractions. \emph{Duke Mathematical Journal} \textbf{4}, 323--340.
\index{continued cotangents}
\index{real numbers!rational approximations to}

{\footnotesize \href{https://mathscinet.ams.org/mathscinet/article?mr=1546053}{MR1546053 (citation only)}}

{\footnotesize Source: \href{https://doi.org/10.1215/S0012-7094-38-00424-7}{doi.org (Project Euclid)}.} 

{\footnotesize Cited in 
\hyperref[Lei1940]{\textsc{Leighton 1940}}\index[ppl]{Leighton, Walter}, 
\hyperref[Bis1944]{\textsc{Bissinger 1944}}\index[ppl]{Bissinger, B. H.}, 
\hyperref[Spi1974]{\textsc{Spira and Scheeline 1974}}\index[ppl]{Spira, Robert}\index[ppl]{Scheeline, Alexander}, 
\hyperref[Sha1976]{\textsc{Shallit 1976}}\index[ppl]{Shallit, Jeffrey}, 
\hyperref[Sch1992]{\textsc{Sch\"{o}nefuss 1992}}\index[ppl]{Schonefuss@Sch\"{o}nefuss, Lutz W.}, 
\hyperref[Riv2007]{\textsc{Rivoal 2007}}\index[ppl]{Rivoal, T.}, and 
\hyperref[Sch2016]{\textsc{Schweiger 2016}}\index[ppl]{Schweiger, Fritz}.}

\vspace{4pt}

The author proposes the ``continued iteration"\index{continued iteration} of a function $f(x,y)$:
\begin{equation}\label{E:Leh1938-1}
f(x_1,f(x_2,f(x_3,\ldots)))\;.
\end{equation}
The $x_i$ are not initially defined, but later are assumed to be integers. Examples of \eqref{E:Leh1938-1} are given for various $f$: 
\begin{equation*}
f(x,y)=
\begin{cases}
x+y &\textrm{yields an infinite sum}\index{infinite sums}\\
xy &\textrm{\quad$^{\prime\prime}$ \quad an infinite product}\index{infinite products}\\
x+\dfrac{1}{y} &\textrm{\quad$^{\prime\prime}$ \quad a regular continued fraction}\index{continued fractions!regular}\\
\phantom{\Biggl |}x+\dfrac{y}{c} &\textrm{\quad$^{\prime\prime}$ \quad a power series of terms $\tfrac{x_i}{c^i}$\;.}\index{power series}\index{series!power} 
\end{cases}
\end{equation*}
If $c=10$ in the last case, one obtains a decimal expansion\index{real numbers!decimal expansions of}. The author's interest is in 
\[f(x,y)=\dfrac{1+xy}{y-x}\;,\]
which can be expressed as $f(x,y)=\cot(\arccot x-\arccot y)$ by using the identity
\[\cot(\alpha-\beta)=\dfrac{1+\cot\alpha\cot\beta}{\cot\alpha-\cot\beta}\;.\]
In the notation of this bibliography, \eqref{E:Leh1938-1} can be expressed as a continued composition\index{continued compositions!of cotangent functions} \eqref{E:ccomp} using $t_i(x)=(a_ix+1)/(x-a_i)=\cot(\arccot a_i-\arccot x)$, which produces 
\[\cot(\arccot a_0-\arccot a_1 +\arccot a_2 -\ldots)\;.\]
Lehmer\index[ppl]{Lehmer, D. H.} uses integer terms ``in order to obtain sequences of rational approximations to a real number."\index{real numbers!rational approximations to} After defining conditions for a ``regular" continued cotangent expansion\index{continued cotangents}, he proves that all infinite regular continued cotangent expansions converge, and, in a certain sense, do so more rapidly than continued fractions. He introduces a constant $\xi\approx 0.59263$\index{Lehmer's constant $(0.59263\ldots)$}, for which the continued cotangent expansion converges least rapidly. Some properties of continued cotangents and continued fractions\index{continued fractions} are analogous: for instance, $x$ is rational or irrational accordingly as its continued cotangent expansion is finite or infinite. 


\item\vspace{9pt} \label{Lei1940} (1940) Walter Leighton\index[ppl]{Leighton, Walter|textbf}, Proper continued fractions. \emph{The American Mathematical Monthly} \textbf{47}(5), 274--280.
\index{continued fractions!regular}
\index{continued cotangents}

{\footnotesize \href{https://mathscinet.ams.org/mathscinet/article?mr=0002567}{MR0002567}}

{\footnotesize Source: \href{https://doi.org/10.1080/00029890.1940.11990975}{doi.org (Taylor \& Francis)}.} 

{\footnotesize Cited in
\hyperref[Bis1944]{\textsc{Bissinger 1944}}\index[ppl]{Bissinger, B. H.}, 
\hyperref[Riv2007]{\textsc{Rivoal 2007}}\index[ppl]{Rivoal, T.}, and 
\hyperref[Sch2016]{\textsc{Schweiger 2016}}\index[ppl]{Schweiger, Fritz}.}
\index{continued fractions!proper}
\index{continued fractions!regular}

\vspace{4pt}

From the introduction: ``This paper generalizes the so-called `regular' continued fraction expansion\index{continued fractions!regular} of a real number. The treatment includes as a special case the `continued cotangent' expansion of Lehmer [\hyperref[Leh1938]{\textsc{Lehmer 1938}}].\index[ppl]{Lehmer, D. H.}"\index{continued cotangents}


\item\vspace{9pt} \label{Mac1941} (1941) Roy MacKay\index[ppl]{MacKay, Roy|textbf}, Problem E 474. \emph{The American Mathematical Monthly} \textbf{48}(5), 337. Solution by Eduardo Gaspar, \emph{ibid.} \textbf{49}(3), 1942, 197.
\index{iterated square roots!of constant nonnegative real terms}

{\footnotesize Solution source: \href{https://doi.org/10.1080/00029890.1942.11991206}{doi.org (Taylor \& Francis)}.}

{\footnotesize Cited in 
the solution to \hyperref[Ogi1949]{\textsc{Ogilvy 1949}}\index[ppl]{Ogilvy, C. S.}.}

\vspace{4pt}

The problem states: ``For $k>1$, define 
\begin{align*}
a_1&=\sqrt{k(k-1)}\\
a_n&=\sqrt{k(k-1)+a_{n-1}}\\
b_1&=\sqrt{k}\\
b_n&=\sqrt{kb_{n-1}}\;.
\end{align*}
Prove that $\lim_{n\to\infty} a_n=\lim_{n\to\infty} b_n=k.$"


\item\vspace{9pt} \label{Bis1944} (1944) B. H. Bissinger\index[ppl]{Bissinger, B. H.|textbf}, A generalization of continued fractions. \emph{Bulletin of the American Mathematical Society} \textbf{50}, 868--876.
\index{f@$f$-expansions}
\index{real numbers!f@$f$-expansions of}
\index{continued fractions!generalizations of}

{\footnotesize \href{https://mathscinet.ams.org/mathscinet/article?mr=0011338}{MR0011338}}

{\footnotesize Source: \href{https://projecteuclid.org/journals/bulletin-of-the-american-mathematical-society/volume-50/issue-12/A-generalization-of-continued-fractions/bams/1183506620.full}{Project Euclid}.}

{\footnotesize Cited in 
\hyperref[Bis1945]{\textsc{Bissinger and Herzog 1945}}\index[ppl]{Bissinger, B. H.}\index[ppl]{Herzog, F.}, 
\hyperref[Her1945]{\textsc{Herzog and Bissinger 1945}}\index[ppl]{Herzog, F.}\index[ppl]{Bissinger, B. H.}, 
\hyperref[Eve1946]{\textsc{Everett 1946}}\index[ppl]{Everett, C. J.}, 
\hyperref[Ren1957]{\textsc{R\'{e}nyi 1957}}\index[ppl]{Renyi@R\'{e}nyi, A.}, 
\hyperref[Swa1960]{\textsc{Swartz and Wendroff 1960}}\index[ppl]{Swartz, B. K.}\index[ppl]{Wendroff, B.},
\hyperref[Thr1961]{\textsc{Thron 1961}}\index[ppl]{Thron, W. J.}, 
\hyperref[Wat1970]{\textsc{Waterman 1970}}\index[ppl]{Waterman, Michael S.}, 
\hyperref[Sch1992]{\textsc{Sch\"{o}nefuss 1992}}\index[ppl]{Schonefuss@Sch\"{o}nefuss, Lutz W.}, 
\hyperref[Mar2004]{\textsc{Martin 2004}}\index[ppl]{Martin, Greg}, and 
\hyperref[Sch2016]{\textsc{Schweiger 2016}}\index[ppl]{Schweiger, Fritz}.}

\vspace{4pt}

Introduces the \emph{$f$-expansion} of a real number $x\in(0,1)$,
\begin{equation}\label{E:Bis1944-1}
x=f(a_1+f(a_2+f(\ldots)))\;,
\end{equation}
where the $a_i$ are positive integers and $f$ is monotone, strictly decreasing, defined on $(1,\infty)$ with $f(1)=1$ and $f(\infty)=0$, and is such that the chord defined by any two points on the graph of $y=f(t)$ has a slope bounded away from $-1$ for $t>f(2)+1$. See also \hyperref[Kak1924]{\textsc{Kakeya 1924}}\index[ppl]{Kakeya, S\^{o}ichi}.


\item\vspace{9pt} \label{Her1945} (1945) F. Herzog\index[ppl]{Herzog, F.|textbf} and B. H. Bissinger\index[ppl]{Bissinger, B. H.|textbf}, A generalization of Borel's and F. Bernstein's theorems on continued fractions. \emph{Duke Mathematical Journal} \textbf{12}, 325--334.
\index{f@$f$-expansions}
\index{continued fractions!generalizations of}

{\footnotesize \href{https://mathscinet.ams.org/mathscinet/article?mr=0012346}{MR0012346}}

{\footnotesize Source: \href{https://doi.org/10.1215/S0012-7094-45-01227-0}{doi.org (Project Euclid)}.} 

{\footnotesize Cited in
\hyperref[Bis1945]{\textsc{Bissinger and Herzog 1945}}\index[ppl]{Bissinger, B. H.|textbf}\index[ppl]{Herzog, F.}.} 

\vspace{4pt}

Regarding the $f$-expansion \eqref{E:Bis1944-1} above, the authors write, ``Borel\index[ppl]{Borel, \'{E}.}$\,\ldots\,$and F. Bernstein\index[ppl]{Bernstein, F.}$\,\ldots\,$have shown in the case $f(t)= t^{-1}$ (simple continued fractions\index{continued fractions!simple}) that (a) the set of all $x, 0 < x < 1$, for which $a_n(x)\le k_n$ for all $n$, where the $k_n$ are positive integers, is of measure zero if and only if $\sum 1/k_n$ diverges, (b) the set of all $x, 0 < x < 1$, for which $a_n(x) > 1$ for all $n$ is of measure zero. F. Bernstein showed these statements to be true also when only a fixed subsequence of the $a_n(x)$ is considered. Certain consequences of statements (a) and (b), concerning the boundedness of the sequence $a_n(x)$ as well as the occurrence of finitely many ones in the sequence $a_n(x)$, were derived in these two papers.

``This paper is concerned with the problem of generalizing Borel's and Bemstein's results to [$f$-expansions]. The first results in this direction were obtained by Bissinger\index[ppl]{Bissinger, B. H.}$\,\ldots\,$who proved these statements to be true when $f$ belongs to a certain class of polygonal functions\index{polygonal functions}. The proofs of these statements are based on the fact that the $a_n(x)$ are statistically independent for polygonal functions. Consequently these proofs cannot be used to re-establish Borel's and Bernstein's results since, as already pointed out by them, the $a_n(x)$ are not statistically independent for $f(t) = t^{-1}$.

``This paper gives a solution of the problem which remained --- to characterize analytically a [class of functions] which does include $f(t) = t^{-1}$ and for which the statements (a) and (b) hold." 


\item\vspace{9pt} \label{Bis1945} (1945) B. H. Bissinger\index[ppl]{Bissinger, B. H.|textbf} and F. Herzog\index[ppl]{Herzog, F.|textbf}, An extension of some previous results on generalized continued fractions. \emph{Duke Mathematical Journal} \textbf{12}, 655--662.
\index{f@$f$-expansions}

{\footnotesize \href{https://mathscinet.ams.org/mathscinet/article?mr=0015534}{MR0015534}}

{\footnotesize Source: \href{https://doi.org/10.1215/S0012-7094-45-01257-9}{doi.org (Project Euclid)}.}

{\footnotesize Cited in
\hyperref[Sch1992]{\textsc{Sch\"{o}nefuss 1992}}\index[ppl]{Schonefuss@Sch\"{o}nefuss, Lutz W.}.}

\vspace{4pt}

This paper builds on \hyperref[Bis1944]{\textsc{Bissinger 1944}}\index[ppl]{Bissinger, B. H.} and \hyperref[Her1945]{\textsc{Herzog and Bissinger 1945}}\index[ppl]{Herzog, F.} largely by loosening the restrictions on functions used to create $f$-expansions. 


\item\vspace{9pt} \label{Eve1946} (1946) C. J. Everett\index[ppl]{Everett, C. J.|textbf}, Representations for real numbers. \emph{Bulletin of the American Mathematical Society} \textbf{52}, 861--869.
\index{f@$f$-expansions}
\index{real numbers!f@$f$-expansions of}

{\footnotesize \href{https://mathscinet.ams.org/mathscinet/article?mr=0018221}{MR0018221}}

{\footnotesize Source: \href{https://projecteuclid.org/journals/bulletin-of-the-american-mathematical-society/volume-52/issue-10/Representations-for-real-numbers/bams/1183509722.full}{Project Euclid}.} 

{\footnotesize Cited in 
\hyperref[Rec1950]{\textsc{Rechard 1950}}\index[ppl]{Rechard, O. W.},
\hyperref[Ren1957]{\textsc{R\'{e}nyi 1957}}\index[ppl]{Renyi@R\'{e}nyi, A.}, 
\hyperref[Swa1960]{\textsc{Swartz and Wendroff 1960}}\index[ppl]{Swartz, B. K.}\index[ppl]{Wendroff, B.},
\hyperref[Thr1961]{\textsc{Thron 1961}}\index[ppl]{Thron, W. J.}, 
\hyperref[Wat1970]{\textsc{Waterman 1970}}\index[ppl]{Waterman, Michael S.}, and 
\hyperref[Sch2016]{\textsc{Schweiger 2016}}\index[ppl]{Schweiger, Fritz}.}

\vspace{4pt}

Considers Bissinger's\index[ppl]{Bissinger, B. H.} $f$-expansions using continuous strictly \emph{increasing} functions $f(t)$ on $[0, p]$ for which $f(0) = 0, f(p) = 1$. Where Bissinger generalizes continued fraction\index{continued fractions} expansions of real numbers, Everett generalizes decimal expansions\index{real numbers!decimal expansions of}; see \hyperref[Kak1924]{\textsc{Kakeya 1924}}\index[ppl]{Kakeya, S\^{o}ichi} for an earlier treatment of both topics.


\item\vspace{9pt} \label{Ogi1949} (1949) C. S. Ogilvy\index[ppl]{Ogilvy, C. S.|textbf}, Problem E 874. \emph{The American Mathematical Monthly} \textbf{56}(6), 404. Solution by C. W. Trigg, \emph{ibid.} \textbf{57}(3), 1950, 186.
\index{continued square roots!of terms $a_n=c$}

{\footnotesize Solution source: \href{https://doi.org/10.2307/2304433}{doi.org (JSTOR)}.} 

{\footnotesize Cited in 
the solution to \hyperref[Her1957]{\textsc{Herschfeld 1957}},\index[ppl]{Herschfeld, Aaron} and \hyperref[Jon2008]{\textsc{Jones 2008}}\index[ppl]{Jones, Dixon J.}.}

\vspace{4pt}

The problem essentially asks for proof that the iterated composition
\[\cdots\sqrt{x+\sqrt{x+\sqrt{x}}}\]
has an integral limit if and only if $x$ is of the form $n(n-1)$, in which case the limit is $n$. The solution generalizes this to $r$th roots, and gives $x=(n-1)n(n+1)$ for $r=3$. An editor's note points to \hyperref[Mac1941]{\textsc{MacKay 1941}} as a similar problem.


\item\vspace{9pt} \label{Rec1950} (1950) O. W. Rechard\index[ppl]{Rechard, O. W.|textbf}, The representation of real numbers, \emph{Proceedings of the American Mathematical Society} \textbf{1}, 674--681.
\index{real numbers!f@$f$-expansions of}

{\footnotesize \href{https://mathscinet.ams.org/mathscinet/article?mr=0042459}{MR0042459}}

{\footnotesize Source: \href{https://doi.org/10.1090/S0002-9939-1950-0042459-6}{doi.org (AMS)}.} 

{\footnotesize Cited in 
\hyperref[Swa1960]{\textsc{Swartz and Wendroff 1960}}\index[ppl]{Swartz, B. K.}\index[ppl]{Wendroff, B.} and 
\hyperref[Sch2016]{\textsc{Schweiger 2016}}\index[ppl]{Schweiger, Fritz}.}

\vspace{4pt}

From the introduction: ``Let $p\ge2$ be a positive integer and denote by $E_p$ the class of all continuous, strictly increasing functions $f(x)$ on the interval $0\le x\le p$ with $f(0)=0$ and $f(p)=1$. In a generalization of the decimal representation\index{real numbers!decimal expansions of}, [\hyperref[Eve1946]{\textsc{Everett 1946}}\index[ppl]{Everett, C. J.}] has associated with every real number $0\le\gamma<1$ a sequence of integers [by means of an algorithm involving] $f(x)$$\,\ldots\,$an arbitrary function in the class $E_p$.

``Some functions (for example, $f(x)=x/p$, which leads to the representation of a number as a decimal to the base $p$) when employed in [Everett's] algorithm$\,\ldots\,$yield one-one correspondences between real numbers and sequences of integers$\mod p$. On the other hand, any function, for example, whose graph has more than one point in common with any of the straight line segments connecting the points $(j, 0)$ and $(j+1,1)$, $j=0, 1, \ldots, p-1$, will obviously lead to a correspondence which is many-one. We shall denote by $E_p*$ the subclass of $E_p$ consisting of those functions which in [Everett's algorithm] give[s] rise to one-one correspondences.

``The present paper contains very simple characterizations of those correspondences between real numbers and sequences of integers mod $p$ which can be obtained by applying [Everett's algorithm] with functions from the classes $E_p*$ and $E_p-E_p*$ respectively. By means of these characterizations it is possible to settle two of the problems raised by Everett and to give an answer (albeit not a completely satisfactory one) to a third, namely that of characterizing the class $E_p*$ itself."


\item\vspace{9pt} \label{Kom1951} (1951) Karl Kommerell\index[ppl]{Kommerell, Karl}, Berechnung der trigonometrischen und zyklometrischen Funktionen durch Kettenwurzeln. \emph{Mathematisch-physikalische Semesterberichte zur Pflege des Zusammenhangs von Schule und Universit\"{a}t} \textbf{2}, 126--134.
\index{Kettenwurzeln@\emph{Kettenwurzeln}}
\index{Vi\`{e}te's formula for $\tfrac{2}{\pi}$}
\index{pi@$\pi\;(3.14159\ldots)$!continued square root expressions for}\index{constants, named!pi@$\pi\;(3.14159\ldots)$}
\index{continued square roots!of terms $a_n=\pm2$}
\index{continued square roots!and trigonometric functions}

\vspace{4pt}

The paper begins: ``In the \emph{Zeitschrift f\"{u}r mathematischen und naturwissenschaftlichen Unterricht}, Vol. 41, 1910, K. Bochow\index[ppl]{Bochow, Karl} developed formulas for the calculation of trigonometric and inverse trigonometric functions$\ldots$ Bochow has only infinite continued square roots\index{continued square roots!and trigonometric functions}; he is therefore not in a position to estimate the error which one incurs when an infinite continued square root is terminated at a certain point; furthermore, the description and the argumentation are quite cumbersome. Finally, he does not have the courage to carry out the calculation. [His] formulas $\ldots$ can now be proved very easily and can be used to calculate the functions quickly and reliably and to estimate the degree of accuracy of the calculation." The \hyperref[Cat1842]{\textsc{Catalan 1842}} formula for $\pi$ (equation \eqref{E:Cat1842-1} above) and Vi\`{e}te's\index[ppl]{Vi\`{e}te, Fran\c{c}ois} formula for $2/\pi$ (equation \eqref{E:Vie1593-1} above) are discussed.

(In a few internet links, this article is incorrectly attributed to ``Karl Koerli," perhaps due to inaccurate optical character recognition of scanned pages.) 


\item\vspace{9pt} \label{Han1955} (1955) Herman Hanisch\index[ppl]{Hanisch, Herman|textbf}, Problem 4620. \emph{The American Mathematical Monthly} \textbf{62}(1), 45. Solution by N. J. Fine\index[ppl]{Fine, N. J.|textbf}, \emph{ibid.} \textbf{63}(3), 1956, 194--195.
\index{continued square roots!Ramanujan's}

{\footnotesize Solution source: \href{https://doi.org/10.2307/2306675}{doi.org (JSTOR)}.}

This is identical to the first problem in \hyperref[Ram1911]{\textsc{Ramanujan 1911}}, as noted in an editorial comment to the solution.\index[ppl]{Ramanujan, Srinivasa}\index{continued square roots!Ramanujan's}


\item\vspace{9pt} \label{Wol1956} (1956) G. N. Wollan\index[ppl]{Wollan, G. N.|textbf} and D. M. Mesner\index[ppl]{Mesner, D. M.|textbf}, On a function defined by means of an infinite radical. Paper given at the 33rd annual meeting of the Indiana Section of the Mathematical Association of America, May 5, 1956. Abstract printed in \emph{The American Mathematical Monthly} \textbf{63}(8), 614.
\index{infinite radicals}
\index{continued square roots!of terms $a_n=\pm c$}

{\footnotesize Source: \href{https://www.jstor.org/stable/2310232}{JSTOR}.}

{\footnotesize Cited in
\hyperref[Wol1957]{\textsc{Wollan and Mesner 1957}}\index[ppl]{Wollan, G. N.}\index[ppl]{Mesner, D. M.}.}

\vspace{4pt}

From the abstract: ``If $0<x\le 1$, then $x$ has a non-terminating binary representation\index{real numbers!binary representations of} $x=a_1a_2\cdots a_n\cdots$. Let $\alpha_n=(-1)^{a_n}$ and \index{continued square roots!of terms $a_n=\pm c$}
\[f_n(x)=\sqrt{k+\alpha_1\sqrt{k+\alpha_2\sqrt{k+\cdots+\alpha_n\sqrt{k}}}}\;,\]
$n=1,2,3,\ldots$, and let $I$ denote the interval $0<x\le 1$. When $k>2+\sqrt{2}$, then $\lim_{n\to\infty}f_n(x)$ exists and is real for every $x$ in $I$ and this defines a function $f(x)=\lim_{n\to\infty} f_n(x)$. This function is discontinuous at every point in $I$ having a terminating binary representation\index{real numbers!binary representations of} and is continuous elsewhere in $I$. The function is not monotone in any sub-interval on $I$ and yet has a derivative equal to zero at each point of a dense set of points of $I$ and has a left derivative equal to zero at every point of discontinuity."


\item\vspace{9pt} \label{Her1957} (1957) Aaron Herschfeld\index[ppl]{Herschfeld, Aaron|textbf}, Problem E 1258. \emph{The American Mathematical Monthly} \textbf{64}(3), 197. Solution by D. A. Freedman, \emph{ibid.} \textbf{64}(9), 1957, 673. 

{\footnotesize Solution source: \href{https://doi.org/10.1080/00029890.1957.11988959}{doi.org (Taylor \& Francis)}.}

\vspace{4pt}

``Prove that a necessary and sufficient condition for the rationality of\index{continued cube roots}
\[R=\sqrt[3]{a+\sqrt[3]{a+\cdots}}\;,\]
where $a$ is a positive integer, is that $a=N(N+1)(N+2)$, the product of three consecutive integers. In that case find $R$." $R$ is shown to equal $N+1$. This problem was essentially solved in \hyperref[Ogi1949]{\textsc{Ogilvy 1949}}\index[ppl]{Ogilvy, C. S.}.


\item\vspace{9pt} \label{Ren1957} (1957) A. R\'{e}nyi\index[ppl]{Renyi@R\'{e}nyi, A.|textbf}, Representations for real numbers and their ergodic properties. \emph{Acta mathematica Academiae Scientiarum Hungaricae} \textbf{8}, 477--493.
\index{ergodic theory}
\index{f@$f$-expansions}
\index{beta@$\beta$-expansions}
\index{real numbers!beta@$\beta$-expansions of}
\index{real numbers!f@$f$-expansions of}
\index{real numbers!representations as continued $r$th roots}


{\footnotesize \href{https://mathscinet.ams.org/mathscinet/article?mr=0097374}{MR0097374}}

{\footnotesize Source: \href{https://doi.org/10.1007/BF02020331}{doi.org (Springer)}.} 

{\footnotesize Cited in 
\hyperref[Thr1961]{\textsc{Thron 1961}}\index[ppl]{Thron, W. J.}, 
\hyperref[Wat1970]{\textsc{Waterman 1970}} and \hyperref[Wat1975]{\textsc{1975}}\index[ppl]{Waterman, Michael S.}, and
\hyperref[Sch2016]{\textsc{Schweiger 2016}}\index[ppl]{Schweiger, Fritz}\index{f@$f$-expansions}.}

\vspace{4pt}

Puts the $f$-expansion results of \hyperref[Bis1944]{\textsc{Bissinger 1944}}\index[ppl]{Bissinger, B. H.} and \hyperref[Eve1946]{\textsc{Everett 1946}}\index[ppl]{Everett, C. J.}, as well as known theorems about regular continued fraction\index{continued fractions!regular} and decimal expansions\index{real numbers!decimal expansions of}, in a general context. An example is given expanding a real number $x\in[0,2^n-1]$ as a continued $m$th root\index{continued robb roots@continued $r$th roots}. The $\beta$-expansion of a real number is introduced. A summary of this work is given in Chapter 10 of \hyperref[Sch2016]{\textsc{Schweiger 2016}}\index[ppl]{Schweiger, Fritz}.


\item\vspace{9pt} \label{Wol1957} (1957) G. N. Wollan\index[ppl]{Wollan, G. N.|textbf} and D. M. Mesner\index[ppl]{Mesner, D. M.|textbf}, Some additional remarks on a function defined by means of an infinite radical. Paper given at the 34th annual meeting of the Indiana Section of the Mathematical Association of America, May 11, 1957. Abstract printed in \emph{The American Mathematical Monthly} \textbf{64}(8), 621. 
\index{continued square roots!of terms $a_n=\pm c$}
\index{infinite radicals}

{\footnotesize Source: \href{https://www.jstor.org/stable/2308868}{JSTOR}.} 

\vspace{4pt}

From the abstract: ``This paper presents some additional properties of the function $f(x)$ defined on $0<x\le1$ by the relation $f(x)=\lim_{n\to\infty} f_n(x)$ where
\begin{gather*}
f_1(x) = \sqrt{k + \alpha_1\sqrt{k}}\;,\quad f_2(x)=\sqrt{k+\alpha_1\sqrt{k+\alpha_2\sqrt{k}}}\;,\\
f_n(x) =\sqrt{k + \alpha_1\sqrt{k+\cdots+\alpha_n\sqrt{k}}}\quad \textrm{with $n$ nested root signs,}
\end{gather*}
$n= 1, 2,\ldots$ and $\alpha_n= (-1)^{a_n}$ where $a_n$ is the $n$th digit in the nonterminating binary representation of $x$\index{real numbers!binary representations of}. $\ldots$The authors show that when $k > 2 + \sqrt{2}$, although the function has a denumerably infinite set of discontinuities and is not monotone in any subinterval, it is of bounded variation; although it has a value at each point of the interval with $f(x_1)\ne f(x_2)$ when $x_1\ne x_2$, yet the set of values of the function is of measure zero. Furthermore the derivative exists almost everywhere and whenever it exists its value is zero, but there is a nondenumerable set of points (of measure zero) at which the derivative does not exist."


\item\vspace{9pt} \label{Myr1958} (1958) P. J. Myrberg\index[ppl]{Myrberg, P. J.|textbf}, Iteration von Quadratwurzeloperationen. \emph{Annales Academiae Scientiarum Fennicae. Series A. I, Mathematica} \textbf{259}, 16 pp.
\index{continued square roots!and trigonometric functions}

{\footnotesize \href{https://mathscinet.ams.org/mathscinet/article?mr=0108662}{MR0108662 (in German)}}

{\footnotesize Cited in
\hyperref[Sch1962]{\textsc{Schuske and Thron 1962}}\index[ppl]{Schuske, Georgellen}\index[ppl]{Thron, W. J.} and
\hyperref[Jon2008]{\textsc{Jones 2008}}\index[ppl]{Jones, Dixon J.}.}

\vspace{4pt}

The paper addresses convergence questions for the expression\index{continued square roots!of terms $a_n=\pm c$}
\[z=\pm\sqrt{p\pm\sqrt{p\pm\sqrt{p+z}}}\]
and derives formulas involving trigonometric functions of special angles. Compare \hyperref[Luc1878]{\textsc{Lucas 1878}}\index[ppl]{Lucas, Edouard}, \hyperref[Cip1908]{\textsc{Cipolla 1908}}\index[ppl]{Cipolla, Michele}, and the works of Bochow\index[ppl]{Bochow, Karl}, Wiernsberger\index[ppl]{Wiernsberger, Paul}, and Pincherle\index[ppl]{Pincherle, Salvatore}.


\item\vspace{9pt} \label{Swa1960} (1960) B. K. Swartz\index[ppl]{Swartz, B. K.|textbf} and B. Wendroff\index[ppl]{Wendroff, B.|textbf}, Continued function expansions of real numbers. \emph{Proceedings of the American Mathematical Society} \textbf{11}(4), 634--639. 
\index{f@$f$-expansions}
\index{real numbers!f@$f$-expansions of}

{\footnotesize \href{https://mathscinet.ams.org/mathscinet/article?mr=0122753}{MR0122753}}

{\footnotesize Source: \href{https://doi.org/10.1090/S0002-9939-1960-0122753-7}{doi.org (AMS)}.} 

From the introduction: ``We present a theory of continued function expansions of numbers which contains the generalized continued fractions of [\hyperref[Bis1944]{\textsc{Bissinger 1944}}]\index[ppl]{Bissinger, B. H.} and the generalized decimal representations of [\hyperref[Eve1946]{\textsc{Everett 1946}}\index[ppl]{Everett, C. J.}]$\ldots$ We generalize [the latter] by admitting a wider class of functions than those of the form $f^{-1}(x-n)$. [\hyperref[Rec1950]{\textsc{Rechard 1950}}] gave a necessary and sufficient condition that the correspondence between numbers and sequences resulting from Everett's algorithm be $1-1$. This condition appears in our theory as a simple functional relation similar to one considered [in Schreier\index[ppl]{Schreier, J.} and Ulam\index[ppl]{Ulam, S.}, Eine Bemerkung \"{u}ber die Gruppe der topologischen Abbildungen der Kreislinie auf sich selbst, \emph{Studia Mathematica} \textbf{5} (1934), 155--159]."


\item\vspace{9pt} \label{Sch1961} (1961) Georgellen Schuske\index[ppl]{Schuske, Georgellen|textbf} and W. J. Thron\index[ppl]{Thron, W. J.|textbf}, Infinite radicals in the complex plane. \emph{Proceedings of the American Mathematical Society} \textbf{12}(4), 527--532. 
\index{continued square roots!of arbitrary complex terms}
\index{infinite radicals}

{\footnotesize \href{https://mathscinet.ams.org/mathscinet/article?mr=0151586}{MR0151586}}

{\footnotesize Source: \href{https://doi.org/10.1090/S0002-9939-1961-0151586-1}{doi.org (AMS)}.} 

{\footnotesize Cited in
\hyperref[Thr1961]{\textsc{Thron 1961}}\index[ppl]{Thron, W. J.}, 
\hyperref[Sch1962]{\textsc{Schuske and Thron 1962}}\index[ppl]{Schuske, Georgellen}\index[ppl]{Thron, W. J.}, and
\hyperref[Lor1995]{\textsc{Lorentzen 1995}}\index[ppl]{Lorentzen, Lisa} and \hyperref[Lor1998]{\textsc{1998}}.}


The following theorem is proved about continued square roots with complex terms: Let $\theta$ be a fixed number, $0 < \theta < \pi$. Define $g(\theta)$ by
\[
g(\theta)=
\begin{cases}
\pi-\theta/2 &\text{when $0 < \theta \le 2\pi/3$,}\\
2(\pi-\theta) &\text{when $2\pi/3 \le \theta < \pi$.}
\end{cases}
\]
Let $\epsilon$ be an arbitrarily small fixed positive number, $0 < \epsilon < \min [\theta, g(\theta)]$. Let $a\in A_\epsilon$ if and only if
\[ |a| > 0\quad\text{and}\quad-\theta + \epsilon < \arg a < g(\theta) - \epsilon\;. \]
Then the sequence of continued square root\index{continued square roots} approximants $\tK_{i=0}^n \sqrt{a_i}$ converges if
\[a_n\in A_\epsilon \quad \text{for every $n$,}\]
and
\[\limsup_{n\to\infty}|a_n|^{2^{-n}}<\infty\;.\]
The proof makes use of the bound \eqref{E:Her1935-1} above, from \hyperref[Her1935]{\textsc{Herschfeld 1935}}\index[ppl]{Herschfeld, Aaron}\index{continued square roots!convergence conditions for}. 


\item\vspace{9pt} \label{Thr1961} (1961) W. J. Thron\index[ppl]{Thron, W. J.|textbf}, Convergence regions for continued fractions and other infinite processes. \emph{The American Mathematical Monthly} \textbf{68}, 734--750.
\index{infinite processes}
\index{continued compositions!of complex-valued functions}

{\footnotesize \href{https://mathscinet.ams.org/mathscinet/article?mr=0133444}{MR0133444}}

{\footnotesize Source: \href{https://doi.org/10.1080/00029890.1961.11989761}{doi.org (Taylor \& Francis)}.} 

{\footnotesize Cited in
\hyperref[Lor1995]{\textsc{Lorentzen 1995}} and \hyperref[Lor1998]{\textsc{1998}}\index[ppl]{Lorentzen, Lisa}, and 
\hyperref[Jon2015]{\textsc{Jones 2015}}\index[ppl]{Jones, Dixon J.}.}

\vspace{4pt}

The author begins with an arbitrary complex valued generating function
\[f(a_n^{(1)},\ldots,a_n^{(k)},z)\]
of $k+1$ complex variables and calls it $t_n(z)$. He then inductively defines an ``infinite process" as a sequence $\{T_n(z)\}$ such that
\[T_1(z) = t_1(z);\quad T_n(z) = T_{n-1}(t_n(z))\;,n > 1\;.\]
He writes, ``Through suitable choice of $z=c$ [where $c$ is the initial value], many well-known infinite processes can be obtained. Some of these are $\ldots$ infinite series\index{infinite sums}, infinite products\index{infinite products}, infinite radicals\index{infinite radicals}, infinite exponentials\index{infinite exponentials}, [and] two different kinds of continued fractions\index{continued fractions}. [W]e are concerned mainly with describing methods $\ldots$ for obtaining convergence-region criteria for certain of these infinite processes. The methods have been successfully applied to all except the first two types$\ldots$ That no nontrivial convergence regions exist for infinite series and infinite products will become clear after we have defined what a convergence region is." 


\item\vspace{9pt} \label{Sch1962} (1962) Georgellen Schuske\index[ppl]{Schuske, Georgellen|textbf} and W. J. Thron\index[ppl]{Thron, W. J.|textbf}, On periodic infinite radicals. \emph{Annales Academiae Scientiarum Fennicae. Series A. I, Mathematica} \textbf{307}, 8 pp.
\index{continued square roots!of constant complex terms}
\index{infinite radicals!of periodic complex terms}

{\footnotesize \href{https://mathscinet.ams.org/mathscinet/article?mr=0137937}{MR0137937}}

{\footnotesize Source: \href{https://www.acadsci.fi/mathematica/1962/no307pp01-08.pdf}{Academiae Scientiarum Fennicae}.}

This paper appears to have been the first to prove that\index{continued square roots!of constant complex terms}
\[\sqrt{a+\sqrt{a+\sqrt{a+\cdots}}}\]
converges to
\[\dfrac{1+\sqrt{1+4a}}{2}\]
when $a$ is any complex number, and where $\sqrt{1+4a}$ is taken in accordance with this definition: $\sqrt{x}$ is the square root whose real part is positive or 0, and is defined to lie on the positive imaginary axis if $x$ is negative. An independent rediscovery of this result is \hyperref[Siz1996]{\textsc{Sizer and Wiredu 1996}}. Compare \hyperref[Ber1692]{\textsc{Bernoulli 1692}}\index[ppl]{Bernoulli, Jacob}, \hyperref[Dop1832b]{\textsc{Doppler 1832}}\index[ppl]{Doppler, Christian}, and \hyperref[Gin1916]{\textsc{Ginsburg 1916}}\index[ppl]{Ginsburg, J. J.} for the case where $a$ is a positive real number.


\item\vspace{9pt} \label{Shk1962} (1962) D. O. Shklarsky\index[ppl]{Shklarsky, D. O.|textbf}, N. N. Chentzov\index[ppl]{Chentzov, N. N.|textbf}, and I. M. Yaglom\index[ppl]{Yaglom, I. M.|textbf}, \emph{The USSR Olympiad Problem Book: Selected Problems and Theorems of Elementary Mathematics.} W. H. Freeman, San Francisco, ISBN 978-0716704126. Reprinted by Dover Publications, Inc., Mineola, New York, 1993,  ISBN-13: 978-0486277097.
\index{continued square roots!and trigonometric functions}
\index{continued square roots!of terms $a_n=\pm2$}

{\footnotesize Source: \href{https://books.google.com/books?id=XuHCAgAAQBAJ&newbks=1&newbks_redir=0&lpg=PR1&dq=The\%20USSR\%20Olympiad\%20Problem\%20Book\%3A\%20Selected\%20Problems\%20and\%20Theorems\%20of\%20Elementary\%20Mathematics\&pg=PA44\#v=onepage\&q\&f=false}{Google Books}.}

{\footnotesize Cited in 
\hyperref[Eft2012]{\textsc{Efthimiou 2012}} and \hyperref[Eft2013]{\textsc{2013}}\index[ppl]{Efthimiou, Costas J.}.}

\vspace{4pt}

Problem 195 reads: ``Let some of the numbers $a_1$, $a_2,\ldots$, $a_n$ be $+1$ and the rest be $-1$. Prove that\index{continued square roots!and trigonometric functions}
\begin{align*}
&2\sin\left(a_1+\dfrac{a_1a_2}{2}+\dfrac{a_1a_2a_3}{4}+\cdots+\dfrac{a_1a_2\cdots a_n}{2^{n-1}}\right)45^{\circ} \\
&=a_1\sqrt{2+a_2\sqrt{2+a_3\sqrt{2+\cdots a_n\sqrt{2}}}}\;.
\end{align*}
For example, let $a_1=a_2=$$\cdots=a_n=1$:
\begin{align*}
2\sin\left(1+\dfrac{1}{2}+\dfrac{1}{4}+\cdots+\dfrac{1}{2^{n-1}}\right)45^{\circ}&=2\cos\dfrac{45^{\circ}}{2^{n-1}}\\
&=\sqrt{2+\sqrt{2+\cdots+\sqrt{2}}}\;."
\end{align*}

The problem, which is unattributed, is essentially the same as Problem 183 in \hyperref[Pol1925]{\textsc{P\'{o}lya and Szeg\"{o} 1925}}\index[ppl]{Polya@P\'{o}lya, G.}; see equation \eqref{E:Pol1925-1} above. Compare also \hyperref[Boc1899]{\textsc{Bochow 1899}}\index[ppl]{Bochow, Karl} and \hyperref[Cip1908]{\textsc{Cipolla 1908}}\index[ppl]{Cipolla, Michele}. 


\item\vspace{9pt} \label{McK1967} (1967) J. H. McKay, The William Lowell Putnam Mathematical Competition. \emph{The American Mathematical Monthly} \textbf{74}(7), 771--777.
\index{continued square roots!Ramanujan's}

{\footnotesize Source: \href{https://doi.org/10.1080/00029890.1967.12000030}{doi.org (Taylor \& Francis)}.}

{\footnotesize Cited in
\hyperref[Bor1991]{\textsc{Borwein and de Barra 1991}}\index[ppl]{Borwein, Jonathan M.}\index[ppl]{de Barra, G.}.}

\vspace{4pt}

Putnam Exam problem A-6 asks for ``justification" of the identity
\begin{equation}\label{E:McK1967-1}
3=\sqrt{1+2\sqrt{1+3\sqrt{1+\cdots}}}\;.
\end{equation}
The original source, \hyperref[Ram1911]{\textsc{Ramanujan 1911}}, is not cited in the exam or in the published solution.


\item\vspace{9pt} \label{Ogi1970} (1970) C. S. Ogilvy\index[ppl]{Ogilvy, C. S.|textbf}, To what limits do complex iterated radicals converge? \emph{The American Mathematical Monthly} \textbf{77}(4), 388--389. 
\index{continued square roots!of constant complex terms}

{\footnotesize \href{https://mathscinet.ams.org/mathscinet/article?mr=1535864}{MR1535864 (citation only)}}

{\footnotesize Source: \href{https://doi.org/10.1080/00029890.1970.11992497}{doi.org (Taylor \& Francis)}.} 

{\footnotesize Cited in 
\hyperref[Ger1973]{\textsc{Gerber 1973}}\index[ppl]{Gerber, Leon}, 
\hyperref[Roh1974]{\textsc{Rohde 1974}}\index[ppl]{Rohde, Hanns-Walter}, 
\hyperref[Wal1983]{\textsc{Walker 1983}}\index[ppl]{Walker, Peter L.}, 
\hyperref[Jon1991]{\textsc{Jones 1991}}\index[ppl]{Jones, Dixon J.}, and
\hyperref[Lau1999]{\textsc{Laugwitz and Sch\"{o}nefuss 1999}}\index[ppl]{Laugwitz, Detlef}\index[ppl]{Schonefuss@Sch\"{o}nefuss, Lutz W.}.}

\vspace{4pt}

This short note, printed in the old ``Research Problems" section of the \emph{Monthly}, asked a number of questions concerning the limiting value of $\sqrt{a+\sqrt{a+\cdots}}$ for complex and negative values of $a$. Some answers to these questions had already been given in \hyperref[Sch1961]{\textsc{Schuske and Thron 1961}} and \hyperref[Sch1962]{\textsc{1962}}\index[ppl]{Thron, W. J.}\index[ppl]{Schuske, Georgellen}.


\item\vspace{9pt} \label{Wat1970} (1970) Michael S. Waterman\index[ppl]{Waterman, Michael S.|textbf}, Some ergodic properties of multi-dimensional $f$-expansions. \emph{Zeitschrift f\"{u}r Wahrscheinlichkeitstheorie und Verwandte Gebiete} \textbf{16}, 77--103.
\index{f@$f$-expansions!multidimensional}
\index{ergodic theory}

{\footnotesize \href{https://mathscinet.ams.org/mathscinet/article?mr=0282939}{MR0282939}}

{\footnotesize Source: \href{https://doi.org/10.1007/BF00535691}{doi.org (Springer)}.}

{\footnotesize Cited in 
\hyperref[Sch2016]{\textsc{Schweiger 2016}}\index[ppl]{Schweiger, Fritz}.}

\vspace{4pt}

From the introduction: ``This paper is concerned with probabilistic aspects of the expansion of points in $n$-dimensional Euclidean space. The expansions we consider need not converge although previous work has required convergence$\ldots$ In 1869 Jacobi presented an extension of the continued fraction\index{continued fractions!multidimensional} to two dimensions. Perron extended Jacobi's work to $n$-dimensions. In 1964 Schweiger\index[ppl]{Schweiger, Fritz} began an examination of the measure theoretic properties of Jacobi's algorithm. It was this work which motivated our paper. However Schweiger has $\ldots$ published some results which also concern general $F$-expansions for $n$-dimensions. The class of algorithms he considers does not include the Jacobi algorithm and is a natural generalization of R\'{e}nyi. Our results generalize most of Schweiger's work and have the Jacobi algorithm as an example. We also include a central limit theorem and a law of the iterated logarithm."


\item\vspace{9pt} \label{Leh1971} (1971) D. H. Lehmer\index[ppl]{Lehmer, D. H.|textbf}, On the compounding of certain means. \emph{Journal of Mathematical Analysis and Applications} \textbf{36}, 183--200.
\index{continued square roots!of terms $a_n=\pm2$}

{\footnotesize \href{https://mathscinet.ams.org/mathscinet/article?mr=0281696}{MR0281696}}

{\footnotesize Source: \href{https://doi.org/10.1016/0022-247X(71)90029-1}{doi.org (ScienceDirect)}.}

From the introduction: ``By the compound of two means $M(u,v)$ and $M'(u,v)$ is meant the function $M\times M'(u,v)$ defined as follows. Let the two sequences $\{a_n\}$ and $\{b_n\}$ be constructed recursively by\index{means!recursive construction of}
\begin{align*}
a_0=u\;,&\quad b_0=v\\
a_{i+1}=M(a_i,b_i)\;,&\quad b_{i+1}=M'(a_i,b_i)\;.
\end{align*}
Then $M\times M'(u,v)=\lim_{n\to\infty}a_n=\lim_{n\to\infty}b_n\;,$ provided both limits exist and are equal." Lehmer\index[ppl]{Lehmer, D. H.} defines two classes of means, $\mu$ and $M$ (both generalized from $A,G$, and $H$, the arithmetic\index{means!arithmetic}, geometric\index{means!geometric} and harmonic\index{means!harmonic} means, respectively) and considers compounds\index{means!compounding of} of means from each class. He develops a Taylor expansion\index{Taylor series}\index{series!Taylor} for the mean $R(u,v)=M_1\times M_2(u,v)$ (where both are in $M$, and $M_1=A$), and shows that the series coefficients $g_n$ are even integers which ``satisfy recurrences vaguely like [those for] the Bernoulli numbers\index{Bernoulli numbers}\index{numbers!Bernoulli}." Using polynomials $P_n(x)$ defined recursively by $P_0(x)=1$ and $P_{n+1}(x)=P_n^2(x)+2x^{2^n}$, he gives an alternate formulation for the $g_n$. For calculating the radius of convergence of $R$, he defines an additional polynomial $Q_n(x)=x^{2^{n-1}}P_n(1/x)$, and employs the following lemma: ``The $2^n$ roots of $Q_{n+1}(x)$ are the $2^n$ values of the [continued square root]\index{continued square roots!of terms $a_n=c$}
\[-2\pm\sqrt{-2\pm\sqrt{-2\pm\ldots\pm\sqrt{-2}}}\;,\]
where $n$ square roots are given signs in all $2^n$ ways."


\item\vspace{9pt} \label{Ger1973} (1973) Leon Gerber\index[ppl]{Gerber, Leon|textbf}, Complex iterated radicals. \emph{Proceedings of the American Mathematical Society} \textbf{41}(1), 205--210.
\index{continued square roots!of constant complex terms}
\index{iterated square roots!of constant complex terms}

{\footnotesize \href{https://mathscinet.ams.org/mathscinet/article?mr=0318721}{MR0318721}}

{\footnotesize Source: \href{https://doi.org/10.2307/2038842}{doi.org (JSTOR)}.}

{\footnotesize Cited in 
\hyperref[Wal1983]{\textsc{Walker 1983}}\index[ppl]{Walker, Peter L.}.} 

\vspace{4pt}

This is the first of at least three papers addressing the problems posed in \hyperref[Ogi1970]{\textsc{Ogilvy 1970}}\index[ppl]{Ogilvy, C. S.}, all apparently unaware of the joint papers of Schuske and Thron\index[ppl]{Thron, W. J.} from 1961 and 1962. From the abstract: ``We prove the convergence of the sequence $S$ defined by $z_{n+1} = (z_n-c)^{1/2}$, $c$ real, for any choice of $z_0$. Let $k=|\tfrac{1}{4}-c|^{1/2}$. If $c<0$ or $c=\tfrac{1}{4}$, $S$ has only one fixed point $w=\tfrac{1}{2}+k$ and converges to $w$ for any $z_0$. If $0\le c<\tfrac{1}{4}$, $S$ has the fixed points $w_1=\tfrac{1}{2}+k$ and $w_2=\tfrac{1}{2}-k$, and for any $z_0\ne w_2$, $S$ converges to $w_1$. If $c>\tfrac{1}{4}$, $S$ has the fixed points $w_1=\tfrac{1}{2}+ik$ and $w_2=\tfrac{1}{2}-ik$ and converges to $w_1$ if $\operatorname{Re}(z_0)\ge0$ and to $w_2$ otherwise. We show that convergence is strictly monotone when the neighborhood system is the pencil of coaxial circles with $w_1$ and $w_2$ as limiting points, and give rates of convergence." \hyperref[Ogi1970]{\textsc{Ogilvy 1970}}\index[ppl]{Ogilvy, C. S.} is the only reference cited.


\item\vspace{9pt} \label{Guy1973} (1973) R. K. Guy\index[ppl]{Guy, R. K.|textbf} and J. L. Selfridge\index[ppl]{Selfridge, J. L.|textbf}, The nesting and roosting habits of the laddered parenthesis. \emph{The American Mathematical Monthly} \textbf{80}(8), 868--876.
\index{associativity of function compositions}

{\footnotesize \href{https://mathscinet.ams.org/mathscinet/article?mr=0347625}{MR0347625}}

{\footnotesize Source: \href{https://doi.org/10.1080/00029890.1973.11993395}{doi.org (Taylor \& Francis)}.}

This paper considers combinatorial aspects of the associativity of the continued exponential\index{continued exponentials}
\[a^{a^{\cdot^{\cdot^{\cdot^{a}}}}}\]
with $k$ $a$'s, which, as the authors say, ``is ambiguous until the order of the $k-1$ operations has been indicated, say by the insertion of $k-2$ pairs of parentheses." (Although continued exponentials have been excluded from this bibliography, this reference is included because the nature of the underlying function used in this paper seems of less consequence than the issues of associativity that accompany compositions of multiple functions.)


\item\vspace{9pt} \label{Roh1974} (1974) Hanns-Walter Rohde\index[ppl]{Rohde, Hanns-Walter|textbf}, Complex iterated radicals. \emph{The American Mathematical Monthly} \textbf{81}(1), 14--21.
\index{continued square roots!of constant complex terms}
\index{iterated square roots!of constant complex terms}
\index{continued rocc roots@continued $r_i$th roots!of arbitrary complex terms}

{\footnotesize \href{https://mathscinet.ams.org/mathscinet/article?mr=0338333}{MR0338333}}

{\footnotesize Source: \href{https://doi.org/10.1080/00029890.1974.11993490}{doi.org (Taylor \& Francis)}.}

{\footnotesize Cited in 
\hyperref[Wal1983]{\textsc{Walker 1983}}\index[ppl]{Walker, Peter L.}.}

\vspace{4pt}

This is the second paper, essentially contemporaneous with \hyperref[Ger1973]{\textsc{Gerber 1973}}\index[ppl]{Gerber, Leon}, to address the questions in \hyperref[Ogi1970]{\textsc{Ogilvy 1970}}\index[ppl]{Ogilvy, C. S.}. The author defines a sequence $\{f_n\}$ inductively by $f_0(z)=0$ and $f_n(z)=(a_n(z)+f_{n-1}(z))^{r_n}$, $n=1,2,3,\ldots$, where the $a_n(z)$ are functions holomorphic in a simply connected domain $D$ of the complex plane, and the $r_n$ are complex numbers. (Note that these $f_n$ are defined as iterated compositions; see Section \ref{S:assoc} above). Convergence of the $f_n$ is established ``in a simple way by means of a classical theorem due to Vitali\index{Vitali's theorem} $\ldots$ and the Monodromy theorem\index{monodromy theorem}." The questions in \hyperref[Ogi1970]{\textsc{Ogilvy 1970}}\index[ppl]{Ogilvy, C. S.} are addressed by letting $r_n=\tfrac{1}{2}$ and $a_n(z)=z$ for $n=1,2,3,\ldots$ \hyperref[Her1935]{\textsc{Herschfeld 1935}}\index[ppl]{Herschfeld, Aaron} is cited, but not the papers of Schuske and Thron\index[ppl]{Thron, W. J.}\index[ppl]{Schuske, Georgellen} from the early 1960s.


\item\vspace{9pt} \label{Spi1974} (1974) Robert Spira\index[ppl]{Spira, Robert|textbf} and Alexander Scheeline\index[ppl]{Scheeline, Alexander|textbf}, Table errata: ``A cotangent analogue of continued fractions" (Duke Mathematical Journal 4 (1938), 323--340) by D. H. Lehmer\index[ppl]{Lehmer, D. H.}. \emph{Mathematics of Computation} \textbf{28}, 677.
\index{continued cotangents}

{\footnotesize \href{https://mathscinet.ams.org/mathscinet/article?mr=0345897}{MR0345897}}

{\footnotesize Source: \href{https://doi.org/10.1090/S0025-5718-1974-0345897-1}{doi.org (AMS)}.}

The authors write: ``A new calculation of D. H. Lehmer's constant $\xi$\index{Lehmer's constant $(0.59263\ldots)$} to 1092D has revealed that the 80D value recorded on p. 334 of the original paper [\emph{Duke Mathematical Journal} \textbf{4} (1938), 323--340; Zbl \textbf{19}, 9] is correct to only 71D, corresponding to the eighth convergent to the infinite simple continued fraction appearing in equation (36) on the same page.''


\item\vspace{9pt} \label{Wat1975} (1975) M. S. Waterman\index[ppl]{Waterman, Michael S.|textbf}, $F$-expansions of rationals. \emph{Aequationes Mathematicae} \textbf{13}, 263--268.
\index{f@$f$-expansions}

{\footnotesize \href{https://mathscinet.ams.org/mathscinet/article?mr=0396457}{MR0396457}} 

{\footnotesize Source: \href{https://doi.org/10.1007/BF01836529}{doi.org (Springer)}.}

From the abstract: ``The ergodic theorem has been used to deduce results about the $F$-expansions of almost all $x$ in $(0, 1)^n$. A simple lemma from measure theory yields some corresponding statements about the expansions of the rationals, a set of measure zero."


\item\vspace{9pt} \label{Sha1976} (1976) Jeffrey Shallit\index[ppl]{Shallit, Jeffrey|textbf}, Predictable regular continued cotangent expansions. \emph{Journal of Research of the National Bureau of Standards, Section B} \textbf{80B}(2), 285--290.
\index{continued cotangents}

{\footnotesize \href{https://mathscinet.ams.org/mathscinet/article?mr=0429723}{MR0429723}}

{\footnotesize Source: \href{https://archive.org/details/jresv80Bn2p285}{Internet Archive}.}

{\footnotesize Cited in
\hyperref[Riv2007]{\textsc{Rivoal 2007}}\index[ppl]{Rivoal, T.}.}

\vspace{4pt}

From the abstract: ``Expansions of the form
\[
x = \cot(\arccot n_0 - \arccot n_1 + \arccot n_2 - \cdots)
\]
are discussed. It is shown that if $x$ is of the form $\tfrac{1}{2}(c + \sqrt{c^2 + 4})$, then the $n$s are predictable by a simple recurrence\index{recurrence relation}. Continued fractions\index{continued fractions} derived from the expansion of $x$ are also given."


\item\vspace{9pt} \label{Sto1980} (1980) Kenneth B. Stolarsky, Mapping properties, growth, and uniqueness of Vieta (infinite cosine) products. \emph{Pacific Journal of Mathematics} \textbf{89}(1), 209--227.
\index{Vi\`{e}te's formula for $\tfrac{2}{\pi}$!generalizations}
\index[ppl]{Vi\`{e}te, Fran\c{c}ois}
\index{infinite products!of cosine functions}

{\footnotesize \href{https://mathscinet.ams.org/mathscinet/article?mr=0596932}{MR0596932}}

{\footnotesize Source: \href{http://dx.doi.org/10.2140/pjm.1980.89.209}{dx.doi.org (msp.org)}.}

{\footnotesize Cited in
\hyperref[Lev2005]{\textsc{Levin 2005}}\index[ppl]{Levin, Aaron}.}

\vspace{4pt}

From the abstract: ``The natural logarithm of $z$ can be written as an infinite product involving iterated square roots of $z.$ A \emph{Vieta product} is defined to be a more general infinite product involving $z$ raised to arbitrary fractional powers. Restricted to the unit circle, Vieta products generalize infinite cosine products studied$\,\ldots\,$in connection with [Pisot-Vijayaraghavan] numbers\index{Pisot-Vijayaraghavan numbers}\index{numbers!Pisot-Vijayaraghavan}. Vieta products are shown to have conformal mapping, monotonicity, and growth properties very similar to those of the natural logarithm. By using certain properties of Eulerian polynomials\index{Eulerian polynomials}\index{polynomials!Eulerian}, the exponents of $z$ in a Vieta product are shown to be unique in a strong sense." (Pisot-Vijayaraghavan or PV-numbers are real algebraic integers greater than 1 whose Galois conjugates\index{Galois conjugates} have absolute value less than 1.)


\item\vspace{9pt} \label{Den1983} (1983) Thomas P. Dence\index[ppl]{Dence, Thomas P.|textbf}, Problem 1174. \emph{Mathematics Magazine} \textbf{56}, 178. Solution by the Chico Problem Group, \emph{ibid.} \textbf{57}, 1984, 299--300.
\index{continued square roots!of terms $a_n=\pm c$}


{\footnotesize Solution source: \href{https://doi.org/10.1080/0025570X.1984.11977132}{doi.org (Taylor \& Francis)}.}

{\footnotesize Cited in
\hyperref[Siz1986]{\textsc{Sizer 1986}}\index[ppl]{Sizer, Walter S.},
\hyperref[Bor1991]{\textsc{Borwein and de Barra 1991}}\index[ppl]{Borwein, Jonathan M.}\index[ppl]{de Barra, G.},
\hyperref[Jon2008]{\textsc{Jones 2008}}\index[ppl]{Jones, Dixon J.}, and 
\hyperref[Eft2012]{\textsc{Efthimiou 2012}}\index[ppl]{Efthimiou, Costas J.}.}

\vspace{4pt}

The problem asks for real numbers $A$ and $k$ such that the sequence\index{continued square roots!of periodic real terms}
\[\sqrt{k},\;\sqrt{k-\sqrt{k}},\;\sqrt{k-\sqrt{k+\sqrt{k}}},\;\sqrt{k-\sqrt{k+\sqrt{k-\sqrt{k}}}},\ldots\]
converges to $A$, and to write $k$ explicitly in terms of $A$. The solution, adapted by the problems section editor, is as follows: ``If $k=0$, the sequence converges to $0$. Otherwise, for all $k\ge k_0\approx 1.7548777$ (where $k_0$ is the unique positive solution to $k^3-2k^2+k-1=0$), the sequence converges to $A=\sqrt{k-\tfrac{3}{4}}-\tfrac{1}{2}$, whence $k=A^2+A+1$. Moreover, for each $A\ge A_0=\sqrt{k_0-\tfrac{3}{4}}-\tfrac{1}{2}\approx 0.5024359$, as well as for $A=0$, but for no other value of $A$, there is a (unique) value of $k$ for which the sequence has limit $A$." Much the same question is asked, and the value $1.7548777$ independently derived, in \hyperref[Zim2008a]{\textsc{Zimmerman and Ho 2008}}\index[ppl]{Zimmerman, Seth}\index[ppl]{Ho, Chungwu}.


\item\vspace{9pt} \label{Wal1983} (1983) Peter L. Walker\index[ppl]{Walker, Peter L.|textbf}, Iterated complex radicals. \emph{The Mathematical Gazette} \textbf{67}, 269--273. 
\index{continued square roots!of constant complex terms}
\index{iterated square roots!of constant complex terms}

{\footnotesize Source: \href{https://doi.org/10.2307/3617262}{doi.org (CambridgeCore)}.}

{\footnotesize Cited in
\hyperref[Jon1991]{\textsc{Jones 1991}}\index[ppl]{Jones, Dixon J.}.}

\vspace{4pt}

From the introduction: ``The question of the convergence of the sequence defined recursively by
\[z_{n+1}= \sqrt{c + z_n}\quad (n\ge0)\]
for a given complex value of $c$, and a fixed determination of the square root was posed [in \hyperref[Ogi1970]{\textsc{Ogilvy 1970}}\index[ppl]{Ogilvy, C. S.}] and partial solutions have been given [in \hyperref[Roh1974]{\textsc{Rohde 1974}} and \hyperref[Ger1973]{\textsc{Gerber 1973}}\index[ppl]{Gerber, Leon}]$\ldots\,$ It turns out that if $c$ and $z_n$, $(n\ge 1)$ are chosen to have non-negative imaginary parts, then the sequence converges for any initial $z_0$. This note gives two solutions to the problem: the first uses some complex analytic machinery (Schwarz'[s] Lemma\index{Schwartz's Lemma} and the Riemann Mapping Theorem\index{Riemann Mapping Theorem}) while the second is entirely elementary, using only the properties of bilinear mappings\index{bilinear mappings}."


\item\vspace{9pt} \label{All1985} (1985) Edward L. Allen\index[ppl]{Allen, Edward L.|textbf}, Continued radicals. \emph{The Mathematical Gazette} \textbf{69}, 261--263. 
\index{continued robb roots@continued $r$th roots!of terms $a_n=c$}

{\footnotesize Source: \href{https://doi.org/10.2307/3617569}{doi.org (CambridgeCore)}.}

{\footnotesize Cited in
\hyperref[Rao2005]{\textsc{Rao and Vanden Berghe 2005}}\index[ppl]{Rao, K. Srinivasa}\index[ppl]{Vanden Berghe, G.} and 
\hyperref[Cla2014]{\textsc{Clark and Richmond 2014}}\index[ppl]{Clark, Tyler}\index[ppl]{Richmond, Tom}.}

\vspace{4pt}

Shows that, for any integer $r>1$,
\[\sqrt[r]{a+\sqrt[r]{a+\sqrt[r]{\cdots}}}\]
converges to the unique positive root of $x^r-x-a=0$; mentions that this root is $k$ if $a=k(k-1)$ for $r=2$ (compare \hyperref[Ogi1949]{\textsc{Ogilvy 1949}}\index[ppl]{Ogilvy, C. S.}); and proves that the partial products of the reciprocal of Vi\`{e}te's\index[ppl]{Vi\`{e}te, Fran\c{c}ois} formula for $\tfrac{2}{\pi}$, equation \eqref{E:Vie1593-1} above, are increasing and bounded. No literature is cited.
\index{Vi\`{e}te's formula for $\tfrac{2}{\pi}$}


\item\vspace{9pt} \label{And1985} (1985) R. L. Andrushkiw\index[ppl]{Andrushkiw, R. L.|textbf}, On the convergence of continued radicals with applications to polynomial equations. \emph{Journal of the Franklin Institute} \textbf{319}, 391--396.
\index{continued rocc roots@continued $r_i$th roots!as solutions to polynomial equations}

{\footnotesize \href{https://mathscinet.ams.org/mathscinet/article?mr=0790965}{MR0790965}}

{\footnotesize Source: \href{https://doi.org/10.1016/0016-0032(85)90007-9}{doi.org (ScienceDirect)}.}

{\footnotesize Cited in
\hyperref[Joh2008]{\textsc{Johnson and Richmond 2008}}\index[ppl]{Johnson, Jamie}\index[ppl]{Richmond, Tom}, 
\hyperref[Eft2012]{\textsc{Efthimiou 2012}}\index[ppl]{Efthimiou, Costas J.}, 
\hyperref[Cla2014]{\textsc{Clark and Richmond 2014}}\index[ppl]{Clark, Tyler}\index[ppl]{Richmond, Tom}, and 
\hyperref[Lyn2014]{\textsc{Lynd 2014}}\index[ppl]{Lynd, Chris D.}.}

\vspace{4pt}

The paper begins with an independent derivation of the result in \hyperref[All1985]{\textsc{Allen 1985}}\index[ppl]{Allen, Edward L.} for the convergence of $\sqrt[r]{a+\sqrt[r]{a+\sqrt[r]{\cdots}}}$\index{continued robb roots@continued $r$th roots!of terms $a_n=1$} in the special case $a=1$. It then develops the following generalization of \hyperref[Pol1916]{\textsc{P\'{o}lya 1916}}\index[ppl]{Polya@P\'{o}lya, G.}: Given $a_n\ge 0$, integers $r_n>1$, let
\[\alpha=\limsup_{a_n\in S}\left(\dfrac{\log\log a_n}{n}\right), \quad S=\{a_n\;|\;a_n>1, n\in\mathbb{Z}^+\}\;,\]
\[r=\liminf_{n\to\infty} r_n\;,\quad\textrm{and} \quad R=\limsup_{n\to\infty}r_n\;.\]
Then\index{continued rocc roots@continued $r_i$th roots!convergence conditions for}
\[\sqrt[r_1]{a_1+\sqrt[r_2]{a_2+\sqrt[r_3]{a_3+\cdots}}}\]
converges if $\alpha<\log r$; diverges if $\alpha>\log R$; may converge or diverge if $\alpha\in[\log r, \log R]$; and converges if $a_n\in[0,1]$ for all $n\in\mathbb{Z}^+$. (Compare \hyperref[Her1935]{\textsc{Herschfeld 1935, Theorem III}}.)\index[ppl]{Herschfeld, Aaron}\index{continued rocc roots@continued $r_i$th roots!convergence conditions for} The author gives an error estimate for the approximants, and shows how roots of some polynomial equations\index{continued rocc roots@continued $r_i$th roots!as solutions to polynomial equations} can be approximated by continued $r_i$th roots. The only reference given is \hyperref[Pol1925]{\textsc{P\'{o}lya and Szeg\"{o} 1925}}\index[ppl]{Polya@P\'{o}lya, G.}.


\item\vspace{9pt} \label{Siz1986} (1986) Walter S. Sizer\index[ppl]{Sizer, Walter S.|textbf}, Continued roots. \emph{Mathematics Magazine} \textbf{59}, 23--27.
\index{continued square roots!convergence conditions for}

{\footnotesize \href{https://mathscinet.ams.org/mathscinet/article?mr=0828417}{MR0828417}}

{\footnotesize Source: \href{https://doi.org/10.1080/0025570X.1986.11977215}{doi.org (Taylor \& Francis)}.}

{\footnotesize Cited in 
\hyperref[Lau1990]{\textsc{Laugwitz 1990}}\index[ppl]{Laugwitz, Detlef}, 
\hyperref[Bor1991]{\textsc{Borwein and de Barra 1991}}\index[ppl]{Borwein, Jonathan M.}\index[ppl]{de Barra, G.}, 
\hyperref[Pic1991]{\textsc{Pickover and Lakhtakia 1991}}\index[ppl]{Pickover, C. A.}\index[ppl]{Lakhtakia, A.}, 
\hyperref[Jon1991]{\textsc{Jones 1991}}\index[ppl]{Jones, Dixon J.}, 
\hyperref[Sch1992]{\textsc{Sch\"{o}nefuss 1992}}\index[ppl]{Schonefuss@Sch\"{o}nefuss, Lutz W.}, 
\hyperref[Jon1995]{\textsc{Jones 1995}}\index[ppl]{Jones, Dixon J.}, 
\hyperref[Siz1996]{\textsc{Sizer and Wiredu 1996}}\index[ppl]{Sizer, Walter S.}\index[ppl]{Wiredu, E. K.}, 
\hyperref[Lau1999]{\textsc{Laugwitz and Sch\"{o}nefuss 1999}}\index[ppl]{Laugwitz, Detlef}\index[ppl]{Schonefuss@Sch\"{o}nefuss, Lutz W.}, 
\hyperref[Joh2008]{\textsc{Johnson and Richmond 2008}}\index[ppl]{Johnson, Jamie}\index[ppl]{Richmond, Tom}, 
\hyperref[Zim2008a]{\textsc{Zimmerman and Ho 2008a}}\index[ppl]{Zimmerman, Seth}\index[ppl]{Ho, Chungwu}, 
\hyperref[Lim2010]{\textsc{Lim 2010}}\index[ppl]{Lim, Teik-Cheng}, 
\hyperref[Eft2012]{\textsc{Efthimiou 2012}}\index[ppl]{Efthimiou, Costas J.}, 
\hyperref[Glu2012]{\textsc{Gluzman and Yukalov 2012}}\index[ppl]{Gluzman, S.}\index[ppl]{Yukalov, V. I.},
\hyperref[Muk2013]{\textsc{Mukherjee 2013}}\index[ppl]{Mukherjee, Soumendu Sundar}, 
\hyperref[Cla2014]{\textsc{Clark and Richmond 2014}}\index[ppl]{Clark, Tyler}\index[ppl]{Richmond, Tom}, 
\hyperref[Lyn2014]{\textsc{Lynd 2014}}\index[ppl]{Lynd, Chris D.}, 
\hyperref[Les2016]{\textsc{Lesher and Lynd 2016}}\index[ppl]{Lesher, Devyn A.}\index[ppl]{Lynd, Chris D.}, 
\hyperref[Vel2016c]{\textsc{Vellucci and Bersani 2016c}}\index[ppl]{Vellucci, Pierluigi}\index[ppl]{Bersani, Alberto Maria}, and
\hyperref[Wei]{\textsc{Weisstein n.d.}\index[ppl]{Weisstein, Eric W.}}}

\vspace{4pt}

Presents an independent rediscovery of the continued square root\index{continued square roots} convergence condition given in \hyperref[Wie1904a]{\textsc{Wiernsberger 1904a}} (equation \eqref{E:Wie1904a-3} above)\index{continued square roots!convergence conditions for}, and proves that any nonnegative real number can be represented as a continued square root with terms $\{a_i\}_{i=0}^\infty$ where, for $i\ge 1$, $a_i$ is 0, 1, or 2.\index{real numbers!representations as continued square roots}


\item\vspace{9pt} \label{Par1987} (1987) R. B. Paris\index[ppl]{Paris, R. B.|textbf}, An asymptotic approximation connected with the golden number. \emph{American Mathematical Monthly} \textbf{94}, 272--278.
\index{continued square roots!of terms $a_n=1$}
\index{continued square roots!of terms $a_n=c$}

{\footnotesize \href{https://mathscinet.ams.org/mathscinet/article?mr=883295}{MR883295}}

{\footnotesize Source: \href{https://doi.org/10.1080/00029890.1987.12000627}{doi.org (Taylor \& Francis)}.}

{\footnotesize Cited in 
\hyperref[Fin2003]{\textsc{Finch 2003}}\index[ppl]{Finch, Steven R.}.}

\vspace{4pt}

The author writes, ``The subject of this note is the repeated square root sequence $\{u_n\}$ defined by
\begin{equation}\label{E:Par1987-1}
u_1=1,\quad u_n=(1+u_{n-1})^{1/2}\;,\;\;(n\ge2)
\end{equation}
the $n$th term of which can be represented in the form of $n-1$ nested square roots as
\[u_n=\sqrt{1+\cdots+\sqrt{1+\sqrt{1+1}}}\;.\]
[\ldots]Here, we shall not be concerned with the actual limiting value of \eqref{E:Par1987-1} but with the considerably more subtle problem of the determination of the manner in which $u_n$ approaches the limit $\theta\,[=(1+\sqrt{5})/2]$."\index{golden ratio ($\tfrac{1+\sqrt{5}}{2}=1.61803\ldots$)}\index{constants, named!golden ratio ($\tfrac{1+\sqrt{5}}{2}=1.61803\ldots$)} It is shown that
\[\theta-u_n \sim \dfrac{2K}{(2\theta)^n}\]
where $K$ is a constant (now known as the Paris constant\index{Paris constant ($1.09863\ldots$)}) with approximate value $1.098630$. In section 3, the sequence \eqref{E:Par1987-1} is generalized ``by taking $u_1=x$, where $x$ denotes a real variable, and defining an iterated function $u_n(x)$ by
\[u_1=x,\quad u_n(x)=\{\lambda+\mu u_{n-1}(x)\}^{1/2}\;,\;\;(n\ge2)\]
where $\lambda$ and $\mu$ are real constants.\,\ldots\,If $\mu^2+4\lambda\ge0$, the limit of the sequence, when it exists, is real and independent of $x$ and is given by
\[u_\infty=\dfrac{1}{2}\mu+\dfrac{1}{2}(\mu^2+4\lambda)^{1/2}."\]
The paper concludes with two special cases, one of which has $\lambda=2\mu^2$ and $\mu>0$, leading to a hypergeometric function which can be expressed in terms of elementary functions. In the particular case $\lambda=\mu=\tfrac{1}{2}$, an algebraic form of Euler's product \eqref{E:Eul1744-2} is derived as a function of $x$, which in turn yields Vi\`{e}te's formula \eqref{E:Vie1593-1} when $x=0$.\index{Vi\`{e}te's formula for $\tfrac{2}{\pi}$}


\item\vspace{9pt} \label{Cas1988} (1988) Dario Castellanos\index[ppl]{Castellanos, Dario|textbf}, The ubiquitous $\pi$, part I. \emph{Mathematics Magazine} \textbf{61}, 67--98.
\index[ppl]{Vi\`{e}te, Fran\c{c}ois}
\index{Vi\`{e}te's formula for $\tfrac{2}{\pi}$}

{\footnotesize \href{https://mathscinet.ams.org/mathscinet/article?mr=0934824}{MR0934824}}

{\footnotesize Source: \href{https://doi.org/10.1080/0025570X.1988.11977363}{doi.org (Taylor \& Francis)}.}

{\footnotesize Cited in
\hyperref[Ber1994]{\textsc{Berndt 1994}}\index[ppl]{Berndt, Bruce C.}.}

\vspace{4pt}

Section 2 derives Vi\`{e}te's formula for $\tfrac{2}{\pi}$ (equation \eqref{E:Vie1593-1} above), starting with a finite product formula for $\sin\theta/\theta$, which becomes Euler's formula \eqref{E:Eul1763-1} above (given in \hyperref[Rud1891]{\textsc{Rudio 1891}}) as the number of terms increase without bound. Vi\`{e}te's formula is the case $\theta=\tfrac{\pi}{2}$.


\item\vspace{9pt} \label {Ber1989} (1989) Bruce C. Berndt\index[ppl]{Berndt, Bruce C.|textbf}, \emph{Ramanujan's Notebooks, Part II}, Springer-Verlag New York, 107--112. eBook: ISBN 978-1-4612-4530-8. Softcover: ISBN 978-1-4612-8865-7. Hardcover: ISBN 978-0-387-96794-3.
\index[ppl]{Ramanujan, Srinivasa}
\index{continued square roots!Ramanujan's}

{\footnotesize \href{https://mathscinet.ams.org/mathscinet/article?mr=0970033}{MR0970033}}

{\footnotesize Cited in
\hyperref[Ber1993]{\textsc{Berndt and Bhargava 1993}}\index[ppl]{Berndt, Bruce C.}\index[ppl]{Bhargava, S.}, 
\hyperref[Lor1998]{\textsc{Lorentzen 1998}}\index[ppl]{Lorentzen, Lisa}, 
\hyperref[Ber1999]{\textsc{Berndt, Choi, and Kang 1999}}\index[ppl]{Berndt, Bruce C.}\index[ppl]{Choi, Youn-Seo}\index[ppl]{Kang, Soon-Yi}, 
\hyperref[Lev2005]{\textsc{Levin 2005}}\index[ppl]{Levin, Aaron}, 
\hyperref[Rao2005]{\textsc{Rao and Vanden Berghe 2005}}\index[ppl]{Rao, K. Srinivasa}\index[ppl]{Vanden Berghe, G.}, and
\hyperref[Vel2016c]{\textsc{Vellucci and Bersani 2016c}}\index[ppl]{Vellucci, Pierluigi}\index[ppl]{Bersani, Alberto Maria}.}

\vspace{4pt}

Entries 3, 4, 5 (given in several parts), and 6 in Chapter 12 of Ramanujan's second notebook are identities and formulas related to continued square roots\index{continued square roots!Ramanujan's}\index{continued square roots}. Entries 3 and 4 are general formal identities, generated by successive substitution\index{successive substitution!20th century examples}, of which the two formulas in \hyperref[Ram1911]{\textsc{Ramanujan 1911}} are examples.


\item\vspace{9pt} \label{Lau1990} (1990) Detlef Laugwitz\index[ppl]{Laugwitz, Detlef|textbf}, Kettenwurzeln und Kettenoperationen. \emph{Elemente der Mathematik} \textbf{44/5}, 89--98.
\index{Kettenwurzeln@\emph{Kettenwurzeln}}
\index{Kettenoperationen@\emph{Kettenoperationen}}
\index{continued reciprocal powers}
\index{continued reciprocal roots}

{\footnotesize \href{https://mathscinet.ams.org/mathscinet/article?mr=1059548}{MR1059548}}

{\footnotesize Source: \href{https://doi.org/10.5169/seals-42415}{doi.org (ETHz\"{u}rich)}.}

{\footnotesize Cited in 
\hyperref[Sch1992]{\textsc{Sch\"{o}nefuss 1992}}\index[ppl]{Schonefuss@Sch\"{o}nefuss, Lutz W.}, 
\hyperref[Lau1999]{\textsc{Laugwitz and Sch\"{o}nefuss 1999}}\index[ppl]{Laugwitz, Detlef}\index[ppl]{Schonefuss@Sch\"{o}nefuss, Lutz W.}, 
\hyperref[Joh2008]{\textsc{Johnson and Richmond 2008}}\index[ppl]{Johnson, Jamie}\index[ppl]{Richmond, Tom}, 
\hyperref[Cla2014]{\textsc{Clark and Richmond 2014}}\index[ppl]{Clark, Tyler}\index[ppl]{Richmond, Tom}, 
and \hyperref[Jon2015]{\textsc{Jones 2015}}\index[ppl]{Jones, Dixon J.}.}

\vspace{4pt}

Generalizing from the examples of continued fractions\index{continued fractions} and continued square roots\index{continued square roots}, the author defines a \emph{Kettenoperation}\index{Kettenoperationen@\emph{Kettenoperationen}}, or continued operation, as the sequence $\{P_n\}$ defined by
\begin{equation*} 
P_1=f(p_1), \; P_2=f(p_1+f(p_2)), \; P_3=f(p_1+f(p_2+f(p_3))),\ldots
\end{equation*}
where the $p_i$ are positive real numbers and $f(x)>0$ for $x>0$. He then explores convergence criteria for $\{P_n\}$ when $f(x)=x^\alpha$ for various real $\alpha$.

He rediscovers the convergence condition (\eqref{E:Wie1904a-3} and \eqref{E:Her1935-1} above) for continued square roots\index{continued square roots!convergence conditions for}, but also offers this generalization: Let $f(x)$ be positive and monotone nondecreasing for $x>0$, and $f(ax)\le a^\alpha f(x)$ for all $a\ge1$, $x \ge x_0$, with $0<\alpha<1$. Moreover, let $f(x)$ be strictly monotone increasing for $x\ge x'$. Then a necessary and sufficient condition for the convergence of $\tK_{n=0}^\infty f(p_n)$ is that $\limsup f^n(p_n) < \infty$.\index{continued compositions!of arbitrary real-valued functions}

Laugwitz gives some rough convergence conditions for the cases $f(x)=x^{-\alpha}$, where $0<\alpha<1$ and $1<\alpha$,\index{continued reciprocal powers}\index{continued reciprocal roots} but notes that his approach does not include continued fractions\index{continued fractions}, for which $\alpha=1$. Finally, he gives examples in which decreasing $f$ yield divergent sequences.\index{continued compositions!of functions $a_i+x^p$}


\item\vspace{9pt} \label{Bor1991} (1991) Jonathan M. Borwein\index[ppl]{Borwein, Jonathan M.|textbf} and G. de Barra\index[ppl]{de Barra, G.|textbf}, Nested radicals. \emph{The American Mathematical Monthly} \textbf{98}(8), 735--739.
\index{continued square roots!Ramanujan's}
\index{continued square roots!of constant nonnegative real terms}

{\footnotesize \href{https://mathscinet.ams.org/mathscinet/article?mr=1130684}{MR1130684}}

{\footnotesize Source: \href{https://doi.org/10.1080/00029890.1991.11995783}{doi.org (Taylor \& Francis)}.} 

{\footnotesize Cited in
\hyperref[Ber1994]{\textsc{Berndt 1994}}\index[ppl]{Berndt, Bruce C.}, 
\hyperref[Ber1999]{\textsc{Berndt, Choi, and Kang 1999}}\index[ppl]{Berndt, Bruce C.}\index[ppl]{Choi, Youn-Seo}\index[ppl]{Kang, Soon-Yi},
\hyperref[Fin2003]{\textsc{Finch 2003}}\index[ppl]{Finch, Steven R.},
\hyperref[Hum2007]{\textsc{Humphries 2007}}\index[ppl]{Humphries, Peter J.}, 
\hyperref[Joh2008]{\textsc{Johnson and Richmond 2008}}\index[ppl]{Johnson, Jamie}\index[ppl]{Richmond, Tom}, 
\hyperref[Zim2008a]{\textsc{Zimmerman and Ho 2008a}}\index[ppl]{Zimmerman, Seth}\index[ppl]{Ho, Chungwu}, 
\hyperref[Lim2010]{\textsc{Lim 2010}}\index[ppl]{Lim, Teik-Cheng}, 
\hyperref[Eft2012]{\textsc{Efthimiou 2012}}\index[ppl]{Efthimiou, Costas J.}, 
\hyperref[Glu2012]{\textsc{Gluzman and Yukalov 2012}}\index[ppl]{Gluzman, S.}\index[ppl]{Yukalov, V. I.},
\hyperref[Muk2013]{\textsc{Mukherjee 2013}}\index[ppl]{Mukherjee, Soumendu Sundar}, 
\hyperref[Cla2014]{\textsc{Clark and Richmond 2014}}\index[ppl]{Clark, Tyler}\index[ppl]{Richmond, Tom}, 
\hyperref[Lyn2014]{\textsc{Lynd 2014}}\index[ppl]{Lynd, Chris D.}, 
\hyperref[Vel2016c]{\textsc{Vellucci and Bersani 2016c}}\index[ppl]{Vellucci, Pierluigi}\index[ppl]{Bersani, Alberto Maria}, and
\hyperref[Wei]{\textsc{Weisstein n.d.}}\index[ppl]{Weisstein, Eric W.}}

\vspace{4pt}

The article begins by citing the identities \eqref{E:Ram1911-1} and \eqref{E:Ram1911-2} above, from \hyperref[Ram1911]{\textsc{Ramanujan 1911}}\index[ppl]{Ramanujan, Srinivasa}\index{continued square roots!Ramanujan's}, and the continued square root of constants \eqref{E:Ber1692-1} (but traced back only to the solution of \hyperref[Gin1916]{\textsc{Ginsburg 1916}}\index[ppl]{Ginsburg, J. J.})\index{continued square roots!of constant nonnegative real terms}. The authors continue: ``In this note we consider a general class of nested radicals with arithmetic progressions which includes [\eqref{E:Ber1692-1}, \eqref{E:Ram1911-1}, and \eqref{E:Ram1911-2} above] as special cases. We deal with $p$th roots, though more general monotone functions could be considered$\ldots$ We look for a finite nonnegative solution to the equation\index{continued robb roots@continued $r$th roots!of terms in arithmetic progression}
\begin{equation}\label{E:Bor1991-1}
f(k) = \sqrt[p]{g(k) + k^\alpha f(k + a)}\;,
\end{equation}
where $k, \alpha, a \ge 0$, $p > 1$, and $g: [0, \infty]\to[0, \infty]$ with $y = \inf g > 0$. We solve \eqref{E:Bor1991-1} iteratively. Take $f_0, 0 < f_0 < g^{1/p}$ and define $f_{n+1}$ by
\[f_{n+1}^p(k) = g(k) + k^\alpha f_n(k + a)\;.\]
Then $f_0 \le g^{1/p} \le f_1$. Hence, inductively, if $f_{n+1}\ge f_n$, 
\[f_{n+2}^p(k) - f_{n+1}^p(k) = k^\alpha[f_{n+1}(k + a) - f_n(k+a)]\ge 0\;,\]
and $f_n$ increases to $f^*$ which takes values in $(0,\infty]$ and which solves \eqref{E:Bor1991-1}." A ``growth condition" is then used to prove that $f^*$ has the desired form. The authors give several examples, then extend their method to a functional equation in two variables, and finally consider solutions to \eqref{E:Bor1991-1} which fail the growth condition.


\item\vspace{9pt} \label{Gil1991} (1991) John Gill\index[ppl]{Gill, John|textbf}, Inner composition of analytic mappings on the unit disk. \emph{International Journal of Mathematics and Mathematical Sciences} \textbf{14}(2), 221--226.
\index{continued compositions!of complex-valued functions}

{\footnotesize \href{https://mathscinet.ams.org/mathscinet/article?mr=1096859}{MR1096859}}

{\footnotesize Source: \href{https://doi.org/10.1155/S0161171291000236}{doi.org (Wiley)}.}

{\footnotesize Cited in
\hyperref[Lor1995]{\textsc{Lorentzen 1995}}\index[ppl]{Lorentzen, Lisa}.}

\vspace{4pt}

What we are here calling \emph{continued compositions}, the author of this paper calls \emph{inner compositions}\index{inner compositions}. He writes, ``In the present paper the following question is posed, and, to some extent, answered: Suppose each member of the sequence $\{f_n\}$ is analytic on Int$(D)$ and continuous on $D$ with Int$D\supseteq f(D)$ (it is not assumed that $f_n\to f$). Under what conditions does $F(z)=f_1\circ \cdots \circ f_n(z)\to\lambda$, a constant, for all $z\in D$, as $n\to\infty$? Thus we are considering `inner' compositions of essentially random sequences of functions mapping the unit disk into itself. Although our approach focuses on mappings of $D$ into $D$, more general results are possible$\ldots$ We shall present several theorems describing conditions on the $f_n$s that imply $F(D)\to\lambda$."


\item\vspace{9pt} \label{Jon1991} (1991) Dixon J. Jones\index[ppl]{Jones, Dixon J.|textbf}, Continued powers and roots. \emph{The Fibonacci Quarterly} \textbf{29}(1), 37--46.
\index{continued pobb powers@continued $p$th powers}
\index{continued powers}
\index{continued robb roots@continued $r$th roots}

{\footnotesize \href{https://mathscinet.ams.org/mathscinet/article?mr=1089518}{MR1089518}}

{\footnotesize Source: \href{https://www.fq.math.ca/Scanned/29-1/jones.pdf}{Fibonacci Quarterly}.}

{\footnotesize Cited in 
\hyperref[Jon1995]{\textsc{Jones 1995}}\index[ppl]{Jones, Dixon J.}, 
\hyperref[Lau1999]{\textsc{Laugwitz and Sch\"{o}nefuss 1999}}\index[ppl]{Laugwitz, Detlef}\index[ppl]{Schonefuss@Sch\"{o}nefuss, Lutz W.}, 
\hyperref[Muk2013]{\textsc{Mukherjee 2013}}\index[ppl]{Mukherjee, Soumendu Sundar}, and 
\hyperref[Jon2015]{\textsc{Jones 2015}}\index[ppl]{Jones, Dixon J.}.}

\vspace{4pt}

Gives some convergence conditions for the continued $p$th power\index{continued powers}\index{continued robb roots@continued $r$th roots} \eqref{E:contpower} above, when $p$ and the terms $a_i$ are nonnegative. It is shown that the continued $p$th power with nonnegative constant terms $a_i=a\ge 0$ converges if and only if\index{continued powers!of constant nonnegative real terms}
\begin{align*}
a\ge0\quad &\textrm{for } 0 < p < 1\;;\\
a=0\quad &\textrm{for }p=1\;; \textrm{ and}\\
0\le a \le R\quad &\textrm{for }p >1\;,
\end{align*}
where 
\begin{equation*}
R=\sqrt[p-1]{\dfrac{(p-1)^{p-1}}{p^p}}\;.
\end{equation*}
(These were given independently, and nearly simultaneously, in \hyperref[Sch1992]{\textsc{Sch\"{o}nefuss 1992}}\index[ppl]{Schonefuss@Sch\"{o}nefuss, Lutz W.}; similar calculations were made in \hyperref[Hey1894b]{\textsc{Heymann 1894b}}\index[ppl]{Heymann, W.}.) For $0<p<1$ and arbitrary nonnegative terms, a minor extension of Herschfeld's\index[ppl]{Herschfeld, Aaron} convergence condition \eqref{E:Her1935-1} (actually a special case of Herschfeld's Theorem III) for continued $r$th roots\index{continued robb roots@continued $r$th roots!convergence conditions for} is cited and proved. For $p>1$, the continued $p$th power\index{continued powers!convergence conditions for} is shown to converge if $\limsup_{i\to\infty}a_i<R$, and to diverge if $p\ge 1$ and $\liminf_{i\to\infty} a_i>R$. A necessary and sufficient condition for convergence of 
\[b+^2\!(a+^2\!(b+^2\!(a+\cdots)))\]
is given. Finally it is proved that, for $p>1$, a continued $p$th power of positive terms $a_i$ converges if
\[\dfrac{a_{i+1}^p}{a_i}\le \dfrac{(p-1)^{p-1}}{p^p}\]
for all sufficiently large values of $i$.\index{continued powers!convergence conditions for}


\item\vspace{9pt} \label{Pic1991} (1991) C. A. Pickover\index[ppl]{Pickover, C. A.|textbf} and A. Lakhtakia\index[ppl]{Lakhtakia, A.|textbf}, Continued roots in the complex plane. \emph{Journal of Recreational Mathematics} \textbf{23}, 198--202.
\index{continued square roots!of constant complex terms}

{\footnotesize Source: \href{https://www.researchgate.net/publication/264919078_Continued_roots_in_the_complex_plane}{ResearchGate}.}

The authors write, ``[T]he striking beauty and complexity of the patterns resulting from $\ldots$ iterative calculations has only recently been explored in detail, due in part to advances in computer graphics$\ldots$ The simple iterated functions discussed in the literature naturally lead to curiosity about the behavior of more complicated functions. We have found particularly visually interesting stability plots for 
\[z\to\sin(z+f^n(z))\]
where $z$ is a complex variable, and $f^n(z)$ is a continued radical, and $n$ is the number of terms used in the nesting$\ldots$ In particular, we focus on the function defined recursively by $f^{n+1}(z)=\sqrt{z+f^n(z)}\;,\;n>1$, with the initiator $f^1(z)=\sqrt{z}$." Julia sets\index{Julia sets} of these iterated functions are depicted. 


\item\vspace{9pt} \label{Sch1992} (1992) Lutz W. Sch\"{o}nefuss\index[ppl]{Schonefuss@Sch\"{o}nefuss, Lutz W.|textbf}, \emph{Nichtautonome Differenzengleichungen und Kettenoperationen}, Mitteilungen aus dem mathematischen Seminar der Universit\"{a}t Giessen(207), 110 pp. ISSN 0373-8221.
\index{Kettenoperationen@\emph{Kettenoperationen}}
\index{continued compositions!of arbitrary real-valued functions}

{\footnotesize \href{https://mathscinet.ams.org/mathscinet/article?mr=1156689}{MR1156689}}

{\footnotesize Cited in 
\hyperref[Lau1999]{\textsc{Laugwitz and Sch\"{o}nefuss 1999}}\index[ppl]{Laugwitz, Detlef}\index[ppl]{Schonefuss@Sch\"{o}nefuss, Lutz W.}, and
\hyperref[Jon2015]{\textsc{Jones 2015}}\index[ppl]{Jones, Dixon J.}.}

\vspace{4pt}

From the introduction: ``An infinite nested expression of the form 
\[f(a_0 +f(a_1+f(a_2+\ldots))),\]
which we call a continued operation\index{continued compositions!of arbitrary real-valued functions}, denoted by the symbol $\tK_{n=0}^\infty f(a_n)$, and its sequence of partial operations $f(a_0)$, $f(a_0 + f(a_1))$, $f(a_0+f(a_1+f(a_2)))$, etc., can be viewed as a generalization of infinite sums and continued fractions$\ldots\,$\index{continued fractions!generalizations of} [In] the first chapter of this work$\,\ldots\,$we study the convergence of continued operations with monotonic functions $f$, in particular the family of functions $f(x)=x^{\alpha}, \alpha \in \mathbb{R}$.\index{continued compositions!of functions $a_i+x^p$} Besides numerous convergence results we obtain a complete survey of the convergence of $\tK_{n=0}^\infty a^\alpha, \alpha \in \mathbb{R}$, for constant terms $a > 0\ldots\,$ Here we also make clear the special r\^{o}le of infinite sums\index{infinite sums}.

``With a trick, the convergence question can be answered. Denoting the `$k$th tail' of a convergent continued operation by $x_k$ for $k\in\mathbb{N}$:
\[ x_k=\dK_{n=k}^\infty f(a_n)\;, \]
the numbers $x_k$ satisfy the nonautonomous difference equation\index{difference equations!nonautonomous} 
\[ x_k=f(a_k+x_{k+1}),\quad k\in\mathbb{N}\;. \]
Depending upon conditions on $f$, the convergence of $\tK_{n=0}^\infty f(a_n)$ is equivalent to the existence of a solution --- or of a unique solution --- to this equation$\ldots$ [I]n certain cases a special asymptotic property distinguishes the one solution $(x_n)$ --- for which $x_0=\dK_{n=0}^\infty f(a_n)$ --- from the other possible solutions of that equation: its minimality$\ldots$

``Thus one is lead in a completely natural way to the study of the general nonautonomous difference equation of first order
\[ x_{n+1}=F(n,x_n),\quad n\in\mathbb{N}\;.\]
This happens in the third chapter, the heart of this work$\ldots$

``In order to get on the track of the asymptotic behavior of solutions of the equation, we examine first the limiting set of such a solution. Now the special case of an autonomous equation\index{difference equations!autonomous}
\[x_{n+1}=F(x_n),\quad n\in\mathbb{N} \]
admits of three invariance characteristics of limiting sets, which permit a more exact localization$\ldots$ We can transfer these in an appropriate way to the general equation$\ldots$ Here we find our way using the corresponding approach from the theory of ordinary differential equations\index{differential equations}$\ldots$ 

``The transfer makes it necessary to define limit functions for the function $F:\mathbb{N}\times D\to\mathbb{R},\;D\subseteq\mathbb{R}$. The totality of these we call the limiting set $\Lambda(F)$. It is true that the formulation of the invariance results must be restricted to limit equations formed by limit functions, but nevertheless a beautiful structure develops, analogous to the limiting set $L(x_0)$ of the solution sequence $(x_n)$ --- the limiting set $\Lambda(F)$ of the function $F$. This becomes more transparent if we use various facts from the theory of dynamical systems and results of the Ukrainian mathematician Sarkovskii$\ldots$\index[ppl]{Sharkovskii, O. M.}

``In the fourth chapter we translate results about continued operations into the language of chapter 3$\ldots$ We explain how one can interpret the formation of a continued operation as a generalization of a summation procedure for linear difference equations$\ldots$

``[T]he theory of the third chapter presupposes to a large extent so-called compactness characteristics of the function $F$ and the solution $(x_n)$. This makes more difficult the parallel development of this theory with the difference equations of the continued operations, because the interesting cases of the continued operations\index{continued compositions!of arbitrary real-valued functions} are usually those which lead to unbounded solutions and non-compact functions. But the remarks of the fourth chapter point also to a possible way out, namely stability analyses of the solutions.

``The appendix consists of two parts: A supremely strange continued operation --- the infinite exponential --- is presented\index{continued exponentials}; then two continued radicals of Ramanujan are examined more precisely."\index[ppl]{Ramanujan, Srinivasa}\index{continued square roots!Ramanujan's}


\item\vspace{9pt} \label{Abi1993} (1993) Alexander Abian\index[ppl]{Abian, Alexander|textbf} and Sergei Sverchkov\index[ppl]{Sverchkov, Sergei|textbf}, Solutions and representations by iterated radicals. \emph{International Journal of Mathematical Education in Science and Technology} \textbf{24}(3), 449--455.
\index{iterated radicals}
\index{continued robb roots@continued $r$th roots!as solutions to polynomial equations}
\index{continued square roots!representation of functions as}

{\footnotesize \href{https://mathscinet.ams.org/mathscinet/article?mr=1229874}{MR1229874 (citation only)}}

{\footnotesize Source: \href{https://doi.org/10.1080/0020739930240316}{doi.org (Taylor \& Francis)}.}

{\footnotesize Cited in
\hyperref[Rao2005]{\textsc{Rao and Vanden Berghe 2005}}\index[ppl]{Rao, K. Srinivasa}\index[ppl]{Vanden Berghe, G.}.}

\vspace{4pt}

From the abstract: ``We give solutions of polynomial equations as well as some remarkable representations of functions as infinite iterations of radicals." \index{continued robb roots@continued $r$th roots!as solutions to polynomial equations}\index{continued square roots!representation of functions as} Many of the results are independent rediscoveries; for instance, the authors give continued $r$th root representations of roots of trinomial equations\index{trinomial equations}, and exhibit (without attribution) the Ramanujan identity \eqref{E:Ram1911-1} above.\index[ppl]{Ramanujan, Srinivasa} The main theorems involve an infinitely differentiable function $f$ having nonnegative $m$th derivatives for every $m=0,1,2,\ldots$, and a function $g_n$ defined for every $n=0,1,2,\ldots$ by
\begin{equation*}
g_0(x)=f(x)\quad\textrm{and}\quad g_n(x)=
\begin{cases}
\dfrac{g_{n-1}^2(x)-g_{n-1}^2(0)}{x} &\textrm{for $x\ne 0$}\\
(g_{n-1}^2(0))' &\textrm{for $x=0$.}
\end{cases}
\end{equation*}
With $a_n=g_{n-1}^2(0)$, it is proved that
\begin{equation}\label{E:Abian1}
f(x)=\sqrt{a_1+x\sqrt{a_2+x\sqrt{a_3+x\sqrt{\cdots+x\sqrt{a_n+xg_n(x)}}}}}\;.
\end{equation}
If it is further assumed that the $m$th derivatives of $f$ are nonnegative on $(-\tfrac{1}{2}, \tfrac{1}{2})$ and $1\le g_n'(x)\le g_n(x)$ for $n=0,1,2,\ldots$, then
\begin{align}\label{E:Abian2}
f(x)&=\lim_{n\to\infty}\sqrt{a_1+x\sqrt{a_2+x\sqrt{a_3+x\sqrt{\cdots+x\sqrt{a_n+x}}}}}\notag\\
&=\sqrt{a_1+x\sqrt{a_2+x\sqrt{a_3+x\sqrt{\cdots}}}}\;.
\end{align}


\item\vspace{9pt} \label{Ber1993} (1993) Bruce C. Berndt\index[ppl]{Berndt, Bruce C.|textbf} and S. Bhargava\index[ppl]{Bhargava, S.|textbf}, Ramanujan --- for lowbrows. \emph{The American Mathematical Monthly} \textbf{100}(7), 644--656.
\index[ppl]{Ramanujan, Srinivasa}
\index{continued square roots!Ramanujan's}

{\footnotesize \href{https://mathscinet.ams.org/mathscinet/article?mr=1237220}{MR1237220}}

{\footnotesize Source: \href{https://doi.org/10.1080/00029890.1993.11990465}{doi.org (Taylor \& Francis)}.}

{\footnotesize Cited in
\hyperref[Ber1994]{\textsc{Berndt 1994}}\index[ppl]{Berndt, Bruce C.} and
\hyperref[Ber1999]{\textsc{Berndt, Choi, and Kang 1999}}\index[ppl]{Berndt, Bruce C.}\index[ppl]{Choi, Youn-Seo}\index[ppl]{Kang, Soon-Yi}.} 

\vspace{4pt}

In Section 3, the authors discuss the problems proposed in \hyperref[Ram1914]{\textsc{Ramanujan 1914}} and \hyperref[Ram1915]{\textsc{Ramanujan 1915}}, including entries in Ramanujan's notebooks pertaining to these and similar problems.


\item\vspace{9pt} \label{Ber1994} (1994) Bruce C. Berndt\index[ppl]{Berndt, Bruce C.|textbf}, \emph{Ramanujan's Notebooks, Part IV}, Springer-Verlag New York, 10--20. eBook: ISBN 978-1-4612-0879-2. Hardcover: ISBN 978-0-387-94109-7. Softcover: ISBN 978-1-4612-6932-8.
\index[ppl]{Ramanujan, Srinivasa}

{\footnotesize \href{https://mathscinet.ams.org/mathscinet/article?mr=1261634}{MR1261634}}

{\footnotesize Cited in 
\hyperref[Ber1993]{\textsc{Berndt and Bhargava 1993}}\index[ppl]{Berndt, Bruce C.}\index[ppl]{Bhargava, S.}, 
\hyperref[Ber1999]{\textsc{Berndt, Choi, and Kang 1999}}\index[ppl]{Berndt, Bruce C.}\index[ppl]{Choi, Youn-Seo}\index[ppl]{Kang, Soon-Yi}, 
\hyperref[Xi2013]{\textsc{Xi and Qi 2013}}\index[ppl]{Xi, Bo-Yan}\index[ppl]{Qi, Feng}, and
\hyperref[Wei]{\textsc{Weisstein n.d.}}\index[ppl]{Weisstein, Eric W.}}

\vspace{4pt}

The publisher's blurb about this book states: ``This is the first of two volumes devoted to proving the results found in the unorganized portions of [Ramanujan's] second notebook and in the third notebook." Entries 4 and 5 involve continued square roots.\index{continued square roots!Ramanujan's} Entry 4 is a portion of the problem proposal \hyperref[Ram1914]{\textsc{Ramanujan 1914}}, which is solved by finding six nontrivial zeroes of a polynomial of degree 8; pages 11--17 of this book are devoted to rigorously calculating these zeroes. Entry 5 exhibits continued square roots that converge to these zeroes. The following two pages make use of two computer algebra packages to obtain closed form expressions.


\item\vspace{9pt} \label{Ber1995} (1995) Bruce C. Berndt\index[ppl]{Berndt, Bruce C.|textbf} and Robert A. Rankin\index[ppl]{Rankin, Robert A.|textbf}, \emph{Ramanujan: Letters and Commentary}, American Mathematical Society, Providence; London Mathematical Society, London. ISBN 978-0821802878.
\index[ppl]{Ramanujan, Srinivasa}

{\footnotesize \href{https://mathscinet.ams.org/mathscinet/article?mr=1353909}{MR1353909}}

{\footnotesize Cited in
\hyperref[Ber1994]{\textsc{Berndt 1994}}\index[ppl]{Berndt, Bruce C.} and
\hyperref[Ber1999]{\textsc{Berndt, Choi, and Kang 1999}}\index[ppl]{Berndt, Bruce C.}\index[ppl]{Choi, Youn-Seo}\index[ppl]{Kang, Soon-Yi}.}

\vspace{4pt}

This book reproduces the text of a letter, dated 4 January 1928, from T. Vijayaraghavan\index[ppl]{Vijayaraghavan, T.|textbf} to B. M. Wilson\index[ppl]{Wilson, B. M.} (one of the editors of \hyperref[Ram1927]{\textsc{Ramanujan 1927}})\index[ppl]{Ramanujan, Srinivasa}, ``[r]egarding the justification of the formal processes in Ramanujan's solution" of \hyperref[Ram1911]{\textsc{Ramanujan 1911}}.\index{continued square roots!Ramanujan's} Vijayaraghavan continues, ``I now do not remember whether I justified Ramanujan's solution or simply established the criterion of convergency" \eqref{E:Ram1927-1} above. He then notes that ``to prove that [equation \eqref{E:Ram1911-1} above] has the value 3 some additional remarks are necessary$\,\ldots\,$the only proof I can think of is unconscionably long and tedious$\ldots$" The remainder of the letter sketches this proof.\index{continued square roots!convergence conditions for} In a commentary following the letter, the editors clarify and correct a few of Vijayaraghavan's steps, and give a thumbnail biography of the man. 


\item\vspace{9pt} \label{Jon1995} (1995) Dixon J. Jones\index[ppl]{Jones, Dixon J.|textbf}, Continued powers and a sufficient condition for their convergence. \emph{Mathematics Magazine} \textbf{68}(5), 387--392.
\index{continued powers!of arbitrary nonnegative real terms}

{\footnotesize \href{https://mathscinet.ams.org/mathscinet/article?mr=1365650}{MR1365650}}

{\footnotesize Source: \href{https://doi.org/10.1080/0025570x.1995.11996363}{doi.org (Taylor \& Francis)}.}

{\footnotesize Cited in 
\hyperref[Zim2008a]{\textsc{Zimmerman and Ho 2008a}}\index[ppl]{Zimmerman, Seth}\index[ppl]{Ho, Chungwu}, 
\hyperref[Glu2012]{\textsc{Gluzman and Yukalov 2012}}\index[ppl]{Gluzman, S.}\index[ppl]{Yukalov, V. I.}, 
\hyperref[Muk2013]{\textsc{Mukherjee 2013}}\index[ppl]{Mukherjee, Soumendu Sundar}, 
\hyperref[Xi2013]{\textsc{Xi and Qi 2013}}\index[ppl]{Xi, Bo-Yan}\index[ppl]{Qi, Feng}, 
\hyperref[Kef2014]{\textsc{Kefalas 2014}}\index[ppl]{Kefalas, Kyriakos}, 
\hyperref[Lyn2014]{\textsc{Lynd 2014}}\index[ppl]{Lynd, Chris D.}, and 
\hyperref[Jon2015]{\textsc{Jones 2015}}\index[ppl]{Jones, Dixon J.}.}

\vspace{4pt}

It is proved that, for real $p>1$, the continued $p$th power\index{continued powers!of arbitrary nonnegative real terms} $\tK_{i=0}^\infty a_i^p$ converges if $(a_n/R)^{p^n}$ is bounded, where
\[R=\dfrac{p-1}{p^{\frac{p-1}{p}}}=\sqrt[p-1]{\dfrac{(p-1)^{p-1}}{p^p}}\;.\]
(See \hyperref[Jon1991]{\textsc{Jones 1991}} for an earlier appearance of the second form of $R$ above.) An error in Example III of this paper, noted by J. Nichols-Barrer\index[ppl]{Nichols-Barrer, J.} in a Letter to the Editor, \textbf{69}(3), 238, is corrected in a \href{https://doi.org/10.1080/0025570X.1996.11996468}{Letter to the Editor, \textbf{69}(4), 316}


\item\vspace{9pt} \label{Lor1995} (1995) Lisa Lorentzen\index[ppl]{Lorentzen, Lisa|textbf}, A convergence question inspired by Stieltjes\index[ppl]{Stieltjes, T. J} and by value sets in continued fraction\index{continued fractions} theory. In \emph{Proceedings of the International Conference on Orthogonality, Moment Problems and Continued Fractions}, Delft, 1994. \emph{Journal of Computational and Applied Mathematics} \textbf{65} nos. 1--3, 233--251.
\index{continued compositions!of complex-valued functions}

{\footnotesize \href{https://mathscinet.ams.org/mathscinet/article?mr=1379134}{MR1379134}}

{\footnotesize Source: \href{https://doi.org/10.1016/0377-0427(95)00113-1}{doi.org (ScienceDirect)}.}

{\footnotesize Cited in
\hyperref[Lor1998]{\textsc{Lorentzen 1998}}\index[ppl]{Lorentzen, Lisa}.}

\vspace{4pt}

From the abstract: ``Let $V$ be a subset of the complex plane $\mathbb{C}$. Let $\{f_n\}_{n=1}^\infty$ be a sequence of self-mappings of $V$; i.e., $f_n(V)\in V$. The question is then: Under what conditions will the sequence
\[F_n(w):=f_1\circ f_2 \circ\cdots\circ f_n(w), \quad n=1,2,3\ldots\]
of composite maps converge to a constant function in $V$? In this paper we give a survey of some of the answers and open problems connected with this question." \index{continued compositions!of complex-valued functions}




\item\vspace{9pt} \label{Siz1996} (1996) Walter S. Sizer\index[ppl]{Sizer, Walter S.|textbf} and E. K. Wiredu\index[ppl]{Wiredu, E. K.|textbf}, Convergence of some continued radicals in the complex plane. \emph{Bulletin of the Malaysian Mathematical Society} \textbf{19}(1), 1--7.
\index{continued square roots!of constant complex terms}
\index{iterated square roots!of constant complex terms}

{\footnotesize \href{https://mathscinet.ams.org/mathscinet/article?mr=1465806}{MR1465806}} 

{\footnotesize Cited in 
\hyperref[Joh2008]{\textsc{Johnson and Richmond 2008}}\index[ppl]{Johnson, Jamie}\index[ppl]{Richmond, Tom} and 
\hyperref[Muk2013]{\textsc{Mukherjee 2013}}\index[ppl]{Mukherjee, Soumendu Sundar}.}

\vspace{4pt}

Independently rediscovers the result from \hyperref[Sch1962]{\textsc{Schuske and Thron 1962}}\index[ppl]{Thron, W. J.}\index[ppl]{Schuske, Georgellen} that the continued square root \eqref{E:Ber1692-2} above converges for all complex constants $a$, where the principal value of the square root function is used. The authors note that, using a different branch of the square root function, there are values of $a$ for which \eqref{E:Ber1692-2} does not converge. The MathSciNet review MR1465806 mistakenly displays the general continued square root \eqref{E:croot} above, instead of \eqref{E:Ber1692-2}.


\item\vspace{9pt} \label{Tou1996} (1996) Dominique Tourn\`{e}s\index[ppl]{Tourn\`{e}s, Dominique|textbf}, \emph{L'int\'{e}gration approch\'{e}e des \'{e}quations diff\'{e}rentielles ordinaires (1671-1914), th\`{e}se de doctorat de l'universit\'{e} Paris 7--Denis Diderot}, June 1996. Reprinted by Presses Universitaires du Septentrion, Villeneuve d'Ascq, 1997. 
\index{successive substitution}
\index{differential equations}
\index{history!of successive substitution}

{\footnotesize Source: \href{https://hal.science/tel-03948671v1/document}{HAL Open Science}.}

Chapter 4 of this doctoral thesis traces the history of ``the method of successive approximations,"\index{successive substitution} from ancient Greek times through the early 20th century. Subsection 2.1.1. gives a valuable appraisal of \hyperref[Sch1821]{\textsc{Schmidten 1821}} in its historical context. Tourn\`{e}s writes, ``In a dissertation published in 1821, Henri Gerner Schmidten\index[ppl]{Schmidten, Henri Gerner} gives us a panorama of the method of successive approximations and its possible applications in the various branches of analysis$\ldots$ [The] Schmidten text$\,\ldots\,$shows that at the beginning of the 19th century the empirical technique of successive substitutions had already assumed the aspect of a very general fixed point problem perfectly theorized in its formal aspect." 


\item\vspace{9pt} \label{Lor1998} (1998) Lisa Lorentzen\index[ppl]{Lorentzen, Lisa|textbf}, Continued compositions of self-mappings. \emph{Proceedings of the Second Asian Mathematical Conference 1995 (Nakhon Ratchasima)}, 34--38, World Sci. Publ., River Edge, NJ.
\index{continued compositions!of complex-valued functions}
\index{continued compositions!of linear fractional transformations}

{\footnotesize \href{https://mathscinet.ams.org/mathscinet/article?mr=1660544}{MR1660544 (citation only)}}

{\footnotesize Source: \href{https://books.google.com/books?id=CSDQDgAAQBAJ&lpg=PA34&ots=wY9ZoMiqZL&dq=\%22Continued\%20compositions\%20of\%20self-mappings\%22\&pg=PA34\#v=onepage\&q\&f=false}{Google Books}.}

From the abstract: ``Let $\{f_n\}_{n=1}^\infty$ be a sequence of self-mappings of a set $V$; i.e., $f_n(V)\subset V$. Under what conditions will the sequence $\{F_n(w)\}_{n=1}^\infty$ given by $F_n(w):= f_1\circ f_2\circ\cdots \circ f_n(w)$ converge to a constant function in $V$? Answers to this question have applications in dynamical systems, Schur analysis\index{Schur analysis}, continued fractions\index{continued fractions} and other similar structures like towers of exponentials\index{continued exponentials} and infinite radicals\index{infinite radicals}." The author mentions the Denjoy-Wolff Theorem\index{Denjoy-Wolff Theorem}; notes the $f_i$ which produce continued square roots\index{continued square roots!of arbitrary complex terms} and continued exponentials\index{continued exponentials}; states other results that pertain when the $f_i$ are linear fractional transformations; and closes with a convergence theorem from a paper on continued exponentials\index{continued exponentials} (I. N. Baker\index[ppl]{Baker, I. N.} and P. J. Rippon\index[ppl]{Rippon, P. J.}, Towers of exponents and other composite maps, \emph{Complex Variables, Theory and Application: An International Journal} \textbf{12}(1--4), 1989, 181--200).


\item\vspace{9pt} \label{Bez1999} (1999) Mih\'{a}ly Bencze\index[ppl]{Bencze, Mih\'{a}ly|textbf}, About Ramanujan's nested radicals. \emph{Octogon Mathematics Magazine} \textbf{7}(2), 15--21.
\index[ppl]{Ramanujan, Srinivasa}
\index{continued square roots!Ramanujan's}
\index{nested radicals}

{\footnotesize \href{https://mathscinet.ams.org/mathscinet/article?mr=1730807}{MR1730807 (citation only)}}

The note begins by quoting equation \eqref{E:Pol1925-1} above, citing \hyperref[Pol1925]{\textsc{P\'{o}lya and Szeg\"{o} 1925}}\index[ppl]{Polya@P\'{o}lya, G.}; it then give proof sketches of several common continued square root expressions\index{continued square roots}, including \eqref{E:Vie1593-1}, \eqref{E:Ram1911-1}, \eqref{E:Ram1911-2} above, and a version of \eqref{E:McG-1} below. The only other citation is the name ``Kreyszig," mentioned in connection with iteration of a real-valued function $g(x)$ and the condition that such iteration converges if $|g'(x)|<1$; this may refer to a result in a text by Erwin Kreyszig\index[ppl]{Kreyszig, Erwin}, but the result was known long before Kreyszig's birth in 1922.


\item\vspace{9pt} \label{Ben1999} (1999) Carl M. Bender\index[ppl]{Bender, Carl M.|textbf} and Steven A. Orszag\index[ppl]{Orszag, Steven A.|textbf}, \emph{Advanced Mathematical Methods for Scientists and Engineers I}, Springer-Verlag, New York, ISBN 978-0387989310. Originally published by McGraw-Hill, New York, 1978.

{\footnotesize \href{https://mathscinet.ams.org/mathscinet/article?mr=1721985}{MR1721985}}

{\footnotesize Cited in
\hyperref[Luo2003]{\textsc{Luo, Qi, Barnett, and Dragomir 2003}}\index[ppl]{Luo, Qiu-Ming}\index[ppl]{Qi, Feng}\index[ppl]{Barnett, Neil S.}\index[ppl]{Dragomir, Sever S.}.}

\vspace{4pt}

In Section 8.4, on continued fractions\index{continued fractions} and Pad\'{e} approximants\index{Pad\'{e} approximants}\index{approximants!Pad\'{e}}, the authors include a brief subsection on continued function representations\index{continued compositions!of complex-valued functions}. Examples are given of continued exponentials\index{continued exponentials}, continued square roots\index{continued square roots!of arbitrary complex terms}, and continued logarithms\index{continued logarithms}. Problems 8.32 to 8.36 develop convergence properties of these examples, with particular emphasis on continued exponentials.


\item\vspace{9pt} \label{Ber1999} (1999) Bruce C. Berndt\index[ppl]{Berndt, Bruce C.|textbf}, Youn-Seo Choi\index[ppl]{Choi, Youn-Seo|textbf}, and Soon-Yi Kang\index[ppl]{Kang, Soon-Yi|textbf}, The problems submitted by Ramanujan to the Journal of the Indian Mathematical Society. In \emph{Continued fractions: from analytic number theory to constructive approximation. A volume in honor of L. J. Lange\index[ppl]{Lange, L. J.}. Papers from the International Conference held at the University of Missouri, Columbia, MO, May 20--23, 1998.} Edited by Bruce C. Berndt and Fritz Gesztesy\index[ppl]{Gesztesy, Fritz}. Contemporary Mathematics, 236. American Mathematical Society, Providence, RI, ISBN 978-0821812006.
\index[ppl]{Ramanujan, Srinivasa}

{\footnotesize \href{https://mathscinet.ams.org/mathscinet/article?mr=1665361}{MR1665361}} 

\vspace{4pt}

Gives a detailed overview of the 58 problems Ramanujan submitted to the \emph{Journal of the Indian Mathematical Society} between 1911 and 1919, including \hyperref[Ram1911]{\textsc{Ramanujan 1911}}, \hyperref[Ram1914]{\textsc{Ramanujan 1914}}, \hyperref[Ram1915]{\textsc{Ramanujan 1915}} listed above\index{continued square roots!Ramanujan's}.


\item\vspace{9pt} \label{Lau1999} (1999) Detlef Laugwitz\index[ppl]{Laugwitz, Detlef|textbf} and Lutz Sch\"{o}nefuss\index[ppl]{Schonefuss@Sch\"{o}nefuss, Lutz W.|textbf}, Convergence of continued operations. In \emph{Iteration Theory (ECIT '96), Proceedings of the European Conference on Iteration Theory, Urbino, Italy, September 8--14, 1996}, edited by L. Gardini\index[ppl]{Gardini, L.|textbf}, et al., \emph{Grazer Mathematische Berichte} \textbf{339}, 243--250.
\index{continued operations}

{\footnotesize \href{https://mathscinet.ams.org/mathscinet/article?mr=1748827}{MR1748827}}

{\footnotesize Cited in 
\hyperref[Jon2015]{\textsc{Jones 2015}}\index[ppl]{Jones, Dixon J.}.}

\vspace{4pt}

From the summary: ``For a monotonic bijection of $\prls$ and a sequence $a_n > 0$ we consider both the continued operation\index{continued operations}
\[f(a_0 +f(a_1 +f(a_2 + \cdots)))\]
and the difference equation\index{difference equations!nonautonomous} $x_n = f(a_n +x_{n+1})$. In Section 2, the convergence of the continued operation is discussed for $a_n = c$, and is completely settled for $f(x) = x^r$, $r\in\rls$, by means of fixed point iterations\index{continued compositions!of functions $a_i+x^p$}\index{fixed-point iteration}\index{iteration!fixed-point}. The general case is considered in Section 3 (increasing $f$) and Section 4 (decreasing $f$), and results on the convergence of the continued operation are obtained in terms of positive solutions of the difference equation. Some conclusions can be drawn on the existence and uniqueness of positive solutions of this nonautonomous equation."


\item\vspace{9pt} \label{Osl1999} (1999) Thomas J. Osler, The union of Vieta's and Wallis's products for pi. \emph{The American Mathematical Monthly} \textbf{106}(8), 774--776.
\index{Vi\`{e}te's formula for $\tfrac{2}{\pi}$!generalizations}
\index{Wallis's product for $\tfrac{2}{\pi}$}

{\footnotesize \href{https://mathscinet.ams.org/mathscinet/article?mr=1718586}{MR1718586 (citation only)}} 

{\footnotesize Source: \href{https://doi.org/10.1080/00029890.1999.12005119}{doi.org (Taylor \& Francis)}.}

{\footnotesize Cited in
\hyperref[Lev2005]{\textsc{Levin 2005}}\index[ppl]{Levin, Aaron}.}

\vspace{4pt}

Wallis's\index[ppl]{Wallis, John} product is
\begin{equation}\label{E:Osl1999-1}
\dfrac{2}{\pi}=\dfrac{1\cdot 3}{2\cdot 2}\cdot\dfrac{3\cdot 5}{4\cdot 4}\cdot\dfrac{5\cdot 7}{6\cdot 6}\cdots\;.
\end{equation}
The author shows that this expansion and Vi\`{e}te's\index[ppl]{Vi\`{e}te, Fran\c{c}ois} formula for $\tfrac{2}{\pi}$, equation \eqref{E:Vie1593-1} above, are special cases of the following product:
\begin{align}\label{E:Osl1999-2}
\dfrac{2}{\pi}&=\prod_{n=1}^p\overbrace{\sqrt{\dfrac{1}{2}+\dfrac{1}{2}\sqrt{\dfrac{1}{2}+\dfrac{1}{2}\sqrt{\dfrac{1}{2}+\cdots+\dfrac{1}{2}\sqrt{\dfrac{1}{2}}}}}}^{n\textrm{ square roots}}\notag\\
&\times \prod_{n=1}^\infty\dfrac{2^{p+1}n-1}{2^{p+1}n}\cdot\dfrac{2^{p+1}n+1}{2^{p+1}n}\;.
\end{align}


\item\vspace{9pt} \label{Cor2000a} (2000) Juan Carlos Cort\'{e}s L\'{o}pez\index[ppl]{Cort\'{e}s L\'{o}pez, Juan Carlos|textbf} and Juan \'{A}ngel Aledo S\'{a}nchez\index[ppl]{Aledo S\'{a}nchez, Juan \'{A}ngel|textbf}, C\'{a}lculo geom\'{e}trico del l\'{i}mite de sucesiones trigonom\'{e}tricas. \emph{SUMA} \textbf{34}, 53--58. 

{\footnotesize Source: \href{https://riunet.upv.es/entities/publication/c2e7f47b-7d96-4637-9d30-bc1d3cdefb41}{Universitat Polit\`{e}cnica de Val\`{e}ncia}.}

{\footnotesize Cited in
\hyperref[Cor2000b]{\textsc{Cort\'{e}s L\'{o}pez 2000}}\index[ppl]{Cort\'{e}s L\'{o}pez, Juan Carlos}.}

\vspace{4pt}

From the abstract: ``At present, the calculation of the limit of a sequence, both in high school and in the first university courses, has been carried out from an exclusively analytical-algebraic approach. In this article we propose a rich geometric approach to initiate this theme in the classroom, particularizing our proposal for the trigonometric sequences." By ``trigonometric sequences" the authors mean not only the continued square roots arising from trig functions\index{continued square roots!and trigonometric functions}\index{continued square roots!geometric interpretations of}, but also limits like $\lim_{x\to0}(\sin x)/x$, all of which receive geometric interpretations. 


\item\vspace{9pt} \label{Cor2000b} (2000) Juan Carlos Cort\'{e}s L\'{o}pez\index[ppl]{Cort\'{e}s L\'{o}pez, Juan Carlos|textbf}, Algunas representaciones radicales infinitas de los n\'{u}meros naturales. \emph{Bolet\'{i}n Sociedad ``Puig Adam"} \textbf{54}, 29--38. 
\index[ppl]{Ramanujan, Srinivasa}
\index{radicales infinitas@\emph{radicales infinitas}} 

{\footnotesize \href{https://mathscinet.ams.org/mathscinet/article?mr=1751288}{MR1751288}}

{\footnotesize Source: \href{https://www.ucm.es/data/cont/media/www/pag-89521/Boletin 54 de Soc PUIG ADAM.pdf}{Universidad Complutense Madrid}.}

\vspace{4pt}

The author independently rediscovers the two continued square root expansions of \hyperref[Ram1911]{\textsc{Ramanujan 1911}} and their underlying structure\index{continued square roots!Ramanujan's}, then extends this to cube roots\index{continued cube roots} and $i$th roots\index{continued robb roots@continued $r$th roots!representation of integers as}, and ultimately shows how any positive integer may be represented\index{integers!representation by continued $r$th roots}\index{continued robb roots@continued $r$th roots!representation of integers as} The only reference is to a 1999 preprint later published as \hyperref[Cor2000a]{\textsc{Cort\'{e}s L\'{o}pez and Aledo S\'{a}nchez 2000}}.


\item\vspace{9pt} \label{Swe2000} (2000) N. M. Swerdlow\index[ppl]{Swerdlow, N. M.|textbf}, Kepler's iterative solution to Kepler's equation. \emph{Journal for the History of Astronomy}, \textbf{31}(4), 339--341.
\index{Kepler problem}

{\footnotesize \href{https://mathscinet.ams.org/mathscinet/article?mr=1802586}{MR1802586}}

{\footnotesize Source: \href{https://doi.org/10.1177/002182860003100404}{doi.org (SageJournals)}.}

Gives an account of the geometric reasoning behind Kepler's iterative solution to equation \eqref{E:Kep1621}. The author writes, ``In \emph{Astronomia nova 60}, Kepler\index[ppl]{Kepler, Johannes} set out demonstrations and equations for finding the position of a planet moving in either an ellipse or a circle such that a line from the Sun to the planet describes areas proportional to time. One of them, written as $M = E + e \sin E$ [compare \eqref{E:Gau1809-0}], known as `Kepler's Equation'\index{Kepler problem}, has a well known difficulty, which Kepler himself was the first to recognize, in that it can easily be solved in one direction, for $M$, but not in the other, for $E$\ldots. The problem, as Kepler suspected and as Newton\index[ppl]{Newton, Isaac}, among others, was later to show, is insoluble, but there are various ways of finding $E$ from $M$ to any degree of precision by iteration or interpolation. In \emph{Epitome 5.2.4} Kepler gives an iteration for solving Kepler's equation, which must surely be the first\ldots. Kepler's explanation is purely numerical, with no statement of the procedure, but with a really bizarre error --- the eccentricity of Mars is given in seconds of arc with two digits inverted, $11910''$ for $19110''$, and so used --- as a result of which the whole series of computations is incorrect." See also \hyperref[Gau1809]{\textsc{Gauss 1809}}.


\item\vspace{9pt} \label{Gal2001} \dag(2001) Massimo Galuzzi\index[ppl]{Galuzzi, Massimo|textbf}, Galois' note on the approximative solution of numerical equations (1830). \emph{Archive for History of Exact Sciences} \textbf{56}, 29--37.
\index{iteration!of functions}

{\footnotesize \href{https://mathscinet.ams.org/mathscinet/article?mr=1868932}{MR1868932}}

{\footnotesize Source: \href{https://www.jstor.org/stable/41134129}{JSTOR}.}

This paper succinctly explains both \hyperref[Leg1816]{\textsc{Legendre 1816}}\index[ppl]{Legendre, Adrien-Marie} and Galois' improvements to it in 1830. The author is keen to emphasize, via \hyperref[Gal1830]{\textsc{Galois 1830}}\index[ppl]{Galois, \'{E}variste}, that Galois was interested in areas of mathematics other than the purely abstract ones for which he is best known.


\item\vspace{9pt} \label{Kub2002} (2002) G. Kuba\index[ppl]{Kuba, G.|textbf} and J. Schoissengeier\index[ppl]{Schoissengeier, J.|textbf}, Dyadische Kettenwurzeln. \emph{Wissenschaftliche Nachrichten. Herausgegeben vom Bundesministerium f\"{u}r Bildung, Wissenschaft und Kultur} \textbf{119}, 23--24. 
\index{Kettenwurzeln@\emph{Kettenwurzeln}}
\index{continued square roots!of terms $a_n=\pm2$}

{\footnotesize Source: \href{https://web.archive.org/web/20160409121329/https://www.bmbf.gv.at/schulen/sb/wina/wina119.pdf?4rxn45}{Internet Archive}.}

Citing two of Pincherle's papers from 1917--1918 and \hyperref[Pol1925]{\textsc{P\'{o}lya and Szeg\"{o} 1925}}\index[ppl]{Polya@P\'{o}lya, G.}\index[ppl]{Szeg\"{o}, G.}, the authors write, ``The aim of this note is to prove that every real number $z$ can be represented in the interval $[0,2]$ by means of dyadic continued square roots."\index{continued square roots!of terms $a_n=\pm2$} They take an unusual, symbolic approach to these representations: ``Let $\Sigma = \{+, -\}^N$ be the set of all sequences whose $n$th term $\sigma(n)$ is either a plus or a minus symbol. For any sequence $\sigma\in\Sigma$ let $W(\sigma)$ be given by
\[W(\sigma):=\sqrt{2\sigma(1)\sqrt{2\sigma(2)\sqrt{2\sigma(3)\sqrt{2\sigma(4)\sqrt{2\ldots}}}}}\;,\]
where we first assume that the expression $W(\sigma)$ represents a well-defined real number for each $\sigma\in\Sigma$." 


\item\vspace{9pt} \label{Ser2002} (2002) L. D. Servi\index[ppl]{Servi, L. D.|textbf}, Problem 10973. \emph{The American Mathematical Monthly} \textbf{109} (9), 854. 
Solution by Richard Stong\index[ppl]{Stong, Richard}, \emph{ibid.} \textbf{111}(7), 627--628.
\index{continued square roots!of terms $a_n=\pm2$}

{\footnotesize Solution source: \href{https://doi.org/10.1080/00029890.2004.11920123}{doi.org (Taylor \& Francis)}.}

\vspace{4pt}

The problem states: ``With $R_k(n)$ defined as below, prove that $\lim_{k\to\infty}$ $R_k(2)/$$R_k(3)=3/2$.\index{continued square roots!of terms $a_n=\pm2$}
\[R_k(n)=\underbrace{\sqrt{2-\sqrt{2+\sqrt{2+\sqrt{2+\cdots+\sqrt{2+\sqrt{n}}}}}}}_{k\textrm{ square roots}}\;."\]


\item\vspace{9pt} \label{Fin2003} (2003) Steven R. Finch\index[ppl]{Finch, Steven R.|textbf}, \emph{Mathematical Constants}. Cambridge University Press, Cambridge, England, 7--9. ISBN 0-521-81805-2. 

{\footnotesize \href{https://mathscinet.ams.org/mathscinet/article?mr=2003519}{MR2003519}}

\vspace{4pt}

Section 1.2 mentions the continued square root having terms $a_n=1$\index{continued square roots!of terms $a_n=1$}; subsection 1.2.1 analyzes its convergence (as an iterated function composition, using methods given in \hyperref[Par1987]{\textsc{Paris 1987}}) to the golden ratio\index{golden ratio ($\tfrac{1+\sqrt{5}}{2}=1.61803\ldots$)}\index{constants, named!golden ratio ($\tfrac{1+\sqrt{5}}{2}=1.61803\ldots$)} $\tfrac{1}{2}(1+\sqrt{5})$; the subsection also alludes to the continued square root having terms $a_i=i$ (whose limit is called the Kasner number\index{Kasner number $(1.75793\ldots)$}\index{constants, named!Kasner number $(1.75793\ldots)$} in \hyperref[Her1935]{\textsc{Herschfeld 1935}}\index[ppl]{Herschfeld, Aaron}). Subsection 1.2.2 analyzes the continued cube root of terms $a_i=1$\index{continued cube roots} (whose limit is called the plastic constant\index{plastic constant ($1.32471\ldots$)}\index{constants, named!plastic constant ($1.32471\ldots$)} in \hyperref[Lim2010]{\textsc{Lim 2010}}).


\item\vspace{9pt} \label{Hau2003} (2003) Clemens Hauser\index[ppl]{Hauser, Clemens|textbf}, Pi, e und Kettenwurzeln. \emph{Der mathematische und naturwissenschaftliche Unterricht} \textbf{56}(4), 201--203.
\index{Kettenwurzeln@\emph{Kettenwurzeln}}
\index{continued square roots!and trigonometric functions}
\index{continued square roots!of terms $a_n=\pm2$}

\vspace{4pt}

From the abstract: ``If we approach the circumference of the unit circle by means of a regular polygon\index{polygons, regular}, we obtain attractive representations of $\pi$ in the form of continued square roots\index{continued square roots}. Only the mathematics of the intermediate level are needed for this. Further investigations of continued square roots can be envisaged in the advanced level. Here, a connection to Euler's number $e$ can be established." The author derives the long-established\index{pi@$\pi\;(3.14159\ldots)$!continued square root expressions for}\index{constants, named!pi@$\pi\;(3.14159\ldots)$}
\[\pi=\lim_{n\to\infty}2^n\underbrace{\sqrt{2-\sqrt{2+\sqrt{2+\sqrt{2+\cdots+\sqrt{2}}}}}}_{n\textrm{ square roots}}\;,\]
found at least as far back as \hyperref[Cat1842]{\textsc{Catalan 1842}}\index[ppl]{Catalan, E.}. The number $e$\index{e@$e\;(2.71828\ldots)$}\index{constants, named!e@$e\;(2.71828\ldots)$} arises in connection with the expression
\[\ln 2=\lim_{n\to\infty}2^n\underbrace{\sqrt{\sqrt{2+\ldots+\sqrt{2+2.5}}-2}}_{n+1\textrm{ square roots}}\;,\]
which is obtained by iterating hyperbolic trigonometric functions.\index{hyperbolic trigonometric functions}


\item\vspace{9pt} \label{Luo2003} (2003) Qiu-Ming Luo\index[ppl]{Luo, Qiu-Ming|textbf}, Feng Qi\index[ppl]{Qi, Feng|textbf}, Neil S. Barnett\index[ppl]{Barnett, Neil S.|textbf}, Sever S. Dragomir\index[ppl]{Dragomir, Sever S.|textbf}, Inequalities involving the sequence $\sqrt[3]{a+\sqrt[3]{a+\cdots+\sqrt[3]{a}}}$. \emph{Mathematical Inequalities and Applications} \textbf{6}(3), 413--419. 

{\footnotesize \href{https://mathscinet.ams.org/mathscinet/article?mr=1992481}{MR1992481}}

{\footnotesize Source: \href{https://dx.doi.org/10.7153/mia-06-38}{dx.doi.org (Ele-Math)}.}

{\footnotesize Cited in
\hyperref[Xi2013]{\textsc{Xi and Qi 2013}}\index[ppl]{Xi, Bo-Yan}\index[ppl]{Qi, Feng}.}

\vspace{4pt}

From the abstract: ``The convergence of the sequence [with terms]
\[\underbrace{\sqrt[3]{a+\sqrt[3]{a+\sqrt[3]{a+\cdots+\sqrt[3]{a}}}}}_{n\textrm{ radicals}}\]
is proved, and some inequalities involving this sequence are established for $a>0$. As by-product[s], two identities involving irrational numbers are obtained. Two open problems are proposed." The authors define
\[S_{n,t}(a)=\underbrace{\sqrt[t]{a+\sqrt[t]{a+\sqrt[t]{a+\cdots+\sqrt[t]{a}}}}}_{n\textrm{ radicals}}\]
and
\[f_{n,t}(a)=\dfrac{a-S_{n+1,t}}{a-S_{n,t}}\;,\]
and derive upper and lower bounds for $f_{3,t}(a)$ for ranges of positive $a$. 


\item\vspace{9pt} \label{Ser2003} (2003) L. D. Servi\index[ppl]{Servi, L. D.|textbf}, Nested square roots of 2. \emph{The American Mathematical Monthly} \textbf{110}(4), 326--330.
\index{continued square roots!and trigonometric functions}
\index{continued square roots!of terms $a_n=\pm2$}

{\footnotesize \href{https://mathscinet.ams.org/mathscinet/article?mr=1984573}{MR1984573}}

{\footnotesize Source: \href{https://doi.org/10.1080/00029890.2003.11919968}{doi.org (Taylor \& Francis)}.}

{\footnotesize Cited in 
\hyperref[Nyb2005]{\textsc{Nyblom 2005}}\index[ppl]{Nyblom, M. A.}, 
\hyperref[Jon2008]{\textsc{Jones 2008}}\index[ppl]{Jones, Dixon J.}, 
\hyperref[Zim2008a]{\textsc{Zimmerman and Ho 2008a}}\index[ppl]{Zimmerman, Seth}\index[ppl]{Ho, Chungwu}, 
\hyperref[Eft2012]{\textsc{Efthimiou 2012}}\index[ppl]{Efthimiou, Costas J.}, 
\hyperref[Nyb2012]{\textsc{Nyblom 2012}}\index[ppl]{Nyblom, M. A.}, 
\hyperref[Mor2012]{\textsc{Moreno and Garc\'{i}a-Caballero 2012}}, \hyperref[Mor2013a]{\textsc{2013a}}, and \hyperref[Mor2013b]{\textsc{2013b}}\index[ppl]{Moreno, Samuel G\'{o}mez}\index[ppl]{Garc\'{i}a-Caballero, Esther M.}, 
\hyperref[Sen2013]{\textsc{Senadheera 2013}}\index[ppl]{Senadheera, Jayantha}, 
\hyperref[Gar2014a]{\textsc{Garc\'{i}a-Caballero, Moreno, and Prophet 2014a}}, \hyperref[Gar2014b]{\textsc{2014b}}, and \hyperref[Gar2014c]{\textsc{2014c}}\index[ppl]{Garc\'{i}a-Caballero, Esther M.}\index[ppl]{Moreno, Samuel G\'{o}mez}\index[ppl]{Prophet, Michael P.}, 
\hyperref[Lyn2014]{\textsc{Lynd 2014}}\index[ppl]{Lynd, Chris D.}, and 
\hyperref[Les2016]{\textsc{Lesher and Lynd 2016}}\index[ppl]{Lesher, Devyn A.}\index[ppl]{Lynd, Chris D.}, and 
\hyperref[Vel2016c]{\textsc{Vellucci and Bersani 2016c}}\index[ppl]{Vellucci, Pierluigi}\index[ppl]{Bersani, Alberto Maria}.}

\vspace{4pt}

Another independent rediscovery of continued square roots derived from iterated trigonometric functions.\index{continued square roots}\index{continued square roots!and trigonometric functions} Compare \hyperref[Wie1904b]{\textsc{Wiernsberger 1904b}}, \hyperref[Cip1908]{\textsc{Cipolla 1908}}\index[ppl]{Cipolla, Michele}, \hyperref[Pol1925]{\textsc{P\'{o}lya and Szeg\"{o} 1925}}\index[ppl]{Polya@P\'{o}lya, G.}, \hyperref[Myr1958]{\textsc{Myrberg 1958}}, and others above.


\item\vspace{9pt} \label{Ber2004} (2004) Adriana Berechet\index[ppl]{Berechet, Adriana|textbf}, Solving a conjecture about certain $f$-expansions. \emph{Proceedings of the Romanian Academy -- Series A: Mathematics, Physics, Technical Sciences, Information Science} \textbf{5}(3), 231--235. 
\index{f@$f$-expansions}

{\footnotesize \href{https://mathscinet.ams.org/mathscinet/article?mr=2122307}{MR2122307}}

{\footnotesize Source: \href{https://acad.ro/sectii2002/proceedings/doc3_2004/01_Berechet.pdf}{Romanian Academy}.}

From the abstract: ``The conjecture asserts that the equivalence of the label sequence of the regular continued fraction expansion\index{continued fractions!regular!label sequence of} to the sequence $(\xi_n)_{n\in N+}$ associated with it $\ldots$ still holds for the label sequence of any $f$-expansion satisfying [two conditions]. We prove that [the first] condition and a strengthening of a Lipschitz condition\index{Lipschitz condition} $\ldots$ are sufficient to ensure a necessary and sufficient condition under which the asserted equivalence holds. The proof involves processes on several probability spaces\index{probability spaces} and some associated dynamical systems\index{dynamical systems} relating the $f$-expansion $\ldots$ to [random variables] on the probability space used in the concluding theorem."


\item\vspace{9pt} \label{Mar2004} (2004) Greg Martin\index[ppl]{Martin, Greg|textbf}, The unreasonable effectualness of continued function expansions. \emph{Journal of the Australian Mathematical Society} \textbf{77}(3), 305--319.
\index{f@$f$-expansions}
\index{continued compositions}

{\footnotesize \href{https://mathscinet.ams.org/mathscinet/article?mr=2099803}{MR2099803}}

{\footnotesize Source: \href{https://doi.org/10.1017/S1446788700014452}{doi.org (CambridgeCore)}.}

{\footnotesize Cited in 
\hyperref[Jon2015]{\textsc{Jones 2015}}\index[ppl]{Jones, Dixon J.}.}

\vspace{4pt}

From the introduction: ``[W]e focus on the $f$-expansions introduced [in \hyperref[Bis1944]{\textsc{Bissinger 1944}}]$\ldots$ The purpose of this paper is to demonstrate that the function $f$ can be chosen so that the expansions of prescribed real numbers can have essentially any desired behavior. The following results, listed in roughly increasing order of unlikeliness, are representative of what we can prove.

``\textbf{Theorem 1.} For any two real numbers $x,y \in (0,1)$, there exists a function $f$ such that the $f$-expansion of $x$ is the same as the usual continued fraction expansion of $y$.

``\textbf{Theorem 2.} There exists a function $f$ such that the $f$-expansion of any rational or quadratic irrational terminates.

``\textbf{Theorem 3.} There exists a function $f$ such that the $f$-expansion of a real number $x$ is periodic if and only if $x$ is a cubic irrational number.

``\textbf{Theorem 4.} There exists a function $f$ such that, simultaneously for every integer $d\ge 1$, a real number $x$ is algebraic of degree $d$ if and only if the $f$-expansion of $x$ terminates with the integer $d+1$.

``[W]e should confess what the reader might already suspect, that the functions giving the nice behaviors of Theorems 1--4 are infeasible for actual computations. Indeed, the existence of such functions is essentially a consequence of the existence of continuous functions on the interval $(0,1)$ with certain properties."

This confession notwithstanding, the paper concludes with some computations involving functions $f_\alpha$ defined by
\[a_0 + (a_1 + (a_2 + (a_3 + (a_4 + \cdots )^{-\alpha})^{-\alpha})^{-\alpha})^{-\alpha}\;,\]
where $\alpha>0$. (For $\alpha=\tfrac{1}{2}$, this has been called a continued reciprocal square root\index{continued reciprocal roots}.) Questions and conjectures are posed about the termination of such expansions.


\item\vspace{9pt} \label{Lev2005} (2005) Aaron Levin\index[ppl]{Levin, Aaron|textbf}, A new class of infinite products generalizing Vi\`{e}te's\index[ppl]{Vi\`{e}te, Fran\c{c}ois} product formula for $\pi$. \emph{The Ramanujan Journal} \textbf{10}(3), 305--324. 
\index{Vi\`{e}te's formula for $\tfrac{2}{\pi}$!generalizations}
\index[ppl]{Vi\`{e}te, Fran\c{c}ois}

{\footnotesize \href{https://mathscinet.ams.org/mathscinet/article?mr=2193382}{MR2193382}}

{\footnotesize Source: \href{https://doi.org/10.1007/s11139-005-4852-z}{doi.org (Springer)}.}

{\footnotesize Cited in 
\hyperref[Lev2006]{\textsc{Levin 2006}}\index[ppl]{Levin, Aaron}, 
\hyperref[Mor2013a]{\textsc{Moreno and Garc\'{i}a-Caballero 2013a}} and \hyperref[Mor2013b]{\textsc{2013b}}\index[ppl]{Moreno, Samuel G\'{o}mez}\index[ppl]{Garc\'{i}a-Caballero, Esther M.}, 
\hyperref[Sen2013]{\textsc{Senadheera 2013}}\index[ppl]{Senadheera, Jayantha}, 
\hyperref[Gar2014a]{\textsc{Garc\'{i}a-Caballero, Moreno, and Prophet 2014a}}\index[ppl]{Garc\'{i}a-Caballero, Esther M.}\index[ppl]{Moreno, Samuel G\'{o}mez}\index[ppl]{Prophet, Michael P.}, 
\hyperref[Nis2015b]{\textsc{Nishimura 2015b}}\index[ppl]{Nishimura, Ryo}, 
\hyperref[Nis2016]{\textsc{Nishimura 2016}}\index[ppl]{Nishimura, Ryo}, and 
\hyperref[Osl2016a]{\textsc{Osler 2016}}\index[ppl]{Osler, Thomas J.}.}

\vspace{4pt}

From the abstract: ``We show how functions $F(z)$ which satisfy an identity of the form $F(\alpha z) = g(F(z))$ for some complex number $\alpha$ and some function $g(z)$ give rise to infinite product formulas\index{infinite products!of continued square roots} that generalize Vi\`{e}te's\index[ppl]{Vi\`{e}te, Fran\c{c}ois} product formula for $\pi$. Specifically, using elliptic\index{infinite products!of elliptic functions} and trigonometric functions\index{infinite products!of trigonometric functions} we derive closed form expressions for some of these infinite products. By evaluating the expressions at certain points we obtain formulas expressing infinite products involving nested radicals\index{nested radicals!infinite products of}\index{infinite products!of nested radicals} in terms of well-known constants. In particular, simple infinite products for $\pi$\index{pi@$\pi\;(3.14159\ldots)$!infinite products for}\index{constants, named!pi@$\pi\;(3.14159\ldots)$} and the lemniscate constant\index{lemniscate!constant\;$(2.62205\ldots)$} are obtained." 


\item\vspace{9pt} \label{Nyb2005} (2005) M. A. Nyblom\index[ppl]{Nyblom, M. A.|textbf}, More nested square roots of 2. \emph{The American Mathematical Monthly} \textbf{112}(9), 822--825. 

{\footnotesize \href{https://mathscinet.ams.org/mathscinet/article?mr=2179862}{MR2179862}}

{\footnotesize Source: \href{https://doi.org/10.1080/00029890.2005.11920256}{doi.org (Taylor \& Francis)}.} 

{\footnotesize Cited in
\hyperref[Zim2008a]{\textsc{Zimmerman and Ho 2008a}}\index[ppl]{Zimmerman, Seth}\index[ppl]{Ho, Chungwu}, 
\hyperref[Eft2012]{\textsc{Efthimiou 2012}}\index[ppl]{Efthimiou, Costas J.}, 
\hyperref[Nyb2012]{\textsc{Nyblom 2012}}\index[ppl]{Nyblom, M. A.}, 
\hyperref[Mor2012]{\textsc{Moreno and Garc\'{i}a-Caballero 2012}} and \hyperref[Mor2013b]{\textsc{2013b}}\index[ppl]{Moreno, Samuel G\'{o}mez}\index[ppl]{Garc\'{i}a-Caballero, Esther M.}, 
\hyperref[Sen2013]{\textsc{Senadheera 2013}}\index[ppl]{Senadheera, Jayantha}, 
\hyperref[Gar2014a]{\textsc{Garc\'{i}a-Caballero, Moreno, and Prophet 2014a}} and \hyperref[Gar2014c]{\textsc{2014c}}\index[ppl]{Garc\'{i}a-Caballero, Esther M.}\index[ppl]{Moreno, Samuel G\'{o}mez}\index[ppl]{Prophet, Michael P.}, 
\hyperref[Lyn2014]{\textsc{Lynd 2014}}\index[ppl]{Lynd, Chris D.}, and
\hyperref[Vel2016c]{\textsc{Vellucci and Bersani 2016c}}\index[ppl]{Vellucci, Pierluigi}\index[ppl]{Bersani, Alberto Maria}.}

\vspace{4pt}

This note's principle result is: If $x\ge 2$ and $k$ is a positive integer, then\index{continued square roots!of terms $a_n=\pm2$}
\begin{align*}
R_k(x)&=\sqrt{2+\sqrt{2+\cdots+\sqrt{2+x}}}\\
&=\left(\dfrac{x+\sqrt{x^2-4}}{2}\right)^{1/2^k}+\left(\dfrac{x+\sqrt{x^2-4}}{2}\right)^{-1/2^k}\;,
\end{align*}
and
\begin{align*}
\tilde{R}_k(x)&=\sqrt{-2+\sqrt{2+\cdots+\sqrt{2+x}}}\\
&=\left(\dfrac{x+\sqrt{x^2-4}}{2}\right)^{1/2^k}-\left(\dfrac{x+\sqrt{x^2-4}}{2}\right)^{-1/2^k}\;,
\end{align*}
where there are $k$ square roots in each of $R_k$ and $\tilde{R}_k$. 


\item\vspace{9pt} \label{Rao2005} (2005) K. Srinivasa Rao\index[ppl]{Rao, K. Srinivasa|textbf} and G. Vanden Berghe\index[ppl]{Vanden Berghe, G.|textbf}, On an entry of Ramanujan in his notebooks: a nested square root expansion. \emph{Journal of Computational and Applied Mathematics} \textbf{173}(2), 371--378. 
\index[ppl]{Ramanujan, Srinivasa}
\index{continued square roots!Ramanujan's}

{\footnotesize \href{https://mathscinet.ams.org/mathscinet/article?mr=2102903}{MR2102903}}

{\footnotesize Source: \href{https://doi.org/10.1016/j.cam.2004.04.009}{doi.org (ScienceDirect)}.}

From the abstract: ``In this letter, the elementary result of Ramanujan for nested roots, also called continued or infinite radicals\index{continued square roots!Ramanujan's}\index{infinite radicals}, for a given integer $N$, expressed by him as a simple sum of three parts $(N=x+n+a)$ [see equation \eqref{E:McG-1} below] is shown to give rise to two distinguishably different expansion formulas. One of these is due to Ramanujan and surprisingly, it is this other formula, not given by Ramanujan, which is more rapidly convergent!"


\item\vspace{9pt} \label{Lev2006} (2006) Aaron Levin\index[ppl]{Levin, Aaron|textbf}, A geometric interpretation of an infinite product for the lemniscate constant\index{lemniscate!constant\;$(2.62205\ldots)$}. \emph{The American Mathematical Monthly} \textbf{113}(6), 510--520.
\index{Vi\`{e}te's formula for $\tfrac{2}{\pi}$!generalizations}
\index[ppl]{Vi\`{e}te, Fran\c{c}ois}

{\footnotesize \href{https://mathscinet.ams.org/mathscinet/article?mr=2231136}{MR2231136}}

{\footnotesize Source: \href{https://doi.org/10.1080/00029890.2006.11920331}{doi.org (Taylor \& Francis)}.} 

{\footnotesize Cited in
\hyperref[Mor2013a]{\textsc{Moreno and Garc\'{i}a-Caballero 2013a}}\index[ppl]{Moreno, Samuel G\'{o}mez}\index[ppl]{Garc\'{i}a-Caballero, Esther M.}, 
\hyperref[Sen2013]{\textsc{Senadheera 2013}}\index[ppl]{Senadheera, Jayantha},
\hyperref[Nis2015b]{\textsc{Nishimura 2015b}}\index[ppl]{Nishimura, Ryo}, 
\hyperref[Nis2016]{\textsc{Nishimura 2016}}\index[ppl]{Nishimura, Ryo}, and 
\hyperref[Osl2016a]{\textsc{Osler 2016}}\index[ppl]{Osler, Thomas J.}.}

\vspace{4pt}

The author begins by comparing the Vi\`{e}te formula for $\tfrac{2}{\pi}$ (equation \eqref{E:Vie1593-1} above) \index{Vi\`{e}te's formula for $\tfrac{2}{\pi}$} with 
\begin{equation}\label{E:Lev2006-1}
\dfrac{2}{L}=\sqrt{\tfrac{1}{2}}\cdot\sqrt{\tfrac{1}{2}+\dfrac{\tfrac{1}{2}}{\sqrt{\tfrac{1}{2}}}}\cdot\sqrt{\tfrac{1}{2}+\dfrac{\tfrac{1}{2}}{\sqrt{\tfrac{1}{2}+\dfrac{\tfrac{1}{2}}{\sqrt{\tfrac{1}{2}}}}}}\cdots
\end{equation}
where
\begin{equation*}
L=\dfrac{B\left(\tfrac{1}{4},\tfrac{1}{4}\right)}{2\sqrt{2}}=\dfrac{\Gamma\left(\tfrac{1}{4}\right)^2}{2\sqrt{2\pi}}=2.6220575542\ldots
\end{equation*}
is the lemniscate constant\index{lemniscate!constant\;$(2.62205\ldots)$}, $B(x,y)$ is the beta function\index{beta function}, and $\Gamma(z)$ is the gamma function\index{gamma function}. The author continues, ``We will see that the area enclosed by the curve $C_4$ defined by $x^4+y^4=1$ is $L\sqrt{2}$, and this will allow us to give a geometric meaning to the product formula \eqref{E:Lev2006-1}$\,\ldots$ We will show that the similarity between equations \eqref{E:Vie1593-1} and \eqref{E:Lev2006-1} goes beyond mere typographical appearances. We will see that \eqref{E:Lev2006-1} is related to the curve $C_4$ in much the same way that Vi\`{e}te's\index[ppl]{Vi\`{e}te, Fran\c{c}ois} product is related to the circle." 


\item\vspace{9pt} \label{Osl2006} (2006) Thomas J. Osler, Interesting finite and infinite products from simple algebraic identities. \emph{The Mathematical Gazette} \textbf{90}(517), 90--93.
\index{Vi\`{e}te's formula for $\tfrac{2}{\pi}$!generalizations}

{\footnotesize Source: \href{https://www.jstor.org/stable/3621420}{JSTOR}.} 

{\footnotesize Cited in 
\hyperref[Lev2005]{\textsc{Levin 2005}}\index[ppl]{Levin, Aaron}.}

\vspace{4pt}

From the introduction: ``The difference of two squares, $x^2 - y^2 = (x + y)$ $(x - y)$ and its immediate generalisations $x^3 - y^3 = (x - y)(x2 + xy + y)$,\ldots, are among the most elementary identities. It is the purpose of this note to show that these simple formulas are at the heart of a number of advanced products, some finite, and some infinite. The reader may find some of these surprising, especially those involving trigonometric functions." Vi\`{e}te's\index[ppl]{Vi\`{e}te, Fran\c{c}ois} formula \eqref{E:Vie1593-1} above\index{Vi\`{e}te's formula for $\tfrac{2}{\pi}$} is one of the products derived.


\item\vspace{9pt} \label{Hum2007} (2007) Peter J. Humphries\index[ppl]{Humphries, Peter J.|textbf}, Nesting polynomials in infinite radicals. \emph{Bulletin of the Korean Mathematical Society} \textbf{44}(2), 331--336.
\index{infinite radicals}
\index{continued square roots!of terms which are polynomials}

{\footnotesize \href{https://mathscinet.ams.org/mathscinet/article?mr=2325034}{MR2325034}}

{\footnotesize Source: \href{https://doi.org/10.4134/BKMS.2007.44.2.331}{doi.org (KoreaScience)}.}

Let $\mathbb{R}[x]$ be the ring of polynomials\index{polynomials!ring of} in $x$ with real coefficients. The author writes: ``We consider infinite nested radicals in which the arguments are positive polynomial sequences\index{continued square roots!of terms which are polynomials}. It is shown that the evaluation of such a nesting is always finite, and we prove necessary and sufficient conditions for the evaluation to be a finite polynomial$\ldots$ In Section 2, we characterize when an infinite nested radical involving polynomials from $\mathbb{R}[x]$ has a simple closed form as another polynomial in $\mathbb{R}[x]$." In Section 3 it is proved that, for all positive polynomial sequences $a_n, b_n$, the limit
\[\lim_{n\to\infty}\sqrt{a_1+b_1\sqrt{a_2+\cdots+b_{n-1}\sqrt{a_n+b_n}}}\]
exists and is finite. 


\item\vspace{9pt} \label{Lim2007a} (2007a) Teik-Cheng Lim\index[ppl]{Lim, Teik-Cheng|textbf}, Significance of an infinite nested radical number and its application in van der Waals potential functions\index{van der Waals potential functions}. \emph{MATCH Communications in Mathematical and in Computer Chemistry} \textbf{57}(3), 549--556.
\index{nested radicals}
\index{continued robb roots@continued $r$th roots!of terms $a_n=1$}

{\footnotesize \href{https://mathscinet.ams.org/mathscinet/article?mr=2337287}{MR2337287 (citation only)}}

{\footnotesize Source: \href{https://match.pmf.kg.ac.rs/electronic_versions/Match57/n3/match57n3_549-556.pdf}{MATCH}.}

From the abstract: ``This paper demonstrates the relations between a new mathematical constant and its relevance to the molecular potential energy function\index{molecular potential energy function} commonly adopted in computational chemistry\index{computational chemistry} softwares. This mathematical constant, $1.7767750401$ (correct up to $12$ decimal places), which fulfills the following infinite nested radical equation\index{continued robb roots@continued $r$th roots!of terms $a_n=1$}
\begin{align*}
n&=\sqrt[n]{1+\sqrt[n]{1+\sqrt[n]{1+\sqrt[n]{1+\cdots}}}}\\
&=\sqrt[n]{n+\tfrac{1}{n}\times\sqrt[n]{n+\tfrac{1}{n}\times\sqrt[n]{n+\tfrac{1}{n}\times\sqrt[n]{n+\tfrac{1}{n}\times\cdots}}}}\;,
\end{align*}
is shown to be applicable as the indices of a generalized Lennard-Jones potential energy function\index{Lennard-Jones potential energy function}. This new potential function demonstrates very good agreement with the various versions of specific Lennard-Jones potential energy functions and the Buckingham potential function\index{Buckingham potential function} converted from the Lennard-Jones (12-6) function when the indices are positive integer multiples of the mathematical constant or when the indices are raised to the first four positive integer powers." 


\item\vspace{9pt} \label{Lim2007b} (2007b) Teik-Cheng Lim\index[ppl]{Lim, Teik-Cheng|textbf}, Infinite nested radicals: a mathematical poem. \emph{The American Mathematical Monthly} \textbf{114}(3), 182.
\index{continued robb roots@continued $r$th roots}

{\footnotesize Source: \href{https://www.jstor.org/stable/27642164}{JSTOR}.}

{\footnotesize Cited in 
\hyperref[Lim2010]{\textsc{Lim 2010}}\index[ppl]{Lim, Teik-Cheng}.}

\vspace{4pt}

The ``poem" consists of 9 identities, in as many lines, of the form
\[\sqrt[n]{n\times ((n+1)^n-1)+\sqrt[n]{n\times ((n+1)^n-1)+\cdots}}\]
\[=\sqrt[n]{n+((n+1)^n-1)\times\sqrt[n]{n+((n+1)^n-1)\times\cdots}}\]
for $n=2,\ldots,10$.


\item\vspace{9pt} \label{Nyb2007} (2007) M. A. Nyblom\index[ppl]{Nyblom, M. A.|textbf}, On the evaluation of a definite integral involving nested square root functions. \emph{Rocky Mountain Journal of Mathematics} \textbf{37}, number 4, 1301--1304. 
\index{nested square roots}
\index{continued square roots!finite products of}
\index{continued square roots!of terms $a_n=2$}
\index{continued square roots!definite integrals of}

{\footnotesize \href{https://mathscinet.ams.org/mathscinet/article?mr=2360300}{MR2360300}}

{\footnotesize Source: \href{https://www.jstor.org/stable/44239675}{JSTOR}.}

{\footnotesize Cited in 
\hyperref[Nyb2012]{\textsc{Nyblom 2012}}\index[ppl]{Nyblom, M. A.}.}

\vspace{4pt}

From the introduction: ``[W]e shall in this paper take advantage of a structural feature of a class of nested square root functions, in order that we may evaluate their corresponding definite integrals which on first acquaintance appear rather intractable. In particular the functions in question will be composed of a finite product of reciprocals of the form
\index{continued square roots!finite products of}
\index{continued square roots!of terms $a_n=2$}
\index{continued square roots!definite integrals of}
\[R_N(x)=\sqrt{2+\sqrt{2+\cdots+\sqrt{2+x}}}\;,\]
in which $R_N(x)$ consists of $N$ nested square roots\index{nested square roots}. By restricting the intervals of integration to finite subintervals of $[2,\infty)$, we shall see that a simple application of a hyperbolic function\index{hyperbolic trigonometric functions} substitution together with some standard identities will result in closed-form expressions for these definite integrals." The main theorem is that, for $\xi\ge2$ and $N$ a positive integer,
\begin{align*}
\dfrac{1}{2^N}\int_2^\xi\dfrac{dx}{\sqrt{2+x}\cdot\sqrt{2+\sqrt{2+x}}\cdots\sqrt{2+\sqrt{2+\cdots+\sqrt{2+x}}}}&\\
=\left(\dfrac{\xi+\sqrt{\xi^2-4}}{2}\right)^{1/2^N}+\left(\dfrac{\xi+\sqrt{\xi^2-4}}{2}\right)^{-1/2^N}&-2\;.
\end{align*}
where the last expression in the denominator of the integrand is $R_N(x)$. \index{continued square roots!and hyperbolic trigonometric functions} 


\item\vspace{9pt} \label{Osl2007} (2007) Thomas J. Osler\index[ppl]{Osler, Thomas J.|textbf}, Vieta-like products of nested radicals with Fibonacci and Lucas numbers. \emph{The Fibonacci Quarterly} \textbf{45}(3), 202--204.
\index{nested radicals!infinite products of}

{\footnotesize \href{https://mathscinet.ams.org/mathscinet/article?mr=2437033}{MR2437033}}

{\footnotesize Source: \href{https://www.fq.math.ca/Papers1/45-3/osler.pdf}{Fibonacci Quarterly}.}

{\footnotesize Cited in 
\hyperref[Gar2014a]{\textsc{Garc\'{i}a-Caballero, Moreno, and Prophet 2014a}} and \hyperref[Gar2014c]{\textsc{2014c}}\index[ppl]{Garc\'{i}a-Caballero, Esther M.}\index[ppl]{Moreno, Samuel G\'{o}mez}\index[ppl]{Prophet, Michael P.}.}

\vspace{4pt}

From the abstract: ``We present two infinite products of nested radicals\index{infinite products!of continued square roots}\index{nested radicals!infinite products of} involving Fibonacci\index{numbers!Fibonacci}\index{Fibonacci!numbers} and Lucas numbers\index{Lucas!numbers}\index{numbers!Lucas}. These products resemble Vieta's\index[ppl]{Vi\`{e}te, Fran\c{c}ois} classical product of nested radicals for $\tfrac{2}{\pi}$ [equation \eqref{E:Vie1593-2} above]\index{Vi\`{e}te's formula for $\tfrac{2}{\pi}$!generalizations}\index{nested radicals!infinite products of}. A modern derivation of Vieta's product involves trigonometric functions, while our product involves similar manipulations involving hyperbolic functions\index{hyperbolic trigonometric functions}." 


\item\vspace{9pt} \label{Riv2007} (2007) T. Rivoal\index[ppl]{Rivoal, T.|textbf}, Propri\'{e}t\'{e}s diophantiennes du d\'{e}veloppement en cotangente continue de Lehmer\index[ppl]{Lehmer, D. H.}. \emph{Monatshefte f\"{u}r Mathematik} \textbf{150}(1), 49--71. 
\index{continued cotangents} 


{\footnotesize \href{https://mathscinet.ams.org/mathscinet/article?mr=2297253}{MR2297253}}

{\footnotesize Source: \href{https://doi.org/10.1007/s00605-006-0415-7}{doi.org (Springer)}.} 

{\footnotesize Cited in 
\hyperref[Sch2016]{\textsc{Schweiger 2016}}.\index[ppl]{Schweiger, Fritz}}

\vspace{4pt}

From the English abstract: ``This article deals with an algorithm which enables us to write any real positive number as the sum of an alternating series of cotangents of integers $n_\nu, \nu \ge 0$, in a unique way. We continue the work begun by Lehmer [\hyperref[Leh1938]{\textsc{Lehmer 1938}}] and continued by Shallit [\hyperref[Sha1976]{\textsc{Shallit 1976}}]\index[ppl]{Shallit, Jeffrey}: amongst other things, we give explicitly the link between the rational approximations of a given real number coming from this algorithm and the usual convergents\index{convergents} of the same real number and we produce a quasi-optimal bound for the growth of the sequence $(n_\nu)_{\nu\ge0}$ associated to an algebraic number. We also determine the regular continued fractions\index{continued fractions!regular} of an exceptional class of continued cotangent developments, which enables us to produce optimal irrationality measures of these expansions."


\item\vspace{9pt} \label{Son2007} (2007) Jonathan Sondow\index[ppl]{Sondow, Jonathan|textbf} and Petros Hadjicostas\index[ppl]{Hadjicostas, Petros|textbf}, The generalized-Euler-constant function $\gamma(z)$ and a generalization of Somos's quadratic recurrence constant. \emph{Journal of Mathematical Analysis and Applications} \textbf{332}(1), 292--314. 

{\footnotesize \href{https://mathscinet.ams.org/mathscinet/article?mr=2319662}{MR2319662}}

{\footnotesize Source: \href{https://doi.org/10.1016/j.jmaa.2006.09.081}{doi.org (ScienceDirect)}.}

\vspace{4pt}

The authors define the generalized-Euler-constant function as
\[\gamma(z)=\sum_{n=1}^\infty z^{n-1}\left(\frac{1}{n}-\log\frac{n+1}{n}\right)\;, \quad |z|\le1\;,\]
which includes Euler's constant $\gamma=\gamma(1)$\index{Euler's constant!$\gamma$ ($0.57721\ldots$)}\index{constants, named!Euler's $\gamma$ ($0.57721\ldots$)} and the ``alternating Euler's constant'' $\log\tfrac{4}{\pi}=\gamma(-1)$\index{Euler's constant!alternating ($\log\tfrac{4}{\pi}=0.241564\ldots)$}\index{constants, named!Euler's alternating ($\log\tfrac{4}{\pi}=0.241564\ldots)$}. In Section 3, the authors define the generalized Somos\index[ppl]{Somos, Michael} constant\index{Somos's constant ($1.66168\ldots$)}\index{constants, named!Somos's constant ($1.66168\ldots$)} $\sigma_t$, $t > 1$:
\[\sigma_t = \sqrt[t]{1\sqrt[t]{2\sqrt[t]3\sqrt[t]{\cdots}}}=1^{1/t}2^{1/t^2}3^{1/t^3}\cdots=\prod_{n=1}^\infty n^{1/t^n}\;.\]
This is shown to be related to the \hyperref[Ram1911]{\textsc{Ramanujan 1911}} continued square root \eqref{E:Ram1911-1} above.\index[ppl]{Ramanujan, Srinivasa}


\item\vspace{9pt} \label{Joh2008} (2008) Jamie Johnson\index[ppl]{Johnson, Jamie|textbf} and Tom Richmond\index[ppl]{Richmond, Tom|textbf}, Continued radicals. \emph{The Ramanujan Journal} \textbf{15}(2), 259--273. 
\index{continued radicals}
\index{continued square roots!of finitely many nonnegative real terms}

{\footnotesize \href{https://mathscinet.ams.org/mathscinet/article?mr=2377579}{MR2377579}}

{\footnotesize Source: \href{https://doi.org/10.1007/s11139-007-9076-y}{doi.org (Springer)}.}

{\footnotesize Cited in
\hyperref[Eft2012]{\textsc{Efthimiou 2012}}\index[ppl]{Efthimiou, Costas J.}, 
\hyperref[Cla2014]{\textsc{Clark and Richmond 2014}}\index[ppl]{Clark, Tyler}\index[ppl]{Richmond, Tom}, and 
\hyperref[Lyn2014]{\textsc{Lynd 2014}}\index[ppl]{Lynd, Chris D.}, and
\hyperref[Vel2016c]{\textsc{Vellucci and Bersani 2016c}}\index[ppl]{Vellucci, Pierluigi}\index[ppl]{Bersani, Alberto Maria}.}

\vspace{4pt}

From the abstract: ``We consider the set of real numbers $S(M)$ representable as a continued radical whose [non-negative] terms $a_1,a_2,\ldots$ are all from a finite set $M$. We give conditions on the set $M$ for $S(M)$ to be (a) an interval, and (b) homeomorphic to the Cantor set.\index{Cantor set}" The authors also derive an upper bound on the derivative of a finite approximant to a continued $r_i$th root.\index{continued rocc roots@continued $r_i$th roots!approximants of}


\item\vspace{9pt} \label{Jon2008} (2008) Dixon J. Jones\index[ppl]{Jones, Dixon J.|textbf}, Letter to the Editor. \emph{Mathematics Magazine} \textbf{81}(3), 230.
\index{history!of continued square roots}
\index{letters to editors}

{\footnotesize Source: \href{https://doi.org/10.1080/0025570X.2008.11953557}{doi.org (Taylor \& Francis)}.}

\vspace{4pt}

Provides some additional references for \hyperref[Zim2008a]{\textsc{Zimmerman and Ho 2008a}}\index[ppl]{Zimmerman, Seth}\index[ppl]{Ho, Chungwu}, and points out that the continued square roots derived from trigonometric identities\index{continued square roots!and trigonometric functions} have been independently rediscovered many times over the past century.


\item\vspace{9pt} \label{Lim2008} (2008) Teik-Cheng Lim\index[ppl]{Lim, Teik-Cheng|textbf}, Two infinite nested radical constants. \emph{The Mathematical Gazette}, \textbf{92}(523), 96--97. 
\index{nested radicals}

{\footnotesize Source: \href{https://www.jstor.org/stable/27821728}{JSTOR}.} 

\vspace{4pt}

Computes constants $m=0.4758608124$ and $n=2.398384383$ (correct to 10 significant figures) satisfying\index{continued robb roots@continued $r$th roots!of constant nonnegative real terms}
\[m=\sqrt[m]{m+m\cdot\sqrt[m]{m+m\cdot\sqrt[m]{m+m\sqrt[m]{\cdots}}}}\;,\]
and
\[n=\sqrt[n]{n+n\cdot\sqrt[n]{n+n\cdot\sqrt[n]{\cdots}}}=\sqrt[n]{n\cdot n+\sqrt[n]{n\cdot n+\sqrt[n]{\cdots}}}\;.\]


\item\vspace{9pt} \label{Zim2008a} (2008a) Seth Zimmerman\index[ppl]{Zimmerman, Seth|textbf} and Chungwu Ho\index[ppl]{Ho, Chungwu|textbf}, On infinitely nested radicals. \emph{Mathematics Magazine} \textbf{81}(1), 3--15. 
\index{nested radicals}
\index{continued square roots}

{\footnotesize \href{https://mathscinet.ams.org/mathscinet/article?mr=2380054}{MR2380054}}

{\footnotesize Source: \href{https://doi.org/10.1080/0025570X.2008.11953522}{doi.org (Taylor \& Francis)}.} 

{\footnotesize Cited in 
\hyperref[Eft2012]{\textsc{Efthimiou 2012}}\index[ppl]{Efthimiou, Costas J.}, 
\hyperref[Eft2013]{\textsc{Efthimiou 2013}}\index[ppl]{Efthimiou, Costas J.}, 
\hyperref[Muk2013]{\textsc{Mukherjee 2013}}\index[ppl]{Mukherjee, Soumendu Sundar}, 
\hyperref[Sen2013]{\textsc{Senadheera 2013}}\index[ppl]{Senadheera, Jayantha}, 
\hyperref[Gar2014b]{\textsc{Garc\'{i}a-Caballero, Moreno, and Prophet 2014b}\index[ppl]{Garc\'{i}a-Caballero, Esther M.}}\index[ppl]{Moreno, Samuel G\'{o}mez}\index[ppl]{Prophet, Michael P.}, 
\hyperref[Lyn2014]{\textsc{Lynd 2014}}\index[ppl]{Lynd, Chris D.}, 
\hyperref[Les2016]{\textsc{Lesher and Lynd 2016}}\index[ppl]{Lesher, Devyn A.}\index[ppl]{Lynd, Chris D.}, and
\hyperref[Vel2016c]{\textsc{Vellucci and Bersani 2016c}}\index[ppl]{Vellucci, Pierluigi}\index[ppl]{Bersani, Alberto Maria}.}

\vspace{4pt}

The authors write: ``[I]s it possible to write any arbitrary integer, rational number, or indeed $\pi$\index{pi@$\pi\;(3.14159\ldots)$!continued square root expressions for}\index{constants, named!pi@$\pi\;(3.14159\ldots)$} or $e$\index{e@$e\;(2.71828\ldots)$}\index{constants, named!e@$e\;(2.71828\ldots)$} as the limit of some sequence of nested radicals? And if an integer $k$ is such a limit, how many different sequences of radicals will converge to $k$? Although there seems to be some revived interest in this topic$\,\ldots\,$previous research has not considered these questions$\ldots$ In this paper we will make a systematic study of nested radicals\index{nested radicals}, answering many such questions and suggesting further lines of research for the interested reader." Eight references are cited. 


\item\vspace{9pt} \label{Zim2008b} (2008b) Seth Zimmerman\index[ppl]{Zimmerman, Seth|textbf} and Chungwu Ho\index[ppl]{Ho, Chungwu|textbf}, Erratum: On Infinitely Nested Radicals. \emph{Mathematics Magazine} \textbf{81}(3), 190. 

{\footnotesize \href{https://mathscinet.ams.org/mathscinet/article?mr=2422950}{MR2422950}} 

{\footnotesize Source: \href{https://www.jstor.org/stable/27643105}{JSTOR}.} 

\vspace{4pt}

The journal editors write: ``Authors Chungwu Ho and Seth Zimmerman have written to point out an error on page 14 of their paper `On Infinitely Nested Radicals,' this \textsc{Magazine}, Vol. 81, February 2008: The gaps mentioned for the set $S_2$ do not exist. Gaps exist only for sets $S_a$ with $a\ge 3$." 


\item\vspace{9pt} \label{Kos2010} (2010) Thomas Koshy\index[ppl]{Koshy, Thomas|textbf}, Generalized nested Pell radical sums. \emph{Bulletin of the Calcutta Mathematical Society} \textbf{102}(1), 37--42. 

{\footnotesize \href{https://mathscinet.ams.org/mathscinet/article?mr=2680839}{MR2680839}} 

\vspace{4pt}

From the abstract: ``This article presents an extended Pell family of polynomial functions\index{Pell!polynomials} $g_n(x)$ which includes the well known Fibonacci\index{polynomials!Fibonacci}\index{Fibonacci!polynomials}, Lucas\index{polynomials!Lucas}\index{Lucas!polynomials}, Pell\index{polynomials!Pell}, and Pell-Lucas\index{polynomials!Pell-Lucas} polynomials $f_n(x)$, $\ell_n(x)$, $p_n(x)$, and $q_n(x)$, respectively. It investigates the convergence of the sequence $\{S_n(x)\}$ of nested radical sums [continued square root approximants], where\index{continued square roots!of terms which are polynomials}
\[S_n(x) = \sqrt{g_1(x)+\sqrt{g_2(x)+\sqrt{g_3(x)+\cdots+\sqrt{g_n(x)}}}}\;;\]
it shows that the sequence converges and $\lim_{n\to\infty}S_n(x)<\lambda_(\gamma)$, where $\gamma=\gamma(x)=x+\sqrt{x^2+1}$, $\lambda=\lambda(x)=(b\gamma+a)/\sqrt{5}$, $a=a(x)= g_0(x)$, and $b=b(x)= g_1(x)$." 


\item\vspace{9pt} \label{Lim2010} (2010) Teik-Cheng Lim\index[ppl]{Lim, Teik-Cheng|textbf}, Continued nested radical fractions. \emph{Mathematical Spectrum} \textbf{42} (2009-2010), 59--63.
\index{continued reciprocal roots}

{\footnotesize Source: \href{https://drive.google.com/file/d/17ixtCnUnzhtWdbm1RdTPaHr1zuZ_mVHH/view}{Google Drive (Mathematical Spectrum)}.}

{\footnotesize Cited in
\hyperref[Lim2011]{\textsc{Lim 2011}}\index[ppl]{Lim, Teik-Cheng}.}

\vspace{4pt} 

From the abstract: ``We define a continued nested radical fraction (CNRF)\index{continued nested radical fractions} as a hybrid of a continued fraction and a nested radical:
\begin{equation*}
\textrm{CNRF}=\sqrt[m]{a_0+\cfrac{b_1}{\sqrt[n]{a_1+\cfrac{b_2}{\sqrt[n]{a_2+\ldots}}}}}\,."
\end{equation*}
The note presents CNRF representations of the golden ratio\index{golden ratio ($\tfrac{1+\sqrt{5}}{2}=1.61803\ldots$)}\index{constants, named!golden ratio ($\tfrac{1+\sqrt{5}}{2}=1.61803\ldots$)}\index{constants, named!golden ratio ($\tfrac{1+\sqrt{5}}{2}=1.61803\ldots$)} using Fibonacci numbers\index{Fibonacci!numbers}\index{numbers!Fibonacci} as terms; of the ``silver ratio"\index{silver ratio\;$(1+\sqrt{2}=2.41421\ldots$)}\index{constants, named!silver ratio\;$(1+\sqrt{2}=2.41421\ldots$)} $1+\sqrt{2}$ using Pell numbers\index{Pell!numbers}\index{numbers!Pell} as terms; and of the ``plastic constant"\index{plastic constant ($1.32471\ldots$)}\index{constants, named!plastic constant ($1.32471\ldots$)} 1.32471957\ldots with all terms equal to 1. The manipulations are purely formal; no rigorous convergence criteria are mentioned. Tables of values are given to support the claim of  the CNRF's superior convergence versus continued fraction\index{continued fractions} and continued radical approximations. (The author's ``CNRF" is a form of the continued reciprocal root\index{continued reciprocal roots} discussed in \hyperref[Jon2015]{\textsc{Jones 2015}}\index[ppl]{Jones, Dixon J.}; see also \hyperref[Gun1880]{\textsc{G\"{u}nther 1880}}\index[ppl]{Gunther@G\"{u}nther, Siegmund}, \hyperref[Lau1990]{\textsc{Laugwitz 1990}}\index[ppl]{Laugwitz, Detlef}, \hyperref[Sch1992]{\textsc{Sch\"{o}nefuss 1992}}, and \hyperref[Lau1999]{\textsc{Laugwitz and Sch\"{o}nefuss 1999}}.) 


\item\vspace{9pt} \label{Oht2010} (2010) Hideyuki Ohtsuka\index[ppl]{Ohtsuka, Hideyuki|textbf}, Problem B-1066. \emph{The Fibonacci Quarterly} \textbf{48}(2), 182. Solution by Kenneth B. Davenport\index[ppl]{Davenport, Kenneth B.}, \emph{ibid.} \textbf{49}(2), 183.

{\footnotesize Solution source: \href{https://www.fq.math.ca/Problems/ElemProbMay2011.pdf}{Fibonacci Quarterly}.}

\vspace{4pt}

The problem states: ``Determine the value of 
\[\sqrt{1+F_2\sqrt{1+F_4\sqrt{1+F_6\sqrt{\cdots\sqrt{1+F_{2n}\sqrt{\cdots}}}}}}\;."\]
Here $F_1=F_2=1$ and $F_n$ is the $n$th Fibonacci number\index{Fibonacci!numbers}\index{numbers!Fibonacci}. The value is shown to be $F_4=3$.\index{continued square roots!of Fibonacci numbers}


\item\vspace{9pt} \label{Sza2010} (2010) P\'{e}ter G\'{a}bor Szab\'{o}\index[ppl]{Szab\'{o}, P\'{e}ter G\'{a}bor|textbf}, On the roots of the trinomial equation. \emph{Central European Journal of Operations Research}, \textbf{18}, 97--104.
\index{trinomial equations}
\index{history!of the ``Bolyai algorithm''}

{\footnotesize \href{https://mathscinet.ams.org/mathscinet/article?mr=2593126}{MR2593126}}

{\footnotesize Source: \href{https://doi.org/10.1007/s10100-009-0130-2}{doi.org (Springer)}.}

\vspace{4pt}

The article summarizes the work of Jan\H{o}s Egerv\'{a}ry\index[ppl]{Egerv\'{a}ry, Jan\H{o}s} (1891--1958) concerning solutions of trinomial equations. Section 2 presents historical background centering on the ``Bolyai algorithm" (\hyperref[Bol1832]{\textsc{Bolyai 1832}}\index[ppl]{Bolyai, Farkas (\emph{aka} Wolfgang)}), which is the method of successive substitution: for $m\in\nats, m>1, a>0$, the equation $x^m=a+x$ is approximately solved by finite continued $m$th roots $\sqrt[m]{a+\sqrt[m]{a+\sqrt[m]{\cdots+\sqrt[m]{a}}}}$. Approximations to $x^m=a+bx$ were given by Gyul\'{a}t\'{o}l Farkas (\hyperref[Far1881]{\textsc{Farkas 1881}}\index[ppl]{Farkas, Gyul\'{a}t\'{o}l (\emph{aka} Gyula or Jules)}), who also investigated generalizations of this method (\hyperref[Far1884]{\textsc{Farkas 1884}}\index[ppl]{Farkas, Gyul\'{a}t\'{o}l (\emph{aka} Gyula or Jules)}).


\item\vspace{9pt} \label{Gil2011} (2011) John Gill\index[ppl]{Gill, John|textbf}, A mathematical note: convergence of infinite compositions of complex functions.
\index{continued compositions!of complex-valued functions}

{\footnotesize Source: \href{https://www.coloradomesa.edu/math-stat/documents/JohnGillResearchnoteInfiniteCompositions2.pdf}{Colorado Mesa University}.}

\vspace{4pt}

From the abstract: ``\emph{Inner Composition}\index{inner compositions} of analytic functions $(f_1\circ f_2\circ\cdots\circ f_n(z))$ and \emph{Outer Composition}\index{outer compositions} of analytic functions $(f_n\circ f_{n-1}\circ \cdots \circ f_1(z))$ are variations on simple iteration, and their convergence behaviors as $n$ becomes infinite may reflect that of simple iterations of \emph{contraction mappings} ($\phi$ defined on a simply-connected domain $S$ with $\phi(S)\subset\Omega\subset S$, $\Omega$ compact). Several theorems are combined to give a summary of work in this area. In addition, recent results by the author and others provide convergence information about such compositions that involve functions that are not contractive, and in some cases, neither analytic nor meromorphic."


\item\vspace{9pt} \label{Lim2011} (2011) Teik-Cheng Lim\index[ppl]{Lim, Teik-Cheng|textbf}, Alternate continued nested radical fractions. \emph{Mathematical Spectrum} \textbf{43} (2010-2011), 55--59. 

{\footnotesize Source: \href{https://drive.google.com/file/d/1sAvjVOkcrsKT2gjpUfKxxXz1-FuHJRko/view}{Google Drive (Mathematical Spectrum)}.}

\vspace{4pt}

Following the conventions established in \hyperref[Lim2010]{\textsc{Lim 2010}}, the author presents expansions of the golden ratio\index{golden ratio ($\tfrac{1+\sqrt{5}}{2}=1.61803\ldots$)}\index{constants, named!golden ratio ($\tfrac{1+\sqrt{5}}{2}=1.61803\ldots$)}, silver ratio\index{silver ratio\;$(1+\sqrt{2}=2.41421\ldots$)}\index{constants, named!silver ratio\;$(1+\sqrt{2}=2.41421\ldots$)}, plastic constant\index{plastic constant ($1.32471\ldots$)}\index{constants, named!plastic constant ($1.32471\ldots$)}, and other real numbers using continued nested radical fractions\index{continued nested radical fractions} of constant or periodic terms. The manipulations are again formal, with no convergence conditions given. 


\item\vspace{9pt} \label{Mor2011} (2011) Samuel G\'{o}mez Moreno\index[ppl]{Moreno, Samuel G\'{o}mez|textbf}, Proposed Mayhem Problem M487. \emph{Crux Mathematicorum with Mathematical Mayhem} \textbf{37}(3), 137. Solution by Florencio Cano Vargas\index[ppl]{Vargas, Florencio Cano}, \emph{ibid.} \textbf{38}(2), 50.

\vspace{4pt}

The problem states: ``Let $m$ be a positive integer. Find all real solutions to the equation
\begin{equation*}
m+\sqrt{m+\sqrt{m+\cdots\sqrt{m+\sqrt{m+\sqrt{x}}}}} = x\;,
\end{equation*}
in which the integer $m$ occurs $n$ times." It is shown that the only solution is $x_0=\tfrac{1}{2}(1+2m+\sqrt{1+4m})$. \index{continued square roots!of constant nonnegative real terms}


%
%



\item\vspace{9pt} \label{Eft2012} (2012) Costas J. Efthimiou\index[ppl]{Efthimiou, Costas J.|textbf}, A class of periodic continued radicals. \emph{The American Mathematical Monthly} \textbf{119}(1), 52--58.
\index{continued square roots!of periodic real terms}

{\footnotesize \href{https://mathscinet.ams.org/mathscinet/article?mr=2877666}{MR2877666}}

{\footnotesize Source: \href{https://doi.org/10.4169/amer.math.monthly.119.01.052}{doi.org (Taylor \& Francis)}.} 

{\footnotesize Cited in 
\hyperref[Eft2013]{\textsc{Efthimiou 2013}}\index[ppl]{Efthimiou, Costas J.}, 
\hyperref[Mor2013b]{\textsc{Moreno and Garc\'{i}a-Caballero 2013b}}\index[ppl]{Moreno, Samuel G\'{o}mez}\index[ppl]{Garc\'{i}a-Caballero, Esther M.}, 
\hyperref[Muk2013]{\textsc{Mukherjee 2013}}\index[ppl]{Mukherjee, Soumendu Sundar}, 
\hyperref[Sen2013]{\textsc{Senadheera 2013}}\index[ppl]{Senadheera, Jayantha}, 
\hyperref[Cla2014]{\textsc{Clark and Richmond 2014}}\index[ppl]{Clark, Tyler}\index[ppl]{Richmond, Tom}, 
\hyperref[Lyn2014]{\textsc{Lynd 2014}}\index[ppl]{Lynd, Chris D.}, and
\hyperref[Vel2016c]{\textsc{Vellucci and Bersani 2016c}}\index[ppl]{Vellucci, Pierluigi}\index[ppl]{Bersani, Alberto Maria}.}

\vspace{4pt}

The author writes: ``In this brief article we find the values for a class of periodic continued radicals of the form
\[a_0\sqrt{2+a_1\sqrt{2+a_2\sqrt{2+a_3\sqrt{2+\cdots}}}}\;,\]
where for some positive integer $n$,
\[a_{n+k}=a_k\;, k=0, 1, 2, \ldots\;,\]
and\index{continued square roots!of periodic real terms}\index{continued square roots!of terms $a_n=\pm2$}
\[a_k\in\{-1, +1\}\;, k=0, 1, 2, \ldots\;, n-1\;."\]
Chebyshev polynomials\index{polynomials!Chebyshev}\index{Chebyshev polynomials} are invoked; compare the virtually simultaneous \hyperref[Mor2012]{\textsc{Moreno and Garc\'{i}a-Caballero 2012}}\index[ppl]{Moreno, Samuel G\'{o}mez} below. For continued square roots of the form studied here, compare \hyperref[Boc1899]{\textsc{Bochow 1899}}\index[ppl]{Bochow, Karl}, \hyperref[Cip1908]{\textsc{Cipolla 1908}}\index[ppl]{Cipolla, Michele}, \hyperref[Pol1925]{\textsc{P\'{o}lya and Szeg\"{o} 1925}}\index[ppl]{Polya@P\'{o}lya, G.}\index[ppl]{Szeg\"{o}, G.}, and \hyperref[Wol1956]{\textsc{Wollan and Mesner 1956}}, \hyperref[Wol1957]{\textsc{1957}}. 


\item\vspace{9pt} \label{Glu2012} (2012) S. Gluzman\index[ppl]{Gluzman, S.|textbf}, V. I. Yukalov\index[ppl]{Yukalov, V. I.|textbf}, Self-similar continued root approximants. \emph{Physics Letters A} \textbf{377}(1-2), 124--128. 
\index{continued robb roots@continued $r$th roots!self-similar}

{\footnotesize \href{https://mathscinet.ams.org/mathscinet/article?mr=2998997}{MR2998997}}

{\footnotesize Source: \href{https://doi.org/10.1016/j.physleta.2012.11.005}{doi.org (ScienceDirect)}.}

\vspace{4pt}

From the introduction: ``There exists a more general approach for extrapolating asymptotic series\index{asymptotic series}\index{series!asymptotic} in powers of a small parameter, or a variable, to finite and even infinite values of such variables. This approach is based on the self-similar approximation theory\index{approximation theory}$\ldots$ In the frame of this theory, we have developed the methods of extrapolating asymptotic series by using several types of self-similar approximants\index{approximants!self-similar}, such as optimized approximants\index{approximants!optimized}, nested exponentials\index{continued exponentials}, nested roots\index{nested square roots}, iterated roots\index{iterated square roots}, and factor approximants$\ldots$

``In the present Letter, we advance a novel type of self-similar approximants that may be called \emph{self-similar continued radical approximants}, or, for short, \emph{self-similar continued radicals}. In a particular case, these continued radicals reduce to continued fractions\index{continued fractions}$\ldots\,$and, respectively, to Pad\'{e} approximants\index{approximants!Pad\'{e}}\index{Pad\'{e} approximants}. But, generally, their form is different and not reduceable [sic] to continued fractions. The self-similar continued radicals could be transformed into expressions of the type of the numerical nested radicals$\ldots\,$, which, however, is not convenient for the extrapolation procedure applied to functions.

``In Section 2, we explain how the self-similar continued radicals arise in the process of the self-similar renormalization\index{renormalization} of asymptotic series and prove the convergence of these root approximants. In Section 3, we demonstrate, by several examples from condensed-matter physics\index{condensed-matter physics}, that the continued radicals [continued $r$th roots]\index{continued robb roots@continued $r$th roots!asymptotic series extrapolation for} can be employed as approximants extrapolating asymptotic series and providing good accuracy. Possible generalizations for the continued radical approximants are also mentioned." 

The authors investigate the continued $p$th power
\[(1+A_1x(1+A_2x(\cdots(1+A_kx)^p\cdots)^p)^p)^p\;,\]
noting that the case $p=-1$ yields continued fractions; ultimately the restriction $|p|<1$ is imposed. The continued $p_i$th power\index{continued pocc powers@continued $p_i$th powers}
\[(1+A_1x(1+A_2x(\cdots(1+A_kx)^{p_k}\cdots)^{p_2})^{p_1})^{p_0}\]
is mentioned in passing.


\item\vspace{9pt} \label{Mor2012} (2012) Samuel G. Moreno\index[ppl]{Moreno, Samuel G\'{o}mez|textbf} and Esther M. Garc\'{i}a-Caballero\index[ppl]{Garc\'{i}a-Caballero, Esther M.|textbf}, Chebyshev polynomials and nested square roots, \emph{Journal of Mathematical Analysis and Applications} \textbf{394}(1), 61--73. 

{\footnotesize \href{https://mathscinet.ams.org/mathscinet/article?mr=2926204}{MR2926204}}

{\footnotesize Source: \href{https://doi.org/10.1016/j.jmaa.2012.04.065}{doi.org (ScienceDirect)}.} 

{\footnotesize Cited in 
\hyperref[Mor2013b]{\textsc{Moreno and Garc\'{i}a-Caballero 2013b}}\index[ppl]{Moreno, Samuel G\'{o}mez}\index[ppl]{Garc\'{i}a-Caballero, Esther M.}, 
\hyperref[Gar2014a]{\textsc{Garc\'{i}a-Caballero, Moreno, and Prophet 2014a}} and \hyperref[Gar2014c]{\textsc{2014c}}\index[ppl]{Garc\'{i}a-Caballero, Esther M.}\index[ppl]{Moreno, Samuel G\'{o}mez}\index[ppl]{Prophet, Michael P.}, and
\hyperref[Vel2016c]{\textsc{Vellucci and Bersani 2016c}}\index[ppl]{Vellucci, Pierluigi}\index[ppl]{Bersani, Alberto Maria}.}

\vspace{4pt}

From the abstract: ``The purpose of this note is to report a curious relation between the Chebyshev polynomials\index{Chebyshev polynomials}\index{polynomials!Chebyshev} of degree $2n$ in a complex variable and the nested square roots of depth $n$ of the form
\[\pm\sqrt{2\pm\sqrt{2\pm\ldots\pm\sqrt{2+2\xi}}}\;,\]
$\xi$ being a complex number. \index{continued square roots!of terms $a_n=\pm2$} Our approach leads us to generalize and recover, in a unified manner, the closed-form expressions recently given [in \hyperref[Ser2003]{\textsc{Servi 2003}}], corresponding to the case $\xi \in [-1, 1]$, and [in \hyperref[Nyb2005]{\textsc{Nyblom 2005}}], corresponding to $\xi \in [1,\infty)$."


\item\vspace{9pt} \label{Nyb2012} (2012) M. A. Nyblom\index[ppl]{Nyblom, M. A.|textbf}, Some closed-form evaluations of infinite products involving nested radicals. \emph{Rocky Mountain Journal of Mathematics} \textbf{42}(2), 751--758. 
\index{nested radicals!infinite products of}

{\footnotesize \href{https://mathscinet.ams.org/mathscinet/article?mr=2915517}{MR2915517}}

{\footnotesize Source: \href{https://www.jstor.org/stable/44240073}{JSTOR}.} 

{\footnotesize Cited in 
\hyperref[Mor2013b]{\textsc{Moreno and Garc\'{i}a-Caballero 2013b}}\index[ppl]{Moreno, Samuel G\'{o}mez}\index[ppl]{Garc\'{i}a-Caballero, Esther M.}.}

\vspace{4pt}

From the abstract: ``By applying double and triple angle identities for hyperbolic and trigonometric cosine functions, we obtain closed-form evaluations for two families of infinite products involving nested radicals\index{nested radicals!infinite products of}. The first group of results represents a generalization of the classic Vi\`{e}te\index[ppl]{Vi\`{e}te, Fran\c{c}ois} infinite product expansion for $2/\pi$\index{Vi\`{e}te's formula for $\tfrac{2}{\pi}$}, while the second comprises variations on Vi\`{e}te type infinite products and infinite products involving nested square roots of $2$\index{nested square roots!infinite product of}. In addition, specific examples of Vi\`{e}te type infinite product expansions are presented for such numbers as $\tfrac{3\sqrt{3}}{2\pi}$ and $\tfrac{3}{\pi}$." 


\item\vspace{9pt} \label{Eft2013} (2013) Costas J. Efthimiou\index[ppl]{Efthimiou, Costas J.|textbf}, A class of continued radicals. \emph{The American Mathematical Monthly} \textbf{120}(5), 459--461.
\index{continued radicals}

{\footnotesize \href{https://mathscinet.ams.org/mathscinet/article?mr=3035446}{MR3035446}}

{\footnotesize Source: \href{https://doi.org/10.4169/amer.math.monthly.120.05.459}{doi.org (Taylor \& Francis)}.} 

{\footnotesize Cited in 
\hyperref[Muk2013]{\textsc{Mukherjee 2013}}\index[ppl]{Mukherjee, Soumendu Sundar}, 
\hyperref[Sen2013]{\textsc{Senadheera 2013}}, and
\hyperref[Vel2016c]{\textsc{Vellucci and Bersani 2016c}}\index[ppl]{Vellucci, Pierluigi}\index[ppl]{Bersani, Alberto Maria}.}

\vspace{4pt}

From the introduction: ``In \hyperref[Eft2012]{\textsc{Efthimiou 2012}} the author discussed the values for a class of periodic continued radicals of the form
\begin{equation}\label{E:Eft2013-1}
a_0\sqrt{2+a_1\sqrt{2+a_2\sqrt{2+a_3\sqrt{2+\cdots}}}}\;,
\end{equation}
where for some positive integer $n$, $a_{n+k} = a_k\;, k=0,1,2,\ldots$, $a_k\in\{-1, +1\}$, $k=0,1,\ldots, n-1$. It was also shown that the radicals given by equation \eqref{E:Eft2013-1} have limits two times the fixed points of the Chebycheff polynomials\index{Chebyshev polynomials}\index{polynomials!Chebyshev} $T_{2n}(x)$, thus unveiling an interesting relation between these topics.

``In \hyperref[Zim2008a]{\textsc{Zimmerman and Ho 2008a}}\index[ppl]{Zimmerman, Seth}\index[ppl]{Ho, Chungwu}, the authors defined the set $S_2$ of all continued radicals of the form \eqref{E:Eft2013-1} (with $a_0=1$), and they investigated some of their properties by assuming that the limit of the radicals exists. In particular, they showed that all elements of $S_2$ lie between $0$ and $2$, any two radicals cannot be equal to each other, and $S_2$ is uncountable.

``[M]y previous note partially bridged this gap, but left unanswered the question, \emph{`what are the limits if the radicals are not periodic?'} I answer the question in this note." 


\item\vspace{9pt} \label{Mor2013a} (2013a) Samuel G. Moreno\index[ppl]{Moreno, Samuel G\'{o}mez|textbf} and Esther M. Garc\'{i}a\index[ppl]{Garc\'{i}a-Caballero, Esther M.|textbf}, New infinite products of cosines and Vi\`{e}te\index[ppl]{Vi\`{e}te, Fran\c{c}ois}-like formulae. \emph{Mathematics Magazine} \textbf{86}(1), 15--25.
\index{Vi\`{e}te's formula for $\tfrac{2}{\pi}$!generalizations}

{\footnotesize \href{https://mathscinet.ams.org/mathscinet/article?mr=4878550}{MR4878550 (citation only)}}

{\footnotesize Source: \href{https://doi.org/10.4169/math.mag.86.1.015}{doi.org (Taylor \& Francis)}.}

{\footnotesize Cited in 
\hyperref[Mor2013b]{\textsc{Moreno and Garc\'{i}a-Caballero 2013b}}\index[ppl]{Moreno, Samuel G\'{o}mez}\index[ppl]{Garc\'{i}a-Caballero, Esther M.}, and
\hyperref[Gar2014a]{\textsc{Garc\'{i}a-Caballero, Moreno, and Prophet 2014a}}\index[ppl]{Garc\'{i}a-Caballero, Esther M.}\index[ppl]{Moreno, Samuel G\'{o}mez}\index[ppl]{Prophet, Michael P.} and \hyperref[Gar2014c]{\textsc{2014c}}.}

\vspace{4pt}

The authors write: ``In this note we evaluate infinite products similar to the product [\eqref{E:Vie1593-1} above\index{Vi\`{e}te's formula for $\tfrac{2}{\pi}$!generalizations}, in \hyperref[Vie1593]{\textsc{Vi\`{e}te 1593}}]\index[ppl]{Vi\`{e}te, Fran\c{c}ois}, with the novelty that some of the plus signs are replaced by minus signs. We call these products \emph{Vi\`{e}te-like formulae}. To derive them, first we manipulate a simple trigonometric relation (other than the double angle formula) in order to obtain a family of infinite products of cosines\index{infinite products!of cosine functions}. Next, with the aid of a formula of \hyperref[Ser2003]{\textsc{Servi 2003}} (see also \hyperref[Mor2012]{\textsc{Moreno and Garc\'{i}a-Caballero 2012}}\index[ppl]{Garc\'{i}a-Caballero, Esther M.} for a generalization), we transform these infinite products of cosines into infinite products of nested square roots of 2\index{infinite products!of continued square roots}. These Vi\`{e}te\index[ppl]{Vi\`{e}te, Fran\c{c}ois}-like expressions turn out to represent numbers like $\pi$, $\sqrt{3}$, and $\sqrt{5-\sqrt{5}}$."


\item\vspace{9pt} \label{Mor2013b} (2013b) Samuel G. Moreno\index[ppl]{Moreno, Samuel G\'{o}mez|textbf} and Esther M. Garc\'{i}a-Caballero\index[ppl]{Garc\'{i}a-Caballero, Esther M.|textbf}, On Vi\`{e}te\index[ppl]{Vi\`{e}te, Fran\c{c}ois}-like formulas. \emph{Journal of Approximation Theory} \textbf{174}, 90--112. 
\index{nested square roots!infinite product of}
\index{Vi\`{e}te's formula for $\tfrac{2}{\pi}$!generalizations}

{\footnotesize \href{https://mathscinet.ams.org/mathscinet/article?mr=3090772}{MR3090772}}

{\footnotesize Source: \href{https://doi.org/10.1016/j.jat.2013.06.006}{doi.org (ScienceDirect)}.}

{\footnotesize Cited in 
\hyperref[Mor2013a]{\textsc{Moreno and Garc\'{i}a-Caballero 2013a}}\index[ppl]{Moreno, Samuel G\'{o}mez}\index[ppl]{Garc\'{i}a-Caballero, Esther M.}, 
\hyperref[Gar2014c]{\textsc{Garc\'{i}a-Caballero, Moreno, and Prophet 2014c}}\index[ppl]{Garc\'{i}a-Caballero, Esther M.}\index[ppl]{Moreno, Samuel G\'{o}mez}\index[ppl]{Prophet, Michael P.}, and
\hyperref[Vel2016c]{\textsc{Vellucci and Bersani 2016c}}\index[ppl]{Vellucci, Pierluigi}\index[ppl]{Bersani, Alberto Maria}.}

\vspace{4pt}

From the abstract: ``The very first recorded use of an infinite product in mathematics is the so-called Vi\`{e}te's\index[ppl]{Vi\`{e}te, Fran\c{c}ois} formula\index{Vi\`{e}te's formula for $\tfrac{2}{\pi}$}, in which each of its factors contains nested square roots of 2 with plus signs inside [equation \eqref{E:Vie1593-1} above] $\ldots$ and it can be proved by iterating the double angle formula $\sin 2x=2\cos x\sin x$, thus obtaining the infinite product $\tfrac{2}{\pi}=\prod_{n=2}^\infty\cos(\pi/2^n)$.\index{infinite products!of cosine functions}

``This paper focuses, first, on the wide variety of iterations that the identity $\cos x=2\cos((\pi+2x)/4)\cos((\pi-2x)/4)$ admits; next, on the infinite products of cosines derived from these iterations and finally, on how these infinite products of cosines give rise to striking formulas."\index{Vi\`{e}te's formula for $\tfrac{2}{\pi}$!generalizations}

The main theorem proved is the following: Let $N$ be a positive integer and let $(b_1,b_2,\ldots,b_N)$ be a sequence for which $b_i$ equals either $1$ or $-1$ and at least one $b_i$ equals $-1$. Define
\[\alpha=\sum_{i=1}^N\left(\prod_{j=1}^{i-1}b_{N-j+1}\right)2^{N-i}\;,\quad\textrm{where} \prod_{j=1}^0=1\;,\]

\[\sigma=\prod_{j=1}^Nb_{N-j+1}\;,\]

\[\widehat{b}_k=b_{1+(k-1)\bmod N}\;,\quad k=1,2,\ldots\,.\]
For each complex number $z$, define the sequence $\{\rho_i(z)\}$ by
\[\rho_0(z)=z\;,\quad \rho_i(z)=\sqrt{2+\widehat{b}_i\rho_{i-1}}\;,\quad i\ge1\;.\]
Then
\[\dfrac{\sqrt{1-\tfrac{z^2}{4}}}{\cos\left(\tfrac{\alpha\sigma}{2^N\sigma-1}\tfrac{\pi}{2}\right)}=\prod_{j=0}^\infty\left(\prod_{i=jN+1}^{(j+1)N}\rho_i(z)\right)\;.\]

In an appendix, the authors give a lucid explanation of the method used by Fran\c{c}ois Vi\`{e}te\index[ppl]{Vi\`{e}te, Fran\c{c}ois} to obtain his original formula \eqref{E:Vie1593-1} above.


\item\vspace{9pt} \label{Muk2013} (2013) Soumendu Sundar Mukherjee\index[ppl]{Mukherjee, Soumendu Sundar|textbf}, An approximation inequality for continued radicals and power forms. 

{\footnotesize Source: \href{http://arxiv.org/abs/1303.4251}{arxiv.org}.}
\index{continued radicals}
\index{continued rocc roots@continued $r_i$th roots!of arbitrary nonnegative real terms}
\index{continued powers}

\vspace{4pt}

From the abstract: ``In this article we derive an approximation inequality for continued radicals, generalizing an inequality of Herschfeld\index[ppl]{Herschfeld, Aaron} for continued square roots to arbitrary radicals [i.e., continued $r_i$th roots, $r_i\in(1,\infty)$]\index{continued square roots}\index{continued rocc roots@continued $r_i$th roots!convergence conditions for}, which is useful in exploring convergence issues and obtaining convergence rates. In fact, we generalize this inequality further to encompass the more general continued power forms.\index{continued powers} We demonstrate the use of this inequality by obtaining estimates for the convergence rates of several continued radicals including the famous Ramanujan radical [equation \eqref{E:Ram1911-1} above]."\index{continued square roots!Ramanujan's}


\item\vspace{9pt} \label{Pau2013} (2013) William Paulsen\index[ppl]{Paulsen, William|textbf}, \emph{Asymptotic analysis and perturbation theory}. CRC Press, Hoboken, ISBN 978-1466515116.

{\footnotesize \href{https://mathscinet.ams.org/mathscinet/article?mr=3185159}{MR3185159}}

\vspace{4pt}

Following Problem 18 in Section 3.5 (pp. 153--154), the author introduces a ``continued-function"\index{continued functions}\index{continued compositions!of arbitrary real-valued functions} representation of a function $f(x)$ as follows: Let $h(x)$ be a simple function such that $h(0)=1$. Then inductively define
\begin{align*}
f_0(x)&=f(x)\\
g_n(x)&=\textrm{the leading term of $f_n(x)$}\\
f_{n+1}(x)&=h^{-1}\left(\dfrac{f_n(x)}{g_n(x)}\right)\;.
\end{align*}
Then
\[f(x)\sim g_0(x)\cdot h\bigl(g_1(x)\cdot h\bigl(g_2(x)\cdot h\bigl(\cdots\;,\]
where the relation $\sim$ (read ``is similar to") means that the ratio of $f$ and its expansion approaches $1$ as $n$ increases without bound. Problems 19--22 ask the reader to compute continued square root\index{continued square roots!of arbitrary nonnegative real terms} expansions (using $h(x)=\sqrt{1+x}$) and continued exponential\index{continued exponentials} expansions ($h(x)=e^x$) for functions similar to divergent power series\index{power series}\index{series!power}.


\item\vspace{9pt} \label{Sen2013} (2013) Jayantha Senadheera\index[ppl]{Senadheera, Jayantha|textbf}, On the periodic continued radicals of 2 and generalization for Vieta product. 
\index{Vi\`{e}te's formula for $\tfrac{2}{\pi}$!generalizations}
\index{continued square roots!of terms $a_n=\pm2$}
\index{continued square roots!and trigonometric functions}
\index{continued radicals}

{\footnotesize Source: \href{http://arxiv.org/abs/1304.5659}{arxiv.org}.}

\vspace{4pt}

From the abstract: ``In this paper we study periodic continued radicals of 2\index{continued square roots!of terms $a_n=\pm2$}\index{continued square roots!and trigonometric functions}. We show that any periodic continued radicals of 2 converg[e] to $2\sin(q\pi)$, for some rational number $q$ depend[ing] on the continued radical. Furthermore we show that if $r_n$ is a periodic nested radicals [sic] of 2, which has $n$ nested roots, then the limit points of the sequence $2^n(2\sin(q\pi)-r_n)$ have the form $\alpha\pi$, where $\alpha$ is an algebraic number. This result give[s] a set of sub sequences converg[ing] to $\alpha\pi$, for each $\alpha$. Also we show that limit of these sub sequences can be represented as Vieta\index[ppl]{Vi\`{e}te, Fran\c{c}ois} like nested radical products\index{nested radicals!infinite products of}.\index{Vi\`{e}te's formula for $\tfrac{2}{\pi}$!generalizations} Hence this result generalizes the Vieta product for $\pi$. Several interesting examples are illustrated." 


\item\vspace{9pt} \label{Tiw2013} (2013) Gaurav Tiwari\index[ppl]{Tiwari, Gaurav|textbf}, Complete elementary analysis of nested radicals. 
\index{continued square roots!Ramanujan's}
\index{nested radicals}

{\footnotesize Source: \href{http://gauravtiwari.org/2013/10/08/projects-ramanujan-nested-radical/}{gauravtiwari.org}.}

\vspace{4pt}

A web site, initiated in 2011, by a freelance web designer and blogger to ``collect and expand what Ramanujan did with Nested Radicals and summariz[e] all important facts [in] one article." \index[ppl]{Ramanujan, Srinivasa}\index{continued square roots!Ramanujan's}


\item\vspace{9pt} \label{Xi2013} (2013) Bo-Yan Xi\index[ppl]{Xi, Bo-Yan|textbf}, Feng Qi\index[ppl]{Qi, Feng|textbf}, Convergence, monotonicity, and inequalities of sequences involving continued powers. \emph{Analysis (Berlin)} \textbf{33}(3), 235--242.
\index{continued robb roots@continued $r$th roots!of constant nonnegative real terms}
\index{continued pobb powers@continued $p$th powers}

{\footnotesize \href{https://mathscinet.ams.org/mathscinet/article?mr=3118425}{MR3118425}}

{\footnotesize Source: \href{https://doi.org/10.1524/anly.2013.1191}{doi.org (De Gruyter)}.}

\vspace{4pt}

From the summary: ``$\ldots$[T]he convergence, monotonicity, and inequalities of sequences involving continued powers $\sqrt[t]{a-\sqrt[t]{a-\cdots-\sqrt[t]{a}}}$ are investigated and established." The authors consider the cases $a\in(1,\infty)$ and $t\in(1,\infty)$; $a\in(0, 1)$ and $t\in(0,1)$; $a\in(-\infty,-1)$ and $t$ odd; and $a\in(-1,0)$ and $\tfrac{1}{t}$ odd.


\item\vspace{9pt} \label{Bag2014} (2014) Nikos Bagis\index[ppl]{Bagis, Nikos|textbf}, Solution of polynomial equations with nested radicals.
\index{continued robb roots@continued $r$th roots!of periodic real terms}
\index{nested radicals}

{\footnotesize Source: \href{http://arxiv.org/abs/1406.1948}{arxiv.org}.}

From the abstract: ``In this note we present solutions of arbitrary polynomial equations in nested periodic radicals." To obtain expressions that resemble continued $r$th roots\index{continued robb roots@continued $r$th roots!as solutions to polynomial equations}, the author uses the trick of setting $\sqrt[d]{x}:= x^{1/d}$, where $d$ is a non-zero rational number (used also in \hyperref[Ise1897]{\textsc{Isenkrahe 1897}}).\index{continued robb roots@continued $r$th roots!of periodic real terms}\index{continued robb roots@continued $r$th roots!as solutions to polynomial equations} Compare \hyperref[And1985]{\textsc{Andrushkiw 1985}}\index[ppl]{Andrushkiw, R. L.}.


\item\vspace{9pt} \label{Cla2014} (2014) Tyler Clark\index[ppl]{Clark, Tyler|textbf} and Tom Richmond\index[ppl]{Richmond, Tom|textbf}, Cantor sets arising from continued radicals, \emph{The Ramanujan Journal} \textbf{33}(3), 315--327. 
\index{Cantor set}
\index{continued radicals}

{\footnotesize \href{https://mathscinet.ams.org/mathscinet/article?mr=3182536}{MR3182536}}

{\footnotesize Source: \href{https://doi.org/10.1007/s11139-012-9457-8}{doi.org (Springer)}.}

\vspace{4pt}

From the abstract: `` If $a_1$, $a_2$, $a_3, \ldots$ are nonnegative real numbers and $f_j(x) = \sqrt{a_j + x}$, then $\lim_{n\to\infty} f_1\circ f_2\circ\cdots\circ f_n(0)$ is a continued radical\index{continued radical} [continued square root]\index{continued square roots!of arbitrary nonnegative real terms} with terms $a_1$, $a_2$, $a_3, \ldots$ The set of real numbers representable as a continued radical whose terms $a_i$ are all from a set $S = \{a,b\}$ of two natural numbers is a Cantor set. We investigate the thickness, measure, and sums of such Cantor sets." 


\item\vspace{9pt} \label{Gar2014a} (2014a) Esther M. Garc\'{i}a-Caballero\index[ppl]{Garc\'{i}a-Caballero, Esther M.|textbf}, Samuel G. Moreno\index[ppl]{Moreno, Samuel G\'{o}mez|textbf}, and Michael P. Prophet\index[ppl]{Prophet, Michael P.|textbf}, New Vi\`{e}te\index[ppl]{Vi\`{e}te, Fran\c{c}ois}-like infinite products of nested radicals with Fibonacci\index{Fibonacci!numbers}\index{numbers!Fibonacci} and Lucas numbers\index{Lucas!numbers}\index{numbers!Lucas}. \emph{The Fibonacci Quarterly} \textbf{52}, 27--31. 

{\footnotesize \href{https://mathscinet.ams.org/mathscinet/article?mr=3181093}{MR3181093}}

{\footnotesize Source: \href{https://www.fq.math.ca/Papers1/52-1/CaballeroMorenoProphet.pdf}{Fibonacci Quarterly}.}

{\footnotesize Cited in
\hyperref[Mor2013b]{\textsc{Moreno and Garc\'{i}a-Caballero 2013b}}\index[ppl]{Moreno, Samuel G\'{o}mez}\index[ppl]{Garc\'{i}a-Caballero, Esther M.}, and
\hyperref[Gar2014c]{\textsc{Garc\'{i}a-Caballero, Moreno, and Prophet 2014c}}\index[ppl]{Garc\'{i}a-Caballero, Esther M.}\index[ppl]{Moreno, Samuel G\'{o}mez}\index[ppl]{Prophet, Michael P.}.}
\index{Vi\`{e}te's formula for $\tfrac{2}{\pi}$!generalizations}
\index{nested radicals!infinite products of}

\vspace{4pt}

From the abstract: ``In a 2007 contribution by Osler\index[ppl]{Osler, Thomas J.} in this \emph{Quarterly}, the so-named Vi\`{e}te\index[ppl]{Vi\`{e}te, Fran\c{c}ois}-like products were introduced as two eye-catching formulas representing either the $n$th Fibonacci number\index{Fibonacci!numbers}\index{numbers!Fibonacci} in terms of a product of nested radicals with the $n$th Lucas number\index{Lucas!numbers}\index{numbers!Lucas} inside, or vice-versa. As [in] the original and famous Vi\`{e}te's\index[ppl]{Vi\`{e}te, Fran\c{c}ois} infinite product, Osler's infinite products\index{infinite products!of continued square roots} have plus signs inside the nested radicals\index{nested radicals!infinite products of}. In this paper we explore infinite products of nested square roots with Fibonacci and Lucas numbers with the novelty that inside the radical symbols there are minus signs instead of plus signs." 


\item\vspace{9pt} \label{Gar2014b} (2014b) Esther M. Garc\'{i}a-Caballero\index[ppl]{Garc\'{i}a-Caballero, Esther M.|textbf}, Samuel G. Moreno\index[ppl]{Moreno, Samuel G\'{o}mez|textbf}, and Michael P. Prophet\index[ppl]{Prophet, Michael P.|textbf}, The Golden Ratio and Vi\`{e}te's\index[ppl]{Vi\`{e}te, Fran\c{c}ois} Formula. \emph{Teaching Mathematics and Computer Science} \textbf{12}, 43--54. 
\index{Vi\`{e}te's formula for $\tfrac{2}{\pi}$!generalizations}
\index{golden ratio ($\tfrac{1+\sqrt{5}}{2}=1.61803\ldots$)}\index{constants, named!golden ratio ($\tfrac{1+\sqrt{5}}{2}=1.61803\ldots$)}

{\footnotesize Source: \href{https://doi.org/10.5485/TMCS.2014.0351}{doi.org (University of Debrecen)}.}

\vspace{4pt}

From the abstract: ``Vi\`{e}te's\index[ppl]{Vi\`{e}te, Fran\c{c}ois} formula uses an infinite product to express $\pi$. In this paper we find a strikingly similar representation for the Golden Ratio."


\item\vspace{9pt} \label{Gar2014c} (2014c) Esther M. Garc\'{i}a-Caballero\index[ppl]{Garc\'{i}a-Caballero, Esther M.|textbf}, Samuel G. Moreno\index[ppl]{Moreno, Samuel G\'{o}mez|textbf}, and Michael P. Prophet\index[ppl]{Prophet, Michael P.|textbf}, A complete view of Vi\`{e}te\index[ppl]{Vi\`{e}te, Fran\c{c}ois}-like infinite products with Fibonacci and Lucas numbers. \emph{Applied Mathematics and Computation} \textbf{247}, 703--711.
\index{Vi\`{e}te's formula for $\tfrac{2}{\pi}$!generalizations}

{\footnotesize \href{https://mathscinet.ams.org/mathscinet/article?mr=3270876}{MR3270876 (citation only)}}

{\footnotesize Source: \href{https://doi.org/10.1016/j.amc.2014.09.018}{doi.org (ScienceDirect)}.} 

\vspace{4pt}

From the abstract: ``The main goal of this paper is to link the $n$th Fibonacci\index{Fibonacci!numbers}\index{numbers!Fibonacci} and Lucas\index{Lucas!numbers}\index{numbers!Lucas} numbers through certain infinite products of nested radicals\index{nested radicals!infinite products of}. This work relies on recent results on Vi\`{e}te\index[ppl]{Vi\`{e}te, Fran\c{c}ois}-like infinite products [which] appeared in \hyperref[Mor2013a]{\textsc{Moreno and Garc\'{i}a-Caballero 2013}}. We will analyze in detail one particular case of these formulas and we will show how our treatment covers and extends previous results in the literature."


\item\vspace{9pt} \label{Kef2014} (2014) Kyriakos Kefalas\index[ppl]{Kefalas, Kyriakos|textbf}, On smooth solutions of non linear dynamical systems, $f_{n+1} = u(f_n)$, part I. \emph{Physics International} \textbf{5}(1), 112--127. 

{\footnotesize Source: \href{https://doi.org/10.3844/pisp.2014.112.127}{doi.org (Science Publications)}.}

\vspace{4pt}

From the abstract: ``We consider the dynamical system\index{dynamical systems}, $f_{n+1} = u(f_n)$, $(1)$ (where usually $n$ is time) defined by a continuous map $u$. Our target is to find a flow of the system for each initial state $f_0$, i.e., we seek continuous solutions of $(1)$, with the same smoothness degree as $u$. We start with the introduction of continued forms which are a generalization of continued fractions\index{continued fractions!generalizations of}. With the use of continued forms and a modulator function (i.e., weight function) $m$, we construct a sequence of smooth functions, which come arbitrarily close to a smooth flow of $(1)$. The limit of this sequence is a functional transform\index{functional transform}, $K_m[u]$, of $u$, with respect to $m$. The functional transform is a solution of $(1)$, in the sense that $K_m[u](y + c)$, is a flow of $(1)$ for each translation constant $c$. Here we present the first part of our work where we consider a subclass of dissipative dynamical systems in the sence [sic] that they have wandering sets of positive measure. In particular we consider strictly increasing real univariate maps, $u: D\to D$, $D = (a+\infty)$, where $a\ge 0$, or $a=-\infty$, with the property $u(x)-x\ge \epsilon>0$, which implies that $u$ has no real fixed points. We briefly give some mathematical and physical applications and we discuss some open problems. We demonstrate the method on the simple non-linear dynamical system $f_{n+1} = u(f_n)+1.$" Ten references are cited.


\item\vspace{9pt} \label{Lyn2014} (2014) Chris D. Lynd\index[ppl]{Lynd, Chris D.|textbf}, Using difference equations to generalize results for periodic nested radicals. \emph{The American Mathematical Monthly} \textbf{121}(1), 45--59. 
\index{continued rocc roots@continued $r_i$th roots}
\index{continued rocc roots@continued $r_i$th roots!of periodic real terms}
\index{nested radicals!periodic}

{\footnotesize \href{https://mathscinet.ams.org/mathscinet/article?mr=3139581}{MR3139581}}

{\footnotesize Source: \href{https://doi.org/10.4169/amer.math.monthly.121.01.045}{doi.org (Taylor \& Francis)}.}

{\footnotesize Cited in
\hyperref[Les2016]{\textsc{Lesher and Lynd 2016}}\index[ppl]{Lesher, Devyn A.}\index[ppl]{Lynd, Chris D.}, and
\hyperref[Vel2016c]{\textsc{Vellucci and Bersani 2016c}}\index[ppl]{Vellucci, Pierluigi}\index[ppl]{Bersani, Alberto Maria}.}

\vspace{4pt}

From the abstract: ``We investigate sequences of nested radicals [meaning here continued $r_i$th roots]\index{continued rocc roots@continued $r_i$th roots} where the indices, the coefficients, and the radicands are periodic sequences of real numbers. We show that one can determine the end behavior of a periodic nested radical by analyzing the basin of attraction\index{basin of attraction} of each equilibrium point\index{equilibrium point}, and each period-2 point, of the corresponding difference equation\index{difference equations}. Using this method of analysis, we prove a few theorems about the end behavior of nested radicals of this form. These theorems extend previous results on this topic because they apply to large classes of nested radicals\index{nested radicals} that contain arbitrary indices, negative radicands, and periodic parameters with arbitrary periods. In addition, we demonstrate how to construct a periodic nested radical, of a general form, that converges to a predetermined limit\index{continued rocc roots@continued $r_i$th roots!convergence to a predetermined limit}; and we demonstrate how to construct a nested radical\index{nested radicals} that converges asymptotically to a periodic sequence." 


\item\vspace{9pt} \label{Osl2014} (2014) Thomas J. Osler\index[ppl]{Osler, Thomas J.|textbf} and Sky Waterpeace\index[ppl]{Waterpeace, Sky|textbf}, Vieta's\index[ppl]{Vi\`{e}te, Fran\c{c}ois} product for pi from the Archimedean algorithm. \emph{The Mathematical Gazette} \textbf{98}, 429--431.
\index{Archimedean algorithm}
\index{pi@$\pi\;(3.14159\ldots)$!Archimedean algorithm for}\index{constants, named!pi@$\pi\;(3.14159\ldots)$}

{\footnotesize Source: \href{https://www.jstor.org/stable/24496528}{JSTOR}.}

{\footnotesize Cited in 
\hyperref[Osl2016a]{\textsc{Osler 2016}}\index[ppl]{Osler, Thomas J.}.}

\vspace{4pt}

From the abstract: ``In this paper we show how to derive [Vi\`{e}te's]\index[ppl]{Vi\`{e}te, Fran\c{c}ois} famous product of nested radicals\index{nested radicals!infinite products of} for $\pi$ from the Archimedean iterative algorithm for $\pi$. Only simple algebraic manipulations are needed."


\item\vspace{9pt} \label{Jon2015} (2015) Dixon J. Jones\index[ppl]{Jones, Dixon J.|textbf}, Continued reciprocal roots. \emph{The Ramanujan Journal} \textbf{38}, 435--454.
\index{continued reciprocal roots}

{\footnotesize \href{https://mathscinet.ams.org/mathscinet/article?mr=3414500}{MR3414500}} 

{\footnotesize Source: \href{https://doi.org/10.1007/s11139-014-9594-3}{doi.org (Springer)}.} 

{\footnotesize Cited in
\hyperref[Vel2016c]{\textsc{Vellucci and Bersani 2016c}}\index[ppl]{Vellucci, Pierluigi}\index[ppl]{Bersani, Alberto Maria}.}

\vspace{4pt} 

Proves that the continued reciprocal $r$th root
\begin{equation}\label{E:crsr}
a_0+\cfrac{1}{\sqrt[r]{a_1+\cfrac{1}{\sqrt[r]{a_2+\cfrac{1}{\sqrt[r]{\ddots}}}}}}
\end{equation}
diverges if, and only if,
\[\limsup_{i\to\infty}\;a_i^{p^{i}}<1\;,\]
where $r=1/p$, $0<p<1$, and $a_i>0$, $i=0,1,2,\dotsc$.


\item\vspace{9pt} \label{Neu2015} (2015) J\"{o}rg Neunh\"{a}userer\index[ppl]{Neunh\"{a}userer, J\"{o}rg|textbf}, Continued logarithm representation of real numbers. 
\index{continued logarithms}

{\footnotesize \href{https://mathscinet.ams.org/mathscinet/article?mr=3816432}{$\langle${MR3816432}$\rangle$}}

{\footnotesize Source: \href{https://doi.org/10.13140/RG.2.1.1763.1441}{doi.org (ResearchGate)}.}

From the abstract: ``We introduce the continued logarithm representation of real numbers and prove results on the occurrence and frequency of digits with respect to this representation."


\item\vspace{9pt} \label{Nis2015a} (2015a) Ryo Nishimura\index[ppl]{Nishimura, Ryo|textbf}, New properties of the lemniscate function and its transformation. \emph{Journal of Mathematical Analysis and Applications} \textbf{427}, 460--468.
\index{Vi\`{e}te's formula for $\tfrac{2}{\pi}$!generalizations}
\index[ppl]{Vi\`{e}te, Fran\c{c}ois}

{\footnotesize \href{https://mathscinet.ams.org/mathscinet/article?mr=3318209}{MR3318209}}

{\footnotesize Source: \href{https://doi.org/10.1016/j.jmaa.2015.02.066}{doi.org (ScienceDirect)}.}

{\footnotesize Cited in
\hyperref[Nis2016]{\textsc{Nishimura 2016}}\index[ppl]{Nishimura, Ryo} and
\hyperref[Osl2016b]{\textsc{Osler, Jacob, and Nishimura 2016}}\index[ppl]{Osler, Thomas J.}\index[ppl]{Jacob, Walter}\index[ppl]{Nishimura, Ryo}.}

\vspace{4pt}

The paper's main theorem is the following generalization of Vi\`{e}te's formula \eqref{E:Vie1593-1} above for $\tfrac{2}{\pi}$:
\begin{equation*}
\dfrac{\textrm{sl}x}{x}=\sqrt{\dfrac{1+\textrm{sl}'x}{2}}\cdot
\sqrt{\dfrac{1+\sqrt{\frac{2}{1+\textrm{sl}'x}}}{2}}\cdot
\sqrt{\dfrac{1+\sqrt{\frac{2}{1+\sqrt{\frac{2}{1+\textrm{sl}'x}}}}}{2}}\ldots
\end{equation*}
where sl$\,x$ is the lemniscate sine\index{lemniscate!sine}.


\item\vspace{9pt} \label{Nis2015b} (2015b) Ryo Nishimura\index[ppl]{Nishimura, Ryo|textbf}, New inequalities and infinite product formulas for the trigonometric and the lemniscate functions. \emph{Mathematical Inequalities and Applications} \textbf{18}(2), 529--540.
\index{Vi\`{e}te's formula for $\tfrac{2}{\pi}$!generalizations}
\index{infinite products}

{\footnotesize \href{https://mathscinet.ams.org/mathscinet/article?mr=3338872}{MR3338872}}

{\footnotesize Source: \href{http://dx.doi.org/10.7153/mia-18-39}{dx.doi.org (Ele-Math)}.}

{\footnotesize Cited in
\hyperref[Osl2016b]{\textsc{Osler, Jacob, and Nishimura 2016}}\index[ppl]{Osler, Thomas J.}\index[ppl]{Jacob, Walter}\index[ppl]{Nishimura, Ryo},
\hyperref[Nis2016]{\textsc{Nishimura 2016}}\index[ppl]{Nishimura, Ryo}, and
\hyperref[Osl2016c]{\textsc{Osler and Kosior 2016c}}\index[ppl]{Osler, Thomas J.}\index[ppl]{Kosior, Jesse M.}.}

\vspace{4pt}

From the abstract: ``In this paper, we give a new approach to prove inequalities for the Schwab-Borchardt mean\index{means!Schwab-Borchardt}, the lemniscatic mean\index{means!lemniscatic} and the arithmetic geometric mean\index{means!arithmetic-geometric}. Additionally, we apply these means to inequalities for trigonometric functions or the lemniscate functions by considering several functional inequalities. One of these applications includes infinite product formulas for the lemniscate function\index{infinite products!for the lemniscate function}\index{lemniscate!function} and the arithmetic geometric mean\index{infinite products!for the arithmetic-geometric mean} by considering several functional equations."


\item\vspace{9pt} \label{Osl2015} (2015) Thomas J. Osler\index[ppl]{Osler, Thomas J.|textbf}, A product of nested radicals for the AGM. \emph{The American Mathematical Monthly} \textbf{122}(9), 886--887. 
\index{nested radicals}

{\footnotesize \href{https://mathscinet.ams.org/mathscinet/article?mr=3418211}{MR3418211}}

{\footnotesize Source: \href{https://doi.org/10.4169/amer.math.monthly.122.9.886}{doi.org (Taylor \& Francis)}.}

{\footnotesize Cited in
\hyperref[Osl2016c]{\textsc{Osler and Kosior 2016c}}\index[ppl]{Osler, Thomas J.}\index[ppl]{Kosior, Jesse M.}.}

\vspace{4pt}

From the abstract: ``The arithmetic-geometric mean\index{means!arithmetic-geometric} of two positive numbers $a$ and $b$ is the common limit of two sequences generated by an iterative process\ldots In this paper, we derive an infinite product representation for the AGM. The factors of this product are nested radicals resembling Vieta's\index[ppl]{Vi\`{e}te, Fran\c{c}ois} famous product for pi."\index{Vi\`{e}te's formula for $\tfrac{2}{\pi}$}


\item\vspace{9pt} \label{Aok2016} (2016) Noboru Aoki\index[ppl]{Aoki, Noboru|textbf} and Shota Kojima\index[ppl]{Kojima, Shota|textbf}, Nested square roots and Poincar\'{e} functions. \emph{Tokyo Journal of Mathematics} \textbf{39}(1), 241--269.
\index{continued square roots!and trigonometric functions}
\index{continued square roots!of terms $a_n=\pm2$}
\index{continued square roots!of terms $a_n=\pm c$}

{\footnotesize \href{https://mathscinet.ams.org/mathscinet/article?mr=3543142}{MR3543142}}

{\footnotesize Source: \href{https://doi.org/10.3836/tjm/1471873313}{doi.org (Project Euclid)}.}

\vspace{4pt}

From the introduction: ``Let $c$ be a real number with $c\ge2$ and $\epsilon_1, \epsilon_2,\ldots$ an infinite sequence consisting of $\pm1$. In this paper we are concerned with nested square roots\index{nested square roots} of the form
\begin{equation}\label{E:AoK2016-1}
R_c(\epsilon_1,\epsilon_2,\epsilon_3,\ldots,\epsilon_m)=\epsilon_1\sqrt{c+\epsilon_2\sqrt{c+\epsilon_3\sqrt{c+\cdots+\epsilon_m\sqrt{c}}}}
\end{equation}
and infinite nested square roots
\begin{equation}\label{E:AoK2016-2}
R_c(\epsilon_1,\epsilon_2,\epsilon_3,\ldots):=\lim_{m\to\infty}R_c(\epsilon_1,\epsilon_2,\epsilon_3,\ldots,\epsilon_m)\;.
\end{equation}
[\ldots]\;In the case of $c=2$, it is known that the nested root \eqref{E:AoK2016-1} can be expressed \ldots as:
\begin{equation}\label{E:AoK2016-3}
R_2(\epsilon_1,\epsilon_2,\epsilon_3,\ldots,\epsilon_m)=2\cos\pi\left(\dfrac{a_1}{2}+\dfrac{a_2}{2^2}+\cdots+\dfrac{a_m}{2^m}+\dfrac{1}{2^{m+1}}\right)\;,
\end{equation}
where
\[a_i=\dfrac{1-\epsilon_1\cdots\epsilon_i}{2}=
\begin{cases}
0 &(\textrm{if } \epsilon_1\cdots\epsilon_i=1)\notag\\
1 &(\textrm{if } \epsilon_1\cdots\epsilon_i=-1)\;.\notag
\end{cases}
\]
Taking $\lim_{m\to\infty}$ of \eqref{E:AoK2016-3}, we obtain a simple formula for the infinite nested square root:
\begin{equation}\label{E:AoK2016-4}
R_2(\epsilon_1,\epsilon_2,\epsilon_3,\ldots)=2\cos\alpha\pi\;,
\end{equation}
where $\alpha$ is a real number defined by the $2$-adic expansion
\[\alpha=\dfrac{a_1}{2}+\dfrac{a_2}{2^2}+\dfrac{a_3}{2^3}+\cdots\;.\]
These formulas were proved by Wiernsberger\index[ppl]{Wiernsberger, Paul} \ldots in 1905, and about thirty years later Lebesgue\index[ppl]{Lebesgue, Henri}{\ldots}independently found the same formulas.

``The purpose of this paper is to give a generalization of the formulas \eqref{E:AoK2016-3} and \eqref{E:AoK2016-4} to the case $c\ge2$. To accomplish the task, we need a suitable function which will take the place of $\cos x$. In the proof of the formulas \eqref{E:AoK2016-3} and \eqref{E:AoK2016-4}, the duplication formula
\[2\cos2x=(2\cos x)^2-2\]
was crucial. It is therefore natural to seek for a function $f(x)$ satisfying the functional equation
\begin{equation}\label{E:AoK2016-5}
f(sx)=f(x)^2-c\;,
\end{equation}
where $s$ is a constant depending only on $c$. Such functional equations were studied by Poincar\'{e},\index[ppl]{Poincar\'{e}, Henri} who showed that there exists an entire function $f(x)$ satisfying \eqref{E:AoK2016-5}.''

The authors also use their method to give a generalization of Vi\`{e}te's formula \eqref{E:Vie1593-1} for $\tfrac{2}{\pi}$.\index{Vi\`{e}te's formula for $\tfrac{2}{\pi}$!generalizations}

It should be noted that \eqref{E:AoK2016-3} is to be found in \hyperref[Boc1899]{\textsc{Bochow 1899}}\index[ppl]{Bochow, Karl} (expressed using the sine rather than the cosine function), and foreshadowed in \hyperref[Luc1878]{\textsc{Lucas 1878}}\index[ppl]{Lucas, Edouard}. Among several minor bibliographical errors, the ones most in need of correction are the references to ``M. P. Wiernsberger", which should be to ``P. Wiernsberger." 


\item\vspace{9pt} \label{Cha2016} (2016) Mu-Ling Chang\index[ppl]{Chang, Mu-Ling|textbf} and Chia-Chin (Cristi) Chang\index[ppl]{Chang, Chia-Chin (Cristi)|textbf}, Evaluation of pi by nested radicals. \emph{Mathematics Magazine} \textbf{89}(5), 336--337.
\index{nested radicals}

{\footnotesize \href{https://mathscinet.ams.org/mathscinet/article?mr=3593652}{MR3593652}}

{\footnotesize Source: \href{https://doi.org/10.4169/math.mag.89.5.336}{doi.org (Taylor \& Francis)}.}

\vspace{4pt}

The formula \eqref{E:Cat1842-1} for $\pi$ from \hyperref[Cat1842]{\textsc{Catalan 1842}}\index[ppl]{Catalan, E.} is again proved.\index{pi@$\pi\;(3.14159\ldots)$!continued square root expressions for}\index{constants, named!pi@$\pi\;(3.14159\ldots)$}


\item\vspace{9pt} \label{Les2016} (2016) Devyn A. Lesher\index[ppl]{Lesher, Devyn A.|textbf} and Chris D. Lynd\index[ppl]{Lynd, Chris D.|textbf}, Convergence results for the class of periodic left nested radicals. \emph{Mathematics Magazine} \textbf{89}(5), 319--335.
\index{nested radicals!left}
\index{nested radicals!of periodic real terms}

{\footnotesize \href{https://mathscinet.ams.org/mathscinet/article?mr=3593651}{MR3593651}}

{\footnotesize Source: \href{https://doi.org/10.4169/math.mag.89.5.319}{doi.org (Taylor \& Francis)}.}

\vspace{4pt}

The authors write: ``There are two fundamental questions that are found throughout the research on nested radicals: (1) Given a particular form of nested radical\index{nested radicals}, under what conditions does it converge? (2) Given a particular form of nested radical, which numbers can be expressed as the limit of a nested radical?\footnote{By ``nested radicals" the authors mean both continued and iterated $r_i$th roots.\index{continued rocc roots@continued $r_i$th roots}\index{iterated rocc roots@iterated $r_i$th roots}}

``$[\ldots]$In this paper we address both research questions, as they pertain to periodic \emph{left} nested radicals\index{nested radicals!left}. We provide numerical examples to illustrate each result.

``We also provide a recipe for constructing nested radicals with a predetermined end-behavior. If you choose the form of the nested radical and a limit, the recipe shows how to construct the unique nested radical\index{nested radicals} of the chosen form that converges to the chosen limit.\index{iterated rocc roots@iterated $r_i$th roots!convergence to a predetermined limit} If you choose a periodic limiting sequence, the recipe shows how to construct a nested radical\index{nested radicals} that asymptotically converges to the chosen periodic sequence."


\item\vspace{9pt} \label{Nis2016} (2016) Ryo Nishimura\index[ppl]{Nishimura, Ryo|textbf}, A generalization of Vi\`{e}te's infinite product and new mean iterations. \emph{The Australian Journal of Mathematical Analysis and Applications} \textbf{13}(1), Art. 20, 1--9.
\index{Vi\`{e}te's formula for $\tfrac{2}{\pi}$!generalizations}

{\footnotesize \href{https://mathscinet.ams.org/mathscinet/article?mr=3590882}{MR3590882}}

{\footnotesize Source: \href{https://ajmaa.org/searchroot/files/pdf/v13n1/v13i1p20.pdf}{Australian Journal of Mathematical Analysis and Applications}.}

\vspace{4pt}

From the abstract: ``In this paper, we generalize Vi\`{e}te's\index[ppl]{Vi\`{e}te, Fran\c{c}ois} infinite product formula by use of Chebyshev polynomials\index{polynomials!Chebyshev}\index{Chebyshev polynomials}. Furthermore, the infinite product formula for the lemniscate sine\index{infinite products!for the lemniscate sine}\index{lemniscate!sine} is also generalized. Finally, we obtain new mean iterations by use of these infinite product formulas."\index{infinite products} 
 

\item\vspace{9pt} \label{Oht2016} (2016) Hideyuki Ohtsuka\index[ppl]{Ohtsuka, Hideyuki|textbf}, Problem H-788. \emph{The Fibonacci Quarterly} \textbf{54}(2), 185. Solution by the proposer, \emph{ibid.} \textbf{56}(2), 187.

{\footnotesize Solution source: \href{https://www.fq.math.ca/Problems/AdvProbMay2018.pdf}{$\langle$Fibonacci Quarterly$\rangle$}.}

\vspace{4pt}

The problem states: ``Given $c>0$. Determine 
\[\lim_{n\to\infty}\sqrt{cF_2^2+\sqrt{cF_4^2+\sqrt{cF_8^2+\sqrt{\cdots+\sqrt{cF_{2^n}^2}}}}}\;."\]
Here $F_1=F_2=1$ and $F_n$ is the $n$th Fibonacci number\index{Fibonacci!numbers}\index{numbers!Fibonacci}. The value is shown to be $\tfrac{3+\sqrt{4c+5}}{2}$.\index{continued square roots!of Fibonacci numbers}


\item\vspace{9pt} \label{Osl2016a} (2016) Thomas Osler\index[ppl]{Osler, Thomas J.|textbf}, Iterations for the lemniscate constant\index{lemniscate!constant\;$(2.62205\ldots)$} resembling the Archimedean algorithm for pi. \emph{The American Mathematical Monthly} \textbf{123}(1), 90--93. 
\index{Archimedean algorithm}
\index{Vi\`{e}te's formula for $\tfrac{2}{\pi}$!generalizations}

{\footnotesize \href{https://mathscinet.ams.org/mathscinet/article?mr=3453543}{MR3453543}} 

{\footnotesize Source: \href{https://doi.org/10.4169/amer.math.monthly.123.1.90}{doi.org (Taylor \& Francis)}.} 

\vspace{4pt}

From the abstract: ``We give an iterative algorithm that converges to the lemniscate constant $L$\index{lemniscate!constant\;$(2.62205\ldots)$}. This algorithm resembles the famous Archimedean algorithm\index{pi@$\pi\;(3.14159\ldots)$!Archimedean algorithm for}\index{constants, named!pi@$\pi\;(3.14159\ldots)$} for $\pi$. The derivation is based on the recently discovered product of nested radicals\index{nested radicals!infinite products of} for $\tfrac{2}{L}$ by Aaron Levin\index[ppl]{Levin, Aaron}. Levin's product closely resembles Vieta's\index[ppl]{Vi\`{e}te, Fran\c{c}ois} historic product for $\tfrac{2}{\pi}$."\index{Vi\`{e}te's formula for $\tfrac{2}{\pi}$!generalizations} 


\item\vspace{9pt} \label{Osl2016b} (2016) Thomas J. Osler\index[ppl]{Osler, Thomas J.|textbf}, Walter Jacob\index[ppl]{Jacob, Walter|textbf}, and Ryo Nishimura\index[ppl]{Nishimura, Ryo|textbf}, An infinite product of nested radicals for $\log x$ from the Archimedean algorithm. \emph{The Mathematical Gazette} \textbf{100}(548), 274--278.
\index{infinite products!of continued square roots}
\index{nested radicals!infinite products of}
\index{Archimedean algorithm}

{\footnotesize \href{https://mathscinet.ams.org/mathscinet/article?mr=3520821}{MR3520821}}

{\footnotesize Source: \href{https://www.jstor.org/stable/44161814}{JSTOR}.}

\vspace{4pt}

The authors show that $\log x$ may be written as
\begin{equation*}
\dfrac{x-1}
{\sqrt{x}\sqrt{\tfrac{1}{2}+\tfrac{1}{2}\bigl(\tfrac{1+x}{2\sqrt{x}}\bigr)}
\sqrt{\tfrac{1}{2}+\tfrac{1}{2}\sqrt{\tfrac{1}{2}+\tfrac{1}{2}\bigl(\tfrac{1+x}{2\sqrt{x}}\bigr)}}
\sqrt{\tfrac{1}{2}+\tfrac{1}{2}\sqrt{\tfrac{1}{2}+\tfrac{1}{2}\sqrt{\tfrac{1}{2}+\tfrac{1}{2}\bigl(\tfrac{1+x}{2\sqrt{x}}\bigr)}}}\cdots}
\end{equation*}
based on a similar, more general expansion generated by the Archimedean algorithm.


\item\vspace{9pt} \label{Osl2016c} (2016) Thomas J. Osler\index[ppl]{Osler, Thomas J.|textbf} and Jesse M. Kosior\index[ppl]{Kosior, Jesse M.|textbf}, A sequence of good approximations for the period of a pendulum with large initial amplitude. \emph{The Mathematical Scientist} \textbf{41}(1), 40--44. 

{\footnotesize \href{https://mathscinet.ams.org/mathscinet/article?mr=3561649}{MR3561649}}

{\footnotesize Source: \href{https://drive.google.com/file/d/1ThZSf2W_YRIzDZ6GmgznaPMh-3q41oI8/view?usp=sharing}{Google Drive.}}

\vspace{4pt}

From the abstract: ``We present three elementary approximate formulas for the period of a pendulum\index{pendulum, period of} which starts at rest from a large angle of displacement. The first of these formulas is known, but the other two may be new. These three formulas result from taking the first three partial products of a new infinite product of nested radicals\index{nested radicals!infinite products of} for the complete elliptic integral\index{elliptic integral} of the first kind that gives the exact period. Thus, more elementary approximations can be obtained from this exact product, but they become increasingly complex. Therefore, we stopped at three. We give a detailed table clearly displaying the accuracy of the approximations over the full range of possible initial angles of displacement. This infinite product of nested radicals\index{nested radicals!infinite products of} is a special case of a new infinite product for the arithmetic-geometric mean\index{means!arithmetic-geometric} that has appeared recently.'' 


\item\vspace{9pt} \label{Sch2016} (2016) Fritz Schweiger\index[ppl]{Schweiger, Fritz|textbf}, \emph{Continued Fractions and their Generalizations: A Short History of $f$-expansions}, Docent Press, Boston, ISBN 978-1942795933.\index{f@$f$-expansions}

{\footnotesize \href{https://mathscinet.ams.org/mathscinet/article?mr=3468648}{MR3468648}}

\vspace{4pt}

From the introduction: ``This book is about the history of $f$-expansions, their theory, their application, and their connection to other parts of mathematics$\ldots$ As a kind of background theory we occasionally use the language of fibred systems\index{fibred systems}." $g$-adic expansions of real numbers\index{real numbers!g@$g$-adic expansions of}, where $f(x)=\tfrac{x}{g}$, and continued fraction expansions\index{real numbers!continued fraction expansions of}, where $f(x)=\tfrac{1}{x}$, are given as ``prominent examples" of $f$-expansions; the latter are reviewed in Chapter 2, the former in Chapter 3. Infinite products\index{infinite products} are briefly covered in Chapter 4. 

Chapter 5 looks at \hyperref[Kak1924]{\textsc{Kakeya 1924}}\index[ppl]{Kakeya, S\^{o}ichi}, which attempted to unify $g$-adic expansions and continued fractions\index{continued fractions} under a common general form --- the precursor by twenty years of Bissinger's\index[ppl]{Bissinger, B. H.} $f$-expansions. (Schweiger amusingly laments subsequent authors' completely coincidental permutations of some of Kakeya's terminology and notation.) Continued cotangents\index{continued cotangents} and \hyperref[Leh1938]{\textsc{Lehmer 1938}}\index[ppl]{Lehmer, D. H.} are the subject of Chapter 6. Chapter 7 discusses \hyperref[Bis1944]{\textsc{Bissinger 1944}}\index[ppl]{Bissinger, B. H.} and \hyperref[Eve1946]{\textsc{Everett 1946}}\index[ppl]{Everett, C. J.}, which introduce $f$-expansions based on, respectively, increasing and decreasing functions $f$; the chapter's closing comments concern \hyperref[Rec1950]{\textsc{Rechards 1950}}.

The Borel-Bernstein Theorem\index{Borel-Bernstein Theorem} is taken up in Chapter 8, while ergodic properties of $f$-expansions\index{f@$f$-expansions!ergodic properties of} are discussed in Chapter 9. Chapter 10 reviews the work of R\'{e}nyi\index[ppl]{Renyi@R\'{e}nyi, A.} and others in extending Bissinger's\index[ppl]{Bissinger, B. H.} definitions, including the distinction between those $f$-expansions\index{f@$f$-expansions} with independent versus dependent digits. Chapter 11 considers Gauss's investigations into functional equations involving tails of continued fraction representations, and Kuzmin's\index[ppl]{Kuzmin, R. O.} general result that proves one of Gauss's\index[ppl]{Gauss, Carl Friedrich} conjectures. This leads to Chapter 12, L\'{e}vy's\index[ppl]{Levy@L\'{e}vy, P.} alternate proof of the Gauss conjecture based on ``the \emph{dual algorithm} or the \emph{natural extension} of continued fractions." 

The work of Gel'fond\index[ppl]{Gel'fond, A. O.}, Cigler\index[ppl]{Cigler, J.}, and Parry\index[ppl]{Parry, W.} in 1959 and 1960 to extend R\'{e}nyi's\index[ppl]{Renyi@R\'{e}nyi, A.} $\beta$-expansions\index{beta@$\beta$-expansions} is dealt with in Chapter 13. Chapter 14 begins by recalling a mostly-forgotten paper by T. E. McKinney\index[ppl]{McKinney, T. E.} from 1907 on $\lambda$-continued fractions\index{continued fractions!lambda-type@$\lambda$-type}, whose underlying map $f_\lambda:[\lambda-1,\lambda]\to[\lambda-1,\lambda]$ is defined by
\[f_\lambda(x)=\dfrac{\epsilon}{x}-\left\lfloor{\dfrac{\epsilon}{x}+1-\lambda}\right\rfloor\]
where $x\ne 0$ and $\epsilon=\textrm{sign}\,x$. After pointing to a special case of this by A. Hurwitz\index[ppl]{Hurwitz, A.} from 1888, Schweiger reviews the investigations of Nakada\index[ppl]{Nakada, H.}, Ito\index[ppl]{Ito, S.}, and Tanaka\index[ppl]{Tanaka, S.} into continued fractions of this type. Chapter 15, on discontinuous groups, looks at D. Rosen's\index[ppl]{Rosen, D.} extension of continued fractions, which derives from linear fractional transformations\index{linear fractional transformations} in the group generated by $S(z)=z+\lambda$ and $T(z)=-\tfrac{1}{z}$. 

Ergodic theory\index{ergodic theory}, and its connection to $f$-expansions, is alluded to throughout the book, but becomes the focus of Chapters 16, 17, and 18, with invariant measure\index{invariant measure} given particular attention in the latter two. C. Shannon's\index[ppl]{Shannon, C.} concept of entropy as a measure of information is viewed in Chapter 19 as ``a kind of average for the cylinders in an $f$-expansion (more general [sic] of a fibred system\index{fibred systems})." This leads naturally into a discussion, in Chapter 20, of Hausdorff dimension\index{Hausdorff dimension} and P. Billingsley's\index[ppl]{Billingsley, P.} generalization, ``which is adapted to problems connected with $f$-expansions." \index{f@$f$-expansions}

Multidimensional generalizations of continued fractions\index{continued fractions!multidimensional} and $f$-expansions\index{f@$f$-expansions!multidimensional} are taken up in Chapter 21. Here Schweiger is keen to point out a ``serious flaw" in one of his own papers, and its subsequent detection by some later authors but perpetuation by others. The book's final chapter is a three-page synopsis of aspects of chaos theory.\index[ppl]{Schweiger, Fritz}\index{chaos theory}

A bibliography of over 200 sources is provided. There is no index.


\item\vspace{9pt} \label{Vel2016a} (2016a) Pierluigi Vellucci\index[ppl]{Vellucci, Pierluigi|textbf} and Alberto Maria Bersani\index[ppl]{Bersani, Alberto Maria|textbf}, The class of Lucas-Lehmer polynomials. \emph{Rendiconti di Matematica e delle sue Applicazioni. Serie VII}, \textbf{37}, 43--62.
\index{polynomials!Lucas-Lehmer}

{\footnotesize \href{https://mathscinet.ams.org/mathscinet/article?mr=3622304}{MR3622304}}

{\footnotesize Source: \href{https://drive.google.com/file/d/1hTXCblBc3zCKzwoZuhpp04cKXEdHEeuv/view?usp=sharing}{Google Drive.}} 

{\footnotesize Cited in
\hyperref[Vel2016b]{\textsc{Vellucci and Bersani 2016b}} and \hyperref[Vel2016c]{\textsc{2016c}}\index[ppl]{Vellucci, Pierluigi}\index[ppl]{Bersani, Alberto Maria}.}

\vspace{4pt}

From the abstract: ``In this paper we introduce a new sequence of polynomials, which follow the same recursive rule of the well-known Lucas-Lehmer\index{polynomials!Lucas-Lehmer}\index{Lucas-Lehmer!polynomials} integer sequence\index{Lucas-Lehmer!numbers}\index{numbers!Lucas-Lehmer}. We show the most important properties of this sequence, relating them to the Chebyshev polynomials\index{polynomials!Chebyshev}\index{Chebyshev polynomials} of the first and second kind." The polynomials are defined recursively by $L_0(x)= x$ and $L_n(x) = L^2_{n-1}(x)-2$. Zeroes of $L_n$ turn out to be finite continued square roots with terms $\pm 2$.\index{continued square roots!of terms $a_n=\pm2$}


\item\vspace{9pt} \label{Vel2016b} (2016b) Pierluigi Vellucci\index[ppl]{Vellucci, Pierluigi|textbf} and Alberto Maria Bersani\index[ppl]{Bersani, Alberto Maria|textbf}, Ordering of nested square roots of 2 according to Gray code. [\emph{Ramanujan J.}, \textbf{45} (2018)(1), 197--210.]

{\footnotesize \href{https://mathscinet.ams.org/mathscinet/article?mr=3745071}{$\langle${MR3745071}$\rangle$}}

{\footnotesize Source: \href{https://doi.org/10.1007/s11139-016-9862-5}{doi.org (Springer)}.}

{\footnotesize Cited in
\hyperref[Vel2016c]{\textsc{Vellucci and Bersani 2016c}}\index[ppl]{Vellucci, Pierluigi}\index[ppl]{Bersani, Alberto Maria}.}
\index{Gray code}
\index{continued square roots!of terms $a_n=\pm2$}

\vspace{4pt}

From the abstract: ``In this paper we discuss some relations between zeros of Lucas-Lehmer polynomials and Gray code\index{Gray code}. We study nested square roots of 2 applying a `binary code' that associates bits 0 and 1 to $\oplus$ and $\ominus$ signs in the nested form. This gives the possibility to obtain an ordering for the zeros of Lucas-Lehmer polynomials, which assume the form of nested square roots of 2."\index{continued square roots!of terms $a_n=\pm2$}


\item\vspace{9pt} \label{Vel2016c} (2016c) Pierluigi Vellucci\index[ppl]{Vellucci, Pierluigi|textbf} and Alberto Maria Bersani\index[ppl]{Bersani, Alberto Maria|textbf}, $\pi$-formulas and Gray code. \emph{Ricerche di Matematica} \textbf{68}, no. 2, 551--569 $\langle$2019$\rangle$.
\index{Gray code}. 
\index{continued square roots!of terms $a_n=\pm2$}
\index{pi@$\pi\;(3.14159\ldots)$!continued square root expressions for}\index{constants, named!pi@$\pi\;(3.14159\ldots)$}

{\footnotesize \href{https://mathscinet.ams.org/mathscinet/article?mr=4029300}{$\langle${MR4029300}$\rangle$}}

{\footnotesize Source: \href{https://doi.org/10.1007/s11587-018-0426-4}{doi.org (Springer)}.}

\vspace{4pt}

From the abstract: ``In previous papers we introduced a class of polynomials which follow the same recursive formula as the Lucas-Lehmer numbers\index{polynomials!Lucas-Lehmer}\index{Lucas-Lehmer!numbers}\index{Lucas-Lehmer!polynomials}\index{polynomials!Lucas-Lehmer}, studying the distribution of their zeros and remarking that this distribution follows a sequence related to the binary Gray code\index{Gray code}. It allowed us to give an order for all the zeros of every polynomial $L_n$. In this paper, the zeros, expressed in terms of nested radicals\index{nested radicals}, are used to obtain two new formulas for $\pi$: the first can be seen as a generalization of the known formula\index{pi@$\pi\;(3.14159\ldots)$!continued square root expressions for}\index{constants, named!pi@$\pi\;(3.14159\ldots)$}
\[\pi=\lim_{n\to\infty}2^{n+1}\sqrt{2-\sqrt{2+\sqrt{2+\sqrt{2+\cdots+\sqrt{2}}}}}\]
[where there are $(n+1)$ 2s under the outer radical], related to the smallest positive zero of $L_n$; the second is an exact formula for $\pi$ achieved thanks to some identities valid for $L_n$."


\item\vspace{9pt} \label{Mcg} (Undated) Michael McGuffin\index[ppl]{McGuffin, Michael|textbf} and Brian Wong\index[ppl]{Wong, Brian|textbf}, The Museum of Nested Radicals. 

{\footnotesize Source: \href{http://www.dgp.toronto.edu/~mjmcguff/math/nestedRadicals.html}{Dynamic Graphics Project, University of Toronto}.}

{\footnotesize Cited in 
\hyperref[Hum2007]{\textsc{Humphries 2007}}\index[ppl]{Humphries, Peter J.}, 
\hyperref[Zim2008a]{\textsc{Zimmerman and Ho 2008a}}\index[ppl]{Zimmerman, Seth}\index[ppl]{Ho, Chungwu}, and 
\hyperref[Wei]{\textsc{Weisstein n.d.}}\index[ppl]{Weisstein, Eric W.}.}

\vspace{4pt}

This web site offers the formula\index{continued robb roots@continued $r$th roots!of constant nonnegative real terms}
\[x=\sqrt[n]{(1-q)x^n+qx^{n-1}\sqrt[n]{(1-q)x^n+qx^{n-1}\sqrt[n]{\cdots}}}\;,\]
with a number of special cases for various substitutions and values of $q$ and $n$. It also attributes to Ramanujan the formula\index[ppl]{Ramanujan, Srinivasa}\index{continued square roots!Ramanujan's}
\begin{equation}\label{E:McG-1}
\begin{split}
x+n+a=&\sqrt{ax +(n+a)^2+}\\
&x\sqrt{a(x+n) +(n+a)^2+}\\
&(x+a)\sqrt{a(x+2n) +(n+a)^2+}\\
&(x+2a)\sqrt{\cdots}\;,
\end{split}
\end{equation}
(where each square root subsumes its successors), again with several special cases listed, including the first problem from \hyperref[Ram1911]{\textsc{Ramanujan 1911}}. A few words of justification for these expressions are given; no sources are listed. 


\item\vspace{9pt} \label{Wei} (Undated) Eric W. Weisstein\index[ppl]{Weisstein, Eric W.|textbf}, Nested Radical. From \emph{MathWorld --- A Wolfram Web Resource}.
\index{Wolfram MathWorld}

{\footnotesize Source: \href{http://mathworld.wolfram.com/NestedRadical.html}{Wolfram MathWorld}.}

\vspace{4pt}

Gives an overview of the results of Herschfeld\index[ppl]{Herschfeld, Aaron}, Ramanujan\index[ppl]{Ramanujan, Srinivasa}, and Vi\`{e}te\index[ppl]{Vi\`{e}te, Fran\c{c}ois}, along with continued square root expressions for sines and cosines of special angles\index{continued square roots!and trigonometric functions}, the golden ratio\index{golden ratio ($\tfrac{1+\sqrt{5}}{2}=1.61803\ldots$)}\index{constants, named!golden ratio ($\tfrac{1+\sqrt{5}}{2}=1.61803\ldots$)}, the silver ratio\index{silver ratio\;$(1+\sqrt{2}=2.41421\ldots$)}\index{constants, named!silver ratio\;$(1+\sqrt{2}=2.41421\ldots$)}, and the plastic constant\index{plastic constant ($1.32471\ldots$)}\index{constants, named!plastic constant ($1.32471\ldots$)}. 

\end{enumerate} 


\appendix

\section{Five biographical sketches}\label{S:bios}

In the course of compiling this bibliography, some of the authors' life stories came into view. As we have seen, continued square roots have drawn the attention of a broad spectrum of writers, including some of the great names in mathematics, professionals forgotten or less well-known today, many math educators, and many amateurs. Here we present brief accounts of five lives that seemed interesting, were not already too-copiously documented, and which could be recounted beyond the basics.

\small




\subsection*{T. S. E. Dixon \normalfont{(1845--1898)}} \hyperref[Dix1878]{\textsc{Dixon 1878.}}\index[ppl]{Dixon, T. S. E.}

Born Theron Solimon Eugene Dixon in Jericho, Vermont, USA, he seems to have lived most of his life in the American Midwest. He graduated in 1868 from Beloit College in the Wisconsin town of that name, and spent the next year apprenticing at a law firm in Janesville, Wisconsin.\footnote{Possibly with attorneys Charles G. Williams (\href{https://www.quillproject.net/person/7362}{quillproject.net}) or John Winans (\href{https://npgallery.nps.gov/NRHP/GetAsset/NRHP/86000205_text}{National Park Service})} This was a time when many lawyers trained to pass the bar by apprenticing, rather than attending law school; Dixon passed the bar in 1869.\footnote{Beloit College Monthly. \href{https://babel.hathitrust.org/cgi/pt?id=wu.89067543017&seq=293&q1=Dixon}{HathiTrust.}} He moved to Chicago in 1874 and began practice as a patent attorney, eventually as a partner in the firm of Hill \& Dixon.

Dixon gained prominence in the 1880s for his role in the ``Telephone Cases'', a long-running legal dispute brought by many claimants against the American Bell Telephone Co. and Alexander Graham Bell\index[ppl]{Bell, Alexander Graham} over priority for the invention of the telephone. Dixon was involved in the U.S. Supreme Court case on behalf of the Overland Telephone Co. of New Jersey and People's Telephone Co. of Quebec, Canada. The Supreme Court ruled in favor of Bell in its decision of October 1887.\footnote{``Cases Adjudged in the Supreme Court of the United States, at October Term, 1887. The Telephone Cases.'' \href{https://tile.loc.gov/storage-services/service/ll/usrep/usrep126/usrep126001/usrep126001.pdf}{Library of Congress.}}

Beyond handling the legalities of patents, Dixon was an inventor as well, obtaining at least seventeen U.S. patents in his own name or jointly with his law partner Lysander Hill\index[ppl]{Hill, Lysander} and others. Starting in 1869 with a ``street-car starter'' (in the era of animal-drawn street cars),\footnote{``[A]n improved device for attachment to horse cars, and other wheeled vehicles, by the use of which the power will be first applied to revolve the wheels of the vehicle, and thus start it with less effort than when the draft is applied directly to the body of the car.'' ---\emph{Scientific American}, November 27, 1869, p. 348.} his inventions ranged from gadgets such as ``a mouth-piece to hold cigars'', to industrial devices like a rail car axle and an early version of the air brake for rail cars.

The polymathic attorney published at least three mathematics papers, including ``A General Algebraic Method for the Solution of Equations'' in \emph{The Analyst} for January 1882,\footnote{\emph{The Analyst}, Vol. 9, No. 1, pp. 1--8. \href{https://www.jstor.org/stable/2636119}{JSTOR.}} in which he claimed to give a method for solving polynomial equations of degree $n$ in terms of $n$th roots of certain constants. The paper is followed by a parenthetical comment from the journal editor, J. E. Hendricks\index[ppl]{Hendricks, J. E. (Joel Evans)}: ``[W]e have not been able to detect any error in the foregoing very elegant discussion{\ldots}we commend it to the careful consideration and critical scrutiny of our readers.'' Two months later, however, the author himself reports, ``I have discovered `the weak link in the chain' of [my paper's] logic, as applied to the solution of the equation of the $5$th degree. It lies in the statement{\ldots}that in the supplemental equations `one letter may be interchanged with another without disturbing any relations'. If, in brief, numerical values be assigned to the letters, an interchange in the order of their arrangement materially disturbs the relations.'' Dixon closes his mea culpa\index{mea culpas} with ``I regret that pressure of professional labors renders it impossible for me to further prosecute the work at present'' and he trusts that ``further research will yet result in the development of a perfect and complete Theory of Equations.''

The work most tied to T. S. E. Dixon's name in internet searches, however, is a book of 460+ pages published in 1895, claiming that Francis Bacon\index[ppl]{Bacon, Francis} wrote the works of Shakespeare\index[ppl]{Shakespeare, William}. Of the many polite reviews of \emph{Francis Bacon and his Shakespeare},\footnote{\href{https://archive.org/details/francisbaconhiss00dixo/mode/2up?ref=ol}{Internet Archive.}} this from \emph{The Review of Reviews} is typical: ``This work begins and ends with the assumption that the plays attributed to Shakespeare were really written by Bacon. Mr. Dixon does not labor to prove this proposition; he assumes its truth at the outset, and seeks to interpret the plays in the light of this explanation of their origin.''\footnote{\emph{The Review of Reviews}, Vol. 12, No. 6, December 1895, p. 753. \href{https://books.googleusercontent.com/books/content?req=AKW5QafVOaaxEbuZ9lVCTeeGU3_ET4eXJsUx6EtvndAyfJjfhebOhd4u1jkujlNzJC_Hsy5MhdnpJtJ04LiE48yQ_nl49T7B03dKMv07RkF4Rt7ITgxedAm79z1M0PDy8YEuY-f73icXGSNinU9XJkUEJ7o7TgV498Kw3xNmUjJB6m2L5v55Z0LNsd6CAUn00hKtSTcuUjIaVu38ohfJ7UtSR_C9TMjccx1EgNNYe_oHb6IkeLlbq4M8TFt8gIny1G63jjdCmoqi}{Google Books.}} Other reviewers were less generous. A scathing assessment in \emph{Self-Culture}, November 1895, amounting to a six-page rebuttal of the book's premise, states, ``There is no question of the thoughtfulness and refinement of Mr. Dixon's treatment of the subject, but in the matter of knowledge, whether of Bacon or of Shakespeare, there is almost none at all, and in the matter of argument the reasonings are in line with those marvelous attempts of the legal profession in the criminal courts which almost compel one to assume that the practice of law tends to absolute lunacy.''\footnote{\emph{Self-Culture}, Vol. 2, No. 2, November 1895, pp. 498--503. \href{https://books.googleusercontent.com/books/content?req=AKW5QaeMme0q3HTtDKHpl_98b42H3Dvsw2n7FnyXVWGkYlojr8qie5YqHg1MIkBEGHM5KcdP_ydMULupW2-IzEyURVcrGr7uhhVkM6cpB-QiaTk8lYbzR90yu7MFk85Pln9PSEFtwxVvZrGH3L0sY9ryrzGaibWCfTu4L_ie4VOalULDnoioffHSnMq7DXDmoNW0gnj-EFRhEAhJ_WAI_FdNcPSPPJLBkqBtIUVkJaczcG6_SV8wGvdK_lkSf_4LyfU-jcr_Yct2}{Google Books.}} 

Dixon died tragically. Page 4 of the September 15, 1898 \emph{New York Times} reports that, on the previous day, ``Henry Lester\index[ppl]{Lester, Henry}, proprietor of the Princess Anne Cottage at Virginia Beach, [Virginia,] and two of his guests, T. S. E. Dixon of Chicago and Arthur McLaughlin\index[ppl]{McLaughlin, Arthur} of Newark, N. J., were drowned{\ldots}while bathing at the beach\ldots. It is supposed that McLaughlin was carried out by a receding wave, and the others were drowned in the effort to reach him.'' Dixon was 53 years old. He was survived by his wife, Bertha I. (Wright) Dixon\index[ppl]{Dixon, Bertha I. (Wright)} (to whom he dedicated his Bacon/Shakespeare book), and at least one child, a daughter, Jessie Dixon\index[ppl]{Dixon, Jessie} (c. 1878--1963).\footnote{\href{https://dcms.beloit.edu/digital/collection/magazine/id/5683}{Beloit College Digital Collections.}}


\subsection*{Christian Friedrich Eichhorn \normalfont{(1804--1836)}} \hyperref[Eic1834]{\textsc{Eichhorn 1834}}\index[ppl]{Eichhorn, Christian Friedrich}.

Eichhorn was born in Osnabr\"{u}ck, a city of northwest Germany midway between Dusseldorf and Bremen. Graduating from the Protestant Gymnasium there in 1822, he entered the University of G\"{o}ttingen\index{University of G\"{o}ttingen} in the spring of 1823, and received his doctorate in philosophy in 1826. After traveling to Paris, he was a private lecturer in Göttingen in pure and applied mathematics and mathematical physics. The Higher Trade School of Hannover\footnote{Now Liebniz Universit\"{a}t Hannover\index{Liebniz Universit\"{a}t Hannover}.} was established in 1831, and Eichhorn was appointed chair of mechanical engineering and applied mathematics there in March of that year, a position he held until his death.\footnote{\emph{Allgemeine Deutsche Biographie}, Band 5, Duncker \& Humblot, Leipzig 1877, p. 729. \href{https://de.wikisource.org/wiki/ADB:Eichhorn,_Christian_Friedrich}{Wikisource.}}

What little else we can learn about Eichhorn comes from his published works and their reception. The first of his traceable efforts is a play dating from 1824 (and thus a product of his student years in G\"{o}ttingen) titled \emph{Chriemhildens Rache. Ein Trauerspiel.} (\emph{Chriemhild's Revenge. A Tragedy.}), based on the Nibelungenlied. In the 1877 book \emph{Die deutsche Sage von den Nibelungen in der deutschen Poesie (The German Legend of the Nibelungs in German Poetry)}, Karl Rehorn\index[ppl]{Rehorn, Karl Wilhelm Arnold}\footnote{Karl Wilhelm Arnold Rehorn (1840--1917) was a teacher and literary historian, and was headmaster of the Elisabeth School in Frankfurt am Main from 1881 to 1900.} writes, ``A truly unique phenomenon is C. F. Eichhorn's ``Chriemhild's Revenge'', so singular that it would be difficult to find anything in dramatic literature that could rival it. During the course of the play, thunder roars twelve times; furthermore, at various points, the following appear: a storm on the Rhine; darkness; an earthquake; a glare; a mountain that sinks, leaving a chasm from which smoke rises; a comet; a fireball; will-o'-the-wisps; and the ghosts of Siegfried, Brunhild, Ute, and others.'' At the end of a five-page synopsis of the production, Rehorn reports, ``Since the entire work strives solely for a sensual effect, everything else recedes behind this tendency. The fable is not understood at all; but we also do not find that the poet even grasped the tragic elements\ldots. Even less could one speak of a characterization of the individual persons. \ldots[A] lush, empty declamation proliferates, riddled with meaningless lyricism; and both are overtaken by a cloying sentimentality that pushes its way forward in a repulsive manner.''\footnote{By way of softening this critique, it should be noted that by the end of the 19th century, when Rehorn was writing, the cultural phenomenon of Romanticism was out of favor; Rehorn himself calls the movement ``degenerate".}

Eichhorn appears in print again in 1827, with \emph{Dissertatio Philosophico-Mathematica de Semiologistica ex Principiis Arithmeticus Deducta, Quam Amplissimi Philosophorum Ordinis in Academia Georgia Augusta}. In the preface to this Latin work, we find the mysterious statement, ``[A]lthough some semiological notions were imperfectly explained in my doctoral dissertation, it seemed to me that my duty arose from this, to pursue the principles further and to render them more clear and solid[.]''\footnote{This bibliography's compiler, not conversant in Latin, could not readily parse the machine-translated text of Eichorn's \emph{Philosophical-Mathematical Dissertation on Semiology Deduced from the Principles of Arithmetic}.}

Only two other bibliographic references have surfaced. The first is to an 1828 book, \emph{Versuch einer Entwicklungskarte der allgemeinen reinen Mathematik (An attempt at a developmental map of general pure mathematics)}, which we have thus far been unable to locate. The second is \hyperref[Eic1834]{\textsc{Eichhorn 1834}}\index[ppl]{Eichhorn, Christian Friedrich}, published two years before Eichhorn's death at the age of 32. Through the following decades the works of ``C. F. Eichhorn'' were sometimes mistakenly attributed to Carl (also spelled Karl) Friedrich Eichhorn\index[ppl]{Eichhorn, Carl Friedrich} (1781--1854), a noted German legal scholar. 


\subsection*{Moses Ensheim \normalfont{(1750--1839)}} \hyperref[Ens1799]{\textsc{Ensheim 1799.}}\index[ppl]{Ensheim, Moses}

Also known as Brisac and Moses Metz, Ensheim was born in Metz, France. Schecter\footnote{Ronald Schechter\index[ppl]{Schechter, Ronald}, \emph{Obstinate Hebrews: Representations of Jews in France, 1715-1815}, University of California Press, 2003. \href{https://www.jstor.org/stable/10.1525/j.ctt1pnvn3}{JSTOR.}} states, ``[He] was born to a poor family, received a traditional Jewish education, found his teachers lacking, and turned to an intensive autodidacticism, teaching himself French, German, Latin, Greek, and Arabic.'' Defying his parents, who wanted him to become a rabbi, he left Metz and led a somewhat peripatetic life traveling in Germany. For a few years beginning in 1782 he was a preceptor or tutor in the Berlin home of German-Jewish philosopher and theologian Moses Mendelssohn\index[ppl]{Mendelssohn, Moses}, especially to Mendelssohn's son Abraham\index[ppl]{Mendelssohn, Abraham} (whose own distinguished children include composers Fanny\index[ppl]{Mendelssohn, Fanny} and Felix Mendelssohn\index[ppl]{Mendelssohn, Felix}). Moses Mendelssohn died in 1786; around this time Ensheim returned to Metz and made a meager living by teaching mathematics. He applied for the position of professor of mathematics at the Ecole Centrale at Metz but, being a Jew, was rejected.\footnote{\href{https://www.jewishencyclopedia.com/articles/5783-ensheim-moses}{jewishencyclopedia.com.}}

Ensheim is reputed to have been personally acquainted with Lagrange\index[ppl]{Lagrange, Joseph-Louis}, Laplace\index[ppl]{Laplace, Pierre-Simon}, and Monge\index[ppl]{Monge, Gaspard}, and with the French chemist and senator Berthollet\index[ppl]{Berthollet, Claude Louis};\footnote{Ibid. We have not been able thus far to independently verify Ensheim's acquaintance with these men, apart from Lagrange.} he was also a personal friend of the French Catholic priest, revolutionary leader, and abolitionist Henri (also known as Abb\'{e}) Gr\'{e}goire\index[ppl]{Gr\'{e}goire, Henri (\emph{aka} Abb\'{e}) Gr\'{e}goire}.\footnote{Schecter, ibid., p. 185}

A devout Jew, Ensheim was active in the school of the Meassefim\index{Meassefim}, a movement whose focus was educational reform and Hebrew language revival. In Hebrew scholarship Ensheim's name is mainly associated with his role as correspondent on the French Revolution for the Berlin Haskalah journal \emph{Ha-measef},\footnote{Shmuel Feiner\index[ppl]{Feiner, Shmuel}, \emph{The Jewish Enlightenment}, University of Pennsylvania Press (2004), p. 269. \href{http://www.jstor.org/stable/j.ctt3fh815.16}{JSTOR.}} and for having written a ``canticle depict[ing] the abolition of the French monarchy and the victories of the early revolutionary wars in biblical language.''\footnote{Ari Joskowicz\index[ppl]{Joskowicz, Ari}, Jewish anticlericalism and the making of modern Jewish politics in Late Enlightenment Prussia and France, \emph{Jewish Social Studies}, Vol. 17, No. 3 (Spring/Summer 2011), pp. 40-77. \href{https://www.jstor.org/stable/10.2979/jewisocistud.17.3.40}{JSTOR.}}

His later years were spent with the family of Abraham Furtado\index[ppl]{Furtado, Abraham}, for whom he had once been a tutor.\footnote{\href{https://www.jewishencyclopedia.com/articles/5783-ensheim-moses}{jewishencyclopedia.com}.}


\subsection*{Eduard Heis} \normalfont{(1806--1877).} \hyperref[Hei1844]{\textsc{Heis 1844.}}\index[ppl]{Heis, Eduard}

Born in 1806 in Cologne, Germany, Heis graduated from the Carmelite school there in 1824.\footnote{Most of the Heis biography is adapted or quoted from the \emph{Catholic Encyclopedia} (1913) entry (\href{https://en.wikisource.org/wiki/Catholic_Encyclopedia_(1913)/Eduard_Heis}{Wikisource}) written by Johann Georg Hagen\index[ppl]{Hagen, Johann Georg}, an Austrian Jesuit priest, astronomer, and student of Heis's; and from Ulrich Sperberg\index[ppl]{Sperberg, Ulrich}, ``Eduard Heis, an early pioneer in meteor research''. \emph{History of Geo- and Space Sciences}, \textbf{12}, no. 2 (2021), 163--170  (\href{https://hgss.copernicus.org/articles/12/163/2021/}{hgss.copernicus.org}).} At the University of Bonn\index{University of Bonn}, he solved two prize questions involving the reconstruction and interpretation of Latin texts by Apollonius Perg{\ae}us\index[ppl]{Apollonius Perg{\ae}us} and Cicero\index[ppl]{Cicero}. Upon graduating in 1827, he taught math and science for ten years at the Friedrich Wilhelm Gymnasium in Cologne, followed by 15 years at the B\"{u}rger- und Provinzial-Gewerbeschule (vocational school) in Aachen. He married Margareta Will\index[ppl]{Heis, Margareta (\emph{n\'{e}e} Will)} in 1833; two of their five children died young.

In 1852, at Alexander von Humbolt's\index[ppl]{Humbolt, Alexander von} request, Heis was appointed by King Frederick William IV\index[ppl]{Frederick William IV, King} to the chair of mathematics and astronomy at the K\"{o}nigliche Theologische und Philosophische Akademie,\footnote{Now the University of M\"{u}nster\index{University of M\"{u}nster}.} was elected rector of that institution in 1869, and remained there until his death (three months before his 50th anniversary as a teacher) from a stroke.

An ``inexhaustible'' astronomical observer for most of his life, Heis ``was favoured with an unusually keen sight and, therefore, could see far more than ordinary observers with unassisted eyes.'' Although a four-inch telescope was mounted on the Akademie's roof, he observed the night sky mostly by naked eye for over four decades. His major work, \emph{Atlas Coelestis} (the result of 27 years of observations, published in Cologne in 1872, and dedicated to Pope Pius IX\index[ppl]{Pius IX, Pope}) comprises ``12 charts, a catalogue of 5,421 stars, and the first true delineation of the Milky Way.'' His \emph{Zodiakal-Beobachtungen (Zodiacal Observations)} (1875) records the zodiacal light seen between 1847 and 1875. Heis was the first to grade galactic luminosity with a much-used 1--5 brightness scale, plotted on graph contours. His records of variable stars from 1840 to 1870 were published in 1903.

Heis wrote\footnote{Die Sternschnuppen, Feuerkugeln und Meteorsteine, \emph{Natur und Offenbarung}, \textbf{1}, 216--227. \href{https://books.google.co.jp/books?id=m0xOzHUcAHIC&hl=ja&source=gbs_book_other_versions_r&cad=3}{Google Books}.} that the year before he graduated from college he saw a penknife whose blade was made of celestial iron; this may have inspired his later interest in meteors and meteorites\index{meteorites}. He was the first observer to make a precise hourly count for the August Perseids meteor shower\index{Perseids meteor shower}, finding a maximum rate of 160 meteors per hour in 1839. Over 43 years he and his students compiled over 15,000 observations of meteors, cataloged in \emph{Sternschnuppen-Beobachtungen (Shooting Star Observations)} (1877).

He was interested in sunspots and the aurora borealis, and he wrote ``treatises on the eclipses during the Peloponnesian War (1834), on Halley's comet (1835),{\ldots}on periodic shooting stars (1849)'', and on the star Mira [actually a binary star system in the constellation Cetus] (1859). He wrote several successful mathematics textbooks in addition to \hyperref[Hei1844]{\textsc{Heis 1844}}, including \emph{Lehrbuch der Geometrie zum Gebrauche an h\"{o}heren Lehranstalten Dritter Theil, Ebene und sph\"{a}rische Trigonometrie} (1867), which also saw many editions.

Heis directed that his ``tombstone [be] prepared in the proportions of the `golden section', with the symbol of the dove and olive-branch from the catacombs.'' A lunar crater in the Mare Imbrium is named for him.


\subsection*{George Winslow Pierce \normalfont{(1841--1917)}} \hyperref[Pie1891]{\textsc{Pierce 1891.}}\index[ppl]{Pierce, George Winslow} 

Pierce was born in Boston, Massachusetts, USA. A copy of an 1891 book called \emph{The Life-Romance of an Algebraist},\footnote{\emph{The Life-Romance of an Algebraist}, J. G. Cupples, Boston, 1891. \href{https://hdl.handle.net/2027/uc2.ark:/13960/t6542wp0v}{HathiTrust.}} which as of this writing is still held in the Wheaton College Library, has two inscriptions on its front fly-leaf. The first is ``To Rev. John Perry Barrett\index[ppl]{Barrett, Rev. John Perry} with the author-brother's constant regards.'' The second reads as follows:

\vspace{9pt}

\begin{quote}
To Wheaton College Library.

A curiosity!

And a most pathetic account of the breaking of what was perhaps one of the brightest minds of America. He graduated with a mark of 97 for the entire course, a rank almost unknown at Harvard. Of his mathematical ability this vol. will bear witness. It shows him also to be a poet of no mean order{\ldots}. The causes of trouble were intense nervous strain{\ldots}and the keen life sorrow of love unreturned{\ldots}.

His friend and classmate

Jno. P. Barrett

Wheaton. Dec. 1903.\footnote{The ellipses are in the Rev. Barrett's handwriting.}

\end{quote}

\vspace{9pt}

The Rev. Barrett's ``author-brother'' was George Winslow Pierce, and the description of Pierce's book as ``a curiosity'' is gracious. One of the first things one notices about \emph{The Life-Romance} is that the text is aligned parallel to the spine, rather than perpendicular to it as is typical. The twenty-one pages of front matter comprise testimonials to Pierce from Harvard administrators and faculty\footnote{The testimonials bearing dates were twenty years old or more when the book was published. One from Harvard president Charles W. Eliot\index[ppl]{Eliot, Charles W.}, addressed to the ``Boston Latin School'' and dated 1870, might have been in support of his employment there; see the obituary below.}  (including a facsimile of a handwritten, undated note by the unrelated mathematician Benjamin Pierce\index[ppl]{Pierce, Benjamin}), and a brief autobiographical sketch which begins to acquire a surreal tone around page XI. The first chapter is titled ``Solutions of the General Algebraic Equation'' and leads off with a polynomial-like equation called ``IT'', in which ``all $x\Sigma$'s whatever and complete continued products (embryo $\Sigma$'s) are constants (given or calculable), and have parts which are or go into \emph{permanent} functions{\ldots}'' and so on for four more lines. The book then rambles for about 150 more pages, alternating between math symbolism, rhymed poetry in Latin and English, and stream-of-consciousness narrative, all heavily punctuated with exclamation marks.

What little can be gleaned of Pierce's life comes from the parts of his twelve self-published books\footnote{\href{https://search.worldcat.org/search?q=au\%3D\%22Pierce\%2C+George+Winslow\%22\&offset=1}{List of works by George Winslow Pierce (WorldCat).} One of Pierce's publications claims to be a proof of Fermat's Last Theorem.} that seem the most calmly sincere, and from this obituary in the \emph{Boston Globe}, evening edition, Friday, November 9, 1917, page 14:\footnote{Pierce died seven months after the United States entered World War I. The \emph{Globe} gave his death notice prominence; it stands at the top of the page, with the boldface headline ``GEORGE WINSLOW PIERCE IS DEAD'' in the same type treatment as ``COAL CONSERVATION AN ABSOLUTE NEED'' and ``OHIO PROHIBITION SITUATION IS SAME'', suggesting that Pierce may have been well-known.}

\vspace{9pt}

\begin{quote}

George Winslow Pierce, author of many eccentric and enigmatic writings, published by himself, and first scholar of the class of 1864 at Harvard, died yesterday at the Massachusetts General Hospital of cerebral hemorrhage. 

He was born in Boston, March 21, 1841, the son of John Winslow\index[ppl]{Pierce, John Winslow} and Lydia Ann (Osborne) Pierce\index[ppl]{Pierce, Lydia Ann (Osborne)}. He fitted himself for college. In Harvard, which he entered in 1860, he distinguished himself as the leading scholar of his class and was called by Pres. Charles W. Eliot, ``preeminent in mathematics."

In the Fall of 1864 he entered the Harvard Law School, where he remained a year, then sail[ed] for a six months' European trip, which included points of special interest, such as the tip top of Mont Blanc. In the Fall and Winter of 1866-7 he was tutor in mathematics at the college, substituting for J. M. [James Mill] Pierce\index[ppl]{Pierce, James Mill}, who was absent. Mr Pierce delighted to call himself ``substitutor," in his case a contraction for substitute tutor.

In 1867 he entered the law office of Chandler, Shattuck \& Thayer and the next year was admitted to the Suffolk Bar. Directly afterward he went to Minnesota in search of health, there learning to shoot and hunt.

After being submaster of English High School from 1868 to 1870 and for the next three years mathematics instructor at Boston Latin School, he opened his own law office, but each year, in protection of his health, he made a long hunting trip, five into the New York woods and 10 to the Ottawa Valley. In one of these he shot a huge moose, said to be one of the largest ever bagged, the head of which hangs over the entrance to Memorial Hall at Harvard.\footnote{A photo of this taxidermy is on unnumbered page 170 of \emph{The Life-Romance}.}

Mr Pierce was a student of literature, as of mathematics, and a learned man, but his prose works, one of which was filled with algebraic symbols, tables and groupings, were eccentric, though many admitted they had the stamp of genius. Original and peculiar methods of expression characterized them.\footnote{On page 35 of \emph{The Life-Romance} we find ``A NEW PRONOUN (ITH)!{\ldots} Ith is dative, nominative (iths is genetive), and might be accusative or ablative \emph{man, woman, or child}; civil or savage; wild or domestic \emph{beast} or \emph{bird}, \&c.'' along with the example ``If I found a fellow-being in distress, and thought I could do ith any good, it wouldn't make any difference to me what ith \emph{had} done{\ldots}.''} He also wrote odd verses and composed songs, setting the book of Ruth to music. ``Weary Walker's Song," ``Life Romance of an Algebraist,'' ``Sour Grapes," ``An Autobiography in Rhyme," ``Select Círcle" and ``My Sound Speed Discovery[:] Wit and Song," were some of his works.\footnote{In another of his books, \emph{The Song of Song Which Is Solomon's}, copyrighted 1915 and comprising more than 400 pages, Pierce rewrites the Biblical ``Song of Solomon'' in rhymed couplets and prints one line per page.}

He professed to have invented a new method of punctuation with free use of italics and parentheses, and thus became a ``professor of punctuation.” Of recent years he had devoted himself to literary pursuits at his home [on] West Concord [S]t. He had not married. 

\end{quote}

\vspace{9pt}

A few details might be added. From genealogical records, it appears that Pierce had a younger sibling, Hezekiah E. Pierce\index[ppl]{Pierce, Hezekiah E.}, born in 1854, but who died seven years later.\footnote{\href{https://www.findagrave.com/memorial/281010065/george-winslow-pierce}{findagrave.com.}} The publication date of \emph{The Life-Romance} suggests that Pierce had become wealthy enough from his law practice to afford the self-publishing of his typographically challenging books. On page X of \emph{The Life-Romance}'s front matter, he describes tutoring Benjamin Pierce's\index[ppl]{Pierce, Benjamin} daughter-in-law's nephew in 1864, his senior year at Harvard. George was under great pressure to fix an error in homework that was due in less than an hour. ``My pupil was snapping a big knife in the window again and again, which took off my attention\ldots. Full of power, and suspecting no weakness, I bent myself, after every snap{\ldots}to the task before me; but the last snap snapped something in my head, behind, and a trembling thrill spread slowly over the cerebellum and settled in a pressure in my \emph{ears!} I had to shirk my lessons---it was May, and a week's vacation helped me---the remaining few weeks of the College Course; but a glib tongue and my reputation carried me through!'' For what sounds like a kind of seizure, Pierce was treated by the distinguished endocrinologist Charles E. Brown-S\'{e}quard\index[ppl]{Brown-S\'{e}quard, Charles E.}, but ``anyone who had stood under my chamber window the night I took [the doctor's medications] would have thought there was a prayer-meeting overhead. I didn't call any human being to my assistance, knowing I should live through it[.]'' The Rev. Barrett's\index[ppl]{Barrett, Rev. John Perry} diagnosis of ``nervous strain'', Pierce's own description of his senior year breakdown, the onset of eccentric authorship 25 years later, and the report of cerebral hemorrhage as the cause of death---one wonders if these indicate a progressive disease or mental illness.




\normalsize
\newpage
\section{Sources}\label{S:sources}

{\small

\bigskip

\href{https://dcms.beloit.edu/}{Beloit College Digital Collections}\index{Beloit College Digital Collections}

\href{https://digital.beic.it/}{Biblioteca Europea}\index{Biblioteca Europea}

\href{https://www.bnf.fr/}{Biblioth\`{e}que nationale de France (BnF)}\index{Biblioth\`{e}que nationale de France (BnF)}

\href{https://www.kbr.be/en/}{Biblioth\`{e}que royale de Belgique (KBR)}\index{Biblioth\`{e}que royale de Belgique (KBR)}

DigiZeitschriften\index{DigiZeitschriften} (defunct as of 1 January 2026)

\href{https://www.ele-math.com}{Ele-Math}\index{Ele-Math}

\href{https://digital.slub-dresden.de}{Dresden State and University Library (SLUB)}\index{Dresden State and University Library (SLUB)}

\href{https://library.uaf.edu/}{Elmer E. Rasmuson Library, University of Alaska Fairbanks}\index{Rasmuson Library, Elmer E.}\index{University of Alaska Fairbanks}

EMIS --- European Mathematical Information Service (defunct)\index{EMIS}

\href{https://scholarlycommons.pacific.edu/euler/}{The Euler Archive}\index{Euler Archive, The}

\href{https://www.findagrave.com}{FindAGrave.com}\index{FindAGrave.com}

\href{https://gallica.bnf.fr/}{Gallica}\index{Gallica}

\href{https://books.google.com/}{Google Books}\index{Google}

\href{https://translate.google.com/}{Google Translate}\index{Google}

\href{https://gdz.sub.uni-goettingen.de/gdz/}{G\"{o}ttinger Digitalisierungzentrum}\index{G\"{o}ttinger Digitalisierungzentrum}

\href{https://www.hathitrust.org/}{HathiTrust Digital Library}\index{Hathitrust Digital Library}

\href{https://archive.org/}{Internet Archive}\index{Internet Archive}

\href{https://www.jewishencyclopedia.com}{Jewish Encyclopedia}

\href{https://www.jstor.org/}{JSTOR}\index{JSTOR}

\href{https://www.jstage.jst.go.jp}{J-STAGE}\index{J-STAGE}

\href{https://en.langenscheidt.com/}{Langenscheidt}\index{Langenscheidt}

\href{https://www.loc.gov/}{Library of Congress}\index{Library of Congress}

\href{https://www.mathdoc.fr/}{Mathdoc}\index{Mathdoc}

\href{https://www.ams.org/mathscinet/}{MathSciNet}\index{MathSciNet}

\href{https://www.digitale-sammlungen.de}{M\"{u}nchener DigitalisierungsZentrum}\index{M\"{u}nchener DigitalisierungsZentrum}

\href{http://www.numdam.org/}{NUMDAM --- Num\'{e}risation de documents anciens math\'{e}matiques}\index{NUMDAM}

\href{https://sol.unibo.it/SebinaOpac/}{OPAC Sebina OpenLibrary}\index{OPAC Sebina OpenLibrary}

\href{https://openlibrary.org/}{OpenLibrary}\index{OpenLibrary}

Poltran (defunct)\index{Poltran}

\href{https://www.reverso.net/}{Reverso}\index{Reverso}

\href{https://library.si.edu/digital-library/}{Smithsonian Libraries Digital Library}\index{Smithsonian Libraries Digital Library}

\href{https://digital.ulb.hhu.de/}{University and State Library of D\"{u}sseldorf}\index{University and State Library of D\"{u}sseldorf}

\href{https://de.wikisource.org/}{Wikisource (Germany)}\index{Wikisource (Germany)}

\href{https://search.worldcat.org/}{WorldCat}\index{WorldCat}

\href{https://www.worldlingo.com/}{WorldLingo}\index{WorldLingo}

\href{https://mathworld.wolfram.com/}{Wolfram MathWorld}\index{Wolfram MathWorld}

}

\end{RaggedRight}


\index{alternating Euler's constant|seeonly{Euler's constant!alternating}}
\index{continued compositions|see{$f$-expansions, infinite processes, inner compositions, \emph{Kettenoperationen}}}
\index{continued fractions!generalizations of|see{$f$-expansions}}
\index{continued functions|see{continued compositions}}
\index{continued function compositions|seeonly{continued compositions}}
\index{continued nested radical fractions|see{continued reciprocal roots}}
\index{continued operations|see{continued compositions, \emph{Kettenoperationen}}}
\index{continued pobb powers@continued $p$th powers!$p\in(-\infty,-1)$|seeonly{continued reciprocal powers}}
\index{continued pobb powers@continued $p$th powers!$p=-1$|seeonly{continued fractions}}
\index{continued pobb powers@continued $p$th powers!$p\in(-1,0)$|seeonly{continued reciprocal roots}}
\index{continued pobb powers@continued $p$th powers!$p\in(0,1)$|seeonly{continued radicals, continued $r$th roots, continued $r_i$th roots, continued square roots}}
\index{continued pobb powers@continued $p$th powers!$p=1$|seeonly{infinite sums}}
\index{continued pobb powers@continued $p$th powers!$p=2$|seeonly{continued squares}}
\index{continued pobb powers@continued $p$th powers!$p\in(1,\infty)$|seeonly{continued powers}}
\index{continued pobb powers@continued $p$th powers!$p=\tfrac{1}{2}$|seeonly{continued square roots}}
\index{continued pobb powers@continued $p$th powers!$p=\tfrac{1}{3}$|seeonly{continued cube roots}}
\index{continued pobb powers@continued $p$th powers!$p=\tfrac{1}{r}, r\in(1,\infty)$|seeonly{continued $r$th roots}}
\index{continued radicals|see{continued square roots, continued cube roots, continued $r$th roots, continued $r_i$th roots}}
\index{continued roaa roots@continued roots|seeonly{continued square roots, continued cube roots, continued $r$th roots, continued $r_i$th roots, nested radicals}}
\index{continued square roots!infinite product of|seeonly{infinite products, Vi\`{e}te's formula for $\tfrac{2}{\pi}$}}
\index{convergents|see{approximants}}
\index{cross division|seeonly{ordered division}}
\index{Fourier division|seeonly{ordered division}}
\index{infinite exponentials|see{continued exponentials}}
\index{infinite radicals|see{continued square roots, continued cube roots, continued $r$th roots, continued $r_i$th roots, nested radicals}}
\index{infinite series|seeonly{infinite sums}}
\index{inner compositions|see{continued compositions, $f$-expansions}}
\index{iterated compositions|see{outer compositions, successive substitution}}
\index{iterated radicals|seeonly{continued square roots, iterated square roots}}
\index{iteration!cyclic orbits|see{periodic points}}
\index{Kepler's equation|seeonly{Kepler problem}}
\index{Kettenfunctionen@\emph{Kettenfunctionen}|see{continued compositions}}
\index{Kettenoperationen@\emph{Kettenoperationen}|see{continued compositions}}
\index{Kettenwurzeln@\emph{Kettenwurzeln}|see{continued radicals, continued square roots, continued cube roots, continued $r$th roots, continued $r_i$th roots, infinite radicals, nested radicals}}
\index{nested radicals|see{continued square roots, continued cube roots, continued $r$th roots, continued $r_i$th roots, infinite radicals}}
\index{nested square roots|see{continued square roots}}
\index{numbers!real|seeonly{real numbers}}
\index{outer compositions|see{iterated compositions}}
\index{series|see{infinite sums}}
\index[ppl]{Snell, Willebrord|seeonly{Snellius, Willebrordus}}
\index{tower|seeonly{continued exponentials}}

\indexprologue{\small Indices $i$ and $n$ are assumed to be nonnegative integers.}
{\small\printindex}
\indexprologue{\small \textbf{Bold} type indicates a work by the named person. \emph{aka}: also known as.}
{\small\printindex[ppl]}

\end{document}